\documentclass{article}

\usepackage[francais]{babel}
\usepackage[utf8]{inputenc}
\usepackage{amsmath,amsthm,amsfonts,stmaryrd,fullpage,bbold,txfonts,graphicx}

\newtheorem{prout}{Théorème} [subsection]
\newtheorem*{etheo}{Théorème}
\newtheorem{prop}[prout]{Proposition}
\newtheorem{lemme}[prout]{Lemme}

\newtheorem{cor}[prout]{Corollaire}
\newtheorem{defin}[prout]{Définition}
\newtheorem{defins}[prout]{Définitions}

\theoremstyle{remark}

\newtheorem{exemples}[prout]{Exemples}
\newtheorem{exemple}[prout]{Exemple}

\DeclareMathAlphabet{\mathpzc}{OT1}{pzc}{m}{it}

\newcommand{\fex}[1]{ \begin{exemple} #1 \end{exemple} }
\newcommand{\fexs}[1]{ \begin{exemples}\ \begin{enumerate} #1\end{enumerate} \end{exemples} }
\newcommand{\fprop}[1]{ \begin{prop} #1 \end{prop} }
\newcommand{\fdef}[1]{ \begin{defin} #1 \end{defin} }
\newcommand{\flemme}[1]{ \begin{lemme} #1 \end{lemme} }
\newcommand{\fcor}[1]{ \begin{cor} #1 \end{cor} }
\newcommand{\remas}[1]{\textit{Remarques:}\ \begin{itemize} #1 \end{itemize} }
\newcommand{\fprops}[2][\ ]{ \begin{prop} #1 \begin{enumerate}#2 \end{enumerate}\end{prop} }
\newcommand{\flemmes}[2][\ ]{ \begin{lemme} #1 \begin{enumerate}#2 \end{enumerate}\end{lemme} }
\newcommand{\fdefs}[2][\ ]{ \begin{defins} #1 \begin{itemize}#2 \end{itemize}\end{defins} }
\newcommand{\demos}[1]{ \textit{Démonstration: } \begin{enumerate} #1 \end{enumerate}\hfill $\square$ \vspace{1cm}}

\newcommand{\liste}[1]{\begin{itemize} #1 \end{itemize} }
\newcommand{\enum}[1]{\begin{enumerate} #1 \end{enumerate} }

\newcommand{\fix}{\mathrm{Fix}}
\newcommand{\stab}{\mathrm{Stab}}
\newcommand{\aut}{\mathrm{Aut}}
\newcommand{\card}{\mathrm{Card}}

\newcommand{\vect}{\mathrm{Vect}}
\newcommand{\aff}{\mathrm{Aff}}

\newcommand{\demo}{\textit{Démonstration: }}
\newcommand{\pv}{\textit{Preuve du lemme: }}
\newcommand{\rema}{\textit{Remarque: }}

\newcommand{\fl}[1]{\overrightarrow{#1}}

\newcommand{\ital}[1]{\textit{#1}}

\def\phii{\varphi}
\def\maj{\mathpzc}
\def\tens{\otimes}

\def\barre{\overline}
\def\parallele{\varparallel}

\newcommand{\infini}{\infty}
\newcommand{\vide}{\emptyset}

\newcommand{\fleche}{\rightarrow}
\newcommand{\implique}{\:\Rightarrow\:}
\def\ssi{\:\Leftrightarrow\:}

\newcommand\ent[2] {\llbracket #1,#2 \rrbracket}
\newcommand{\tq}{\text{ tq }}
\newcommand{\et}{\text{ et }}
\newcommand {\inv} {^{-1}}
\newcommand{\cqfd}{\hfill $\square$ \vspace{.6cm}}

\def\inte{\mathring}

\newcommand{\calc}[1]{  \begin{eqnarray*}#1 \end{eqnarray*}  }

\newcommand{\N}{\mathbb N}
\def\K{\mathbb K}
\newcommand{\R}{\mathbb R}
\newcommand{\Z}{\mathbb Z}
\newcommand{\Q}{\mathbb Q}

\def\L{\mathcal L}

\newcommand{\ens}[2]{ \left\{ #1 \; | \; #2 \right\} }
\newcommand{\eng}[1]{\left\langle\:  #1\: \right\rangle } 
\newcommand{\engens}[2]{\left\langle\: \left\{ #1 \; | \; #2 \right\}\: \right\rangle }

\newcommand {\fonc}[4] {
\begin{array}{rcl}
 {#1} & \rightarrow & {#2} \\
 {#3} & \mapsto     & {#4}
\end{array}}

\def\I{\mathcal I}
\def\ikv{\vec{\mathcal I}(\K)}

\def\iv{\vec{\mathcal I}}

\def\D{\mathcal D}
\def\F{\mathcal F}

\def\Q{\mathcal Q}
\def\P{\mathcal P}

\def\L{\mathbb L}
\def\M{\mathbb M}
\def\g{\mathfrak g}

\def\point{\; . \;}
\def\gal{\mathpzc{Gal}}
\def\gm{\maj{Germ}}
\def\cl{\text{Cl}}
\def\op{\text{op}}
\def\sph{_\text{sph}}

\newcommand{\die}[1][ ]{_{\sharp #1}}
\newcommand{\bec}[1][ ]{^{\natural #1}}
\newcommand{\becb}[1][ ]{_{\natural #1}}

\def\lien1{2.1}
\def\lien2{2.2}

\def\pinj{($\text{para } inj$)}

\def\psph{($\text{para } sph$)}
\def\pfonc{(fonc)}
\def\pdec{($\text{para } dec$)}
\def\pisom{($\text{para } 5$)}
\def\pinter{($\text{para } 6$)}

\def\plien{(para 2.1)}
\def\plieng{(para 2)}
\def\plienn{(para  2.2)}
\def\plienp{($\text{para } 2.1^+$)}
\def\plienpm{($\text{para } 2.1^{+-}$)}
\def\plienpn{($\text{para } 2.2^+$)}

\renewcommand{\abstract}[1] { \begin{center}\textbf{Abstract} \end{center} \hspace{.5cm} #1} 
\newcommand{\resume}[1]  { \begin{center}\textbf{Résumé} \end{center} \hspace{.5cm} #1}

\begin{document}

\title{Structures immobilières pour un groupe de Kac-Moody sur un corps local}
\author{Cyril Charignon}
\maketitle

\abstract{In this study, we try to generalize Bruhat-Tits's theory to the case of a Kac-Moody group, that is to define an affine building for a Kac-Moody group over a local field. Actually, we will obtain a geometric space wich lacks some of the incidence properties of a building, so that it is called a hovel, following Guy Rousseau's terminology. Hovels have already been obtained for split Kac-Moody groups by Guy Rousseau and Stéphane Gaussent; here we define them for any group with a (generalized) valuated root datum, a situation wich contains the Kac-Moody groups over local fields, split and nearly split.
}

\resume{Le but de ce travail est de généraliser la théorie de Bruhat-Tits au cas des groupes de Kac-Moody, c'est-à-dire de construire un immeuble affine pour un groupe de Kac-Moody sur un corps local. L'objet obtenu ne sera en fait pas un immeuble, car il ne vérifie pas toutes les conditions d'incidence nécessaires, il s'agira plutôt d'une "masure", selon la terminologie de Guy Rousseau. Des masures ont déjà été définies par Guy Rousseau et Stéphane Gaussent pour des groupes de Kac-Moody déployés; ici on propose une construction valable pour n'importe quel groupe muni d'une donnée radicielle valuée (généralisée pour comprendre le cas d'un système de racines infini), ce qui fournira en particulier des masures pour les groupes de Kac-Moody déployés, mais aussi presque déployés, sur un corps local.
}

\tableofcontents

\section{Introduction}

Lorsque $G$ est un groupe réductif sur un corps local $\K$, la théorie de François Bruhat et Jacques Tits permet de lui associer, outre son immeuble vectoriel, ou sphérique, ou "de Tits", un immeuble affine, dit "de Bruhat-Tits". Cet immeuble est plus précis que l'immeuble vectoriel, au sens où la connaissance du premier permet de reconstruire le second. Comme tout immeuble, il s'agit d'une réunion d'appartements, et on le qualifie d'"affine" car ces appartements sont des espaces affines, dont les espaces vectoriels directeurs peuvent en fait être vus comme les appartements de l'immeuble vectoriel.

Maintenant, si $G$ est un groupe de Kac-Moody, déployé ou presque déployé, on sait construire son immeuble vectoriel. Il s'agit en fait de la réunion de deux immeubles jumelés, voir par exemple \cite{remy}. Le but du travail présent est, dans le cas où $G$ est un groupe de Kac-Moody sur un corps local, de définir un objet "affine", similaire à l'immeuble de Bruhat-Tits du cas réductif.

 Ceci a déjà été fait dans \cite{hovels} par  Stéphane Gaussent et Guy Rousseau, pour un groupe déployé sur un corps $\K$ dont le corps résiduel contient $\mathbb C$. L'objet défini est appelé "masure" car il ne satisfait pas toutes les conditions demandées à un immeuble habituel. Pour un corps plus général, mais toujours un groupe déployé, Guy Rousseau (\cite{microaffine}) a défini des "immeubles microaffines", que l'on peut voir comme une partie du bord à l'infini d'une masure, semblable au bord rajouté  à un immeuble affine lorsqu'on définit sa compactification polygonale. Il a également défini des masures, dans \cite{masures2}.
 
 On se propose ici de se placer dans le cadre général des groupes munis d'une donnée radicielle valuée, cadre qui inclut les groupes de Kac-Moody déployés mais aussi presque déployés sur un corps muni d'une valuation réelle quelconque. Nous construirons simultanément une "masure" et son bord (qui contiendra deux "immeubles microaffines") car celui-ci sera utile à l'étude des propriétés géométriques de la masure. L'objet obtenu sera appelé une masure bordée.
 
 On se rend compte que plusieurs choix peuvent être faits, menant à différentes masures bordées. Pour rentrer un peu dans le détail,  la construction d'un appartement $A$ ne présente pas de difficulté (partie \ref{soussection:appart}) , et on cherche donc ensuite, selon le procédé habituel de Bruhat et Tits, à construire un objet immobilier comme quotient de $G\times A$ par la relation d'équivalence déterminée par le choix des sous-groupes de $G$ qui seront les fixateurs des points de $A$. Ces sous-groupes sont appelés, comme dans le cas réductif, des "sous-groupes parahoriques", et contrairement au cas réductif, leur définition n'est pas évidente: plusieurs possibilités sont envisageables. Ainsi, dans \cite{hovels}, Stéphane Gaussent et Guy Rousseau définissent les familles de groupes parahoriques $\P^{min}$, $\P^{pm}$, $\P^{nm}$ dont la définition est quelque peu indirecte, et nécessite l'emploi de "complétions" du groupe de Kac-Moody $G$ considéré.
 
 Toujours par soucis de généralité, nous optons pour une approche axiomatique, c'est-à-dire que nous définissons abstraitement la notion de "famille de parahoriques", et nous étudions les objets immobiliers que l'ont peut construire, en fonction des propriétés vérifiées par une telle famille. C'est la partie \ref{section:construction generale}. On étudiera plus particulièrement deux familles de parahoriques: la "minimale" puis la "maximale".
 Au \ref{section:descente}, on essaie de descendre les structures précédentes (valuation et famille de parahoriques) à un sous-groupe.  
 Dans la partie \ref{section:KM} enfin, on étudie le cas d'un groupe de Kac-Moody. Pour un groupe déployé, la partie \ref{section:construction generale} s'applique directement, mais pour un groupe presque-déployé, il faut d'abord utiliser les résultats de la partie \ref{section:descente} pour obtenir une valuation et une famille de parahoriques. Le résultat final est le suivant:
 
 \begin{etheo}
Soit $G$ un groupe de Kac-Moody presque déployé sur un corps $\K$, déployé sur la clôture séparable de $\K$. On suppose $\K$ muni d'une valuation réelle discrète non triviale, telle que son corps résiduel soit parfait.

 Alors il existe une masure bordée $\I_\K$ pour $G(\K)$, qui provient d'une valuation $\phii_\K$ et d'une bonne famille de parahoriques $Q_\K$ vérifiant \plienp(sph). Pour toute facette sphérique $\vec f_\K$ de $\iv(\K)$, la façade $\I_{\K,\vec f_\K}$ s'injecte dans la façade $\I_{\L,\vec f}$, de la masure bordée $\I_\L$ pour $G(\L)$, où $\vec f$ est la facette de $\iv(\L)$ contenant un ouvert de $\vec f_\K$.
\end{etheo}

%à faire:
% 
%  - Topologie
%  - fonction op_A
%  - Lien avec l'article "masures" de Guy
%  -rappels sur les ordres grignotants

\section{Rappels et notations}

Lorsque $\alpha: X\rightarrow \R$ est une fonction sur un ensemble $X$, on notera pour $Y\subset X$, $\alpha(Y)=0$ si $\alpha(Y)=\{0\}$, $\alpha(Y)>0$ si $\alpha(Y)\subset \R^{+*}$, $\alpha(Y)\geq 0$ si $\alpha(Y)\subset \R^+$ etc...\\
Si $A$ est un complexe simplicial, $\F(A)$ sera l'ensemble de ses facettes. Si $f$ est une facette de $A$, $f^*$ est la réunion des facettes bordées par $f$. Ainsi, lorsque $\vec f$ est une facette dans un immeuble $\vec \I$, $\vec f^*$ désignera la réunion des facettes de $\iv$ dont l'adhérence contient $\vec f$. Si $\vec Z$ est un appartement de $\iv$, la réunion des facettes de $\vec Z$ dont l'adhérence contient $\vec f$ sera donc $\vec f^*\cap \vec Z$.% Enfin, par définition d'un complexe simplicial, $\F(A)$ est muni d'un ordre qu'on notera juste $\leq$.\\

\subsection{Immeubles vectoriels}

Ce que nous appelons ici les immeubles vectoriels sont les immeubles décrits dans \cite{remy}. Ce sont des immeubles jumelés, donc en fait la réunion de deux immeubles classiques. Dans la réalisation géométrique de ces immeubles que nous considérons, les appartements sont inclus dans des espaces vectoriels, d'où l'appellation "immeubles vectoriels". Une autre appellation fréquente est "immeubles coniques", car les appartements sont des cônes dans ces espaces vectoriels.\\
Dans cette sous-partie, nous rappelons les principaux résultats concernant ces immeubles, et fixons les notations.\\

\subsubsection{Donnée radicielle}
\label{soussoussection:donnee radicielle}

 Il y a plusieurs définitions possibles pour une donnée radicielle, selon que l'on considère qu'un système de racines est un sous-ensemble d'un espace vectoriel réel (comme dans \cite{remy} 6.2.4) ou un ensemble de demi-complexes de Coxeter (comme dans \cite{remy} 1.4.1). La seconde possibilité est plus générale, la première plus précise, elle permet notamment de distinguer une racine et son double. Un système de racines du premier type sera dit vectoriel, un système du second type sera dit géométrique.\\
 
  Si $\alpha,\beta$ sont deux racines d'un système $\phi$, l'intervalle $[\alpha,\beta]$ est défini de la sorte:
 \liste{
 \item $[\alpha,\beta]= \ens{p\alpha+q\beta}{p,q\in\N\et p\alpha+q\beta\in\phi}$ lorsque $\phi$ est un système de racines vectoriel.
 \item $[\alpha,\beta]= \ens{\gamma\in\phi}{\alpha\cap \beta\subset\gamma}$ lorsque $\phi$ est un système de racines géométrique.
 }
 On définit aussi $]\alpha,\beta[= [\alpha,\beta]\setminus \{\alpha,\beta\}$ ainsi que $]\alpha,\beta]$ et $[\alpha,\beta[$ de la manière évidente.\\
Une partie $\psi$ d'un système de racine est dite close lorsque pour tout $\alpha,\beta\in\psi$, $[\alpha,\beta]\subset \psi$. La partie $\psi$ est dite de plus nilpotente si elle est finie. Enfin, une partie $\psi$ est dite prénilpotente s'il existe un système positif $\phi^+$ de $\phi$ et un élément $w\in W(\phi)$ du groupe de Weyl associé à $\phi$ tel que $\psi\subset \phi^+\cap w(-\phi^+)$. Une partie $\psi$ est nilpotente si et seulement si elle est close et prénilpotente.

La notion de prénilpotence est principalement utilisée pour les paires de racines. Si $\{\alpha,\beta\}$ est une telle paire, il est presque immédiat que l'intervalle $[\alpha,\beta]$ est clos, ainsi $\{\alpha,\beta\}$ est prénilpotente si et seulement si $[\alpha,\beta]$ est fini.\\

Dans la suite, sauf mention du contraire, les systèmes de racines considérés seront toujours de type vectoriel. Ainsi, l'existence d'un système de racines $\phi$ sous-entend l'existence d'un $\R$-espace vectoriel $\vec V$ tel que $\phi\subset \vec V^*$. On suppose de plus $\phi$ à base libre, c'est-à-dire que toute base $\Pi$ du système de racines $\phi$ est aussi une base de l'espace vectoriel sous-jacent $\vec V^*$.

L'ensemble des $\ker(\alpha)$, $\alpha\in\phi$ est alors appelé l'ensemble des murs de $\Vec V$, et pour toute racine $\alpha\in\phi$, il existe une réflexion dans $Gl(\vec V)$, notée $r_\alpha$, d'hyperplan fixe $\ker(\alpha)$ qui préserve l'ensemble des murs de $\vec V$. Le groupe engendré par ces $r_\alpha$ est le groupe de Weyl de $\phi$, noté $W(\phi)$. Il existe une famille $(\alpha^\vee)_{\alpha\in\phi}$ de vecteurs de $\vec V$ telle que pour tout $\alpha\in\phi$, la réflexion $r_\alpha$ est donnée par la formule $r_\alpha(\vec v)= \vec v- \alpha(\vec v).\alpha^\vee$. On note pour tout $\alpha,\beta\in\phi$, $\eng{\alpha,\beta} =  \beta(\alpha^\vee)$. On a clairement $\eng{\alpha,\alpha}=2$.

  Un système de racines $\phi$ est dit réduit si pour tout $\alpha\in\phi$, $\phi\cap \R.\alpha =\{\pm\alpha\}$. Lorsque $\phi$ n'est pas réduit, la seule possibilité est en fait $\phi\cap \R.\alpha =\{\pm\alpha,\pm 2 \alpha\}$ ou $\phi\cap \R.\alpha =\{\pm\alpha,\pm \frac{1}{2} \alpha\}$. On notera alors $\phi_{\text{red}}=\ens{\alpha\in\phi}{1/2 \alpha\not\in\phi}$.

\fdef{ Soit $\phi$ un système de racines, et $\phi^+$ un système positif dans $\phi$. Soit $G$ un groupe et $(U_\alpha)_{\alpha\in\phi}$ une famille de sous-groupes de $G$. On note $T=\bigcap_{\alpha\in\phi} N_G(U_{\alpha})$ l'intersection des normalisateurs des $U_\alpha$, $U^+=\engens{U_\alpha}{\alpha\in\phi^+}$ et $U^-=\engens{U_\alpha}{\alpha\in -\phi^+}$.\\

 Le couple $(G,(U_\alpha)_{\alpha\in\phi})$ est appelé une donnée radicielle de type $\phi$ si:\liste{
\item (DR1): Chaque $U_\alpha$ est un sous-groupe de $G$ non trivial.

\item (DR2): Pour toute paire prénilpotente de racines $\{\alpha,\beta\}$, le groupe $[U_\alpha,U_\gamma]$ des commutateurs de $U_\alpha$ et $U_\gamma$ est inclus dans $\engens{U_\gamma}{\gamma\in]\alpha,\beta[}$.

\item (DR3): Si $\alpha \in\phi$ et $2\alpha\in\phi$, alors $U_{2\alpha}$ est inclus strictement dans $U_\alpha$.

\item (DR4): Pour tout $\alpha\in\phi$, et tout $u\in U_\alpha\setminus\{e\}$, il existe $u',u''\in U_{-\alpha}$ tels que $n(u):=u'uu''$ conjugue chaque $U_\beta$, $\beta\in\phi$ en $U_{r_\alpha.\beta}$. De plus, les différents $n(u)$ peuvent être choisis de sorte que pour tout $u,v\in U_\alpha\setminus\{e\}$, $n(u).T=n(v).T$.

\item (DR5): $T.U^+\cap U^-=\{e\}$.\\
}

Cette donnée radicielle est dite génératrice si de plus:\liste{
\item (DRG): $G$ est engendré par $T$ et les $U_\alpha$.\\

Enfin, lorsque le système de racines $\phi$ est fini, on dira que $\D$ est de type fini.
}
}

\remas{\item  C'est la définition utilisée dans \cite{microaffine}, 1.5, elle équivaut à la définition de "donnée radicielle jumelée entière" de \cite{remy} 6.2.5. Dans la terminologie de \cite{remy}, le qualificatif "entière" sert à indiquer que le système de racines est de type vectoriel. La définition d'une donnée radicielle pour un système de racines géométriques est exactement la même, à ceci près que la définition d'un intervalle de racines utilisée en (DR2) a changé, et que (DR3) devient inutile.

 Le qualificatif "jumelé" sert quand à lui à se rappeler que dans le cas où $\phi$ est infini, cette donnée radicielle mènera à un immeuble jumelé. Il n'a aucune signification formelle, ce qui explique son omission ici.
 
 Signalons enfin  que c'est la notion géométrique de donnée radicielle qui est définie dans \cite{abramenko-brown}.\\
 
 \item Dans \cite{bruhat-tits}, la classe $n(u)T$, $u\in U_\alpha$ est notée $M_\alpha$. La condition (DR4) y est exprimée avec les $M_\alpha$ plutôt que les $n(u)$.\\}

Lorsque $\D=(G,(U_\alpha)_{\alpha\in\phi})$ est une donnée radicielle, on notera toujours  $T=\bigcap_{\alpha\in\phi} N_G(U_{\alpha})$ comme dans la définition, c'est le tore maximal associé à $\D$. On prouve que les éléments $n(u)$ dans (DR4) sont uniques, on peut donc conserver la notation $n(u)$.
On note enfin $N$ le sous-groupe de $G$ engendré par ces éléments et par $T$. On prouve facilement que $n(u\inv) = n(u)\inv$, que pour tout $m\in N$, $n(mum\inv)=mn(u)m\inv$, et que si $u'$ et $u''$ sont tels que $n(u)= u' u u''$, alors $n(u)=n(u')=n(u'')$ (la preuve de ce dernier point sera rappelée en \ref{prop:V5}). 

 On prouve que $N$ est le normalisateur de $T$ dès que la condition "(CENT)" définie dans \cite{remy} 1.2.5 est vérifiée. Cette condition s'exprime ainsi:
\[(CENT):\: \forall \alpha\in\phi,\: Z_{U_\alpha}(T)=\{e\}\: .\]
Cette condition est vérifiée par tous les groupes de Kac-Moody sur un corps de cardinal au moins 4 (voir infra), et nous verrons qu'elle l'est aussi pour tous les groupes munis d'une donnée radicielle "valuée" (voir \ref{soussection:valuation}). On la supposera toujours vraie dans la suite.

   Le groupe quotient $N/T$ s'identifie au groupe de Weyl du système de racine $\phi$ en associant pour tout $u\in U_\alpha\setminus\{e\}$, $n(u).T$ à la réflexion $r_\alpha$.\\

Le groupe $T$ et tous ses conjugués sont appelés les tores maximaux de $G$, avec un abus de langage puisque leur définition dépend en fait de la donnée radicielle $\D$. % Si la condition (CENT) est vérifiée, 
Si $T'=gTg\inv$ est un tore maximal, on note $N(T') = gNg\inv$ son normalisateur. On note aussi $g.\vec V$ et $g\phi\subset (g\vec V)^*$ l'espace vectoriel et le système de racines  abstraitement isomorphes à $\vec V$ et $\phi$ via les applications $\vec v\mapsto g.\vec v$ et $\alpha\mapsto g\alpha$. Si $g'\in G$ est un autre élément tel que $T'=g'Tg'{}\inv$, alors $g\inv g'\in N$ donc $g\inv g'$ agit sur $\vec V$ et sur $\phi$, on peut donc identifier $g'\vec V$ à $g\vec V$ et $g'\phi$ à $g\phi$ via $g'.\vec v\mapsto g.(g\inv g')\vec v$ et $g'\alpha\mapsto g(g\inv g').\vec v$.
 Pour tout $\alpha\in\phi$, on note enfin $U_{g\alpha} = gU_\alpha g\inv$, ceci est compatible à l'identification $g\phi=g'\phi$. Alors $(G,(U_{g\alpha})_{g\alpha\in g\phi})$ est encore une donnée radicielle.\\

 Dans la suite, on évitera de particulariser la donnée radicielle $(G,(U_\alpha)_{\alpha\in\phi})$ (correspondant au tore $T$), on considérera plutôt que $G$ est muni d'une classe d'équivalence, pour la conjugaison, de données radicielles. Pour chaque tore $T$ on notera $\phi(T)\subset \vec V(T)^*$ le système de racine et $(U_\alpha)_{\alpha\in\phi(T)}$ les groupes radiciels correspondants. \\

Si $\D=(G,(U_\alpha)_{\alpha\in\phi})$ est une donnée radicielle, pour toute partie $\psi$ de $\phi$, on notera $G(\psi)=\engens{U_\alpha}{\alpha\in\psi}$. Ainsi, lorsque $\D$ est génératrice, on a $G=T.G(\phi)$.\\

Par définition même, un groupe de Kac-Moody déployé admet une donnée radicielle. Et c'est un des buts de \cite{remy} que de prouver que c'est encore le cas pour une classe de groupes de Kac-Moody plus générale. Le terme employé dans \cite{remy} pour qualifier ces groupes est "presque déployé".
\fprop{Si $G$ est un groupe de Kac-Moody déployé ou presque déployé, alors il admet une donnée radicielle de type un système de racine $\phi$ vectoriel à base libre (comme ci-dessus). Si $\K$ est un corps de cardinal au moins 4, alors $G(\K)$ vérifie la condition (CENT).\\
}

\ital{Référence}: \cite{remy} 8.4.1.\cqfd

\subsubsection{L'immeuble d'une donnée radicielle}
\label{sousoussection:immeuble vectoriel}

Une donnée radicielle permet de définir un immeuble vectoriel. On notera dans la suite $\iv$ l'immeuble obtenu à partir de la donnée radicielle $(G,(U_\alpha)_{\alpha\in\phi})$, ou plutôt sa réalisation géométrique (voir \cite{remy} chapitre 5). Dès que $\phi$ est infini, il s'agit en fait de la réunion de deux immeubles jumelés comme définis dans \cite{abra}.

Ses \ital{appartements} sont en bijection avec les tores maximaux de $G$. L'appartement correspondant à un tore maximal $T$ est inclus dans le $\R$-espace vectoriel $\vec V(T)$ tel que $\phi(T)\subset \vec V(T)^*$. Le choix d'une base $\Pi$ de $\phi(T)$ définit un cône $\vec C=\vec C_\Pi=\ens{x\in \vec V(T)}{ \alpha(x)>0,\; \forall \alpha\in\Pi}$, c'est la chambre positive relative à $\Pi$. Les ensembles obtenus en remplaçant certaines des inégalités $>0$ par des égalités $=0$ dans la définition de $\vec C$ sont les \ital{facettes }de $\vec C$. La réunion des facettes de $\vec C$ est donc l'adhérence $\barre{\vec C}$ de $\vec C$. L'appartement $\vec A(T)$ est alors $W(\phi(T)).\barre{\vec C}\cup W(\phi(T)).(-\barre{\vec C})\subset \vec V(T)$, ses facettes sont les $\pm w\vec f$, pour $w\in W(\phi(T))$ et $\vec f$ une facette de $\vec C$. Les \ital{chambres} sont les facettes de dimension maximales, donc les images de $\pm \vec C$ par un $w\in W(\phi(T))$, et les \textit{cloisons} sont les facettes de codimension $1$. C'est un cône, réunion de deux cônes convexes $\vec A^+(T)= W(\phi(T)).\barre{\vec C}$ et $\vec A^-(T)= - \vec A^+(T)$. Le cône positif $\vec A^+(T)$ est appelé le cône de Tits. Chacun de ces deux cônes convexes, avec sa structure de facettes, est un complexe de Coxeter pour $W(T)$.

 L'intérieur, noté $\vec A\sph$, de $\vec A(T)$ dans $\vec V(T)$ est une réunion de facettes, appelées les \textit{facettes sphériques}. Ce sont précisément les facettes dont le fixateur dans $W(\phi(T))$ est fini. Les chambres et les cloisons sont toujours sphériques, leur fixateur dans $W(\phi(T))$ étant respectivement $\{e\}$ et un groupe d'ordre 2.
 
 On pourra noter le groupe de Weyl $W(\phi(T))$ par  $W(T)$, ou $W(\vec A(T))$ pour rappeler qu'il s'agit du groupe de Weyl "vectoriel".\\

 L'immeuble $\iv$ est obtenu en collant les $\vec A(T)$ pour tous les tores maximaux $T$. Ces appartements sont permutés transitivement par $G$, selon la formule $g.\vec A(T) = \vec A(gTg\inv)$ (\cite{remy} 2.6.2). En conséquence, $N(T)$ est le stabilisateur de $\vec A(T)$. Pour $\alpha\in\phi$, $u\in U_\alpha$, l'élément $n(u)$ agit sur $\vec A(T)$ comme la réflexion selon le mur $\ker(\alpha)$. Le groupe $T$ quand à lui est le fixateur de $\vec A(T)$. L'ensemble des facettes de $\iv$ muni de la relation d'ordre "être dans l'adhérence de" est un complexe simplicial, c'est un immeuble  au sens abstrait.
 
 Pour toute racine $\alpha\in\phi(T)$, $\vec D(\alpha):=\ens{x\in \vec A(T)}{\alpha(x)\geq 0}$ est le \textit{demi-appartement} dirigé par $\alpha$. Le groupe $U_\alpha$ fixe ce demi-appartement, et est simplement transitif sur les appartements le contenant.\\
Si on fixe un tore maximal $T$, et une base de $\phi(T)$, on définit $\iv^+=G.\vec A^+(T)$ et $\iv^-=G.\vec A^-(T)$. Ce sont deux immeubles au sens classique. Lorsque $\phi$ est fini, ils coïncident. Lorsque $\phi$ est infini, leur intersection ne contient que des facettes non sphériques, et contient toujours $\{0\}$. Le découpage $\iv=\iv^+\cup\iv^-$ est indépendant des choix de $T$ et de la base de $\phi(T)$, mais $\iv^+$ et $\iv^-$ sont échangés si on remplace par exemple une base de $\phi(T)$ par son opposée. Ces deux immeubles sont jumelés.\\
 
Pour toute partie $\vec\Omega$ d'un appartement $\vec A$, on note $\cl_{\vec A}(\vec\Omega)$  son \textit{enclos}, il s'agit de l'intersection de tous les demi-appartements de $\vec A$ contenant $\vec \Omega$. Une partie égale à son enclos sera dite \textit{close}. On verra (\ref{props:paraboliques}) que l'enclos d'une partie $\vec\Omega$ est généralement indépendant de l'appartement $\vec A$ la contenant considéré, et on pourra éliminer l'indice $\vec A$ dans la notation $\cl_{\vec A}(\vec\Omega)$.\\
 
 Lorsque $\vec C$ et $\vec D$ sont deux chambres de même signe, on peut définir leur distance (à valeur dans $\N$, pour l'usage qu'on en aura). Lorsque $\vec C$ et $\vec D$ sont de signe opposé, on définit leur codistance. Celle-ci est nulle si et seulement si $\vec C$ et $\vec D$ sont \textit{opposées}, c'est-à-dire si $\vec C=-\vec D$ dans un certain appartement $\vec A$ les contenant. Dans ces conditions, on note $\vec C= \op_{\vec A}(\vec D)$.\\
 
 Lorsque $\vec f$ est une facette sphérique de $\iv$, on peut définir la \textit{projection} sur $\vec f$ (\cite{abra} I.4). On commence par la définir pour les chambres de $\iv$: si $\vec C$ est une chambre de $\iv$, alors sa projection sur $\vec f$, notée $pr_{\vec f}(\vec C)$ est l'unique chambre de $\vec f^*$ qui est à distance minimale (si $\vec C$ et $\vec f$ sont de même signe) ou à codistance maximale (dans le cas contraire) de $\vec C$. Ensuite, si $\vec g$ est une facette quelconque, on pose $pr_{\vec f}(\vec g) = \bigcap_{\vec C\in \vec f^*} \barre{pr_{\vec f}(\vec C)}$. Insistons sur le fait que la projection sur une facette $\vec f$ n'est définie que lorsque $\vec f$ est sphérique. (Lorsque $\vec f$ n'est pas sphérique, elle est en fait quand même définie pour les facettes $\vec g$ de même signe que $\vec f$.)
 
 On prouve que si $\vec f$ et $\vec g$ sont dans un appartement $\vec A$, alors $pr_{\vec f}(\vec g)\subset \vec A$ (voir \cite{abra}, I.4, corollaire 3). On peut alors donner une caractérisation géométrique de $pr_{\vec f}(\vec g)$: 
 \fprop{ Soient $\vec f$ une facette sphérique et $\vec g$ une facette quelconque alors $pr_{\vec f}(\vec g)$ est la plus grande facette de $\vec f^*$ incluse dans $\cl_{\vec A}(\vec f\cup \vec g)$.
 
 De plus, $\barre{pr_{\vec f}(\vec g)}$ est l'intersection de $\vec f^*$ et de tous les murs de $\vec A$ contenant $\vec f\cup \vec g$.}
 
 \demo 
 Pour tout demi-appartement $\vec D$ contenant $\vec f\cup \vec g$, il existe un appartement $\vec B$ tel que $\vec D=\vec A\cap \vec B$. Comme on vient de le rappeler, $pr_{\vec f}(\vec g)\subset \vec A\cap \vec B = \vec D$. On prouve ainsi que $pr_{\vec f}(\vec g) \subset \cl(\vec f\cup \vec g)$.\\
 
 La  première assertion découle alors de la deuxième, montrons celle-ci. Soit $\maj M$ l'ensemble des murs de $\vec A$ contenant $\vec f\cup \vec g$. Soit $\vec M\in \maj M$, soient $\vec C$ et $\vec D$ les deux chambres de $\vec g^*\cap \vec A$ ayant une cloison commune dans $\vec M$. Alors $\vec M$ ne peut séparer $\vec C$ de $pr_{\vec f}(\vec C)$, ni $\vec D$ de $pr_{\vec f}(\vec D)$ sans quoi $\vec C$ serait plus proche (ou à codistance plus grande) de $r_{\vec M}. pr_{\vec f}(\vec C)$ que de $pr_{\vec f}(\vec C)$. Alors $pr_{\vec f}(\vec g)\subset pr_{\vec f}(\vec C)\cap pr_{\vec f}(\vec D) \subset \vec M$.
 
 D'où $pr_{\vec f}(\vec g)\subset \vec f^*\cap \bigcap_{\vec M\in\maj M}\vec M$.\\
 
 Pour l'autre inclusion, soit $\vec M$ un mur contenant $pr_{\vec f}(\vec g)$, montrons que $\vec M\in\maj M$. Il existe deux chambres $\vec C$ et $\vec D$ de $\vec g^*\cap \vec A$ telles que $\vec M$ sépare $pr_{\vec f}(\vec C)$ de $pr_{\vec f}(\vec D)$. Comme précédemment, $\vec M$ ne peux séparer $\vec C$ ni $\vec D$ de leur projections sur $\vec f$ et par conséquent, $\vec M$ sépare $\vec C$ de $\vec D$. Donc $\vec M\supset \barre{\vec C}\cap \barre{\vec D}\supset \vec g$.\cqfd

Pour toute partie $\vec\Omega$ d'un appartement $\vec A(T)$, on note $P(\vec\Omega)$ son fixateur dans $G$, c'est le \ital{sous-groupe parabolique} de $G$ associé à $\vec\Omega$. Rappelons les propriétés essentielles de ces groupes:
\fprops[\label{props:paraboliques}\ ]{ 

 \item Si $\vec f$ et $\vec g$ sont deux facettes d'un appartement $\vec A(T)$, alors $G = P(\vec f).N(T).P(\vec g)$. Lorsque $\vec f$ et $\vec g$ sont de même signe, cette décomposition est dite de Bruhat, sinon de Birkhoff. Ceci entraine que pour toutes facettes $\vec f$ et $\vec g$ de $\iv$, il existe un appartement les contenant.
 
 Si $\vec f=\pm \vec g$, et si $\vec f$ est une chambre, on a une décomposition plus précise: $G=\bigsqcup_{n\in N(T)} U(\vec f)nU(\vec g)$, et pour tout $t\in T$ il y a écriture unique dans la double classe $U(\vec f)tU(\vec g)$.\\

 \item Soit $\vec\Omega$ une partie de $\vec A(T)$ incluse dans $\vec A^+(T)$, incluse dans $\vec A^-(T)$ ou rencontrant $\vec A^+(T)\sph$ et $\vec A^-(T)\sph$.
 Alors $\bigcap_{\vec f\in \vec\Omega} N(T).P(\vec f) = N(T).P(\vec\Omega)$, ce qui signifie qu'entre deux appartements contenant $\vec\Omega$ existe un isomorphisme induit par un élément de $G$ fixant $\vec\Omega$. Autrement dit, le groupe $P(\vec\Omega)$ est transitif sur les appartements contenant $\vec\Omega$. \\
 
 \item Soit $\vec\Omega$ une partie de $\vec A(T)$ incluse dans $\vec A(T)^+$, incluse dans $\vec A(T)^-$, ou rencontrant $\vec A^+(T)\sph$ et $\vec A^-(T)\sph$. Alors $P(\vec\Omega)=P(\cl(\vec\Omega))$.
 
  Au vu du point 2, ceci signifie que l'intersection de deux appartements $\vec A$ et $\vec B$ est une partie close dans chacun des appartements $\vec A$ et $\vec B$ si elle est incluse dans $\vec A^+$, ou dans $\vec A^-$, ou si elle rencontre $\vec A\sph$.
 
 }
 
\demos{

\item Voir \cite{remy} partie 1. On prouve qu'un groupe muni d'une donnée radicielle est également muni d'une "BN-paire raffinée" en 1.5.4, puis qu'un tel groupe vérifie une version plus précise que celle énoncée ici de la décomposition de Bruhat en 1.2.3 et de la décomposition de Birkhoff en 1.2.4.

Déduisons-en la version géométrique: soient $\vec f$ et $\vec g$ deux facettes de $\iv$, soit $\vec A(T)$ un appartement contenant $\vec f$ et $h.\vec A(T)$, avec $h\in G$, un appartement contenant $\vec g$ ($h$ et $\vec A(T)$ existent car $\iv$ est la réunion de ses appartements, et car $G$ permute ces derniers transitivement). On utilise alors la décomposition de Bruhat/Birkhoff dans l'appartement $\vec A(T)$ avec les facettes $h\inv \vec g$ et $\vec f$: il existe $p\in P(h\inv\vec g)$, $n\in N(T)$ et $q\in P(\vec f)$ tels que $h=qnp$. Alors l'appartement $q.\vec A(T)$ contient $\vec f\cup \vec g$.\\

\item La version géométrique de ce résultat est prouvée dans \cite{abramenko-brown}, 6.73 lorsque $\vec\Omega$ est l'adhérence d'une réunion de chambres.
Elle est de plus connue si $\vec \Omega\subset A^+$ ou $\vec \Omega\subset A^-$ (\cite{abramenko-brown} 4.5).

 Étudions le cas général, où $\vec\Omega$ rencontre $\vec A\sph^+$ et $\vec A\sph^-$. Soit $\vec B$ un autre appartement contenant $\vec\Omega$, soient $\vec f$ et $\vec g$ des facettes maximales de $\vec A^+\cap \vec B^+$  et $\vec A^-\cap \vec B^-$, respectivement.
  Soit $\vec C$ une chambre de $\vec A$ contenant $\vec f$ dans son adhérence, et $\vec D$ une chambre de $\vec B$ contenant $\vec g$ dans son adhérence. Il existe un appartement $\vec Z$ contenant $\vec C\cup \vec D$. On sait que $\vec Z^+\cap \vec A^+$ est une partie close (\ref{prop:intersection close}) contenant la chambre $\vec C$, donc c'est l'adhérence d'un ensemble de chambres. Du côté négatif, comme $\vec g$ est une facette sphérique incluse dans $\vec Z^-\cap \vec A^-$, ce dernier ensemble est clos et contient la chambre $pr_{\vec g}(\vec C)$. Il s'agit donc également de l'adhérence d'un ensemble de chambres. Donc $\vec A\cap \vec Z$ est l'adhérence d'un ensemble de chambres, et par \cite{abramenko-brown}, 6.73, il existe $g_1\in G$ tel que $g.\vec A=\vec Z$, et $g$ fixe $\vec C$, donc $\vec f$, et $\vec g$.

De même, il existe $g_2\in G$ tel que $g.\vec Z=\vec B$ et $g_2$ fixe $\vec f\cup \vec g$. Au final, $g_2g_1$ fixe $\vec f\cup \vec g$ et envoie $\vec A$ sur $\vec B$. Par \ref{prop:isom fixant l intersection}, $g_2g_1$ fixe alors $\vec A^+\cap \vec B^+$ et $\vec A^-\cap \vec B^-$. En particulier $g_2g_1$ fixe $\vec\Omega$.

Pour en déduire la version algébrique du résultat, si $g\in \bigcap_{\vec f\in \vec\Omega} N.P(\vec f)$, posons $\vec B=g.\vec A$, et soit $g_2g_1\in P(\vec\Omega)$ comme ci-dessus, alors $g\inv g_2g_1 .\vec A=\vec A$ donc $g\inv g_2g_1\in N$ et $g\in N.P(\vec\Omega)$.\\

\item Encore une fois, ceci est classique lorsque $\vec\Omega\subset \vec A^\pm$ (voir \ref{prop:intersection close}), on peut donc supposer que $\vec\Omega$ contient des points sphériques positifs et négatifs. On peut aussi supposer que $\vec\Omega$ contient $\cl(\vec\Omega\cap\vec A^+)$ et $\cl(\vec\Omega\cap \vec A^-)$, en particulier $\vec\Omega$ est l'adhérence d'une réunion de facettes sphériques.

 Lorsque $\vec\Omega$ est une partie "équilibrée" (c'est-à-dire une réunion finie de facettes sphériques en contenant au moins une positive et une négative), on a $P(\vec\Omega)=T.G(\phi(\vec\Omega))$ (voir \cite{remy}, chapitre 6, ceci sera détaillé un peu au paragraphe suivant). Chaque groupe radiciel $U_\alpha$, $\alpha\in\phi(\vec\Omega)$ fixe un demi-appartement contenant $\vec\Omega$, donc fixe aussi $\cl(\vec\Omega)$.
 
 Si $\vec\Omega$ est l'adhérence d'un ensemble de chambres, alors $\cl(\vec\Omega)$ est la plus petite partie de $\vec A$ fermée et stable par projection sur ses cloisons (\cite{abramenko-brown} 5.193, où une partie convexe est par définition un ensemble de chambres stable par projection sur ses cloisons intérieures). Si $\vec f$ et $\vec g$ sont deux facettes sphériques de signes opposés, alors $\vec f\cup \vec g$ est une partie équilibrée, et le paragraphe précédent montre que $P(\vec f\cup \vec g)$ fixe $\cl(\vec f\cup \vec g)$, qui contient $pr_{\vec f}(\vec g)$. Il est alors évident que $P(\vec\Omega)$ fixe $P(\cl(\vec\Omega))$.\\
 
 Pour traiter le cas général, on reprend la démonstration de \cite{abramenko-brown} 5.193. On commence par le 
 
 \flemme{ Soit $\maj C$ une partie close de $\vec A^+$ ou de $\vec A^-$ contenant un point sphérique. Alors $\maj C$ est l'intersection des demi-appartements contenant $\maj C$ et dont le bord contient $\maj C$ ou une cloison sphérique de $\maj C$. 
 }
 \rema Une cloison de $\maj C$ est par définition une facette de codimension 1 dans $\maj C$.\\
 
 \pv
Supposons par exemple $\maj C\subset \vec A^+$. Soit $\vec f$ une facette de $\vec A^+\setminus \maj C$, montrons qu'il existe un demi-appartement comme dans l'énoncé qui sépare $\vec f$ de $\maj C$. Soit $\vec c$ une chambre de $\maj C$, à distance minimale de $\vec f$ (rappelons qu'une partie close d'un système de Coxeter est un complexe de chambre, d'après \ref{prop:clos implique cx de chb}), $\vec c$ est donc une facette sphérique de $\vec A ^+$.

Supposons dans un premier temps que $pr_{\vec c}(\vec f)\not= \vec c$. Soit $\vec M$ un mur contenant $\vec c$, et donc $\maj C$, mais pas $pr_{\vec c}(\vec f)$, il ne contient alors pas $\vec f$. Le demi-appartement délimité par $\vec M$ et ne contenant pas $\vec f$ est comme requis dans l'énoncé, contient $\maj C$ et pas $\vec f$.

Supposons à présent que $pr_{\vec c}(\vec f)=\vec c$. Il existe alors une unique cloison $\vec m$ de $\vec c$ à distance minimale de $\vec f$. Remarquons que puisque $\maj C$ est fermé, $\vec f\not\subset \barre{\vec c}$, donc $pr_{\vec m}(\vec f)$ est une facette de même dimension que $\vec c$, différente de $\vec c$, et incluse dans $\vect_{\vec A}(\vec c)$. Sachant que $\vec c$ est sphérique donc dans l'intérieur du cône de Tits, tout point de $\vec m$ est dans un segment ouvert reliant un point de $pr_{\vec m}(\vec f)$ et un point de $\vec c$, donc $\vec m$ est dans l'intérieur du cône de Tits, donc $\vec m$ est une cloison sphérique de $\maj C$. Soit $\vec M$ un mur contenant $\vec m$ et pas $\vec c$, soit $\vec D$ le demi-appartement délimité par $\vec M$ contenant $\vec c$, alors $\vec D$ contient $\maj C$, sans quoi $pr_{\vec m}(\vec f)$ serait dans $\maj C$, contredisant la définition de $\vec c$. De plus, $\vec f\not\subset \vec D$.\cqfd

On montre alors la version géométrique du point 3. Soit $\vec B$ un appartement tel que $\vec A\cap \vec B$ contient des points sphériques positifs et négatifs. Soit $\epsilon$ un signe. La partie $\maj C^\epsilon:=\vec A^\epsilon\cap\vec B^\epsilon$ est close dans $\vec A^\epsilon$, c'est l'intersection des demi-appartements la contenant et dont le bord contient $\maj C^\epsilon$ ou une cloison sphérique de $\maj C^\epsilon$. Notons $\D^\epsilon$ cet ensemble de demi-appartements. Soit $\vec D\in\D^\epsilon$, et $\vec m$ une chambre ou une cloison sphérique de $\maj C$ incluse dans $\partial\vec D$. Comme $\vec A\cap \vec B$ est stable par projection, et que la projection sur $\vec m$ est bien définie puisque $\vec m$ est sphérique, $\maj C^\epsilon$ contient tous les $pr_{\vec m}(d)$ pour $d$ une facette de $\vec A^{-\epsilon}\cap \vec B^{-\epsilon}$. Et ceci entraine que $\vec D$ contient $\vec A^{-\epsilon}\cap \vec B^{-\epsilon}$.

On prouve ainsi que $\cl(\vec A\cap \vec B)$ est égal à l'intersection des demi-appartements de $\D^+\cup \D^-$, puis que $ \cl(\vec A\cap \vec B)=\maj C^+\cup \maj C^- \subset \vec A\cap \vec B$.\\

Pour en déduire la version algébrique du résultat, soit $g\in P(\vec \Omega)$, posons $\vec B=g.\vec A$. Par le point 2, il existe $h\in P(\vec A\cap \vec B)$ tel que $h.\vec A=\vec B$, donc $g\inv h\in N(\vec\Omega)$. Or $N(\vec \Omega)=N(\cl(\vec \Omega))$, et $\cl(\vec\Omega)\subset \vec A\cap \vec B$ d'où $h\in P(\cl(\vec\Omega))$. Au final, on obtient bien $g\in P(\cl(\vec \Omega))$.

}

\subsubsection{Décomposition de Lévi}
On se réfère ici à \cite{remy} chapitre 6. On fixe dans ce paragraphe un tore $T$, et on note $\vec A=\vec A(T)$ et $\phi=\phi(T)$.\\

 Soit $\vec\Omega$ une partie de $\vec A$. On note:\begin{itemize}
\item $\phi^u(\vec\Omega)= \ens{\alpha\in\phi} {\alpha(\vec\Omega)>0}$
\item $\phi^m(\vec\Omega)= \ens{\alpha\in\phi} {\alpha(\vec\Omega)=0}$
\item $\phi(\vec \Omega)=\phi^u(\vec\Omega)\sqcup \phi^m(\vec\Omega)=\ens{\alpha\in\phi} {\alpha(\vec\Omega)\geq 0} $.
\end{itemize}

L'ensemble de racines $\phi^m(\vec\Omega)$ est un sous-système de racines de $\phi$. Lorsque $\vec\Omega$ contient un point sphérique, il est fini.\\

On définit ensuite les sous-groupes de $P(\vec\Omega)$ suivants:\\\begin{itemize}

 \item $M_{\vec A}(\vec \Omega)$, le \ital{facteur de Lévi de $P(\vec \Omega)$ par rapport $\vec A$:}\\ Il est définit par $M_{\vec A}(\vec \Omega):= \fix_G(\vect_{\vec A}(\vec\Omega)) = \fix_G(\vec\Omega\cup \op_{\vec A}(\vec\Omega))$. Il est d'après \ref{props:paraboliques} transitif sur les appartements contenant $\vect_{\vec A}(\vec\Omega)$. Si $\vec\Omega$ est une facette ou si elle contient une facette sphérique, alors $M_{\vec A}(\vec \Omega)=\langle T,\ens{U_\alpha}{\alpha\in \phi^m(\vec A)(\vec\Omega)} \rangle = T.G(\phi^m(\vec \Omega))$.
Enfin, le couple $(M(\vec\Omega),(U_\alpha)_{\alpha\in\phi^m(\vec\Omega)})$ est une donnée radicielle de système de racines $\phi^m(\vec\Omega)$, voir \cite{remy} 6.2.3. En particulier, lorsque $\vec\Omega$ contient un point sphérique, $\phi^m(\vec\Omega)$ est un système de racines fini, et $M(\vec\Omega)$ est muni d'une donnée radicielle de type fini.\\

 \item  $U(\vec\Omega)$, le \ital{facteur unipotent de $P(\vec \Omega)$:}\\ C'est le sous-groupe distingué de $P(\vec\Omega)$ engendré par $G(\phi^u(\vec\Omega))$ $=\engens{U_\alpha}{\alpha\in \phi^u(\vec A)(\vec\Omega)}$. Il est donc indépendant de l'appartement $\vec A$ contenant $\vec\Omega$ considéré.  
 %Dès que $\phi^u(\vec\Omega)$ est une partie nilpotente de racines,
 Si $\vec \Omega$ contient un point sphérique positif et un point sphérique négatif, il admet la décomposition avec écriture unique: $U(\vec\Omega)= \prod_{\alpha\in\phi^u(\vec A)(\vec\Omega)} U_\alpha$, quel que soit l'ordre des facteurs. En particulier, on a alors $U(\vec\Omega)=G(\phi^u(\vec\Omega))$. %ce cas se présente notamment dès que $\vec\Omega$ intersecte une facette positive sphérique et une facette négative sphérique.\\
  Lorsque $\vec \Omega$ est une chambre, $U(\vec \Omega)= G(\phi^u(\vec \Omega))$, et lorsque $\vec\Omega$ est une facette, $U(\vec \Omega)$ est l'intersection des $U(\vec C)$ pour $\vec C$ les chambres de $\vec A$ contenant $\vec\Omega$ dans leur adhérence. On prouve enfin que si $\vec\Omega$ est une facette sphérique, $U(\vec\Omega)=U(\vec C)\cap U(\vec D)$ dès que $\vec C$ et $\vec D$ sont deux chambres opposées dans $\vec\Omega^*\cap\vec A$.\\

\end{itemize}

Une partie équilibrée dans $\iv$ est une partie d'appartement qui contient des points positifs et négatifs et qui est recouverte par un nombre fini de facettes sphériques (éventuellement fermées).

Dans le cas où $\vec\Omega$ est soit une facette soit une partie équilibrée de $\vec A$,
  $P(\vec\Omega)$ admet une décomposition de Lévi (\cite{remy}, 6.2.2 et 6.4.1):
\[P(\vec \Omega)=M_{\vec A}(\vec \Omega)\ltimes U(\vec \Omega)\]

 Une extension vectorielle de $\vec\Omega$ est une partie de $\iv$ de la forme $\vect_{\vec B}(\vec\Omega)$ pour un appartement $\vec B$ contenant $\vec\Omega$. La décomposition de Lévi peut aussi s'exprimer en disant que $U(\vec\Omega)$ est simplement transitif sur les extensions vectorielles de $\vec\Omega$.\\

\subsection{Valuation d'une donnée radicielle}
\label{soussection:valuation}

A l'exemple de \cite{bruhat-tits}, on ajoute maintenant une structure supplémentaire à notre donnée radicielle qui permet de rendre compte, dans le cas d'un groupe sur un corps local, de la valuation du corps.

\fdef{\label{def:valuation}Soit $\phi$ un système de racines. Soit $(G,(U_\alpha)_{\alpha\in\phi})$ une donnée radicielle et pour tout $\alpha\in\phi$ soit $\phii_\alpha$ une fonction de $U_\alpha$ dans $\R\cup\{\infini\}$. Pour tout $\lambda\in\R$ on note $U_{\alpha,\lambda} = \phii_\alpha\inv([\lambda,\infini])$.

On dit que la famille $(\phii_\alpha)_{\alpha\in\phi}$ est une valuation de la donnée radicielle $(G,(U_\alpha)_{\alpha\in\phi})$, ou que $(G,(U_\alpha,\phii_\alpha)_{\alpha\in\phi})$ est une donnée radicielle valuée si:\\
\liste {
\item (V0): $\forall \alpha\in\phi$, $\phii_\alpha(U_\alpha)$ a au moins trois éléments.\\% et contient $0$.\\

\item (V1): Pour tout $\alpha\in\phi$ et $\lambda \in\R$, $U_{\alpha,\lambda}$ est un sous-groupe de $U_\alpha$, et $U_{\alpha,\infini} = \{e\}$.\\

\item (V2.1): %Pour tout $\alpha\in\phi$, pour tout $u\in U_{\alpha}\setminus\{e\}$, la fonction $ \fonc {U_{-\alpha}\setminus \{e\}} {\Lambda} {v} {\phii_{-\alpha} (v) - \phii_\alpha(n(u) v n(u)\inv)} $ est constante.
	Pour tout $\alpha,\beta\in\phi$, $u\in U_\alpha\setminus\{e\}$, $v\in U_\beta\setminus\{e\}$,
	\[ \phii_{r_\alpha.\beta}\left(n(u).v.n(u)\inv\right) = \phii_\beta(v)-\langle \alpha,\beta\rangle  \phii_\alpha(u)  \]

\item (V2.2): Pour tout $\alpha\in\phi$, pour tout $t\in T$, $\fonc {U_{\alpha}\setminus \{e\}} {\Lambda} {v} {\phii_{\alpha} (v) - \phii_\alpha(t v t\inv)}$ est constante.\\

\item (V3): Pour toute paire prénilpotente de racines $\{\alpha,\beta\}$, pour tous $\lambda,\mu\in\R$:
\[ [U_{\alpha,\lambda},U_{\beta,\mu}] \subset \engens {U_{p\alpha+q\beta, p\lambda+q\mu} }{p,q\in \N^*\et p\alpha+q\beta\in\phi} \]

\item (V4): Si $\alpha\in\phi$ et $2\alpha\in\phi$ alors $\phii_{2\alpha}$ est la restriction de $\phii_\alpha$ à $U_{2\alpha}$.\\

%\item (V5): Pour tout $\alpha\in\phi$, pour tout $u\in U_\alpha\setminus\{e\}$, soient $u',u''\in U_{-\alpha}$ tels que $n(u)=u'uu''$, alors $-\phii_\alpha(u) = \phii_{-\alpha}(u') = \phii_{-\alpha}(u'')$.
Lorsque  $(G,(U_\alpha,\phii_\alpha)_{\alpha\in\phi})$ est une donnée radicielle valuée, on garde la notation $U_{\alpha,\lambda}$ qu'on vient d'introduire.\\

Si de plus pour tout $\alpha\in\phi_{red}$, $0\in\phii_\alpha(U_\alpha)$, on dit que $\phii$ est une valuation spéciale.\\
}
}

\remas{
\item Soit $\alpha\in\phi$, et $u\in U_\alpha\setminus\{e\}$. Par (V1), on voit que $\phii_\alpha(u)=\phii_\alpha(u\inv)$, et par (V2.1) on obtient $\phii_{-\alpha}(n(u).u.n(u)\inv) = -\phii_\alpha(u)$.

\item Avec (V0) et (V2.1), on voit qu'il existe un sous groupe de $\R$ non trivial $\Lambda$ tel que $\phi_\alpha(U_\alpha)+\Lambda= \phi_\alpha(U_\alpha)$.\\}

L'ensemble des valuation de $\D$ est muni d'une action de $\vec V$: pour tout $\vec v\in\vec V$ et $\phii$ une valuation de $\D$, on définit $\phii+\vec v$ par $\forall \alpha\in\phi$, $u\in U_\alpha$, $(\phii+\vec v)_\alpha(u) = \phii_\alpha(u)+\alpha(\vec v)$. Il est immédiat de vérifier que $\phii+\vec v$ est encore une valuation de $\D$. Deux valuations $\phii$ et $\phii'$ telles qu'il existe $\vec v\in\vec V$ tel que $\phii=\phii'+\vec v$ sont dites équipollentes.
Soit $\phii$ une valuation quelconque. Le fait que $\phi$ engendre $\vec V^*$ entraine que l'action de $\vec V$ sur l'ensemble $\phii+\vec V$ des valuations équipollentes à $\phii$ est simplement transitive, donc $\phii+\vec V$ est un espace affine sous $\vec V$.
 
   Si $\Pi$ est une base de $\phi$, il s'agit aussi d'une base de $\vec V$ et on peut donc trouver $\vec v\in\vec V$ tel que pour tout $\alpha\in\phi$, $0\in\phii_\alpha(U_\alpha)+\alpha(\vec v)$. En utilisant (V2.1), on vérifie alors que cette relation reste vraie pour tout $\alpha\in\phi_{red}$. Ainsi, $\phii+\vec v$ est une valuation spéciale équipollente à $\phii$. On peut donc toujours se ramener à une valuation spéciale, à équipollence près.
   
   Notons enfin que pour tout $g\in G$, la famille de fonction $g.\phii$ définie par $\forall \alpha\in g.\phi$, $\forall u\in U_\alpha$, $(g.\phii)_\alpha(u) = \phii_{g\inv\alpha}(g\inv u g)$ est une valuation de la donnée radicielle $g.\D$.\\

Dans la définition de la valuation d'une donnée radicielle finie de \cite{bruhat-tits}, ou même dans la définition de la valuation d'une donnée radicielle quelconque de \cite{microaffine}, l'axiome (V2) est beaucoup plus faible que celui présenté ici. Par contre on ajoute un cinquième axiome à savoir:\\

(V5): Pour tout $\alpha\in\phi$, pour tout $u\in U_\alpha\setminus\{e\}$, pour tous $u',u''\in U_{-\alpha}$ tels que $n(u)=u'uu''$, alors $-\phii_\alpha(u) = \phii_{-\alpha}(u') = \phii_{-\alpha}(u'')$.\\

On prouvera (en \ref{prop:V5}) que pour une donnée radicielle finie, la définition présente équivaut à celle de \cite{bruhat-tits}, autrement dit que d'une part l'axiome (V2) présenté ici est conséquence de sa version faible et de (V0), (V1), (V3), (V4), (V5), et que d'autre part (V5) découle des (V0)...(V4) ci-dessus. Dans le cas où $\phi$ est infini, l'auteur n'a pas pu se passer de la version forte de (V2).\\

\fprop{ \label{prop:valuation implique cent}Soit $\D=(G,(U_\alpha)_{\alpha\ni\phi})$ une donnée radicielle admettant une valuation $\phii=(\phii_\alpha)_{\alpha\in\phi}$. Alors $\D$ vérifie la condition (CENT).
}
\demo
Soit $\alpha\in\phi$, $u\in U_\alpha\setminus\{e\}$. Nous devons trouver $t\in T$ tel que $utu\inv\not=t$. Par (V0), il existe $v\in U_\alpha\setminus\{e\}$ tel que $\phii_\alpha(u)\not=\phii_\alpha(v)$. Posons $t=n(u)n(v)$, il est clair par (DR4) que $t$ normalise chaque $U_\beta$, donc $t\in T$. Montrons que $t\inv u t\not=u$. On calcule $\phii_\alpha(t\inv u t)$ en utilisant deux fois (V2.1):
\calc{
 \phii_\alpha(n(v)\inv n(u)\inv .u.n(u)n(v)) &=& \phii_{-\alpha}(n(u)\inv un(u)) - \eng{-\alpha,\alpha}\phii_\alpha(v)\\
&=& \phii_\alpha(u)-\eng{\alpha,\alpha}\phii_\alpha(u) +2\phii_\alpha(v) = 2\phii_\alpha(v)-\phii_\alpha(u)\\
&\not=& \phii_\alpha(u)
}
Donc $t\inv u t\not=u$.\cqfd

\fprops{\label{prop:existence de valuation} 
\item Soit $G$ un groupe de Kac-Moody déployé, soit $\K$ un corps muni d'une valuation non triviale $\varpi:\K\fleche \R\cup\{\infini\}$. Alors il existe une valuation $(\phii_\alpha)_{\alpha\in \phi}$ de la donnée radicielle $(G(\K),(U_\alpha(\K))_{\alpha\in\phi})$, elle est définie par $\phii_\alpha(u_\alpha(k))  = \varpi(k)$, où $(u_\alpha)_{\alpha\in\phi}$ est un système de Chevalley pour $\phi$.

\item Soit $\D=(G,(U_\alpha)_{\alpha\in\phi})$ une donnée radicielle  avec $\phi$ un système de racines fini. Soit $\phii=(\phii_\alpha)_{\alpha\in\phi}$ une valuation de $\D$ au sens de \cite{bruhat-tits} 6.2.1. Alors $\phii$ est aussi une valuation de $\D$ au sens présent.
}

\demo 

Pour le cas Kac-Moody, on trouvera dans \cite{microaffine} 2.2 des références qui prouvent toutes les conditions (Vx) sauf (V2). Dans le cas d'un système de racines fini, ces conditions sont les même dans \cite{bruhat-tits} et ici.\\

Il ne reste qu'à étudier (V2). Cette condition découle rapidement de l'existence d'un espace affine muni d'une action convenable de $N$ (un appartement en fait). Expliquons comment sous forme d'un lemme:
\flemme{\label{lemme:pour V2}
Soit $\phi\subset \vec V^*$ un système de racines, et $\D=(G,(U_\alpha)_{\alpha\in\phi})$ une donnée radicielle. Soit $\phii=(\phii_\alpha)_{\alpha\in\phi}$ une famille de fonctions avec $\forall\alpha\in\phi$, $\phii_\alpha:U_\alpha\fleche \R\cup\{\infini\}$ et $\phii_\alpha\inv\{\infini\}=\{e\}$.

On suppose qu'il existe un espace affine $A$ sous $\vec V$, un point $o\in A$ et une action de $N$ sur $A$ par automorphismes affines telle que pour tout $\alpha\in\phi$ et $u\in U_\alpha\setminus\{e\}$, $n(u)$ agit comme la réflexion d'hyperplan $M(\alpha,\phii_\alpha(u)):=\ens{x\in A}{\alpha(\vec {ox})+\phii_\alpha(u)=0}$ dont la direction est la réflexion $r_{\alpha}\in W(\phi)$. Alors la famille $\phii$ vérifie (V2.1). Si de plus $T$ agit sur $A$ par translation, alors $\phii$ vérifie (V2.2).
}
\pv On identifie les $\alpha\in\phi$ à des formes affines sur $A$ par $\alpha(x)=\alpha(\vec{ox})$.
  La valeur $\phii_\alpha(u)$ est caractérisée par l'égalité $\alpha\left( \fix_{\vec V}(n(u)) \right) = \{-\phii_\alpha(u) \}$.

Soient $\alpha,\beta\in\phi$ et $(u,v)\in U_\alpha\setminus\{e\}\times U_\beta\setminus\{e\}$. Pour calculer $\phii_{r_\alpha.\beta}\left(n(u).v.n(u)\inv\right)$, il faut donc étudier la réflexion $n\left( n(u) .v. n(u)\inv \right)$. Par la définition (voir (DR4)), il est clair que $n\left( n(u) .v .n(u)\inv \right) = n(u).n(v).n(u)\inv$, et l'ensemble de ses points fixes est $n(u).\fix_{A}(v)$. Il existe un point $a\in \fix(n(u))$ et un réel $x$ tel que $a+x.\alpha^\vee \in\fix(n(v))$ (si $\alpha$ et $\beta$ sont colinéaires, pour tout $a\in\fix(n(u))$ il existe un tel $x$, sinon il existe $a\in \fix(n(u))\cap \fix(n(v))$, alors $x=0$ convient). Alors $n(u).(a+x.\alpha^\vee) = a-x.\alpha^\vee$, on a donc:
 \calc{
-\phii_\alpha(u) &=& \alpha(a)\; ,\\
-\phii_\beta(v)&=&\beta(a+x.\alpha^\vee)=\beta(a)+x\langle \alpha,\beta\rangle\; ,\\
\text{et: }\: -\phii_{r_\alpha.\beta}(n(u).v.n(u)\inv) &=& 
 (r_\alpha.\beta)(a-x.\alpha^\vee) = (\beta - \langle \alpha,\beta\rangle \alpha)(a-x.\alpha^\vee) \\
 &=& \beta(a) - \langle \alpha,\beta\rangle \alpha(a) - x\langle \alpha,\beta\rangle+2x\langle \alpha,\beta\rangle\\
 & =& \beta(a) - \langle \alpha,\beta\rangle \alpha(a) + x\langle \alpha,\beta\rangle\\
 &=&- \phii_\beta(v) + \langle \alpha,\beta\rangle \phii_\alpha(u)\point
 }
 Ceci prouve (V2.1). On procède de même pour (V2.2): soit $\alpha\in\phi$, $u\in U_\alpha\setminus\{e\}$ et $t\in T$. On vérifie facilement que $n(tut\inv)=tn(u)t\inv$, c'est donc une réflexion d'hyperplan $M(\alpha,\phii_\alpha(u))+\vec v_t$, si $\vec v_t$ désigne le vecteur de la translation induite par $t$ sur $A$. Alors $\phii_\alpha(u)-\phii_\alpha(tut\inv)=\alpha(\vec v_t)$ ne dépend pas de $u$.\cqfd
 
 Pour conclure la preuve de la proposition, dans les deux cas il existe un espace affine $A$ muni d'une action de $N$ comme dans le lemme, voir \cite{microaffine} partie 2 et \cite{bruhat-tits} 6.2.10.\cqfd

%Comme on l'a déjà dit dans le cadre des données radicielles, on préfère s'affranchir du choix d'un tore maximal de référence et considérer que $G$ est muni d'une classe d'équivalence de données radicielles valuées. Soit $T$ le tore maximal correspondant à une donnée radicielle $\D$, et soit $\phii$ est une valuation de $\D$. Si $T'=gTg\inv$ est un autre tore maximal pour $\D$, on définit une valuation $g.\phii$ de la donnée radicielle de la sorte: on pose pour tout $\alpha\in\phi$, $g.\phii_\alpha=\phii_{g\alpha}:\fonc{U_{g\alpha}\setminus\{e\}} {\Lambda} {u} { \phii_\alpha(g\inv ug)}$, et bien sûr toujours $g.\phii_\alpha(e)=\infini$. Si $g'$ est un autre élément de $G$ tel que $T'=g'Tg'{}\inv$, on vérifie immédiatement que la définition $\phii_{g\alpha}=g.\phii_\alpha$ est cohérente avec l'identification $g\phi=g'\phi$.  Il est également immédiat que $g.\phii:=(\phii_{\alpha})_{\alpha\in g.\phi}$ est une valuation de $g.\D$.

 %Ainsi $G$ est muni pour tout tore maximal $T'$ d'une donnée radicielle valuée, et on oubliera le choix particulier de la donnée $\D$ correspondant au tore $T$.
 
 %Remarquons pour finir que ces définitions mènent aux égalités $g U_{\alpha,k}g\inv = U_{g\alpha, k}$, pour tous $g\in G$,$T$ tore maximal, $\alpha\in\phi(T)$ et $k\in\Lambda\cup\{\infini\}$.

\section{Construction générale}\label{section:construction generale}

On fixe pour toute cette partie une donnée radicielle valuée génératrice $\D=(G,(U_\alpha,\phii_\alpha)_{\alpha\in\phi})$, avec $\phii$ spéciale. On rappelle que  la condition (CENT) de \cite{remy} 1.2.5 est alors vérifiée d'après \ref{prop:valuation implique cent}, donc le normalisateur $N(T)$ d'un tore maximal $T$ est le groupe engendré par $T$ et les $n(u)$, $u\in U_\alpha$, $\alpha\in \phi(T)$.

 On construit un objet immobilier $\I$ de manière tout à fait similaire à \cite{bruhat-tits}: on commence par définir pour tout tore maximal $T$ un appartement $A(T)$ muni d'une action de $N(T)$; on définit ensuite les sous-groupes, appelés "parahoriques", qui seront les fixateurs dans $G$ des points de $A$; et on définit finalement $\I$ comme quotient de $G\times A(T)$ par la relation imposant que les sous-groupes parahoriques soient effectivement les fixateurs des points de $A(T)$.

\subsection{L'appartement}
\label{soussection:appartement general}
On fixe un tore maximal $T$ dans $G$. On définit dans ce numéro l'objet $A(T)$ qui sera l'appartement relatif à $T$. Il s'agit d'une réunion disjointe d'espaces affines sous $\vec V(T)$ et certains de ses sous-espace vectoriels. Dans toute cette partie, on notera $\phi=\phi(T)$, $\vec V=\vec V(T)$, $\vec A=\vec A(T)$ et $N=N(T)$.

\label{soussection:appart}
\subsubsection{Façades d'appartement}
\label{sous-sous-section:facades}

Soit $Y(T)$ un espace affine sous $\vec V$. Un cône dans $Y(T)$ est une partie de $Y(T)$ de la forme $f=x+\vec{f}$, où $\vec f$ est un cône de $\vec V(T)$. Le cône vectoriel $\vec f$ est uniquement déterminé, c'est la direction de $x+\vec f$. Deux cônes de même direction sont dits parallèles, on note $g\parallele f$. Lorsque $g$ est parallèle et inclus dans $f$, on dit que c'est un sous-cône parallèle, abrégé en "scp".\\

 On définit une relation d'équivalence $\sim$ (ou $\sim_T$ lorsqu'il faut préciser) sur l'ensemble des cônes convexes de $Y(T)$. Soient $f=x+\vec f$ et $g=y+\vec g$, on pose:\\
\[ f\sim g\ssi (f\parallele g \et f\cap g\not=\vide) \ssi (f\cap g \text{ contient un scp de $f$ et de $g$ })
\ssi (f\parallele g \et \vec{xy}\in \vect(\vec f)) \]

Pour tout cône convexe $\vec f$ on note $Y(T)_{\vec f}$ l'ensemble des cônes dirigés par $\vec f$ quotienté par $\sim$. C'est la \textit{façade} de $Y(T)$ de \textit{direction} $\vec f$. C'est un espace affine isomorphe à $Y(T)/\vect(\vec f)$, et l'action de $\vec V(T)$ sur $A$ passe au quotient sur $Y(T)_{\vec f}$. \\

Soit $\mathcal F(\vec A(T))$ l'ensemble des facettes de $\vec A(T)$, ce sont en particulier des cônes convexes de $\vec V(T)$. L'\textit{appartement (bordé) associé à $T$} est alors:
\[ A(T) := \bigcup_{\vec f\in \mathcal F(\vec A(T))} Y(T)_{\vec f} \]

On notera $A(T)_{\vec f}$ pour désigner la façade $Y(T)_{\vec f}$, et les façades de $Y(T)$ seront appelées les façades de $A(T)$.\\
 Les façades ainsi construites seront appelées façades d'appartement, pour les distinguer des façades d'immeubles dont la définition est à venir. Lorsque $a$ est un point de $A(T)$, on notera $\vec f_a$ la direction de la façade le contenant. L'espace $\vec V$ agit sur l'appartement $A(T)$ par translation et les orbites sont ses façades.\\

On note $A(T)_{sph}$ la réunion des façades sphériques de $A(T)$, c'est-à-dire des façades de type une facette vectorielle sphérique. On note aussi $A(T)^+$, $A(T)^-$ la réunion des façades positives et négatives, $A(T)^+_{sph}$ et $A(T)^-_{sph}$ la réunion des façades sphériques positives et négatives. Dans la terminologie de \cite{microaffine}, $A(T)^+_{sph}$ et $A(T)^-_{sph}$ sont les réalisations de Satake de deux appartements microaffines.\\

On définit une topologie sur $A(T)$, telle que les voisinages d'un point $[x+\vec f]$ sont les: \[\maj V(U,\vec f) =\ens{a\in A(T)}{\text{un représentant de $a$ est inclus dans $U+\vec f$ }}\: ,\] pour tous les voisinages $U$ de $x$ dans l'espace affine $Y(T)$.
Cette topologie induit la topologie classique d'un espace affine de dimension finie sur chaque façade et l'adhérence d'une façade $A_{\vec f}$ est l'union des façades $A_{\vec g}$ pour $\vec f\subset \barre{\vec g}$. Cette topologie est séparée, et si $\vec f$ est une facette sphérique, alors $\barre{A_{\vec f}}=\bigcup_{\vec g\tq \vec f\subset \barre{\vec g} }A_{\vec g}$ est compact (\ref{prop:compactification de lappart}).\\

Comme les bases de $\phi$ sont des bases de $\vec V^*$, la plus petite facette de $\vec A(T)$ est $\{0\}$. La façade $A(T)_{\{0\}}=Y(T)$ est appelée \textit{la façade principale} de $A(T)$, c'est l'intérieur de $A(T)$. On la notera donc $\mathring A(T)$, et on pourra oublier la notation $Y(T)$.\\

Si $\vec f$ et $\vec g$ sont deux facettes de $\vec A$ telles que $\vec f\subset \vect(\vec g)$, l'application $pr_{\vec g}:\fonc{A_{\vec f}}{A_{\vec g}}{[a+\vec f]}{[a+\vec g]}$ est bien définie, c'est la \textit{projection} sur la façade $A_{\vec g}$.\\

Si $\Omega\subset A=A(T)$ est dans une façade $A_{\vec f}$, et si $\vec E$ est un sous-espace vectoriel de $\vec Y(T)$ contenant $\vec f$, alors on définit $\eng{\Omega,\vec E}_A=\barre{pr_{\vec f}\inv(\Omega) +\vec E}$, le \textit{sous-appartement engendré par} $\Omega$ et $\vec E$..

Si $\vec E$ est juste une partie de $\vec A$, la notation $\eng{\Omega,\vec E}$ désignera $\eng{\Omega,\vect_{\vec A}(\vec E)}$.\\

\rema Quels que soient $\Omega$ et $\vec E$, si $\Omega\not=\vide$, alors $\eng{\Omega,\vec E}_A$ coupe toujours la façade principale.\\

Jusqu'à la fin de \ref{soussection:appartement general}, on notera $A=A(T)$.

\subsubsection{Murs et demi-appartements}
 On fixe une origine $o\in \inte A$, et on identifie les formes linéaires sur $\vec V = \vec{\inte A}$ à des formes affines sur $\inte A$ qui s'annulent en $o$, autrement dit on pose pour tout $\alpha\in\phi$ et $a\in A$, $\alpha(a):=\alpha(\vec{oa})$. Soit $\alpha\in\phi$ une racine. Si $\vec f\in\mathcal F$ est une facette telle que $\alpha(\vec f)=0$, alors $\alpha$ définit encore une forme affine sur $A_{\vec f}$. Si $\alpha(\vec f)>0$, on dit que $\alpha$ prend la valeur $\infini$ sur $A_{\vec f}$. Enfin si $\alpha(\vec f)<0$, on dit que $\alpha$ prend la valeur $-\infini$ sur $A_{\vec f}$. De la sorte, $\alpha$ définit une fonction sur $A$, à valeurs dans $\R\cup\{\pm\infini\}$.\\
 
 Ces définitions permettent une caractérisation pratique de la topologie de $A$:
 \flemme{\label{lemme:topo de A par les racines} La topologie de $A$ est engendrée par les demi-espaces ouverts $\ens{x\in A}{\alpha(x)>a}$ pour $a\in\R$ et $\alpha\in\phi$. Autrement dit, une suite $x_n$ tend vers une limite $x$ si et seulement si $\forall \alpha\in\phi$, $\alpha(x_n)\fleche \alpha(x)$.
 }
 \cqfd

Pour tout $\alpha\in\phi$ et $\lambda\in\R$, on pose $M(\alpha,\lambda)=\ens{x\in A}{\alpha(x)+\lambda =0}$ et $D(\alpha,\lambda)=\ens{x\in A}{\alpha(x)+\lambda\geq 0}$. On notera également $D(\alpha,\infini)=A$. L'ensemble des $M(\alpha,\lambda)$ ainsi obtenus pour $\alpha\in\phi$ et $\lambda\in\phii_\alpha(U_\alpha\setminus\{e\})$ est l'ensemble des \textit{murs} de $A$; l'ensemble des $D(\alpha,\lambda)$ correspondants est l'ensemble des \textit{demi-appartements} de $A$.\\
Si $M=M(\alpha,\lambda)$ est un mur de $A(T)$, on notera $\vec M=\ker(\alpha)\subset \vec V$ la \textit{direction} de $M$, c'est un mur de $\vec A$ (ou plutôt sa trace sur $\vec A$ est un mur de $\vec A$). Lorsqu'une intersection de murs est réduite à un seul point, ce point est appelé un \textit{sommet}. Lorsqu'un sommet est inclus dans un mur de chaque direction possible, c'est un \textit{sommet spécial}. Par exemple $o$ est un sommet spécial (ceci est en fait équivalent à la condition "$\phii$ est spéciale").

Notons que si $M$ est un mur contenant $\Omega$ dont la direction contient $\vec E$, alors $\eng{\Omega,\vec E}_A\subset M$.\\

Un isomorphisme affine $\psi$ entre les façades principales de deux appartements $A(T)$ et $A(T')$ dont la partie vectorielle préserve $\mathcal F$ induit une bijection, encore notée $\psi$, entre $A(T)$ et $A(T')$. Si cette bijection préserve l'ensemble des murs, on dit que $\phi$ est un \textit{isomorphisme d'appartements}. Remarquons qu'un isomorphisme d'appartements ainsi défini ne préserve pas forcément les types des facettes de $\vec A(T)$ ni même leur signe.\\

 Pour tout mur $M=M(\alpha,\lambda)$, $r_M$ désigne la réflexion de direction $r_\alpha$ qui fixe $M\cap \inte A$. Elle induit un automorphisme involutif de $A(T)$, qu'on appelle la \textit{réflexion selon} $M$.\\

 Soit $M$ un mur de $A(T)$ et $\vec f\in\mathcal F$. Alors $M\cap A(T)_{\vec f}$  est soit vide, soit un hyperplan de $A(T)_{\vec f}$. Ces hyperplans seront appelés les murs de $A(T)_{\vec f}$.\\

%Si $B$ est l'adhérence d'un sous-espace affine de $\inte A$, muni de la restriction de la forme bilinéaire de $A$ et d'un ensemble $\maj M_B$ d'hyperplans égaux à des traces sur $B$ de murs de $A$, et si pour chaque $M\in \maj M_B$, la réflexion $r^B_M$ orthogonale et d'hyperplan $M$ stabilise $\maj M_B$, on dira que $B$ est un \textit{sous-appartement} de $A$.\\

On notera pour toute partie ou filtre $\Omega$ de $A$, et pour tout $\alpha\in\phi$, $U_\alpha(\Omega)=\ens{u\in U_\alpha}{\Omega\subset \D(\alpha,\phii_\alpha(u))}$. Par exemple, $U_\alpha(\vide)=U_\alpha$ et $U_\alpha(A)=\{e\}$. Pour toute partie $\psi$ de $\phi$, on notera aussi $G(\psi,\Omega)=\engens{U_\alpha(\Omega)}{\alpha\in\psi}$. Enfin, $G(\Omega)$ désignera $G(\phi,\Omega)$.\\\

\rema Dans \cite{bruhat-tits}, $\inte A$ est par définition l'espace affine des valuations équipollentes à $\phii$. Il est isomorphe à celui défini plus haut via $o+\vec v\mapsto \phii+\vec v$. En particulier, la valuation $\phii$ est identifiée au point $o$.\\

	\subsubsection{Parties closes}
\label{soussoussection:parties closes}

\begin{defin} Une partie close de $A(T)$ est une intersection finie de demi-appartements. L'enclos d'une partie ou d'un filtre $E$ de $A(T)$ est le filtre noté $\cl(E)$ engendré par les parties closes de $A(T)$ contenant $E$.\end{defin}

\remas{
\item Avec cette définition, $\cl(\vide) =\vide$.
\item Cette définition de partie close est plus restrictive que celle de \cite{hovels}, elle conduit donc à des enclos plus grands. En effet, dans \cite{hovels}, on autorise des demi-appartements dirigés par des racines imaginaires, et une intersection infinie de demi-appartements est close, pourvu que ces demi-appartements soient dirigés par des racines  distinctes.\\}

Si $\Omega$ est une partie de $A$, on notera $\vec\Omega$ la réunion des directions des façades rencontrées par $\Omega$. Lorsque $\Omega$ est un filtre, $\vec\Omega$ sera la réunion des directions des façades rencontrées par tous les éléments de $\Omega$.\\

\fexs{\label{exemples:direction d enclos}

 \item Soit $\vec C$ une chambre de $\vec A$ et $\Omega= A_{\vec C}$. Alors $\cl(\Omega)$ est le filtre des voisinages de $\Omega$, $\vec\Omega=\vec C$, et $\cl(\vec\Omega) = \fl{\cl(\Omega)} = \barre{\vec C}$. En effet tout élément du filtre $\cl(\Omega)$ contient des points de chaque façade dirigée par une facette de $\barre{\vec C}$, même si $\cl(\Omega)$ ne contient aucun point d'aucune façade $A_{\vec f}$ pour $\vec f\subset \partial \vec C$. \\
 
 \item Soit $\vec m$ une cloison de $\vec A$, prenons $\Omega = A_{\vec m} \cup A_{-\vec m}$. Alors $\vec\Omega=\vec m\cup -\vec m$, et $\cl(\vec\Omega)$ est l'hyperplan contenant $\vec m$. Mais $\cl(\Omega) = A$, donc $\fl{\cl(\Omega)}=\vec A \not= \cl(\vec\Omega)$. (Il suffit même de prendre $\Omega=A_{\vec m} \cup\{x\}$ avec $x$ un point de $A_{-\vec m}$.)\\
 \item Avec encore $\vec m$ une cloison de $\vec A$, en prenant $\Omega=A_{\vec m}$, on obtient $\fl{\cl(\Omega)} = \barre{\vec m^*} \not= \cl(\vec \Omega)=\barre{\vec m}$. Ainsi même en restant dans $A^+$ ou $A^-$ on n'a pas $\fl{\cl(\Omega)} =\cl(\vec \Omega)$.\\

}

\begin{prop}\label{prop:direction d enclos}
Pour toute partie $\Omega$ de $A$, $\fl{\cl(\Omega)}$ est close. En conséquence, $\fl{\cl(\Omega)}\supset \cl(\vec\Omega)$.
\end{prop}

\demo\\
Soit $\vec f\subset \cl(\fl{\cl(\Omega)})$, montrons que $\vec f\subset \fl{\cl(\Omega)}$. Il s'agit de prouver que tout élément du filtre $\cl(\Omega)$ contient un point de $A_{\vec f}$. Soit donc $D_1\cap...\cap D_k \in \cl(\Omega)$ une intersection finie de demi-appartements contenant $\Omega$. Pour tout $i$, $\cl(\Omega)\subset D_i$ donc $\fl{\cl(\Omega)} \subset \vec D_i$, donc $\vec f\subset \vec D_i$. Si $\vec f$ est dans l'intérieur de $\vec D_i$, alors $A_{\vec f}\subset D_i$. On peut donc, quitte à retirer les $D_i$ tels que $\vec f$ est dans l'intérieur de $\vec D_i$, supposer que pour tout $i$, $\vec f\subset \partial \vec D_i$. Alors $D_1\cap D_2...\cap D_k$ est stable par $\vect(\vec f)$ et contient donc $\eng{\Omega,\vec f}_A$, qui contient bien au moins un point de $A_{\vec f}$ si $\Omega\not=\vide$. Le cas $\Omega=\vide$ est trivial..\cqfd

Le résultat suivant fournit une description plus ou moins constructive de la trace de l'enclos d'une partie $\Omega$ dans une façade $A_{\vec f_0}$.\\

\begin{prop}\label{prop:trace dun enclos dans une facade}
Soit $\Omega$ une partie de $A$. Soit $\D$ l'ensemble des facettes de $\fl{\cl(\Omega)})$.  On effectue les opérations suivantes sur $\Omega$:
\begin{enumerate}
\item Pour chaque couple $(a,\vec g)$ tel que $a\in \Omega$, $\vec g\in\D$ et $a$ est dans la façade $A_{\vec f}$ avec $\vec f\subset \vect(\vec g)$, on rajoute $pr_{\vec g}(a)$ à $\Omega$.
\item Pour chaque couple $(b,\vec f)$, avec $b\in \Omega$, $\vec f\in\D$, tels que $b$ est dans la façade $A_{\vec g}$ avec $\vec f\subset \vect(\vec g)$, on choisit $a\in pr_{\vec g}\inv(b)\cap A_{\vec f}$ et on rajoute $a+\vec g$ à $\Omega$.
\end{enumerate}

Appelons $\Omega^1(C_1)$ l'ensemble ainsi obtenu, où $C_1$ représente les choix effectués à chaque opération 2. Si de nouveaux couples $(a,\vec g)$ ou $(b,\vec f)$ vérifiant les conditions ci-dessus sont apparus, on effectue à nouveau les opérations 1 et 2, et on note $\Omega^2(C_2)$ l'ensemble obtenu. On obtient ainsi par récurrence un ensemble $\Omega^n(C_n)$ pour tout $n\in\N$, on note $\Omega^\infini(C)$ la réunion de tous ces ensembles, il dépend de la suite $C$ de tous les choix effectués à chaque opération 2. Notons $\maj C$ l'ensemble de toutes les suites de choix possibles. Alors pour tout$\vec f_0\in\F(\vec A(T))$:

\[\cl(\Omega)\cap A_{\vec f_0}= \bigcap_{C\in\maj C} \cl( \Omega^\infini(C) \cap A_{\vec f_0})\]
Ou plut\^ot, $\cl(\Omega)\cap A_{\vec f_0}$ est le filtre engendré par les $\cl( \Omega^\infini(C) \cap A_{\vec f_0})$, pour $C\in\maj C$.
\end{prop}

\rema Ceci signifie grosso modo que $\cl(\Omega)$ est la clôture de $\Omega$ sous les opérations 1, 2, et "prendre la clôture dans chaque façade". La difficulté de rédaction vient du fait que l'opération 2 n'est pas bien définie puisqu'elle dépend d'un choix.\\

\demo\\
Pour montrer l'inclusion $"\supset "$, il suffit de vérifier que pour tout demi-appartement $D$ contenant $\Omega$, il existe un choix $C\in\maj C$ tel que $D\supset \Omega^\infini(C)$. Ceci revient à vérifier que si $D$ contient une partie $\Theta$, alors il contient toute partie obtenue à partir de $\Theta$ par une opération 1, et que pour chaque couple $(b,\vec f)$ vérifiant les conditions de 2, il existe un choix de $a\in pr_{\vec g}\inv(b)$ tel que la partie obtenue par l'opération 2 à partir de $\Theta$ est encore incluse dans $D$. Ces vérifications sont immédiates.\\

Pour montrer l'autre inclusion, il faut prouver que si $D$ est un demi-appartement, dirigé par une racine $\alpha\in\phi^m(\vec f_0)$, contenant un $\Omega^\infini(C)\cap A_{\vec f_0}$, pour un $C\in\maj C$, alors $D\supset \Omega$. 
Soit donc $D$ un tel demi-appartement, et supposons par l'absurde qu'il existe $\omega\in \Omega\setminus D$. Soit $\vec g$ la direction de la façade contenant $\omega$. Soit $\vec h=pr_{\vec f_0}(\vec g)$, en appliquant l'opération 1 à $\Omega$ avec le couple $(\omega,\vec h)$, on voit que $pr_{\vec h}(\omega)\in \Omega^\infini(C)$. Ensuite, en appliquant l'opération 2 avec le couple $(pr_{\vec h}(\omega), \vec f_0)$, on voit que $\Omega^\infini(C)\cap A_{\vec f_0}$, et donc en particulier $D$, contient un cône de la forme $a+\vec h$. Ceci entraine que $\alpha(\vec h)\geq 0$, d'où $\alpha(\vec g)\geq 0$. D'autre part, $\alpha(\vec g)\leq 0$ sans quoi on aurait $\omega\in A_{\vec g}\subset D$. Ainsi, $\alpha(\vec g)=0$: $\vec g$, $\vec f_0$ et donc aussi $\vec h$ sont dans $\ker(\alpha)$. Donc $\alpha(\omega) = \alpha(pr_{\vec h}(\omega)) = \alpha (a)$. Mais ceci contredit le fait que $\omega\not\in D$ alors que $a\in D$.\cqfd\\

\begin{cor}
Soient $\vec f,\vec g$ deux facettes incluses dans $\fl{\cl(\Omega)}$, telles que $\vec f\subset \vect(\vec g)$. On note $\Omega_{\vec g}=\cl(\Omega)\cap A_{\vec g}$, $\Omega_{\vec f}=\cl(\Omega)\cap A_{\vec f}$. Alors $\Omega_{\vec g}=pr_{\vec g}(\Omega_{\vec f})$.\\
%Si de plus $\vec f\subset \barre{\vec g}$, alors $\Omega_{\vec f}$ contient un voisinage de $\Omega_{\vec g}$ dans $A_{\vec f}\cap \aff(\cl(\Omega))$.
\end{cor}

Pour utiliser le résultat de la proposition \ref{prop:trace dun enclos dans une facade} lorsqu'on ne connait pas précisément les directions des façades rencontrées par $\cl(\Omega)$, on pourra utiliser le lemme suivant:
\begin{lemme}\label{lemme:trace dun enclos dans une facade spherique}
On se place à nouveau dans les conditions de la proposition \ref{prop:trace dun enclos dans une facade}. On suppose en outre que $\vec f_0\subset \cl(\vec\Omega)$. Alors le résultat de la proposition \ref{prop:trace dun enclos dans une facade} est encore valable si on définit les $\Omega^\infini(C)$ de la même manière, mais en n'effectuant les opérations 1 et 2 que lorsque les facettes $\vec f$ ou $\vec g$ concernées sont dans $\cl(\vec\Omega)$.
\end{lemme}

\pv\\
 Les ensembles $\Omega^\infini(C)$ obtenus ici sont plus petits que ceux obtenus en \ref{prop:trace dun enclos dans une facade}, donc l'inclusion\\ $\cl(\Omega)\cap A_{\vec f_0} \supset \bigcap_{C\in\maj C} \cl\left( \Omega^\infini(C) \cap A_{\vec f_0}\right)$ est encore vraie.\\
 Pour l'inclusion réciproque, la preuve de \ref{prop:trace dun enclos dans une facade} est encore vraie puisqu'elle ne passe que par des facettes $\vec g\subset\vec\Omega$ et $\vec h=pr_{\vec f_0}(\vec g)\subset \cl(\vec\Omega)$.\cqfd\\

\subsubsection{Facettes}

\begin{defin}
Soit $x\in A$, soit $\vec f$ la direction de la façade contenant $x$. On note $\vec A_x$ l'espace vectoriel $\vec A_{\vec f}$, muni des directions des murs contenant $x$. C'est donc un complexe de Coxeter, a priori non essentiel, de groupe $W(\vec A_x)= \ens{\vec w\in W(\vec A)}{w.x=x}\subset \fix_{W(\vec A)}(\vec f)$. On y pense comme à l'espace tangent de $A$ en $x$.
\end{defin}

 Soit $x\in A(T)$, soit $A(T)_{\vec f}$ la façade contenant $x$. Soit $\vec F\subset\fl{A(T)_{\vec f}} $ une facette de $\vec A_x$. On note $F(x,\vec F)=\gm_x(x+\vec F)$ le filtre engendré par les parties closes de $A(T)$ contenant un voisinage de $x$ dans $x+\vec F$ (pour la topologie induite). Insistons sur le fait que $F(x,\vec F)$ est engendré uniquement par des parties closes.\\
 L'ensemble de ces filtres est l'ensemble des facettes de $A(T)$. Si $F=F(x,\vec F)$ est une facette de $A(T)$, la facette vectorielle $\vec F$ est uniquement déterminée par $F$, c'est la direction de $F$. Le point $x$ par contre n'est uniquement déterminé que lorsque $\Lambda$ est non discret ou que $x$ est un sommet de $A$.\\
 
 Dans le cas où $\Lambda$ est discret, et où $\phi(T)$ est fini, les facettes sont en fait les filtres associés à des ensembles, et ces ensembles sont les facettes affines fermées habituelles.\\
%mettre des exemples!! 

 \rema Cette définition est identique à celle de \cite{microaffine} pour une facette sphérique, bien que non présentée de la m\^eme manière. Elle diffère cependant de celle de \cite{hovels} pour une facette de $\inte A$.\\

\subsubsection{Action de $N$}
\label{subsubsection:action de N}
 Le normalisateur $N$ du tore $T$ agit sur l'appartement vectoriel $\vec A$, et même sur l'espace $\vec V$. On notera $\vec\nu:N\fleche W(\vec A)=Gl(\vec V)$ cette action. On va définir (suivant l'exemple de \cite{bruhat-tits}) une action affine $\nu$ de $N$ sur $\inte A$, qui s'étendra à $A$, dont la partie vectorielle sera $\vec \nu$, et telle que l'élément $n(u)$, pour $u\in U_\alpha$, $\alpha\in\phi$ agira par réflexion selon le mur $M(\alpha,\phii_\alpha(u))$. Remarquons que puisque le tore $T$ fixe $\vec V$, il devra agir sur $A$ par translation.

\fprop{\label{prop:action de T}Pour tout $t\in T$, il existe un unique vecteur $\vec v_t\in \vec V$ tel que pour tout $\alpha\in \phi$, pour tout $u\in U_\alpha\setminus\{e\}$, 
\[ \D(\alpha,\phii_\alpha(t u t\inv)) = \D(\alpha,\phii_\alpha(u))+\vec v_t \point\]
Ce vecteur est caractérisé par les égalités $\alpha(\vec v_t) = \phii_\alpha(u) -\phii_\alpha(t u t\inv)$, pour tout $\alpha\in \phi$ et $u\in U_\alpha\setminus\{e\}$.

L'application qui à $t$ associe la translation de vecteur $\vec v_t$ est  une action de $T$ sur $Y(T)$.
} 
\demo

On commence par prouver l'unicité. Si $\vec v_t$ est un vecteur convenable, alors pour tout $\alpha\in\phi$ et $u\in U_\alpha\setminus\{e\}$, on a:
\calc{
 \D(\alpha,\phii_\alpha(u))+\vec v_t &=& \ens{ a\in  A}{\alpha(a)+\phii_\alpha(u)\geq 0} +\vec v_t\\
 &=& \ens{ a+\vec v_t \in  A}{\alpha(a)+\phii_\alpha(u)\geq 0}\\
 &=& \ens{ a \in  A}{\alpha(a-\vec v_t)+\phii_\alpha(u)\geq 0}\\
  &=& \ens{ a \in  A}{\alpha(a)-\alpha(\vec v_t)+\phii_\alpha(u)\geq 0} = \D(\alpha,\phii_\alpha(u) - \alpha(\vec v_t))\point
}
 
On en déduit $ \phii_\alpha(t u t\inv)= \phii_\alpha(u) - \alpha(\vec v_t)$ ou encore $\alpha(\vec v_t) = \phii_\alpha(u) -\phii_\alpha(t u t\inv)$. Comme $\phi$ est une famille génératrice de $\vec V^*$, ces conditions pour tous les $\alpha\in\phi$ forcent l'unicité de $\vec v_t$.\\

Passons à l'existence. D'après le calcul précédent, la condition $D(\alpha,\phii_\alpha(t u t\inv)) = D(\alpha,\phii_\alpha(u))+\vec v_t$ équivaut à $\alpha(\vec v_t) = \phii_\alpha(u) -\phii_\alpha(t u t\inv)$.
Soit $\Pi$ un système de racines simples dans $\phi$, il s'agit donc d'une base de $\vec V^*$. Pour chaque $\alpha\in\phi$, on choisit un $u_\alpha\in U_\alpha\setminus\{e\}$. Il existe alors un unique $\vec v_t\in \vec V$ tel que $\forall \alpha\in\Pi$, $\alpha(\vec v_t) = \phii_\alpha(u_\alpha) -\phii_\alpha(t u_\alpha t\inv)$. D'après la condition (V2.2) des valuations de données radicielles, la quantité $\phii_\alpha(u_\alpha) -\phii_\alpha(t u_\alpha t\inv)$ est indépendante du choix de $u_\alpha$, donc l'égalité précédent reste vraie pour tout $u\in U_\alpha\setminus\{e\}$. Il reste à montrer que l'ensemble des racines vérifiant cette propriété est stable par n'importe quelle réflexion $r_\beta$, $\beta\in\Pi$.

C'est une conséquence de (V2.1). Soit $\alpha\in\phi$ une racine vérifiant $\forall u\in U_\alpha\setminus\{e\}$,  $\alpha(\vec v_t) = \phii_\alpha(u_\alpha) -\phii(t u_\alpha t\inv)$. Soit $\beta$ une racine simple. Alors:
\calc{
r_\beta.\alpha(\vec v_t) &=& \alpha(\vec v_t) -\langle \beta,\alpha\rangle \beta(\vec v_t)\\
&=& \phii_\alpha(u_\alpha) -\phii_\alpha(t u_\alpha t\inv) - \langle \beta,\alpha\rangle \left( \phii_\beta(u_\beta) - \phii_\beta(t u_\beta t\inv) \right)\\
&=& \left(\phii_\alpha(u_\alpha) - \langle \beta,\alpha\rangle\phii_\beta(u_\beta) \right) - \left( \phii_\alpha(t u_\alpha t\inv)-\langle \beta,\alpha\rangle\phii_\beta(t u_\beta t\inv) \right)
}
 
 Mais par (V2.1), le premier terme est $\phii_{r_\beta\alpha}(n(u_\beta).u_\alpha.n(u_\beta)\inv)$ et le second est $\phii_{r_\beta\alpha}(n(tu_\beta t\inv) . t u_\alpha t\inv . n(tu_\beta t\inv)\inv))$. Comme $n(tu_\beta t\inv) = tn(u_\beta)t\inv$, ce dernier vaut $\phii_{r_\beta\alpha}(t n(u_\beta).u_\alpha.n(u_\beta) t\inv)$, d'où le résultat.\\
 
Enfin, si $t_1$ et $t_2$ sont deux éléments de $T$, alors pour tout $\alpha\in\phi$ et $u\in U_\alpha\setminus\{e\}$, on a $\D(\alpha,\phii_\alpha(t_1t_2 u t_2\inv t_1\inv)) = \D(\alpha, \phii_\alpha(u)) + \vec v_{t_1 t_2}$ et d'autre part $\D(\alpha,\phii_\alpha(t_1t_2 u t_2\inv t_1\inv)) = \D(\alpha,\phii_\alpha(u)) + \vec v_{t_2} +\vec v_{t_1}$. Par unicité des $\vec v_t$, on obtient $\vec v_{t_1.t_2} = \vec v_{t_1}+\vec v_{t_2}$: on a bien une action de groupe.\cqfd

\fprop{Il existe une unique action $\nu$ de $N$ sur $A$ par automorphismes d'appartement telle que pour tout $\alpha\in\phi$ et $u\in U_\alpha\setminus\{e\}$, $\nu(n(u))$ est la réflexion orthogonale selon le mur $M(\alpha,\phii_\alpha(u))$
 %pour tout $n\in N$, la partie vectorielle de $\nu(n)$ est $\vec \nu(n)$, 
 et pour tout $t\in T$, $\nu(t)$ est la translation de vecteur $\vec v_t$.

Pour tout $n\in N$, on a $\fl{\nu(n)} = \vec\nu(n)$, autrement dit la partie vectorielle de cette action est $\vec\nu$.
Enfin, cette action échange les demi-appartements de $A$ selon la formule:
\[\forall n\in N,\,\alpha\in\phi\et u\in U_\alpha\setminus\{e\}, \:\: \nu(n).\D(\alpha,\phii_\alpha(u)) = \D(\vec\nu(n).\alpha, \;\phii_{\vec\nu(n).\alpha}(nun\inv)) \point\]
} 
 
 \demo\\
 
 Rappelons qu'on a fixé un point $o\in \inte A$ tel que pour tout $\alpha\in\phi$, le mur $M(\alpha,0)$ passe par $o$. Notons $N_o =\engens{ n(u)}{\alpha\in\phi,\, u\in U_\alpha\setminus\{e\} \et \phii_\alpha(u)=0}$. On commence par définir $\nu$ sur $N_o$ en posant que $\forall n\in N_o$, $\nu(n)$ est l'automorphisme affine de $Y(T)$ qui fixe $o$ et dont la partie vectorielle est $\vec \nu(n)$.
 
 On veut ensuite définir $\nu$ sur $T$ en posant pour $t\in T$ que $\nu(t)$ soit la translation de vecteur $\vec v_t$, il faut vérifier que les deux définitions sont compatibles sur $N_o\cap T$. Déjà, $\vec \nu (T)=\{id_{\vec Y}\}$, donc $\nu(N_o\cap T) = \{id_{Y(T)}\}$. Il reste donc à prouver que pour tout $t\in N_o\cap T$, $\vec v_t=\vec 0$. Fixons un tel $t$, au vu de la proposition précédente, il suffit de prouver que pour tout $\alpha\in\phi$ et $u\in U_\alpha\setminus\{e\}$, $\phii_\alpha(u) = \phii_\alpha(tut\inv)$. Fixons de tels $\alpha$ et $u$. Il découle de (V2.1) que pour tout $\beta\in\phi$ et $v\in U_\beta$ tel que $\phii_\beta(v)=0$, $\phii_{r_\beta.\alpha}(n(v).u.n(v)\inv) = \phii_\alpha(u)$. Ceci entraine que pour tout $n\in N_o$, $\phii_{n.\alpha}(nun\inv) = \phii_\alpha(u)$, et en particulier, $\phii_\alpha(u) = \phii_\alpha(tut\inv)$.
 
 On a ainsi défini $\nu$ sur $N_o\cup T$. Montrons que $N=N_o.T$. Pour tout $\alpha\in\phi$ et $u\in U_\alpha\setminus\{e\}$, il existe $u_o\in U_\alpha\setminus\{e\}$ tel que $\phii_\alpha(u_o)=0$. Donc $n(u)n(u_o)\in T$ et $n(u_o)\in N_o$, on prouve ainsi que $N=\eng{N_o,T}$. Mais comme $N_o$ normalise $T$, on obtient bien $N=N_o.T$.
 
  On définit alors $\nu$ sur $N$ par $\nu(n_o t)=\nu(n_o)\circ\nu(t)$, pour tous $n_o\in N_o$ et $t\in T$. Ceci est bien défini car $\nu$ est trivial sur $N_o\cap T$. Comme $\vec\nu(T) = \{id\}$, il est évident que la partie vectorielle de $\nu(n)$ est $\vec\nu(n)$ pour tout $n\in N$. Montrons que $\nu$ est une action de groupe. Soient $n_1,n_2\in N_o$ et $t_1,t_2\in T$, par définition on a $\nu(n_1t_1n_2t_2) = \nu(n_1n_2n_2\inv t_1 n_2t_2) = \nu(n_1)\nu(n_2)\nu(n_2\inv t_1 n_2)\nu(t_2)$. Ceci devrait valoir $\nu(n_1)\nu(t_1)\nu(n_2 )\nu(t_2)$, nous devons donc prouver que $\nu(n_2)\nu(n_2\inv t_1 n_2) = \nu(t_1)\nu(n_2)$ autrement dit que $\vec \nu(n_2\inv).\vec v_{t_1} = \vec v_{n_2\inv t_1 n_2}$. Il suffit de traiter le cas où $\vec\nu(n_2)$ est une réflexion: il existe alors $\beta\in\phi$ et $v\in U_\beta\setminus\{e\}$, avec $\phii_\beta(v)=0$ tel que $n_2=n(v)$. Soit $\alpha\in\phi$ et $u\in U_\alpha\setminus\{e\}$. Alors:
  \calc{
  \alpha( \vec v_{n_2\inv t_1 n_2}) &=&  \phii_\alpha(u) -\phii_\alpha(n_2 \inv t_1 n_2. u. n_2\inv t_1\inv n_2) \\
  &=&  \phii_\alpha(u) - \phii_{r_\beta\alpha}(t_1 n_2. u. n_2\inv t_1\inv) - \eng{\beta,\alpha }\phii_\beta(v) \\
  &=& \phii_\alpha(u)+r_\beta\alpha(\vec v_{t_1}) - \phii_{r_\beta\alpha}(n_2 u n_2\inv) - \eng{\beta,\alpha}\phii_\beta(v)  \\
  &=& \phii_\alpha(u)+\alpha(r_\beta.\vec v_{t_1}) - \phii_\alpha(u) + \eng{\beta,\alpha}\phii_\beta(v) - \eng{\beta,\alpha}\phii_\beta(v) \\
  &=& \alpha(r_\beta.\vec v_{t_1})\\
  &=& \alpha(\vec\nu(n_2\inv).\vec v_{t_1})
  }
Ceci, étant vrai quelque soit $\alpha\in\phi$, prouve bien que $\vec \nu(n_2\inv).\vec v_{t_1} = \vec v_{n_2\inv t_1 n_2}$.\\

Comme la partie vectorielle de cette action est $\vec\nu$, qui préserve l'ensemble des facettes de $\vec A$, elle s'étend à une action sur $A$. Prouvons qu'elle stabilise l'ensemble des demi-appartements de $\inte A$ selon la formule annoncée. Soit $\alpha\in\phi$, $u\in U_\alpha\setminus\{e\}$. La relation $\nu(n).\D(\alpha,\phii_\alpha(u)) = \D(\vec\nu(n).\alpha, \phii_{\vec\nu(n).\alpha}(nun\inv))$ est déjà vraie pour $n\in T$, par la définition des $\vec v_t$. Soit $n\in N_0$, on a déjà vu qu'alors $\phii_{\vec\nu(n).\alpha}(nun\inv) = \phii_\alpha(u)$. Il ne reste plus qu'à calculer:
\calc{
\nu(n).D(\alpha,\:\phii_\alpha(u)) &=& \ens{\nu(n).(o+\vec v) \in  A}{\alpha(o+\vec v)+ \phii_\alpha(u)\geq 0}\\
&=& \ens{o+\vec v\in A}{\alpha\left(\nu(n)\inv(o+\vec v)\right)+ \phii_\alpha(u) \geq 0}\\
&=& \ens{o+\vec v\in  A}{\alpha\left(o+\vec \nu(n)\inv.(\vec v)\right)+ \phii_\alpha(u) \geq 0}\\
&=& \ens{o+\vec v\in  A}{\vec\nu(n).\alpha(o+\vec v)+ \phii_\alpha(u) \geq 0}\\
&=& D\left(\vec\nu(n).\alpha, \:\phii_\alpha(u)\right)\\
&=& D\left(\vec\nu(n).\alpha, \:\phii_{\vec\nu(n).\alpha}(nun\inv)\right)
}

A présent, soit $\alpha\in\phi$, $u\in U_\alpha\setminus\{e\}$, montrons que $\nu(n(u))$ est la réflexion selon le mur $M(\alpha,\phii_\alpha(u))$. On sait déjà que la partie vectorielle de $\nu(n(u))$ est une réflexion selon le mur $\ker(\alpha)$, donc $\nu(n(u))$ est la composée d'une réflexion et d'une translation de vecteur $\vec w\in \ker(\alpha)$.

Alors $\nu(n(u))^2$ est la translation de vecteur $2\vec w=\vec v_{n(u)^2}$. Mais pour tout $\beta\in\phi$, $v\in U_\beta$, on a:
\calc{
\phii_\beta(n(u)^2.v.n(u)^{-2}) &=& \phii_{r_\alpha\beta}(n(u).v.n(u)\inv) - \eng{\alpha,\beta}\phii_\alpha(u)\\
&=& \phii_\beta(v) - \eng{\alpha,r_\alpha \beta}\phii_\alpha(u) - \eng{\alpha,\beta}\phii_\alpha(u)\\
&=& \phii_\beta(v)
}
car $\eng{\alpha,r_\alpha \beta} = \eng{r_\alpha\alpha,\beta} = -\eng{\alpha,\beta}$. Ceci prouve que $\vec w=0$, donc $\nu(n(u))$ est une réflexion.

Maintenant, par le résultat précédent, $\nu(n).M(\alpha,\phii_\alpha(u)) = M(-\alpha,\phii_{-\alpha}(n(u).u.n(u)\inv)$. Mais par (V2.1), $$\phii_{-\alpha}(n(u).u.n(u)\inv)= \phii_\alpha(u) -\eng{\alpha,\alpha}\phii_\alpha(u) = -\phii_\alpha(u)\point$$ Donc $\nu(n).M(\alpha,\phii_\alpha(u)) = M(-\alpha,-\phii_\alpha(u)) = M(\alpha,\phii_\alpha(u))$. 
Sachant que l'hyperplan fixe de $\nu(n(u))$ est parallèle à  $M(\alpha,\phii_\alpha(u))$, ceci entraine que c'est précisément  $M(\alpha,\phii_\alpha(u))$, et donc que $\nu(n(u))$ est la réflexion annoncée.\\

L'unicité de $\nu$ est claire car $N$ est engendré par $T$ et les $n(u)$, pour $u\in U_\alpha\setminus\{e\}$, $\alpha\in\phi$.
\cqfd

 La définition de cette action de $N$ permet de prouver la condition "(V5)" présente dans la définition de la valuation d'une donnée radicielle pour \cite{bruhat-tits} ou  \cite{microaffine}:
\fcor{ \label{prop:V5} Soit $\alpha\in\phi$, $u\in U_\alpha\setminus\{e\}$ et $u',u''\in U_{-\alpha}$ tels que $n(u)=u'uu''$. Alors $\phii_\alpha(u)=-\phii_{-\alpha}(u') = -\phii_{-\alpha}(u'')$.

Lorsque $\phi$ est un système de racines fini, la définition de valuation d'une donnée radicielle donnée en \ref{soussection:valuation} équivaut à celle de \cite{bruhat-tits} 6.2.1 .

}
\demo
Il est classique que $n(u)=n(u')=n(u'')$. Rappelons tout de même la preuve: on a $u'uu'' =uu'' n(u)\inv u' n(u)$. Mais $n(u)\inv u' n(u) \in U_\alpha$, d'où le résultat. Or $n(u)$, $n(u')$, $n(u'')$ agissent respectivement  par les réflexions selon $M(\alpha,\phii_\alpha(u))$, $M(-\alpha,\phii_\alpha(u'))$ et $M(-\alpha,\phii_\alpha(u''))$. L'égalité de ces trois murs entraine bien les égalités annoncées.\\

Comme la condition (V2) de \ref{soussection:valuation} est clairement plus forte que celle de \cite{bruhat-tits}, et comme les autres conditions (  (V0), (V1), (V3), (V4) ) sont inchangée, le point précédent prouve qu'une valuation au sens de \ref{soussection:valuation} est aussi une valuation pour \cite{bruhat-tits}. Nous avons vu l'autre implication en \ref{prop:existence de valuation}.
\cqfd

En conséquence de ce corollaire, pour tout facette sphérique $\vec f$, la famille $(\phii_\alpha)_{\alpha\in\phi^m(\vec f)}$ est une valuation au sens de \cite{bruhat-tits} de la donnée radicielle $(M_{\vec A}(\vec f), (U_\alpha)_{\alpha_in\phi^m(\vec f)})$. La donnée radicielle valuée $(M_{\vec A}(\vec f), (U_\alpha,\phii_\alpha)_{\alpha_in\phi^m(\vec f)})$ sera notée $\D_{\vec f}$.\\

\fdef{
Le fixateur dans $N$ d'un point ou d'une partie $a$ de $A$ ou de $\vec A$ sera noté $N(a)$ (et donc $N(T)(a)$ s'il faut préciser le tore). Le fixateur de $A$ sera noté $H(T)$, ou juste $H$ s'il est inutile de préciser le tore.}

Dans la suite, on omettra souvent de noter $\nu$ et $\vec\nu$ pour l'action d'un élément de $N$ sur un point  de $A$ ou de $\vec A$.\\

\fex{\label{exemple:action de N}
Remarquons tout de suite que pour une partie $\Omega\subset A$, $N(\Omega)\not\subset N(\cl(\Omega))$. Il suffit de choisir deux chambres $\vec c$ et $\vec d$ de $\vec A(T)$, séparées par une cloison $\vec m$, et de trouver $t\in T$ qui induit une translation dont la direction n'est pas incluse dans $\vec m$. Alors $t$ fixe $A_{\vec c}\cup A_{\vec d}$ mais pas $A_{\vec m}$, alors que $\cl(A_{\vec c}\cup A_{\vec d})=A_{\vec c}\cup A_{\vec d}\cup A_{\vec m}$.\\}

Il est par contre clair d'après \ref{props:paraboliques} qu'un $n\in N(\Omega)$ agit comme une translation sur chaque sous-espace\\  $\eng{\omega,\vect(\cl(\vec \Omega))}_A$. La proposition suivante améliore un peu ce résultat, en permettant de remplacer $\cl(\vec\Omega)$ par $\fl{\cl(\Omega)}$. Notons que les différentes translations induites sur chaque $\omega+\vect(\fl{\cl(\Omega)})$ ne sont a priori pas selon le même vecteur.
\fprop{\label{prop:N(Omega) fixe fl(Cl(Omega))} Soit $\Omega$ une partie de $A$, alors $N(\Omega) \subset N(\fl{\cl(\Omega)})$.}
\demo

Soit $n\in N(\Omega)$, soit $\vec E=\fix_{\vec A}(n)$, c'est une partie close de $\vec A$ contenant $\cl(\vec\Omega)$. Supposons $\vec E\not= \fl{\cl(\Omega)}$.
Alors il existe $\alpha\in\phi$ tel que $\vec E\subset \vec D(\alpha)$ et un point $a\in\cl(\Omega)$ tel que $\alpha(a)=-\infini$. Ce point $a$ n'est donc dans aucun demi-appartement dirigé par $\alpha$, et pourtant il est dans $\cl(\Omega)$: il n'existe donc pas de demi-appartement dirigé par $\alpha$ qui contienne $\Omega$. Il existe donc $(\omega)\in\Omega^\N$ tel que $\alpha(\omega_i)\in\R$ $\forall i\in\N$ et $\lim_{i \fleche \infini}\alpha(\omega_i) = -\infini$.\\
 L'ensemble $\vec E\cap \ker\alpha\subset \vec A$ est clos et non vide puisqu'il contient les directions des façades contenant les $\omega_n$. Soit $\vec f$ une facette maximale de $\vec E\cap \ker\alpha$, alors tous les $\omega_n$ se projettent sur $A_{\vec f}$, et ces projetés sont fixes par $n$. Soit $(\omega'_n)_n\in(A_{\vec f})^\N$ la suite ainsi obtenue. Soit $n\in\N$ tel que $\alpha(\fl{\omega'_0\omega'_n})<0$. Alors $n$ fixe la droite contenant $\{\omega'_0,\omega'_n\}$, donc sa direction $\fl{\omega'_0\omega'_n}$, puis la facette $\vec g$ contenant cette direction. Donc $\vec g\subset \vec E$, mais $\alpha(\vec g)=\R^{-*}$ donc $\vec g\not\subset \vec D(\alpha)$, ce qui est impossible.\\
 Donc $\vec E=\fl{\cl(\Omega)}$.\cqfd\\

\rema La définition des façades $A_{\vec f}\simeq A/\vect(\vec f)$ revient à essentialiser $\vec A$ pour le groupe de Coxeter $\fix_{W(\vec A)}(\vec f)$. Cette construction est semblable à la compactification polyhédrale, ou de Satake, d'un appartement d'un immeuble affine. Elle permet d'avoir sur $A$ une topologie séparée, et telle que l'adhérence d'une façade sphérique est compacte. Elle a cependant le défaut de perdre une partie de l'action du tore. En effet, si $t$ induit une translation de direction incluse dans $\vect(\vec f)$ sur $Y(T)$, alors $t$ agit trivialement sur $A_{\vec f}$. Dans la première réalisation d'un appartement microaffine dans \cite{microaffine} par exemple, les façades sont toutes isomorphes comme espaces affines à $Y(T)$, ce qui évite ce souci.\\

\subsubsection{Opposition}

\begin{defin}
Si $a=[x+\vec f]\in A(T)$, le point opposé à $a$ dans $A(T)$ est $op_{A(T)}(a) = [x-\vec f]$.
\end{defin}

L'application $op_{A(T)}$ est une involution qui permute les façades de $A(T)$, préserve l'ensemble des murs et commute à l'action de $W(T)$. Plus généralement, elle commute à tout isomorphisme d'appartements. Cependant, ce n'est pas un automorphisme d'appartement car l'action sur le bord d'une façade n'est pas induite par l'action sur cette façade ($\op_{A(T)}$ n'est pas continue). En fait, $op_{A(T)}$ fixe la façade principale, et le seul automorphisme d'appartement de $A(T)$ fixant la façade principale est $id_{A(T)}$.\\

\subsection{Familles de sous-groupes parahoriques}

  Maintenant que nous disposons des appartements $A(T)$, il faut, pour définir un immeuble selon la méthode usuelle, déterminer quels seront les fixateurs des points de $A(T)$. Ces fixateurs seront appelés des sous-groupes parahoriques de $G$.
  
 Dans cette sous-section, on étudie quelles sont en général les propriétés qu'on peut espérer d'une famille de sous-groupes parahoriques. On étudie également l'exemple le plus simple de telle famille: la "famille minimale de parahoriques".\\
 
 On fixe un tore maximal $T$, et on note $A=A(T)$, $N=N(T)$.\\

\subsubsection{Définition}
\label{sousoussection:definition des parahoriques}

 \begin{defin}
 Soit $Q=(Q(a))_{a\in A}$ une famille de sous-groupes de $G$. Si $\Omega$ est une partie de $A$, on note $Q(\Omega)=\bigcap_{\omega\in\Omega}Q(\omega)$. Si $\Omega$ est un filtre de $A$, on note $Q(\Omega) =\bigcup_{\Omega'\in\Omega}Q(\Omega')$.\\
 
 On dit que $Q$ est une famille de sous-groupes parahoriques pour $\D$ si elle vérifie:
 \liste{
 	\item (para 0.1): Si $a\in A_{\vec f}$, alors $U(\vec f)\subset Q(a) \subset P(\vec f)$.  (compatibilité avec l'immeuble vectoriel)
  \item (para 0.2): $\forall a\in A(T)$, $N(T)_a\subset Q(a)$.  (compatibilité de l'action de $N(T)$)
  \item (para 0.3): $\forall a\in A(T)$, $\forall (\alpha,\lambda)\in \phi(T)\times \R$ tel que $a\in\D(\alpha,\lambda)$, $U_{\alpha,\lambda}\subset Q(a)$.  (points fixes des groupes radiciels)
    \item (para 0.4): $\forall n\in N(T)$, $\forall a\in A(T)$, $n Q(a) n\inv = Q(na)$.\\
}

Si $Q_1$ et $Q_2$ sont deux familles de parahoriques, on dira que $Q_2$ contient $Q_1$ si $\forall a\in A,\, Q_1(a)\subset Q_2(a)$.

On note $\P=(P(a))_{a\in A}$ la famille de sous-groupes parahoriques de $G$ définie par $\forall a\in A$, $P(a) = \eng{U(\vec f_a),N(a),G(a)}$.\\

 On définit encore les condition suivantes sur $Q$, certaines dépendent d'une facette $\vec g\in\F(\vec A(T))$, d'une partie $\Omega\subset A$, d'un point $a\in A$ ou d'une chambre $\vec C$ de $\vec f_a^*\cap \vec A$:\\
  \begin{itemize}
  \item \pinj: $\forall a\in A(T)$, $N(T)_a = Q(a)\cap N(T)$. (inclusion des appartements dans l'immeuble)\\

\item \psph: Pour tout $a\in A_{sph}$, $Q(a)=P(a)$. (valeur sur les points sphériques)\\  

\item \plieng \:(lien entre une façade et son bord) qui s'énonce en plusieurs variantes:
\liste{  
      \item \plien($\vec g$): $\forall \vec f\in\F(\barre{\vec g})$, $\forall a\in A_{\vec f}$, $Q(a)\cap P(\vec g)= Q(\{a,pr_{\vec g}(a)\})$.
    \item \plienn($\vec g$): $\forall \vec f\in\F(\barre{\vec g})$, $\forall a\in A_{\vec f}$, $N(T)Q(a)\cap N(T)P(\vec g)= N(T)Q(\{a,pr_{\vec g}(a)\})$.
    \item \plienp($\vec g$): $\forall \vec f\in\F(\barre{\vec g})$, $\forall a\in A_{\vec f}$, $Q(a)\cap P(\vec g)= Q(\barre{a+\vec g})$.
 %   \item \plienpn($\vec g$): $\forall \vec f\in\F(\barre{\vec g})$, $\forall a\in A_{\vec f}$, $N(T)Q(a)\cap N(T)P(\vec g)= N(T)Q(\barre{a+\vec g})$.
 \\
    }
    
  %\item \pliendeux : Pour toutes facettes $\vec f,\vec g$ telles que $\vec f\subset \barre{\vec g}$, pour tout $b\in A_{\vec g}$ et $q\in Q(b)\cap P(\vec f)$, il existe $a\in pr_{\vec g}\inv(\{b\})\cap A_{\vec f} $ et $n\in N(T)_b$ tel que $nq\in Q(a+\vec g)$. (lien entre une façade et les façades qu'elle borde) \\

\item \pdec($\vec C, a$):  $Q(a)=(Q(a)\cap U(\vec C))\point (Q(a)\cap U(-\vec C)) \point N(a)$. (décomposition de $Q(a)$)\\

\item \pinter($\Omega$): $Q(\Omega) = N_\Omega. Q(\cl(\Omega))$. (intersections d'appartements)\\

\item \pisom($\Omega$):  $\bigcap_{a\in\Omega} (N(T).Q(a)) = N(T). Q(\Omega)$. (isomorphismes entre appartements)\\

\end{itemize}

Lorsque $Q$ vérifie \plien($\vec f$) (ou une de ses variantes) pour toute facette $\vec f$, on dira juste que $Q$ vérifie \plien. Lorsqu'elle vérifie \plien($\vec f$) pour toute facette sphérique $\vec f\in \F(\vec A)$, on dira qu'elle vérifie \plien(sph), lorsqu'elle vérifie \plien($\vec m$) pour toute cloison $\vec m$ de $\vec A$, on dira qu'elle vérifie \plien(cloison), etc...\\

Une famille de parahoriques vérifiant \psph, \pinj\ et \plienn(sph) sera appelée une bonne famille de parahoriques. %Si elle vérifie en outre \plienpn(sph), elle sera appelée une très bonne famille de parahoriques.

 \end{defin}

 \remas{
 \item  Ici, pour une partie $\Omega$ de $A$, $Q(\Omega)$ désigne par définition $\bigcap_{\omega\in\Omega}Q(\omega)$. Dans \cite{hovels} ou \cite{masures2}, on définit directement les valeurs d'une famille de parahoriques sur chaque partie $\Omega$ de $A$, et on prouve ensuite, au moins dans les bons cas, la relation $Q(\Omega) = \bigcap_{\omega\in\Omega}Q(\omega)$.\\

\item On prouvera en \ref{prop:plien et plienn} qu'une bonne famille de parahoriques vérifie automatiquement \plien(sph).\\

\item Il est immédiat que pour toute facette $\vec g$, la conjonction de \plienp($\vec g$) et de \plienn($\vec g$) équivaut à:
 \[ \text{\plienpn}(\vec g):\: \forall \vec f\in\F(\barre{\vec g}), \: \forall a\in A_{\vec f},\quad N(T)Q(a)\cap N(T)P(\vec g)= N(T)Q(\barre{a+\vec g}) \point\]

\item Pour toute famille de parahoriques $Q$, et pour toute partie $\Omega$ de $A$, $G(\Omega)\subset Q(\Omega)$, par (para 0.3). De même pour toute facette $\vec f$, $G(\phi^m(\vec f),\Omega)$ est en quelque sorte le plus petit groupe évidemment inclus dans $Q(\Omega)\cap M(\vec f)$.\\
}

 \fprop{La famille $\P$ est la plus petite famille de parahoriques.
}
\demo Rappelons que $G(a)=\engens{U_{\alpha,k}}{a\in\D(\alpha,k)}$. Par (para 0.1), (para 0.2) et (para 0.3), toute famille de parahorique doit être supérieure à $\P$. Mais il est clair que cette dernière vérifie (para 0).\cqfd

 Toute famille de parahoriques permettra de définir une masure bordée.
 Le but est bien sûr de trouver une famille de parahoriques vérifiant un maximum de conditions (para x), ce qui mènera à une masure possédant un maximum de propriétés semblables à celles d'un immeuble.
 
  Nous verrons que, au moins dans le cas Kac-Moody, la famille $P$ est une bonne famille de parahorique, et nous étudierons les propriétés de la masure bordée que définit une telle famille.
   
 Nous prouverons au \ref{soussection:famille maximale} l'existence d'une bonne famille de parahoriques maximale $\bar P$, de sorte que toute bonne famille de parahoriques sera à chercher entre $P$ et $\bar P$.

\subsubsection{La famille minimale de parahoriques }
\label{soussection:la famille minimale}

Dans ce paragraphe, on étudie le premier exemple de famille de parahoriques disponible: la famille minimale $\P$.\\

Lorsque $\vec f$ est sphérique, la théorie de Bruhat-Tits décrit bien les facteurs $M_{\vec A}(\vec f)\cap P(a)$, permettant de prouver la proposition suivante (qui est la raison d'être de la condition \psph):

\fprops[\label{prop:csq de parasph}Soit $\vec f$ une facette sphérique, et $\Omega\subset A_{\vec f}$.]{
\item $P(\Omega)\cap M_{\vec A}(\vec f)$ est le sous-groupe parahorique de $M_{\vec A}(\vec f)$ associé à la partie $\Omega$ de l'appartement $A_{\vec f}$ pour la donnée radicielle valuée finie $\D_{\vec f}:=(M_{\vec A}(\vec f),(U_\alpha,\phii_\alpha)_{\alpha\in\phi^m(\vec f)})$, au sens de \cite{bruhat-tits}.
\item $P(\Omega) = U(\vec f)\rtimes\left(M_{\vec A}(\vec f)\cap P(\Omega) \right)$ et $M_{\vec A}(\vec f)\cap P(\Omega)=N(\Omega)\point G\left(\phi^m(\vec f),\Omega\right)$.
%\item Pour tout $p\in P(\Omega)$, il existe une partie $\Omega_0\subset \inte A$  et $n\in N(\Omega)$ tels que $pr_{\vec f}(\Omega_0)=\Omega$ et pour tout $\omega\in\Omega_0$, $np$ fixe $\barre{\omega+\vec f}$. 
\item $P(\Omega) = N(\Omega).P(\cl_{A_{\vec f}}(\Omega))$, autrement dit, $\P$ vérifie \pinter\ sur les parties de $A_{\vec f}$.
\item $\bigcap_{\omega\in\Omega}N.P(\omega) = N.P(\Omega)$, autrement dit, $\P$ vérifie \pisom\ sur les parties de $A_{\vec f}$.
\item Pour tout $a\in A_{\vec f}$, et $\vec g\in\vec f^*\cap\vec A$, $P(a)\cap P(\vec g) = P(\barre{a+\vec g})$.
\item Pour toute chambre $\vec C$ de $\vec f^*$, $P(\Omega)=(P(\Omega)\cap U(\vec C))\point(P(\Omega)\cap U(-\vec C))\point N(\Omega)$, autrement dit, $P$ vérifie \pdec($\Omega$).\\
}

\demos{
\item Pour tout $a\in \Omega$, $P(a) = U(\vec f)\rtimes \eng{N(a), \ens{U_\alpha(a)}{\alpha\in\phi^m(\vec f)}}$, d'où $P(a)\cap M_{\vec A}(\vec f) = \eng{N(a), \ens{U_\alpha(a)}{\alpha\in\phi^m(\vec f)}}$. Ceci est précisément la définition du sous groupe parahorique de $M_{\vec A}(\vec f)$ au point $a$. Ensuite, $P(\Omega)\cap M_{\vec A}(\vec f)$ est l'intersection de tous ces groupes pour $a\in\Omega$, c'est bien le sous-groupe parahorique de $M_{\vec A}(\vec f)$ pour la partie $\Omega$.

\item 
Par la décomposition de Lévi de $P(\vec f)$, on a $P(\Omega)\subset P(\vec f)=U(\vec f)\rtimes M_{\vec A}(\vec f)$. Mais $U(\vec f)\subset P(\Omega)$ d'où
$P(\Omega) = U(\vec f)\rtimes ( M(\vec f)\cap P(\Omega))$.
Et par \cite{bruhat-tits}, $M_{\vec A}(\vec f)\cap P(\Omega)=N(\Omega).G(\phi^m(\vec f),\Omega)$.
% Il reste à étudier le facteur $U(\vec f)\cap P(\Omega)$. Tout élément $u$ du groupe $U(\vec f)$ fixe un cône de la forme $\barre{x+\vec C}$ avec $x\in\inte A$ pour toute chambre $\vec C$ de $\vec f^*\cap\vec A$. L'intersection avec $A_{\vec f}$ de la réunion de ces cônes contient une partie dont l'enclos dans $A_{\vec f}$ est $A_{\vec f}$. Ceci entraine que $u$ stabilise $A_{\vec f}$. Mais l'espace affine engendré par les $\barre{x+\vec C}$ fixés par $u$ est $A_{\vec f}$. Donc au final, $u$ fixe $A_{\vec f}$, et $U(\vec f)\cap P(\Omega)=U(\vec f)$.

%\item Comme $N(\Omega)$ normalise $U(\vec f)\cap P(\Omega)$, le point précédent entraine $P(\Omega)=N(\Omega)\point (U(\vec f)\cap P(\Omega))\point G(\phi^m(\vec f),\Omega)$.
% Le groupe $G(\phi^m(\vec f),\Omega)$ fixe $\eng{\Omega,\vec f}_A$, il reste juste à étudier $U(\vec f)\cap P(\Omega)$. Pour tout $a\in\Omega$, la description de $U(\vec f)\cap P(a)$ comme le plus petit groupe contenant $G(\phi^u(\vec f))$ et normalisé par $G(\phi^m(\vec f),a)$ prouve que tout élément de $U(\vec f)\cap P(a)$ fixe un cône de la forme $\omega+\barre{\vec f}$ pour un $\omega\in\inte A$ tel que $pr_{\vec f}(\omega)=a$. D'où le point 2.\\
 
 \item Découle immédiatement de 1, car $G(\phi^m(\vec f),\Omega) = G(\phi^m(\vec f),\cl_{A_{\vec f}}(\Omega))$.

\item On calcule:
\begin{eqnarray*}
\bigcap_{\omega\in\Omega} N.P(\omega) &=& N.\bigcap_{\omega\in\Omega} N(\vec f).P(\omega)\\
&=& N.\bigcap_{\omega\in\Omega} U(\vec f)\point N(\vec f) \point \engens{U_\alpha(a)}{\alpha\in\phi^m(\vec f)}\\
&=&N.U(\vec f)\point \bigcap_{\omega\in\Omega} \point N(\vec f) \point \left( M(\vec f)\cap P(a) \right)\\
&=& N.U(\vec f)\point N(\vec f).\left(M(\vec f)\cap P(\Omega)\right)\\
&=& N.P(\Omega)
\end{eqnarray*} 
 
 La première égalité est vraie car pour tout $\omega\in\Omega$, $P(\omega)\subset P(\vec f)$. La troisième vient de l'unicité de la décomposition de Lévi de $P(\vec f)$, et la quatrième est le résultat classique dans les immeubles, voir \cite{bruhat-tits}.\\
 
\item  Soit $g\in P(a)\cap P(\vec g)=U(\vec f)\point (M(\vec f)\cap P(a))\cap P(\vec g)$. Pour tout $\vec h\in \vec f^*\cap \vec A$, $U(\vec f)\subset P(A_{\vec h})$. En particulier, $U(\vec f)\subset P(\barre{a+\vec g})$, on peut donc supposer $g\in M(\vec f)\cap P(a)\cap P(\vec g)$. Le résultat est alors classique.\\
 
 \item Le groupe $M(\vec f)\cap P(\Omega)$ admet par \cite{bruhat-tits} une telle décomposition, elle s'écrit ici: \[M(\vec f)\cap P(\Omega) = G(\phi(\vec C)\cap \phi^m(\vec f),\Omega)\point G(\phi(\vec C')\cap \phi^m(\vec f),\Omega)\point N(\Omega)\],
  où $\vec C'$ est la chambre opposée à $\vec C$ dans $\vec f^*$. Mais $\phi(\vec C')\cap\phi^m(\vec f) = \phi(-\vec C)\cap \phi^m(\vec f)$, d'où \[M(\vec f)\cap P(\Omega) \subset G(\phi(\vec C),\Omega)\point G(\phi(-\vec C),\Omega)\point N(\Omega) \subset (P(\Omega)\cap U(\vec C))\point (P(\Omega)\cap U(-\vec C))\point N(\Omega)\point\] Comme $U(\vec f)\subset P(\Omega)\cap U(\vec C)$ et $P(\Omega)= U(\vec f).(M(\vec f)\cap P(\Omega)$, on arrive au résultat annoncé.
 
 }

\fcor{\label{cor:csq de parasph} Les résultats de la proposition précédente sont vrais pour n'importe quelle famille $Q$ de parahoriques vérifiant $\psph$.
}

Dans le cas où $\D$ vient d'un groupe de Kac-Moody déployé, Guy Rousseau a prouvé en \cite{masures2} 4.6 que $\P$ vérifie \pdec. Pour le cas d'un groupe de Kac-Moody sur un $\mathbb C((t))$, c'est \cite{hovels} 3.4.1.
\fprop{Si $\D$ est la donnée radicielle valuée issue d'un groupe de Kac-Moody déployé $G$, alors la famille minimale de parahoriques attachée à $\D$ vérifie \pdec.
}
\cqfd

\rema Si on retire la condition $\forall a\in A,\ U(\vec f_a)\subset Q(a)$ dans (para 0), on obtient une famille minimale plus petite que $\P$. Il s'agit de $P_0$, avec $P_0(a)=\eng{N(a),G(a)}$. On peut étudier rapidement cette famille de sous-groupes de $G$.

Soit $a\in A$, et $\vec f$ la direction de la façade de $a$. Pour tout $\alpha\in\phi^u(\vec f)$, on a $U_\alpha(a)=U_\alpha$, et pour $\alpha\in\phi\setminus \phi(\vec f)$, $U_\alpha(a)=\{e\}$. Donc $G(a)=G(\phi(\vec f),a) = \engens{U_\alpha(a)}{\alpha\in\phi(\vec f)}$. Sachant que $P_0(a)\subset P(\vec f)=U(\vec f)\rtimes M_{\vec A}(\vec f)$, le groupe $\eng{ N(a),G(\phi^m(\vec f),a)}=\eng{ N(a),\ens{U_\alpha(a)}{\alpha\in\phi^m(\vec f)}}$, qui est inclus dans $M_{\vec A}(\vec f)$, normalise $U(\vec f)\cap P_0(a)$, et on prouve:
\[ P_0(a) = \left( U(\vec f)\cap P_0(a)\right)\rtimes \eng{ N(a),\: G\left(\phi^m(\vec f),a\right)}\]

Comme $N(a)$ normalise à la fois $U(\vec f)\cap P_0(a)$ et $G\left(\phi^m(\vec f),a\right)$, on a aussi les décompositions:
 
\begin{eqnarray*}
P_0(a) & = & \left( U(\vec f)\cap P_0(a)\right)\rtimes \left( N(a)\point G\left(\phi^m(\vec f),a\right)\right)\\
 & = & N(a)\point \left( U(\vec f)\cap P_0(a)\right) \point G\left(\phi^m(\vec f),a\right)
\end{eqnarray*}

%Le groupe $N(a).G\left(\phi^m(\vec f),a\right)=M_{\vec A}(\vec f)\cap P_0(a)$ est le sous-groupe parahorique au point $a$ de l'appartement\ $\barre{A_{\vec f}}$, pour la famille minimale de parahoriques dans le groupe $M_{\vec A}(\vec f)$ muni de la donnée radicielle valuée $(M_{\vec A}(\vec f), (U_\alpha,\phii_\alpha)_{\alpha\in\phi^m(\vec f)})$.

Le groupe $U(\vec f)\cap P(a)$ est le sous-groupe distingué de $P(a)$ engendré par $G(\phi^u(\vec f))$, et plus précisément le plus petit sous-groupe de $P(a)$ contenant $G(\phi^u(\vec f))$ et normalisé par $G\left(\phi^m(\vec f),a\right)$. Il peut être strictement inclus dans $U(\vec f)$.\\

Lorsque $\Omega$ est une partie de $A_{\vec f}$, par l'unicité dans la décomposition de Lévi $P(\vec f)=U(\vec f)\rtimes M_{\vec A}(\vec f)$, on obtient:
\[P_0(\Omega) = \left( U(\vec f)\cap P_0(\Omega)\right) \rtimes \left(M_{\vec A}(\vec f)\cap P_0(\Omega) \right)\: .\]

Notons que la famille $P_0$ bénéficie d'une propriété de plus que $P$: si $\vec f$ est une facette sphérique et si $\Omega\subset A_{\vec f}$, alors tout $p\in P_0(\Omega)$ fixe une partie de $A$ de la forme $\bigcup_{\omega\in\Omega_0} \barre{\omega_0+\vec f}$, avec $\Omega_0$ une partie de $\inte A$ telle que $pr_{\vec f}(\Omega_0)=\Omega$.

\subsubsection{La condition \pfonc}
\label{soussoussection:definition de fonc}

Si $a$ est un sommet spécial, alors $N(a)$, et donc $Q(a)$ pour n'importe quelle famille $Q$ de parahoriques, contient un système de représentants pour $W(\vec A)$. Donc $N.Q(a)=T.Q(a)$, ceci est un bon point de départ pour prouver par exemple \pinj, \plienn,  ou \plienpn. On est donc souvent capable de prouver ces conditions pour des sommets spéciaux.

 Pour passer à un point $a$ plus général, l'idée retenue ici est de plonger l'appartement $A$ pour la donnée radicielle valuée $\D$ dans un appartement $A^\Delta$ pour une donnée radicielle $\D^\Delta$ qui soit identique à $A$ comme ensemble, mais muni de plus de murs, de sorte que $a$ soit spécial dans $A^\Delta$. Il s'agit donc de trouver une donnée radicielle valuée, sur le même système de racines que $\D$, mais dont le groupe d'arrivée de la valuation soit plus grand que $\Lambda$.
Dans le cas d'un groupe de Kac-Moody déployé, ceci revient juste à considérer une extension (ramifiée) du corps de base.

Ceci est axiomatisé par la condition "\pfonc ":
 
\fdef{La donnée radicielle valuée $\D=(G, (U_\alpha,\phii_\alpha)_{\alpha\in\phi})$ vérifie la condition \pfonc\ si pour tout sous-groupe $\Delta$ de $\R$, il existe une donnée radicielle valuée notée $\maj D^\Delta= (G^\Delta,\Delta, (U^\Delta_\alpha,\phii^\Delta_\alpha)_{\alpha\in\phi})$ telle que, avec les notations évidentes:\enum{

\item $G\subset G^\Delta$ et $T\subset T^\Delta$.
\item  Pour tout $\alpha\in\phi$, $U_\alpha= U^\Delta_\alpha \cap G$.
\item Pour tout $\alpha\in\phi$, $\phii_\alpha= {\phii_\alpha^\Delta} |_{U_\alpha}$.
\item Pour tout $\alpha\in\phi$, $\phii^\Delta_\alpha(U_\alpha)$ est stable par addition avec $\Delta$.
%\item $P(\vec f)=P^\Delta(\vec f)$ pour toute facette $\vec f$ de $\vec A$.
}

Lorsqu'une donnée radicielle $\D=(G,(U_\alpha,\phii_\alpha)_\alpha)$ vérifie \pfonc,  on notera pour tout sous-groupe $\Delta$ de $\R$, $\D^\Delta$, $G^\Delta$, $U_\alpha^\Delta$, $N^\Delta$, $T^\Delta$, $\P^\Delta=(P^\Delta(a))_{a\in A}$, ... tous les objets obtenus grâce à la donnée radicielle valuée $\D^\Delta$.\\

On notera $A^\Delta$ l'appartement pour $\D^\Delta$, muni de ses murs, facettes et son action de $N^\Delta$, défini à partir de l'espace affine $\inte A$ et du point de base $o$ (c'est-à-dire les mêmes que pour $A$).

% Si $\Q$ est une famille de parahoriques pour $\maj D$, on dit que le couple $(\maj D,\maj Q)$ vérifie la condition \pfonc si $\maj D$ la vérifie et si pour tout $\Lambda'$, $\maj D'$ comme avant, il existe une famille de parahoriques $\Q'$ pour $\maj D'$ telle que pour tout $a\in A$, $Q(a)=Q'(a)\cap G$.
}

\fprop{
 Un groupe de Kac-Moody $G$ déployé sur un corps local $\K$ vérifie toujours la condition \pfonc.
}

\demo Lorsque $G$ est un groupe de Kac-Moody déployé sur un corps $\K$, il s'agit en fait de la valeur en $\K$ d'un foncteur $\maj G$ des $\K$-algèbres vers les groupes munis d'une donnée radicielle de système de racine fixé, et pour toute $\K$-algèbre $\K'$, les sous-groupes radiciels de $\maj G(\K')$ sont isomorphes à $(\K',+)$. \'Etant donné $\Delta$, il suffit donc de choisir une extension $\K^\Delta$ de $\K$ ramifiée de sorte que $\varpi(\K^\Delta) =\Delta$, puis de prendre $G^\Delta =\maj G(\K^\Delta)$. Le lemme 8.4.4 de \cite{remy} prouve que $G^\Delta$ est muni d'une donnée radicielle vérifiant le point 2 ci-dessus, les points 1 et 3 sont clairs.\cqfd

\rema Lorsque $G$ n'est que presque déployé, il s'agit encore d'un foncteur des $\K$-algèbres vers les groupes munis d'une donnée radicielle, mais le système de racine varie, et les groupes radiciels ne sont plus isomorphes au groupe additif du corps.\\

\fprops[\label{prop:csq de fonc} Soit $\D$ une donnée radicielle valuée vérifiant \pfonc, soit $\Delta$ un sous-groupe de $\R$. Alors:]
{
\item $T=T^\Delta\cap G$, et $N=N^\Delta\cap G$.

\item L'action de $N$ sur $\vec A$ est la restriction de l'action de $N^\Delta$, c'est-à-dire $\vec\nu=(\vec\nu^\Delta)|_N$. En particulier, pour toute facette $\vec f$ de $\vec A$, $G\cap N^\Delta(\vec f) = N(\vec f)$.

\item Pour tout $\vec f\in\F(\vec A)$, $P^\Delta(\vec f)\cap G=P(\vec f)$, puis $U(\vec f)=U^\Delta(\vec f)\cap G$, $M(\vec f)=M^\Delta(\vec f)\cap G$, et $N^\Delta.P^\Delta(\vec f) \cap G= N.P(\vec f)$.

\item L'immeuble $\iv=\iv(\D)$ s'injecte dans $\iv^\Delta = \iv(\D^\Delta)$, les appartements de ces deux immeubles sont isomorphes.

\item Pour tout $a\in A\sph$, $P(a)=P^\Delta(a)\cap G$.

\item $\nu = {\nu^\Delta}|_{N}$. En particulier, pour tout $a\in A$, $N(a)=N^\Delta(a)\cap G$.
%\item La famille minimale de parahoriques $\P^\Delta$ pour $\D^\Delta$ vérifie $_forall a\in A$, $P(a)= P^\Delta(a)\cap G$.
\item Lorsque $\Delta=\R$, tout point de $A$ est un sommet spécial dans $\I^\Delta$.
}

\demos{
\item On rappelle que par définition $T=\bigcap_{\alpha\in\phi}N_G(U_\alpha)$, où $N_G(U_\alpha)$ est le normalisateur de $U_\alpha$. Soit $g\in G\cap T^\Delta$, alors pour tout $\alpha\in\phi$, $g U_\alpha g\inv \subset G\cap U_\alpha^\Delta= U_\alpha$. Donc $g\in T$.

\'Etudions $N^\Delta\cap G$. Le groupe $N^\Delta$ est engendré par $T^\Delta$ et par un $n(u)$, pour un $u\in U_\alpha^\Delta$ pour chaque $\alpha\in\phi$. On peut choisir $u\in U_\alpha$, il vient alors $N^\Delta=N.T^\Delta$. Pour conclure, $N^\Delta\cap G= N.T^\Delta\cap G= N.(T^\Delta\cap G) = N.T=N$.\\

\item L'action de $N^\Delta$ sur $\vec A$ se fait via $W^\Delta= N^\Delta/T^\Delta$, et celle de $N$ via $W=N/T$. Mais $N^\Delta=N.T^\Delta$ et $T^\Delta\cap N=T$, donc $W^\Delta = N/T = W$.\\

\item 
%Pour commencer, soit $\vec C$ une chambre, montrons que $U^\Delta(\vec C) \cap G= U(\vec C)$. Soit $g\in U^\Delta(\vec C) \cap G$, et soit $g=u_1nu_2$ une écriture de $g$ correspondant à la décomposition de Bruhat $G=\bigsqcup_{n\in N} U(\vec C)nU(\vec C)$. Comme $N\subset N^\Delta$ et $U^\Delta(\vec C)\subset U(\vec C)$, c'est aussi une écriture de $g$ pour la décomposition de Bruhat de $G^\Delta$. Or $g=g\in U^\Delta(\vec C)$ est déjà une telle écriture. Par unicité du facteur dans $N^\Delta$, on en déduit $n=e$, d'où $g\in U(\vec C)$.

%Ensuite, si $\vec f$ est une facette de $\vec A$, alors $U(\vec f)$ est l'intersection des $U(\vec C)$ pour toutes les chambres $\vec C$ de $\vec F^*\cap \vec A$. Il est alors évident de prouver que $U(\vec f)^\Delta\cap G= U(\vec f)$.\\

 Soit $\vec f$ une facette de $\vec A$ et $g\in P^\Delta(\vec f)\cap G$. Soit $\vec C$ une chambre de $\vec f^*\cap \vec A$, on notera $U^+=U(\vec C)$ et $\phi^+=\phi(\vec C)$. On rappelle la décomposition de Bruhat fine:
 \[G=\bigsqcup_{n\in N} U^+.n.\left( U^+\cap n\inv U^-n\right) \point\]
Soit $g=u_1n u_2$ l'unique écriture  de $g$ selon cette décomposition.

La même décomposition pour $P^\Delta(\vec f)$ donne:
 \[P^\Delta(\vec f)= \bigsqcup_{n\in N^\Delta(\vec f)} U^{\Delta +} .n.\left( U^{\Delta +}\cap n\inv U^{\Delta -}n\right) \point\]

Soit $g=u^\Delta_1 n^\Delta u^\Delta_2$ la décomposition correspondante de $g$. Si $n\not= n^\Delta$, les cellules $U^{\Delta +} .n^\Delta.\left( U^{\Delta +}\cap n^{\Delta-1} U^{\Delta -}n^\Delta\right) \subset  U^{\Delta +} n^\Delta U^{\Delta +}$ et $U^+.n.\left( U^+\cap n\inv U^-n\right)\subset U^{\Delta +} n U^{\Delta +}$ sont disjointes.
 Donc $n=n^\Delta\in N\cap N^\Delta(\vec f)=N(\vec f)$. Ensuite, l'unicité de l'écriture dans la cellule $U^{\Delta +} .n.\left( U^{\Delta +}\cap n^{-1} U^{\Delta +}n\right)$ entraine $u_1=u_1^\Delta\in U(\vec C)$ et $u_2=u_2^\Delta\in U(\vec C)$.
Ceci prouve que $g\in U^+N(\vec f)U^+\subset P(\vec f)$.\\

Ensuite, l'unicité de l'écriture dans la décomposition de Lévi $P^\Delta(\vec f)= U^\Delta(\vec f)\rtimes M^\Delta(\vec f)$, et les inclusions $U(\vec f)\subset U^\Delta(\vec f)$ et $M(\vec f)\subset M^\Delta(\vec f)$ entrainent les deux égalités $U(\vec f)=U^\Delta(\vec f)\cap G$, $M(\vec f)=M^\Delta(\vec f)\cap G$.\\

Soit enfin $g\in G\cap N^\Delta.P^\Delta(\vec f)$. Comme $N$ et $N^\Delta$ induisent le même groupe de transformations sur $\vec A$ (c'est juste le groupe de Weyl associé au système de racines $\phi$), il existe $n\in N$ tel que $ng\in P^\Delta(\vec f)$. Alors $ng\in P^\Delta(\vec f)\cap G= P(\vec f)$.

\item Par construction, $\iv=G\times \vec A/{\vec \sim}$ où $\vec\sim$ est la relation d'équivalence $(g,a)\vec\sim (h,b) \ssi \exists n\in N \tq b=na\et g\inv bn\in P(\vec f)$, où $\vec f$ est la facette contenant $a$. Pour $\iv$ et $\iv^\Delta$, les deux appartements de référence sont les m\^emes, et le point précédent permet facilement de vérifier que $\iv$ s'injecte dans $\iv^\Delta$.

\item Soit $\vec f$ une facette sphérique, $a\in A_{\vec f}$ et $p\in P^\Delta(a)\cap G$. Alors $p\in U^\Delta(\vec f)\rtimes M^\Delta(a) \cap G = (U^\Delta(\vec f)\cap G)\rtimes (M^\Delta(a)\cap G)$ grâce au point 2 et à l'unicité dans la décomposition de Lévi de $P^\Delta(\vec f)$. Le groupe $M^\Delta(a)\cap G$ est le fixateur dans $G$ du point $a$ de l'immeuble de Bruhat-Tits de la donnée radicielle $(M^\Delta(\vec f), (U^\Delta_\alpha)_{\alpha\in\phi^m(\vec f)})$. D'après \cite{bruhat-tits} proposition 9.1.17, il s'agit de $M(a)$. D'autre part, il a déjà été vu que $U^\Delta(\vec f)\cap G = U(\vec f)$.

\item Comme le point de base $o$ est le même dans $A$ et $A^\Delta$, les action $\nu$ et $\nu^\Delta$ coïncident sur $N_o$. Il reste à étudier l'action de $T$. Rappelons que pour $t\in T$, $\nu(t)$ est par définition la translation de $Y$ de vecteur $\vec v_t$ tel que $\alpha(\vec v_t) =\phii_\alpha(u)-\phii_\alpha(tut\inv)$ pour tout $\alpha\in\phi$ et $u\in U_\alpha\setminus\{e\}$. L'égalité de $\nu(t)$ et de $\nu^\Delta(t)$ est alors conséquence du troisième point de la définition de la condition \pfonc.\\
%\item D'après \ref{soussection:la famille minimale}, on a $P^\Delta(a)= \left(U^\Delta(a)\cap P^\Delta(a)\right)\rtimes M^\Delta(a)$ où $\vec f$ est la direction de la façade de $a$.
\item Clair.
}

Lorsqu'une donnée radicielle valuée $\D$ vérifie \pfonc, il sera utile, étant donnée une famille de parahoriques $\Q$ pour $\D$ de savoir si $\Q$ se comporte bien vis-à-vis des extensions $\D^\Delta$ données par \pfonc. Ceci justifie d'introduire encore une définition:
\fdefs[Soit $Q$ une famille de parahoriques pour une donnée radicielle valuée $\D$.]{
\item Pour toute partie $\Omega$ de $A$, on dira que $Q$ vérifie la condition \pfonc($\Omega$) si $\D$ vérifie \pfonc\ et si pour tout $\Delta\leq \R$ il existe une famille de parahoriques $Q^\Delta$ pour $\D^\Delta$ telle que $\forall a\in \Omega$, $Q^\Delta(a)\cap G=Q(a)$.\\

\item On dira que $Q$ vérifie un ensemble de conditions (para $x_1,...,x_k$) "fonctoriellement" si $\D$ vérifie \pfonc\ et si pour tout groupe $\Delta\leq \R$, il existe une famille de parahoriques $ Q^\Delta$ pour $\D^\Delta$, contenant $Q$, et vérifiant (para $x_1,...,x_k$).\\

\item On dira que $Q$ vérifie \pfonc($\Omega$) et un ensemble de conditions (para $x_1,...,x_k$) "fonctoriellement" si $\D$ vérifie \pfonc\ et si pour tout groupe $\Delta\leq\R$, il existe une famille de parahoriques $Q^\Delta$ pour $\D^\Delta$, vérifiant (para $x_1,...,x_k$), et telle que $\forall a\in \Omega$, $Q^\Delta(a)\cap G=Q(a)$.\\

\item Les notations \pfonc\ et \pfonc(sph) désigneront respectivement \pfonc($A$) et \pfonc($A\sph$).
}

Nous avons vu par exemple que dès que $\D$ vérifie \pfonc, alors la famille minimale de parahoriques associée vérifie \pfonc(sph).
De plus, toutes les propriétés que nous prouverons être vérifiées par la famille minimale de parahoriques $P$ le seront en fait fonctoriellement. En particulier:
\fprop{\label{prop:fonc pour KM}Si $\D$ est la donnée radicielle valuée attachée à un groupe de Kac-Moody $G$, alors la famille minimale de parahoriques $P$ vérifie \pdec\ fonctoriellement.
}

\subsubsection{Relations directes entre les conditions (para x)}

Nous étudions les relations les plus directes entre les différentes conditions introduites en \ref{sousoussection:definition des parahoriques} et \ref{soussoussection:definition de fonc}. La philosophie est la suivante: les conditions \pisom\ et \pinter\ sont équivalentes aux propriétés d'incidence classiques attendues d'un immeuble. Elles ne sont pas vérifiées en général pour une donnée radicielle sur un système de racine infini, on tâchera cependant de déterminer des parties $\Omega$ pour lesquelles elles sont vraies (dans \cite{hovels} par exemple, on détermine des "parties avec un bon fixateur" qui en particulier vérifient \pisom\ et \pinter).

 Nous nous appuirons plutôt pour construire la masure sur les conditions \pinj, \psph\ et les variantes de \plien. Celles-ci ont déjà une interprétation géométrique, et ceci permettra de les descendre d'un groupe de Kac-Moody déployé à un groupe presque déployé.

Les conditions \pdec\ et \pfonc\ sont plutôt des intermédiaires techniques, vérifiés par les groupes de Kac-Moody déployés et qui entrainent les conditions précédentes. On ne s'en préoccupera a priori plus lors de l'étude d'un groupe presque déployé.\\

On étudiera les relations un peu moins directes entre ces conditions après avoir défini la masure $\I(\D, Q)$, car certaines implications sont plus facilement prouvées grâce à cet outil géométrique.\\

\fprop{\label{prop:pdec fonc implique plien} Soit $Q$ une famille de parahoriques vérifiant \pdec\  fonctoriellement. Alors pour tout point $a$ et toute facette sphérique $\vec f$ de $\vec f_a^*\cap \vec A$, $Q(a)\cap P(\vec f) = Q(a)\cap P(pr_{\vec f}(a))$.

En particulier, $\Q$ vérifie fonctoriellement \psph\ et \plien(sph).
}

\demo
Comme $\vec f$ est sphérique, on a la décomposition de Bruhat de $M_{\vec A}(\vec f)$, d'où:
 \calc {
P(\vec f) &= & U(\vec f)\rtimes M_{\vec A}(\vec f) \\
 &=& U(\vec f)\rtimes ( G(\phi^m(\vec f),a)\point N_{\vec f}\point G(\phi^m(\vec f),a) )\\
 & =& G(\phi^m(\vec f),a) \point  U(\vec f)\point N_{\vec f} \point G(\phi^m(\vec f),a)\\
   &\subset & P(\eng{a,\vec f}_A) \point  U(\vec f) \point N_{\vec f}\point  P(\eng{a,\vec f}_A)\\
}

Soit $g\in Q(a)\cap P(\vec f)$, nous pouvons donc supposer $g\in Q(a)\cap (U(\vec f)\rtimes N_{\vec f})$.\\

Dans un premier temps, supposons que $a$ est un sommet spécial. Alors quitte à multiplier par un élément de $N(a+\vec f)$, on peut supposer $g\in Q(a)\cap (U(\vec f)\rtimes T)$. Soit $u\in U(\vec f),t\in T$ tels que $g=u.t$. Par ailleurs, soit $\vec C$ une chambre de $\vec f^* \cap \vec A$, et soient $u^+\in U(\vec C)\cap Q(a)$, $u^-\in U(-\vec C)\cap Q(a)$ et $n\in N_a$ tels que $g=u^+u^-n$. Alors $g=u^+nn\inv u^-n=ut$, et par l'unicité modulo $T$ du facteur dans $N$ dans la décomposition de Birkhoff pour les chambres $\vec C$ et $-n\inv \vec C$, on obtient que $n.T=T$, donc $g$ est dans la double classe $U^+TU^-$. Enfin, par unicité d'écriture dans cette double classe, on obtient $u^+=u\in U(\vec f)\cap Q(a)$, $n=t\in T\cap N(a)$ et $n\inv u^- n =e$ d'où $u^-=e$.
L'élément $t\in N(a)\cap T$ fixe $a$, sa façade et son adhérence, et $u\in U(\vec f)\cap Q(a)$ fixe $a$ et $A_{\vec f}$ en entier. Donc $g\in Q(a)\cap P(pr_{\vec f}(a))$.\\

Lorsque $a$ n'est pas un sommet spécial de $A$, il l'est dans un certain $A^\Delta$ grâce à \pfonc. Le paragraphe précédent entraine alors que $Q(a)\cap P(\vec f) \subset Q(a)\cap P^\Delta(pr_{\vec f}(a))$. Mais $G\cap P^\Delta(pr_{\vec f}(a)) = P(pr_{\vec f}(a))$ d'après \ref{prop:csq de fonc}.\cqfd

\fprop{ \label{prop:plien implique pinj} Toute famille $Q$ de parahoriques vérifiant fonctoriellement \psph\ et \plien($\vec m$) pour toute cloison $\vec m$ de $\vec A$  vérifie \pinj.
}

\demo

Soit $a\in A$ et $g\in N\cap Q(a)$. Soit $\Delta$ tel que $a$ est un sommet spécial dans $A^\Delta$, soient $\D^\Delta=(G^\Delta,(U^\Delta_\alpha))$, $N^\Delta$, $T^\Delta$... les groupes donnés par la condition \pfonc.

Il existe $n\in N^\Delta(a)$ tel que $ng\in T^\Delta\cap n.Q(a) \subset T^\Delta\cap Q^\Delta(a)$. Pour toute cloison $\vec m$ de $\vec f_a^*\cap \vec A$, $\Q^\Delta$ vérifie \plien($\vec m$) d'où $ng\in T^\Delta\cap Q^\Delta(pr_{\vec m}(a))$. Comme $\vec m$ est sphérique, par \psph\ on obtient $ng\in T^\Delta\cap P^\Delta(pr_{\vec m}(a))$ ce qui vaut $T^\Delta\cap U^\Delta(\vec m)\rtimes (M^\Delta(\vec f)\cap P^\Delta(pr_{\vec m}(a))$ par la proposition \ref{prop:csq de parasph}. %Ici, $M^\Delta(pr_{\vec m}(a))$ est le sous-groupe parahorique au sens de Bruhat-Tits au point $pr_{\vec f}(a)$ pour la donnée radicielle finie $(M^\Delta(\vec m),(U_\alpha^\Delta,\phii^\Delta_\alpha)_{\alpha\in\phi^m(\vec m)})$.
 Comme $T^\Delta\subset M^\Delta(\vec m)$ on obtient $ng\in T^\Delta\cap M^\Delta(pr_{\vec m}(a))$ et par \cite{bruhat-tits}, ceci est le fixateur dans $M^\Delta(\vec m)$ de $A^\Delta_{\vec m}$, donc l'ensemble des $t\in T^\Delta$ induisant une translation de vecteur $\vec v_t\in\vect_{\vec A}(\vec m)$.

Ceci étant pour chaque cloison $\vec m$ de $\vec f_a^*\cap\vec A$, et comme l'intersection de $A^\Delta_{\vec f_a}$ et des murs contenant ces cloisons est triviale (car $\vec A_{\vec f_a}$ est essentiel), on voit que $ng$ induit une translation triviale sur $A^\Delta_{\vec f_a}$, autrement dit, $ng\in T^\Delta\cap N^\Delta(A^\Delta_{\vec f_a})$ d'où $g\in N^\Delta(a)\cap G$. Or ceci vaut $N(a)$ par \ref{prop:csq de fonc}.\cqfd

\rema L'hypothèse la plus précise pour cette proposition est en fait "$\Q$ vérifie fonctoriellement \psph\ et \plien($\vec f$) pour $\vec f$ dans une famille $\F$ de facettes sphériques de $\vec A$ telle que $\bigcap_{\vec f\in\F}\vect(\vec f)=\{0\}$. Par exemple la famille $\F$ des cloisons bordant une chambre donnée convient.
 En général, nous appliquerons ce résultat à des familles vérifiant \plien(sph), ce qui entraine bien \plien($\vec m$) pour toute cloison $\vec m$.\\

En assemblant les deux propositions précédentes, on trouve:

\fcor{\label{cor:csq de pdec pspsh et fonc} Toute famille vérifiant fonctoriellement \pdec\ vérifie fonctoriellement \pinj, \psph\ et \plien(sph).

C'est en particulier le cas, lorsque $\D$ est la donnée radicielle valuée issue d'un groupe de Kac-Moody $G$ déployé, pour la famille minimale de parahoriques $\P$.
}

\subsection{Définition de la masure bordée}

\subsubsection{La relation d'équivalence}
Soit $\Q=(Q(a))_{a\in A}$ une famille de sous groupes de $G$ vérifiant (para 0.2).

\begin{defin}
Soit $\sim_\Q$ la relation sur $G\times A$ définie par:
\[ (g,a)\sim_\Q (h,b) \ssi \exists n\in N(T) \tq b=na \et g\inv hn\in Q(a) \]
\end{defin}

\begin{prop}
 La relation $\sim_\Q$ est une relation d'équivalence si et seulement si $\Q$ vérifie (para 0.4).
\end{prop}

\demo\\
Déjà, $\sim_\Q$ est toujours réflexive.\\

Supposons que $\Q$ vérifie (para 0.4).\\
 Commençons par montrer que $\sim_\Q$ est symétrique.
Soient $(g,a)$ et $(h,b)$ tels que $(g,a)\sim_\Q (h,b)$. Soit $n\in N(T)$ tel que $b= na$ et $g\inv hn\in Q(a)$. Alors $a=n\inv b$ et $h\inv g n\inv= n (g\inv hn)\inv n\inv\in nQ(a)n\inv=Q(b)$.\\
Pour la transitivité, soient $(g,a)$, $(h,b)$ et $(k,c)$ tels que $(g,a)\sim_\Q (h,b)\sim_\Q (k,c)$. Soit $n\in N(T)$ tel que $b=na$ et $g\inv hn\in Q(a)$, et soit $m\in N(T)$ tel que $c= mb$ et $h\inv km\in Q(b)$.\\
 Alors $c=mna$ et $g\inv k mn=  (g\inv h n)n\inv (h\inv km)n\in Q(a). n\inv Q(b)n=Q(a)$.\\

Réciproquement, supposons $\sim_\Q$ une relation d'équivalence. Soit $n\in N(T)$ et $a\in A(T)$. soit $q\in Q(na)$, alors $(1,na)\sim_\Q (q\inv n,a)$. D'où par symétrie, $(q\inv n,a)\sim_\Q (1,na)$, donc il existe $n'\in N(T)$ tel que $na=n'a$ et $n\inv q n'\in Q(a)$. Alors $n\inv n'\in \fix_{N(T)}(a)\subset Q(a)$ par (para 0.2). D'où $n\inv q n\in Q(a)$. Nous avons montré que $n\inv Q(na) n\subset Q(a)$. L'inclusion inverse s'obtient en appliquant ce résultat à $n\inv$.\cqfd\\

\subsubsection{Définition}

On fixe un tore maximal $T_0$, et une famille de parahoriques $Q$ sur $A(T_0)$. On va construire l'objet immobilier $\I(Q)$ en se basant sur l'appartement $A(T_0)$, la construction sera bien sûr indépendante du choix de $T_0$.

\begin{defin}
Soit $\I(Q)= G\times A(T_0)/ \sim_Q$. C'est la masure bordée associée à $(\D,Q)$. Lorsque le contexte sera clair, on notera juste $g.a$ pour la classe de $(g,a)$ dans $\I(Q)$. Sinon on la notera $[g,a]_Q$.\\
On définit une action de $G$ sur $\I$ par $g'.[g,a]=[g'g,a]$.\\

Soit $\iota_{T_0,Q}:$ $\fonc {A(T_0)} {\I_Q} {a}{[1,a]_Q}$. C'est l'injection canonique de $A(T_0)$ dans $\I(Q)$. Son image est l'appartement de $\I(Q)$ associé à $T_0$. Les images de ce dernier par les éléments de $G$ sont les appartements de $\I(Q)$.\\
\end{defin}

\begin{prop}\label{prop:injection des apparts}
 La fonction $\iota_{T_0,Q}$ est $N(T_0)$-équivariante. Elle est de plus injective si et seulement si $Q$ vérifie \pinj.
\end{prop}

\demo\\
La $N(T_0)$-équivariance découle de la définition de $\sim_{Q}$.\\

 Si $Q$ vérifie \pinj, soient $a,b\in A(T_0)$ tels que $(1,a)\sim_Q (1,b)$. Alors il existe $n\in N(T_0)$ tel que $b=na$ et $n\in Q(a)$. Donc $n\in Q(a)\cap N(T_0)=\fix_{N(T_0)}(a)$, donc $a=b$.

Réciproquement, supposons $\iota_Q$ injective. L'inclusion $\fix_{N(T_0)}(a)\subset Q(a)\cap N(T_0)$ est vraie par (para 0.2). Pour l'autre inclusion, soit $q\in Q(a)\cap N(T_0)$. Alors $qa\in A(T_0)$ et $(1,a)\sim_Q (1,qa)$ d'où par injectivité de $\iota_Q$, $a=qa$ et $q\in \fix(a)$.\cqfd

On suppose désormais que $Q$ vérifie \pinj, et on identifie $A(T_0)$ à $\iota_{T_0,Q}(A(T_0))$.\\

\begin{prop}
Le stabilisateur de $A(T_0)$ dans $G$ est $N(T_0)$.
\end{prop}
\demo L'inclusion $N(T_0)\subset \stab_G(A(T_0))$ est vraie par la définition de $\sim_Q$ (ou par la $N(T_0)$-équivariance de l'inclusion de $A(T_0)$ dans $\I$).

 Réciproquement, soit $g\in G$ tel que $g.A(T_0)=A(T_0)$. Soit $\vec f\in\F(\vec A(T_0))$, soit $a\in A_{\vec f}$. Alors $g.a\in A(T)$ et par la définition de $\sim_\Q$, il existe $n\in N(T_0)$ tel que $g.a=n.a$. Donc $g\in n.P(a)\subset N(T_0).P(\vec f)$. Ceci étant valable pour toute facette $\vec f\in\F(\vec A(T_0))$, on obtient $g\in \bigcap_{\vec f\in\F(\vec A(T_0))} N(T_0).P(\vec f) = N(T_0).P(\vec A(T_0))=N(T_0).T_0=N(T_0)$ (\ref{props:paraboliques} pour la première égalité).\cqfd

Les deux propositions précédentes permettent de définir la structure des appartements de $\I(\Q)$:
\begin{defin}
La structure d'appartement sur $A(T_0)$ (\textit{i.e.} les murs et la structure affine sur chaque façade) est celle qui fait de $\iota_{T_0,\Q}$ un isomorphisme d'appartements. % On munit $A(T_0)$ de l'action de $N(T_0)$ qui fait de $\iota_{T_0,\Q}$ un isomorphisme $N(T_0)$-équivariant.\\
 Pour tout autre appartement $g.A(T_0)$, la structure d'appartement sur $g.A(T_0)$ est celle qui fait de $g|_{A_0}$ un isomorphisme d'appartements dont la partie vectorielle est l'application $\fonc{\vec V(T_0)}{g.\vec V(T_0)}{\vec v}{g.\vec v}$ (voir la fin de \ref{soussoussection:donnee radicielle}).% On munit $g.A(T_0)$ de l'action de $N(g T_0 g\inv)=gN(T_0)g\inv$ définie par $n.(g.a)=g.((g\inv n g).a)$ pour $n\in N(gT_0 g\inv)$ et $a\in A(T_0)$. Ces deux définitions sont indépendantes du choix de $g$ d'après la proposition précédente.

Si $A=g.A(T_0)$ est un appartement, l'addition d'un point de $A$ et d'un vecteur de $g.\vec V(T_0)$ sera notée $+_A$.
 
\end{defin}

Par exemple, si $A$ est un appartement, $a$ un point de $A$, $\vec v$ un vecteur de $\vec V(T)$, et $g\in G$, alors $g(a+_A\vec v)=g(a)\; +_{gA} \; g(\vec v)$. Si $n\in N(T)$, alors $n(a+_A\vec v)=n(a)+_A \vec\nu(n).\vec v$.\\

\fprop{ Soit $T$ un tore maximal de $G$. Soit $g\in G$ tel que $T=gT_0 g\inv$. Pour tout $\alpha\in\phi(T_0)$, on pose $g.\phii_\alpha:\fonc{U_{g\alpha}\setminus\{e\}} {\Lambda}{u}{\phii_\alpha(g\inv u g)}$, et $g.\phii_\alpha(e)=\infini$. Alors la famille $(g\phii_\alpha)_{\alpha\in\phi}$ est une valuation de la donnée radicielle $(G,(U_{\alpha})_{\alpha\in g.\phi})$, et l'appartement abstrait $A(T)$ qu'elle définit est isomorphe de manière $N(T)$-équivariante à l'appartement de $g.A(T_0)$ de $\I$.
}

\rema Le système de racine $g.\phi$ et la donnée radicielle $(G,(U_{\alpha})_{\alpha\in g.\phi})$ ne dépendent pas du choix de $g\in G$ tel que $T=gTg\inv$, contrairement à la valuation $g\phi$.\\ % Deux choix différents mèneront à deux valuations équipollentes, pour reprendre le terme de \cite{bruhat-tits}.\\

\demo Il est immédiat que $g.\phii$ est une valuation. Pour celle-ci, on a pour tout $\alpha\in\phi$ et $k\in\Lambda\cup\{\infini\}$, $g.U_{\alpha,k}g\inv= U_{g\alpha,k}$. En conséquence, si $o\in A(T)$ (resp. $o_0\in A(T_0)$) est le point de base pour $g\phii$ (resp. $\phii$), alors le groupe $N(T)_0=\engens{n(u)}{u\in U_{\alpha}, \,\alpha\in g\phi,\et g\phii_{g\inv\alpha}(u)=0}$ qui sert à définir l'action de $N(T)$ sur $A(T)$ dans \ref{subsubsection:action de N} est égal à $g.N(T_0)_0 g\inv$, et fixe le point $g.o_0$.

Alors l'application $\fonc {A(T)}{A(T_0)} {o+\vec v}{g.o_0 +_{T} \vec v}$ est un isomorphisme $N(T)$-équivariant. En effet, elle est clairement $N(T)_0$-équivariante, et concernant l'action de $T$, on vérifie directement avec la proposition \ref{prop:action de T} que pour tout $t\in T$, $\vec v_t=g(\vec v_{g\inv t g})$.\cqfd

Grâce à cette proposition, on identifie désormais pour tout $g\in G$ l'appartement $g.A(T_0)$ de $\I(Q)$  avec l'appartement abstrait $A(gT_0 g\inv)$.\\

Voici quelques propriétés immédiates de $\I$:
\begin{prop}\ \label{prop:immediates de I}
\begin{enumerate}
\item Le fixateur d'un point $x\in A(T_0)$ est $Q(x)$. Plus généralement, pour un appartement $B=g.A(T_0)=A(gT_0 g\inv)$, le fixateur d'un point $x\in B$ est $g Q(g\inv x)g\inv$.
\item Soit $T$ un tore maximal, $a\in A(T)$, $g\in G$. Si $g.a\in A(T)$, alors il existe $n\in N(T)$ tel que $g.a=n.a$.
\item Pour tout $x\in \I$, le groupe $\fix_G(x)$ est transitif sur les appartements contenant $x$.
\end{enumerate}
\end{prop}

\demo \begin{enumerate}
\item Si $(g,x)\sim_\Q x$, alors il existe $n\in N(T_0)$ tel que $nx=x$ et $g\inv n\in Q(x)$. Alors $n\in \fix_{N(T_0)}(x)\subset Q(x)$ par (para 0.2), et donc $g\in Q(x)$. La réciproque est claire, le cas général aussi.
\item Dans $A(T_0)$, c'est la définition de $\sim_\Q$. Le cas plus général s'y ramène car $hN(T_0)h\inv=N(hT_0 h\inv)$.

\item Soient $A$ et $g.A$ deux appartements contenant $x$. On peut supposer $A=A(T_0)$. Alors $g\inv x\in A(T_0)$, et par le point précédent, il existe $n\in N(T_0)$ tel que $g\inv x=nx$. Alors $gn\in \fix_G(x)$, et $g.A(T_0)=gnA(T_0)$.
\end{enumerate}
\cqfd\\

On prolonge naturellement la définition des sous-groupes parahoriques, de manière cohérente avec les notations $Q(x)$, $Q(\Omega)$ déjà introduites:
\begin{defin}
Pour tout $x\in\I$, on note $Q(x)=\fix_G(x)$. Pour toute partie $\Omega$ de $\I$, on note $Q(\Omega)=\fix_G(\Omega)$.
\end{defin}

\remas{
%\item L'action ainsi définie de $N(T)$ sur $A(T)$ coïncide avec l'action induite par celle de $G$ sur $\I$.\\

%\item On a ainsi défini $A(T)$, pour tout tore maximal $T$ de $G$, comme un appartement inclus dans $\I(\Q)$. Mais cet appartement est isomorphe à l'appartement $A(T)$ défini dans la partie précédente, via $g_{|A(T_0)}\circ \iota_{T_0,\Q}$ pour tout $g\in G$ tel que $T=gT_0 g\inv$.\\
\item Si $A$ est un appartement, si $a\in A$ et $\vec v\in\vec A$, alors le point $a+_A \vec v\in A$ est bien défini par la formule $a+_A \vec v= g(g\inv(a)+_{A(T_0)} g\inv(\vec v))$, où $g\in G$ est tel que $g.A(T_0)=A$. Mais ceci dépend a priori de l'appartement $A$ contenant $a$ et tel que $\vec v\in \vec A$ considéré. On ne peut donc pas noter ce point $a+\vec v$. En terme de condition sur la famille $\maj Q$, le point $a+\vec v$ est bien défini si $N.Q(g\inv a)\cap N.P(g\inv \vec v) \subset N.(Q(\{g\inv a,g\inv(a+\vec v)\}) \cap P(g\inv \vec v))$ (où l'addition est faite dans $A(T_0)$). Ceci est une conséquence de \plienpn, voir la section \ref{sousoussection:csq de plienpn}.\\
\item Il en va de même pour les projections: si $a\in A$ et $\vec f\in\F(\vec A)$, on peut noter $pr_{A,\vec f}(a)$ la projection de $a$ sur la façade $A_{\vec f}$ dans l'appartement $A$. Le projeté $pr_{\vec f}(a)$ est bien défini si $N.Q(a)\cap N.P(\vec f)\subset N.Q(\{a,pr_{\vec f}(a)\})$, en particulier si $Q$ vérifie \plienn($\vec f$). Voir \ref{sousoussection:csq de plienp}.\\
\item Pour l'enclos enfin, si $\Omega$ est une partie d'un appartement $A$, on notera $\cl_A(\Omega)$ l'enclos de $\Omega$ dans $A$. Ceci sera indépendant de $A$ si $Q$ vérifie \pisom($\Omega$) et \pinter($\Omega$).\\

}

 %Si $F$ est une facette de $\I$, on s'attend à ce que $Q(F)$ soit transitif sur les appartements contenant $F$. Pour l'instant, on peut prouver le résultat un peu plus faible suivant:
% \flemme{ \label{lemme:isom stabilisant une facette} Soit $A$ un appartement et $F=\gm_a(a+\vec f)$ une chambre de $A$. Soit $B$ un autre appartement contenant $F$. Alors il existe $g\in Q(a)$ tel que $g.A=B$ et $g.F=F$.
 %}
 %\pv
 %Pour commencer il existe $q\in Q(a)$ tel que $q.A=B$.

Enfin, on montre une caractérisation géométrique des sous-groupes radiciels valués $U_{\alpha,k}$.
\begin{prop}\label{prop:points fixes de ualpha}
Soit $T$ un tore maximal, $\alpha\in\phi(T)$, $u\in U_{\alpha}$.  Alors l'ensemble des points fixes de $u$ dans $A(T)$ est précisément $\D(\alpha,\phii_\alpha(u))$.
\end{prop}

\demo\\
Notons $k=\phii_\alpha(u)$. Par (para 0.3) et comme pour tout $a\in A$, $Q(a)=\fix(a)$, $U_{\alpha,k}$ fixe $\D(\alpha,k)$.\\
Supposons que $u$ fixe un point $x\in \D(-\alpha,-k)\setminus M(\alpha,k)$. Soient $u',u''\in U_{-\alpha}$ tels que $n(u)=u'uu''$. D'après \ref{prop:V5}, $\phii_{-\alpha}(u')=\phii_{-\alpha}(u'')=-k$ donc $u'$ et $u''$ fixent $x$. Ainsi $n(u)\in Q(x)$. Or $n(u)$ induit la réflexion selon le mur $M(\alpha,k)$ qui ne contient pas $x$. Donc $n(u)\in (Q(x)\cap N)\setminus N(x)$, ceci contredit \pinj.\cqfd\\

\begin{cor}
Pour $\alpha\in\phi(T)$ et $x\in A(T)$, $U_\alpha(x)=\ens{u\in U_\alpha}{u\text{ fixe } x}$.
\end{cor}

\begin{cor}\label{cor:conjugaison des u(alpha)(x)}
Soit $\alpha\in\phi(T)$, $x\in A(T)$, et $g\in G$. Alors $gU_\alpha(x)g\inv = U_{g\alpha}(g.x)$.
\end{cor}

\fprop{ \label{prop: projection entre immeubles} Soient $Q$ et $R$ deux familles de parahoriques avec $Q\subset R$. Alors il existe une projection naturelle $G$ équivariante $\phi:\I(Q) \twoheadrightarrow \I(R)$, $[g, a]_Q \mapsto [g,a]_{R}$. En particulier, la masure bordée $\I(P)$ correspondant à la famille minimale de parahoriques se projette sur toutes les autres masures bordées de $G$.
}

\demo 
 Si $(g,a)\sim_{Q} (h,b)$, alors il existe $n\in N$ tel que $b=na$ et $g\inv h n\in Q(a)$. Comme $Q(a)\subset R(a)$, ceci entraine $(g,a)\sim_R (h,b)$. Donc $\phi$ est bien définie. Le reste est évident.\cqfd

\subsubsection{Façades d'immeuble}

Rappelons que, par définition, les façades d'appartements dans $\I$ sont toutes les parties de la forme $g.(A_{0 \vec f})$, avec $g\in G$ et $\vec f\in\F(\vec A_0)$.

\flemme{ Soient $\mathfrak f=(hA_0)_{h\vec f}$ et $\mathfrak e=(gA_0)_{g\vec e}$ deux façades d'appartements. Si $\mathfrak f\cap \mathfrak e\not=\vide$, alors $h\vec f = g\vec e$.
}
\pv
On peut supposer $h=e$. Soit $a\in A_{0 \vec f}\cap (g A_0)_{g\vec e}$. Alors $g\inv a\in A_{0\vec e}$ donc il existe $n\in N$ tel que $gn\in Q(a)$. En particulier, $gn\in P(\vec f)$ donc $g\inv \vec f=n\vec f$, et d'autre part $g\inv a=na$ d'où $n.\vec f=\vec e$. Au final, $g\inv\vec f = \vec e$.\cqfd

En conséquence de ceci, la direction d'une façade d'appartement dans $\I$ est bien définie:
\fdef{ Soit $\mathfrak f = g.(A_{0 \vec f})$ une façade d'appartement, alors la facette vectorielle $g.\vec f$ est appelée la direction de $\mathfrak f$.
La réunion de toutes les façades d'appartement dirigées par une facette vectorielle $\vec f$ est appelée la façade de $\I$ de direction $\vec f$. On la note $\I_{\vec f}$.
}

\fprop{ Les façades de $\I$ forment une partition de $\I$. De plus, l'application $\vec f\mapsto \I_{\vec f}$ est une bijection $G$-équivariante entre les facettes de $\iv$ et les façades de $\I$. En conséquence, le stabilisateur dans $G$ de la façade $\I_{\vec f}$ est $P(\vec f)$.
}
\demo Le lemme prouve que deux façades sont disjointes ou égales. De plus, $A_0$ est la réunion de ses façades, donc $\I=G.A_0$ est la réunion de ses façades d'appartement, et donc de ses façades d'immeuble.

L'application $\vec f\mapsto \I_{\vec f}$ est clairement surjective et $G$-équivariante. Si $\I_{\vec f}=\I_{\vec g}$, soit $\mathfrak f$ une façade d'appartement de direction $\vec f$ et $a\in\mathfrak f$. Par hypothèse, il existe $\mathfrak g$ de direction $\vec g$ contenant $a$. Le lemme entraine alors $\vec f=\vec g$.\cqfd

\fprop{Si $\vec f$ est une facette sphérique de $\iv$, et si $Q$ vérifie \psph, la façade $\I_{\vec f}$ est l'immeuble de Bruhat-Tits du groupe $M_{\vec A}(\vec f)$, pour tout appartement $\vec A$ contenant $\vec f$.

Pour une facette $\vec f$ générale, la façade fermée $\barre{\I_{\vec f}} : =\bigcup_{\vec g\in \vec f^*}\I_{\vec g}$ est la masure bordée pour la donnée radicielle valuée $(M(\vec f), (U_\alpha,\phii_\alpha)_{\alpha\in\phi^m(\vec f)})$ avec la famille de parahoriques $Q$ restreinte à $\barre{A_{\vec f}}$.
}
%\rema Le groupe $P(\vec f)= U(\vec f)\rtimes M(\vec f)$, et pas seulement $M(\vec f)$, agit sur $\I_{\vec f}$. Dans le cas où $\vec f$ est sphérique, le facteur $U(\vec f)$ agit trivialement, d'après \ref{prop:csq de parasph}. Dans le cas général, une condition pour que $U(\vec f)$ agisse trivialement sur $A_{\vec f}$ est par exemple que $Q$ vérifie \plienp($\vec C$) pour toute chambre $\vec C$.

\demo Soit $\vec f$ une facette de $\iv$. Soit $\g=B_{\vec f}$ une façade d'appartement de direction $\vec f$. Il existe $p\in P(\vec f)$ tel que $B=p.A$ et donc $\g=p.A_{\vec f}$. Ainsi, $\I_{\vec f}=P(\vec f).A_{\vec f}$. Mais comme $U(\vec f)$ fixe $A_{\vec f}$, on a $\I_{\vec f} = M_{\vec A}(\vec f).A_{\vec f}$.

Supposons $\vec f$ sphérique. Alors $A_{\vec f}$ est un appartement pour la donnée radicielle valuée finie $\D_{\vec A, \vec f}$. Le fait que pour tout $a\in A_{\vec f}$, $M_{\vec A}(\vec f)\cap P(a)$ soit le sous-groupe parahorique pour $\D_{\vec A,\vec f}$ au point $a$ permet de vérifier immédiatement que $\I_{\vec f}$ est l'immeuble de Bruhat-Tits de $\D_{\vec A,\vec f}$.

Dans le cas général, $\barre{A_{\vec f}} = \bigcup_{\vec g\in \vec f^*\cap \vec A} A_{\vec g}$ est un appartement pour $\D_{\vec A,\vec f}$. De plus $\barre{\I_{\vec f}} = M_{\vec A}(\vec f).\barre{A_{\vec f}}$, et on vérifie directement sur la définition que ceci est la masure de $\D_{\vec A,\vec f}$ pour la famille de parahoriques $(M_{\vec A}(\vec f)\cap Q(a))_{a\in \barre{A_{\vec f}}}$.\cqfd

\subsection{Décomposition d'Iwasawa}
\label{soussection:Iwasawa}
 On prouve ici la décomposition d'Iwasawa $G = P(\vec C) \point N(T)  \point G(F)$, valable pour toute chambre $\vec C$ de $\vec A(T)$ et pour toute facette $F\subset A$. Rappelons que selon nos notations, le groupe $G(F)$ est défini par $G(F) =\engens{U_\alpha(F)}{\alpha\in\phi(T)}$, avec $U_\alpha(F):=\ens{u\in U_\alpha}{ F\subset \D(\alpha,\phii_\alpha(u)) }$. Pour toute famille de parahorique $\Q$, on a $G(F)\subset Q(F)$, donc la décomposition d'Iwasawa implique $G=P(\vec f) \point N(T) \point Q(F)$, pour toute facette $\vec f$ de $\vec A(T)$ et $F$ de $A(T)$.\\

\subsubsection{La décomposition}

 \flemme{\label{lemme:une racine a part}Soit $\vec C$ une chambre de $\vec A(T)$, soit $\alpha\in\phi(T)$ qui s'annule sur une cloison $\vec m$ de $\vec C$. Autrement dit, $\alpha$ est une racine simple de $\phi(\vec C)$. Alors :
 \[U(\vec C) = U(\vec m)\rtimes U_\alpha
 %\subset P(\vec m^*)\rtimes U_\alpha \subset P(\vec C)
 \]
 }
 
\rema On rappelle que pour une chambre $\vec C$, $U(\vec C)= G(\phi(\vec C))$. \\
 
\pv\\
 Le groupe $U_\alpha$ fixe la cloison $\vec m$ donc normalise $U(\vec m)$. Soit $u\in U_\alpha \cap U(\vec m)$, alors $u$ fixe les deux chambres de $\vec A$ qui bordent $\vec m$, et $u\in U_\alpha$, ceci entraine $u=e$. Ainsi le groupe engendré par $U_\alpha$ et $U(\vec m)$ est bien un produit semi-direct. Il est inclus dans $U(\vec C)$ car $U_\alpha$ tout comme $U(\vec m)$ le sont ($U(\vec m) = U(\vec C)\cap U(r_{\vec m}.\vec C)$ si $r_{\vec m}$ est la réflexion dans $\vec A$ par rapport à $\vec m$). Enfin, pour toute racine $\beta\in \phi(\vec C)$, on a soit $\beta=\alpha$, soit $\beta\in\phi^u(\vec m)$, et dans les deux cas $U_\beta\subset U(\vec m)\rtimes U_\alpha$. D'où le résultat.\cqfd

\fprop{\label{prop:Iwasawa} (Décomposition d'Iwasawa) Soit $\vec C$ une chambre de $\vec A(T)$ et $F$ une facette de $A(T)$. Alors:
 \[G=U(\vec C)\point N(T) \point G(F) \, .\]
}
\demo\\
La preuve est classique. Il suffit de prouver que l'ensemble $Z:=U(\vec C)\point N(T) \point G(F)$ est stable par multiplication à gauche par n'importe quel élément de $G$. Or $G$ est engendré par $T$, par les $U_\alpha$ avec $\alpha\in\phi(\vec C)$ et par les $U_\beta$ avec $\beta$ une racine simple de $\phi(-\vec C)$.  Déjà, $Z$ est clairement stable par multiplication à gauche par $T$ et les $U_\alpha$. Soit donc $\beta=-\alpha$ une racine simple de $\phi(-\vec C)$, montrons que $Z$ est stable par multiplication à gauche par $U_\beta$.\\

 Par le lemme précédent, $U(\vec C)= U(\vec m)\rtimes U_\alpha$, avec $\vec m $ la cloison de $\vec C$ qui est incluse dans le noyau de $\beta$. Le groupe $U(\vec m)$ est normalisé par $U_\beta$, et $U_\beta.Z\subset U(\vec m).U_\beta.U_\alpha.N(T) .G(F)$. L'ensemble $\{\alpha,\beta\}=\{\alpha,-\alpha\}$ est un système de racine fini, donc par \cite{bruhat-tits}, $U_\beta.U_\alpha\subset \eng{U_\beta,U_\alpha,T} = U_\alpha.T\{e,r_\alpha\} U_\alpha = U_\alpha.T\sqcup U_\alpha U_{-\alpha}Tr_\alpha\subset U_\alpha U_{-\alpha} N$, où $r_\alpha=n(u)$ est la réflexion générée par un élément quelconque $u$ de $U_\alpha$.\\
 Ainsi, $U_\beta.Z \subset U(\vec m)\point U_\alpha\point U_{-\alpha} \point N(T) \point G(F) = U(\vec C) \point U_{-\alpha} \point N(T) \point G(F)$.\\
 
 Étudions maintenant le produit $U_{-\alpha}.N(T).G(F)$. Pour tout $n\in N(T)$, $U_{-\alpha}.n=nn\inv U_{-\alpha}n = nU_{-n\inv.\alpha}$. Dans la donnée radicielle de type finie $(\eng{U_{-n\inv.\alpha}, U_{n\inv.\alpha},T}, ( (U_{-n\inv.\alpha}, \phii_{-n\inv.\alpha},(U_{n\inv.\alpha}, \phii_{n\inv.\alpha})) ) $, en utilisant la décomposition d'Iwasawa ou de Bruhat selon que $n\inv\alpha(F)$ est fini ou  non, il vient $U_{-n\inv.\alpha}\subset U_{n\inv\alpha}.N(T).G(F)$. D'où $U_{-\alpha}.N(T).G(F)\subset U_\alpha.N(T).G(F)$.\\
 
 On a alors obtenu: $U_\beta.Z\subset U(\vec m).U_\alpha\point N(T)\point G(F) \subset U(\vec C)\point N(T)\point G(F)$.\cqfd\\

\fcor{\label{cor:Iwasawa} Pour toutes facettes $\vec f\subset \vec A(T)$ et $F\subset A(T)$, pour toute famille $\Q$ de parahoriques sur $A(T)$, on a 
\[G= P(\vec f).N(T) .Q(F) \]
}
\fcor{ \label{cor:facette et facette vectorielle dans un appart} 
Soit $F$ une facette de $\I$ et $\vec f$ une facette de $\iv$. Pour tout appartement $A$ contenant $F$, il existe $q\in Q(F)$ tel que $\vec f\subset q.\vec A$. En particulier, il existe un appartement contenant $F$ dont l'appartement directeur contient $\vec f$.
%Pour toutes facettes $\vec f$ de $\iv$ et $F$ de $\I(\Q)$, il existe un appartement $A$ tel que $F\subset A$ et $\vec f\subset \vec A$. %De plus, entre deux tels appartements, il existe un isomorphisme induit par un élément $q\in Q(F)$ qui fixe $F$.
}
\demo
 Soit $g\in G$ tel que $\vec f\subset g.\vec A$. Soit $T$ le tore maximal tel que $A=A(T)$. Par le corollaire précédent, $g\in Q(F)N(T) P(g\inv \vec f)$. Soit $g=qnp$ une écriture de $g$ correspondante. Alors $\vec f\subset q.\vec A$.\cqfd

\subsubsection{Unicité}
\label{soussoussection:unicie dans iwasawa}
Lorsque $Q$ vérifie \pdec\ pour une facette $F$, on prouve un résultat d'unicité pour le facteur dans $N$ pour toute décomposition d'Iwasawa faisant intervenir $F$.

 \fprop{\label{prop:unicite dans iwasawa} Soient $\vec C$ une chambre de $\vec A$ et $F$ une facette de $A$, $F'\subset A$ un élément du filtre $F$. Soient $n,n'\in N$ tels que $U(\vec C).n.G(F') \cap U(\vec C).n'.G(F')\not=\vide$. On suppose de plus que $\vec f\subset n\barre{\vec C}$ ou $\vec f\subset n'\barre{\vec C}$, ou $\vec f$ est la direction de la façade de $F$.
 
 Pour tout $a\in F'$, si $Q$ vérifie \pdec($a$), alors $n\inv n'\in N(a)$. Autrement dit, le facteur dans $N$ dans la décomposition de Lévi est unique modulo $N(a)$.
 
 En particulier, si $Q$ vérifie \pdec($a$) pour tout $a\in F'$, alors $n\inv n'\in N(F)$.
 }

 \demo\\
Soit $a\in F'$ tel que $Q$ vérifie \pdec($a$).
 On a $n'\in U(\vec C).n.Q(F')$ , d'où $n\inv n'\in U(n\inv \vec C)\point Q(F')\subset  U(n\inv \vec C)\point Q(a)= U(n\inv \vec C)\point (U(-n\inv \vec C)\cap Q(a))\point N(a)$ par \pdec($a$). Donc il existe $n_a\in N(a)$ tel que $n\inv n' n_a\in U(n\inv \vec C)\point U(-n\inv \vec C)$. Par unicité du facteur de $N$ dans la décomposition de Birkhoff vectorielle, on obtient $n\inv n' n_a=e$.\cqfd
 
 En fait, l'unicité modulo $N(a)$ d'un facteur $n\in N$ dans la décomposition d'Iwasawa $G=U(\vec C).N.Q(a)$, avec $\vec f\subset n\barre{\vec C}$, équivaut à l'égalité: $N\cap U(n\inv\vec C).Q(a)= N(a)$. Nous avons juste vérifié que cette dernière est conséquence de \pdec($a$). Nous verrons plus loin (\ref{prop:unicite dans iwasawa 2}) qu'elle est également vraie lorsque $Q$ est une bonne famille de parahoriques.\\

\fcor{\label{cor:unicite dans iwasawa} Soit $a\in A$. Soit $\Q$ une famille de parahoriques vérifiant \pdec($a$) fonctoriellement.
 Alors pour tout sous-groupe $\Delta$ de $\R$,   $G\cap N^\Delta(\vec f_a).Q^\Delta(a) = N(\vec f_a).(G\cap Q^\Delta(a))$.
 
 Et donc si $\Q$ vérifie en outre \pfonc($a$), alors $G\cap N^\Delta.Q^\Delta(a) = N.Q(a)$.
}
\demo
Soit $g=u.n.p$ une écriture de $g$ dans la décomposition d'Iwasawa $G=U(\vec C)\point N\point Q(a)$, pour $\vec C$ une chambre de $\vec f_a^*\cap \vec A$. Comme $P(a)\subset P^\Delta(a)$, $N\subset N^\Delta$, et $U(\vec C)\subset U^\Delta(\vec C)$, c'est aussi une écriture de $g$ dans $G^\Delta=U^\Delta(\vec C)\point N^\Delta\point Q^\Delta(a)$.

 Par ailleurs, soit $g=n^\Delta.q^\Delta$ une écriture venant de l'hypothèse $g\in N^\Delta(\vec f_a).Q^\Delta(a)$. Par la proposition précédente, $n.N^\Delta(a) = n^\Delta.N^\Delta(a)$. Donc $g\in n.Q^\Delta(a)$, puis $g\in n.(Q^\Delta(a)\cap G)$, ce qui vaut $n.Q(a)$ si \pfonc($a$) est vrai.\cqfd
 
\rema Le même raisonnement permettra la même conclusion lorsqu'on étudiera un groupe presque déployé par descente galoisienne.\\

\subsection{Décomposition de Bruhat/Birkhoff}

\fprop{ Soient $F_1$ et $F_2$ deux facettes d'un appartement $A(T)$. On suppose qu'au moins une des deux est sphérique. Alors $G=U(\vec f_1)G(F_1) N(T) G(F_2)U(\vec f_2)$, où $\vec f_1$, respectivement $\vec f_2$, est la direction de la façade de $F_1$, respectivement $F_2$.
}

\rema En fait, l'hypothèse minimale sur les facettes $F_1$ et $F_2$ pour que la décomposition soit vraie, et prouvée par la preuve à suivre est: si $\vec f_1$ et $\vec f_2$ sont les directions des façades de $F_1$ et $F_2$, alors pour tout $w\in W(\vec A(T))$, $\phi^m(\vec f_1\cup w\vec f_2)$ est fini.\\

\demo

 On fixe $g\in G$. Pour cette preuve, l'appartement par défaut est $A(T)$, c'est-à-dire qu'on notera $N$ pour $N(T)$, $M(f)$ pour $M_A(f)$... .\\

Par la décomposition d'Iwasawa de $G$, $g\in P(\vec f_1)\point N\point G(F_2) = U(\vec f_1).M(\vec f_1)\point N\point G(F_2)$. Comme $U(\vec f_2)$ est normalisé par $G(F_2)$, on peut supposer qu'il existe $n\in N$ tel que $g\in M(\vec f_1)\point n$.\\

On utilise alors la décomposition d'Iwasawa dans $M(\vec f_1)$, avec la facette affine $F_1$ et la facette vectorielle $pr_{\vec f_1}(n\vec f_2)$. Le sous-groupe parabolique de $M(\vec f_1)$ fixant cette facette vectorielle est le groupe engendré par les $U_\alpha$ avec $\alpha\in\phi^m(\vec f_1)\cap \phi(n\vec f_2) = \phi(\vec \Omega)$, en notant $\vec \Omega=\text{conv}(\vec f_1-\vec f_1+n\vec f_2)$.
 Donc $M(\vec f_1)=(G(F_1)\cap M(\vec f_1))\point N(\vec f_1)\point P(\vec \Omega)$. Notons que $\vec\Omega$ est équilibrée, car les deux facettes $pr_{\vec f_1}(\vec f_2)$ et $pr_{-\vec f_1}(\vec f_2)$ sont sphériques, incluses dans $\vec\Omega$, et de signes opposés. Donc $P(\vec \Omega)= M(\vec\Omega)\ltimes U(\vec\Omega)\subset M(\vec\Omega)\ltimes U(n\vec f_2)$. Ainsi, 
 \[g\in G(F_1)\point N(\vec f_1) \point M(\vec\Omega) \point U(n\vec f_2). n = G(F_1)\point N(\vec f_1) \point M(\vec\Omega) \point n.U(\vec f_2)\]
On peut donc supposer qu'il existe $n_1\in N(\vec f_1)$ tel que:
\[g\in n_1. M(\vec\Omega).n \]

Les facettes $nF_2$ et $n_1\inv F_1$ se projettent sur $A_{\vec\Omega}$. Et par la décomposition de Bruhat dans le groupe muni d'une donnée radicielle valuée finie $M(\vec\Omega)$, on a $M(\vec\Omega) = G(\phi^m(\vec\Omega), pr_{\vec\Omega}(n_1\inv F_1))\point N(\vec\Omega)\point G(\phi^m(\vec\Omega), pr_{\vec\Omega}(nF_2))$. Enfin, sachant que pour toute partie $F$ de $A_{\Omega}$, $G(\phi^m(\vec\Omega),F) \subset G(\eng{F,\vect(\vec\Omega)}_A)$,
 on arrive à \[M(\vec\Omega) \subset G(n_1\inv F_1)\point N(\vec\Omega)\point G(nF_2)\]
 , puis:
\[ g \in n_1 G(n_1\inv F_1)\point N(\vec\Omega)\point G(nF_2) n = G(F_1)\point n_1N(\vec\Omega)n\point G(F_2) \subset G(F_1)\point N\point G(F_2)\]
\cqfd

\fcor{Pour toute famille $Q$ de parahoriques, pour toutes facettes $F_1$ et $F_2$ d'un appartement $A(T)$, si l'une des deux est sphérique, alors:
\[G= Q(F_1)\point N(T)\point Q(F_2) \]
}

\fcor{\label{cor:appart contenant deux facettes}Pour toute famille $Q$ de parahoriques, pour toutes facettes $F_1$ et $F_2$ de $\I(\Q)$, si l'une des deux est sphérique, alors il existe un appartement contenant $F_1\cup F_2$.
}
\demo Soit $A(T)$ un appartement contenant $F_1$, soit $g\in G$ tel que $F_2\subset g.A(T)$. Par la décomposition de Bruhat/Birkhoff, $G=Q(F_1)\point N(T)\point Q(g\inv F_2)$. Soit $g=q_1nq_2$ l'écriture correspondante de $g$. Alors $F_1\cup F_2\subset q_1.A(T)$.\cqfd

\subsection{Construction de familles vérifiant \pdec }

\flemme{ Soit $\vec C$ une chambre de $\vec A$, $\vec m$ une cloison de $\vec C$, et $a\in A$ tel que $\vec f_a\subset \barre{\vec m}$. Soit $\alpha\in\phi$ la racine telle que $\alpha(\vec C)>0$ et $\alpha(\vec m)=0$. Pour toute famille $Q$ de parahoriques vérifiant \psph\ et \plien($\vec m$) on a :
\[ Q(a) \cap U(\vec C) = \left( Q(a) \cap U(\vec m)\right) \rtimes U_\alpha(a)\point \]
}
\pv
Soit $g\in Q(a) \cap U(\vec C)$.
On sait par le lemme \ref{lemme:une racine a part} que $U(\vec C) = U(\vec m) \rtimes U_\alpha$, il existe donc une décomposition $g=u.u_\alpha$ avec $u\in U(\vec m)$ et $u_\alpha\in U_\alpha$. Comme $U(\vec C)\subset P(\vec m)$, on a $g\in Q(a)\cap P(\vec m) = Q(\{a,pr_{\vec m}(a)\})$ par \plien($\vec m$). Le facteur $u\in U(\vec m)$ fixe toute la façade $A_{\vec m}$, on en déduit que $u_\alpha\in U_\alpha\cap P(pr_{\vec m}(a))$, c'est-à-dire $u_\alpha\in U_\alpha(a)$. En particulier, $u_\alpha$ fixe $a$, donc $u$ aussi: $u\in U(\vec m)\cap Q(a)$.\cqfd

\fprop{\label{lemme:obtenir pdec} Soit $\epsilon$ un signe, soit $Q$ une famille de parahoriques vérifiant \psph\ et \plien($\vec m$) pour toute cloison $\vec m$ de signe $\epsilon$. On suppose que la famille minimale $P$ vérifie aussi $\plien(\vec m)$ pour toute cloison $\vec m$. Alors pour tout $a\in A^\epsilon$ et toute chambre $\vec C$ de $\vec f_a^*\cap \vec A^\epsilon$, l'ensemble:
 \[R(a)=\left(Q(a)\cap U(\vec C)\right) \point \left(P(a)\cap U(-\vec C)\right) \point N(a) \]
 est un sous-groupe de $Q(a)$ contenant $P(a)$. Il est en outre indépendant de la chambre $\vec C$ de $\vec f_a^*\cap \vec A^\epsilon$ choisie.\\
 
 Si on définit en outre $R(a)=Q(a)$ pour tout $a\in A\setminus A^\epsilon$, on obtient une famille de parahoriques $R$ incluse dans $Q$ et vérifiant \pdec\ pour les chambres de signe $\epsilon$.
 }

\demo 

Supposons $\epsilon=+$. Soit $a\in A$ et $\vec C$ une chambre de $\vec A^+$, notons $R_{\vec C}(a)=\left(Q(a)\cap U(\vec C)\right) \point \left(P(a)\cap U(-\vec C)\right) \point N(a)$.

 Commençons par montrer que $R_{\vec C}(a)$ est indépendant de la chambre $\vec C$ de $\vec f_a^*$: il suffit de prouver que $R_{\vec C}(a) = R_{r_\alpha \vec C}(a)$ pour toute racine simple $\alpha$ de $\phi(\vec C)\cap \phi^m(\vec f_a)$. Soit $\alpha$ une telle racine, et $\vec m_\alpha$ la cloison de $\vec C$ correspondante. D'après le lemme précédent, $Q(a)\cap U(\vec C)=(Q(a)\cap U(\vec m_\alpha))\point U_\alpha(a)$ et de même $P(a)\cap U(-\vec C) = (P(a)\cap U(-\vec m_\alpha)) \rtimes U_{-\alpha}(a)$. De plus, $U_\alpha(a)$ normalise $P(a)\cap U(-\vec m_\alpha)$. Ainsi:
 \calc{
 R_{\vec C}(a) &=& (Q(a)\cap U(\vec m_\alpha))\point U_\alpha(a)\point (P(a)\cap U(-\vec m_\alpha))\point U_{-\alpha}(a)\point N(a) \\
 &=& (Q(a)\cap U(\vec m_\alpha))\point (P(a)\cap U(-\vec m_\alpha))\point U_\alpha(a).U_{-\alpha}(a)\point N(a)
 }

L'ensemble $U_\alpha(a).U_{-\alpha}(a)$ est inclus dans le groupe avec donnée radicielle valuée finie $M(\vec m_\alpha)$. Il est même inclus dans le fixateur $\fix_{M(\vec m_\alpha)}(pr_{\vec m_{\alpha}}(a))$ du point $pr_{\vec m_{\alpha}}(a)$ de  $\I_{\vec m_\alpha}$, l'immeuble de Bruhat-Tits de $M(\vec m_\alpha)$. Or ce fixateur est égal à $U_\alpha(a).U_{-\alpha}(a).(N(a)\cap M(\vec m_\alpha)) = U_{-\alpha}(a).U_{\alpha}(a).(N(a)\cap M(\vec m_\alpha))$ d'après \cite{bruhat-tits}. Donc $U_\alpha(a).U_{-\alpha}(a) \subset U_{-\alpha}(a).U_{\alpha}(a).N(a)$ et finalement:
\calc{
 R_{\vec C}(a) &=&(Q(a)\cap U(\vec m_\alpha))\point (P(a)\cap U(-\vec m_\alpha))\point U_{-\alpha}(a).U_{\alpha}(a).N(a)\\
 &=& (Q(a)\cap U(\vec m_\alpha)).U_{-\alpha}(a) \point (P(a)\cap U(-\vec m_\alpha)).U_\alpha(a)\point N(a)
}

Or $\phi(r_\alpha \vec C) = (\phi(\vec C)\setminus\{\alpha\})\cup\{-\alpha\}$ et similairement pour $\phi(-r_\alpha\vec C)$, donc le lemme précédent indique que $(Q(a)\cap U(\vec m_\alpha)).U_{-\alpha}(a) = Q(a)\cap U(r_\alpha\vec C)$ et $(P(a)\cap U(-\vec m_\alpha)).U_\alpha(a) = P(a)\cap U(r_\alpha \vec C)$. Nous avons bien prouvé que $R_{\vec C}(a) = R_{r_\alpha \vec C}(a)$.\\

On peut maintenant noter $R(a)$ au lieu de $R_{\vec C}(a)$. Pour montrer qu'il s'agit d'un groupe, il suffit de montrer qu'il est stable par multiplication à gauche par $G(a)$ et par $N(a)$. Commençons par la stabilité par multiplication par $G(a)$. Comme $G(a)=U(\vec f_a).\engens{ U_\alpha(a)}{\alpha\in\phi^m(\vec f_a)}$, il suffit de voir la stabilité sous la multiplication par $U(\vec f_a)$ et par les $U_\alpha(a)$, $\alpha\in\phi^m(\vec f_a)$. Pour tout $\alpha\in\phi^m(\vec f_a)$, il existe une chambre $\vec C\subset \vec f_a^*$ positive telle que $\alpha\in\phi(\vec C)$. Alors $U_\alpha(a)\subset Q(a)\cap U(\vec C)$ d'où $U_\alpha(a).R(a)\subset R(a)$. Passons à $U(\vec f_a)$: pour n'importe quelle chambre $\vec C$ de $\vec f_a^*$, on a $U(\vec f_a)\subset Q(a)\cap U(\vec C)$, d'où $U(\vec f_a).R(a)\subset R(a)$.

% Si $\vec f_a$ est négative, il existe une chambre $\vec C$ positive telle que $U(\vec f_a)\subset P(a)\cap U(-\vec C)$, et $Q(a)\subset P(\vec f_a)$ normalise $U(\vec f_a)$. D'où $U(\vec f_a).(Q(a)\cap U(\vec C)).(P(a)\cap U(\vec C)) = (Q(a)\cap U(\vec C)).U(\vec f_a).(P(a)\cap U(\vec C))=(Q(a)\cap U(\vec C)).(P(a)\cap U(\vec C))$.

 Pour finir, soit $n\in N(a)$, fixons $\vec C$ une chambre positive. Alors $n.R(a) = n.R_{\vec C}(a) = n(Q(a)\cap U(\vec C))\point (P(a)\cap U(-\vec C))\point N(a) = (Q(a)\cap U(n\vec C))\point (P(a)\cap U(-n\vec C)) \point nN(a) = R_{n\vec C}(a) = R(a)$.\\

Ceci prouve que $R(a)$ est un groupe. Il est clair que $P(a)\subset R(a) \subset Q(a)$, et ceci entraine immédiatement que $(R(a))_{a\in A}$ est une famille de parahoriques. De plus, pour toute chambre $\vec C$ positive, $Q(a)\cap U(\vec C) = R(a)\cap U(\vec C)$ et $P(a)\cap U(-\vec C)\subset R(a)\cap U(\vec C)$ donc $R$ vérifie \pdec($\vec C$).\cqfd

\rema Ceci et la proposition \ref{prop:pdec fonc implique plien} prouvent que pour la famille minimale de parahoriques $P$, les conditions \plien(sph) fonctorielle (et même \plien\ fonctoriel sur les cloisons) et \pdec\ fonctorielle sont équivalentes, car pour $Q=P$, on obtient $R=P$, quel que soit le signe $\epsilon$ choisi.\\

\fcor{\label{cor:P verifie pdec} S'il existe une famille de parahoriques vérifiant \psph\ et \plien($\vec m$) pour chaque cloison $\vec m$, alors la famille minimale $P$ vérifie \pdec.}

\demo 
Pour toute cloison $\vec m$, le fait que $Q$ vérifie \psph\ et \plien($\vec m$) entraine que $P$ vérifie aussi \plien($\vec m$), nous sommes donc bien dans les conditions d'application de la proposition.

Soit $a\in A$, soit $\vec C$ une chambre de $\vec f_a^*\cap \vec A$. Soit $\epsilon$ le signe de $\vec C$. Soit $R^\epsilon$ la famille de parahoriques construite par la proposition, alors $R^\epsilon(a) = (Q(a)\cap U(\vec C)).(P(a)\cap U(-\vec C)).N(a)$. Prenant l'intersection avec $P(a)$, il vient $P(a)= P(a)\cap R^\epsilon(a) = (P(a)\cap U(\vec C)).(P(a)\cap U(-\vec C)).N(a)$.\cqfd

\subsection{Bonnes familles de parahoriques}
\label{sousoussection:csq de plienp}

\subsubsection{Une condition suffisante pour \plienn}

Nous avons déjà vu au \ref{cor:csq de pdec pspsh et fonc} que la condition \pdec\ fonctorielle implique \pinj, \psph\ et \plien(sph). Nous pouvons maintenant, grâce à la décomposition d'Iwasawa, prouver qu'elle implique également \plienn(sph), autrement dit que toute famille vérifiant \pdec\ fonctoriellement est une bonne famille de parahoriques.

\fprops{ \label{prop:plien et plienn}
\item Soit $Q$ une famille de parahoriques qui vérifie \pdec\ fonctoriellement. Alors $Q$ vérifie \plienn(sph) (fonctoriellement).

\item Soit $\vec f$ une facette telle que $Q$ vérifie \plienn($\vec f$). On suppose de plus que $Q$ vérifie \pinj.
%\psph et \plienn($\vec m$) pour toute cloison $\vec m$.
 Alors $Q$ vérifie \plien($\vec f$). En particulier, toute bonne famille de parahoriques vérifie \plien(sph).
}

\rema Les conditions du deuxième point sont plutôt  plus faibles que celles du premier point, et sa preuve plus simple, on peut donc considérer \plien\ comme plus faible que \plienn.\\

\demos{

\item Rappelons que $Q$ vérifie \psph, \plien\ et \pinj\ par la proposition \ref{prop:pdec fonc implique plien}.
Soit $a\in A$ et $\vec f$ une facette sphérique dans $\vec f_a^*$.
 Soit $g\in N.Q(a)\cap N.P(\vec f)$, on peut supposer $g\in Q(a)\cap N(\vec f_a).P(\vec f)$. Grâce à \pfonc, soit $\Delta$ tel que $a$ est spécial dans $A^\Delta$, puis soit $n\in N^\Delta(a)$ tel que $ng\in P^\Delta(\vec f)$. Comme $Q$ vérifie fonctoriellement \plien($\vec f$) puis \psph, on a $ng\in Q^\Delta(a)\cap P^\Delta(\vec f) = Q^\Delta(\{a,pr_{\vec f}(a)\}) \subset P^\Delta(pr_{\vec f}(a))$. D'où $g\in N^\Delta(a).P^\Delta(pr_{\vec f}(a))\cap G$. Le point $pr_{\vec f}(a)$ étant sphérique, la famille $Q$ vérifie \pfonc($pr_{\vec f}(a)$) (par \ref{prop:csq de fonc}, 5) et \pdec($pr_{\vec f}(a)$) (grâce à \ref{prop:csq de parasph}, 6). Par le corollaire \ref{cor:unicite dans iwasawa}, puis la proposition \ref{prop:csq de parasph}, 1, on obtient $g\in Q(a)\cap  N.P(pr_{\vec f}(a)) = Q(a)\cap N.U(\vec f).G(\phi^m(\vec f),a)$. Comme $G(\phi^m(\vec f),a)\subset Q(\{a,pr_{\vec f}(a)\})$, on peut supposer $g\in Q(a)\cap N U(\vec f)$.

 Soit $\vec C$ une chambre dont l'adhérence contient $\vec f$. Soit $g=n u$ une écriture de $g$ correspondant à $g\in N.U(\vec f)$ et  
 $g=n'u^{-\epsilon}u^\epsilon$ une écriture correspondant à $g\in Q(a) = N(a).(U(-\vec C)\cap Q(a)).(U(\vec C)\cap Q(a))$ (obtenue par \pdec($a$)). Alors $n\inv n'= u (u^\epsilon)\inv (u^{-\epsilon})\inv$. Et $u(u^\epsilon)\inv \in U(\vec C)$. Donc par l'unicité dans la décomposition de Birkhoff vectorielle de $G$, on obtient $u=u^\epsilon\in U(\vec f)\cap Q(a)\subset Q(\{a,pr_{\vec f}(a)\})$. Donc $g=nu\in N.Q(\{a,pr_{\vec f}(a)\})$.\\

\item Soit $g\in Q(a)\cap P(\vec f)$. D'après \plienn$(\vec f)$, on a alors $g\in N.Q(\{a,pr_{\vec f}(a)\})$. Soit $g=nq$ l'écriture correspondante, on a $a=g.a=n.a$, d'où $n\in Q(a)\cap N= N(a)$ grâce à \pinj. De même, $\vec f=g.\vec f= n.\vec f$ d'où $n\in N(\vec f)$. Mais comme $N$ agit de manière affine sur $A$, $n\in N(a)\cap N(\vec f) = N(\barre{a+\vec f})\subset Q(\{a,pr_{\vec f}(a)\})$.
}

\fcor{\label{cor:pdec implique plienn} Toute famille de parahoriques $\Q$ vérifiant fonctoriellement \pdec\ vérifie fonctoriellement \psph, \pinj, \plien(sph) et \plienn(sph), c'est donc une bonne famille de parahoriques.

C'est en particulier le cas pour la famille minimale de parahoriques dans un groupe de Kac-Moody déployé.
}
\demo C'est une conséquence du corollaire \ref{cor:csq de pdec pspsh et fonc} et de la proposition ci-dessus.
\cqfd

\subsubsection{Projections}

La condition \plienn($\vec f$) entraine que la projection $pr_{\vec f}(a)$ est bien définie dans $\I$, dès que $\vec f_a\subset \barre{\vec f}$, indépendamment de l'appartement contenant $a$ et $\vec f$ considéré. Ainsi, dans une bonne famille de parahoriques, la conjonction de \plienn(sph) et de \psph\ permettra souvent de ramener une preuve à une étude dans une façade sphérique, donc dans un immeuble affine, bien connu grâce à \cite{bruhat-tits}.

\fprop{ \label{prop:projection bien def} Si $\Q$ vérifie \plienn($\vec f$), alors la projection $pr_{\vec f}$ est bien définie sur la réunion des $A(T)_{\vec g}$ tels que $\vec f\cup\vec g\subset \vec A(T)$ et $\vec g\subset \barre{\vec f}$. Pour tout $g\in G$, on a $g.pr_{\vec f}(a)=pr_{g\vec f}(g.a)$.
}
\demo \\
Soit $a$ dans une façade $A(T)_{\vec g}$ telle que $\vec f\cup\vec g\subset \vec A(T)$ et $\vec g\subset \barre{\vec f}$. Soit $T_0$ un tore maximal de référence, soit $g\in G$ tel que $g.A(T)=A(T_0)$. Alors $g.\vec f$ et $g.\vec g$ sont deux facettes de $\vec A(T_0)$ telles que $g.\vec g\subset \vect_{\vec A(T_0)}(g.\vec f)$ donc la projection $pr_{g\vec f}(ga)$ est bien définie dans $A(T_0)$. On veut poser $pr_{\vec f}(a) = g\inv(pr_{g\vec f}(ga))$, il s'agit de vérifier que ceci est indépendant de $T$ et de $g$.\\
 Soient donc $T'$ et $g'$ d'autres choix possibles. Alors $g' g\inv$ envoie $ga$ et $g\vec f$ sur un autre point de $A(T_0)$ et une autre facette de $\vec A(T_0)$. Donc $g' g\inv \in N(T_0)Q(ga)\cap N(T_0)P(g\vec f)$. Ceci vaut $N(T_0)Q(\{ga,pr_{g\vec f}(ga)\})$ par \plienn($g\vec f$). Soit $n\in N(T_0)$ tel que $g'g\inv\in n.Q(\{ga,pr_{g\vec f}(ga)\})$. Alors $g'g\inv(pr_{g\vec f}(ga)) =n(pr_{g\vec f}(ga)) = pr_{ng\vec f}(nga) = pr_{g'\vec f}(g' a)$. Ceci prouve que $g\inv(pr_{g\vec f}(ga)) = (g')\inv(pr_{g'\vec f}(g'a))$.\cqfd

\subsubsection{Conséquences variées}

On rassemble ici quelques petits résultats obtenus au moyen des projections pour une bonne famille de parahoriques.\\

On peut dans une certaine mesure caractériser un point de $\I$ par ses projetés dans diverses façades. Par exemple, on prouve le:
\flemme{\label{lemme:action de u sur A} Soit $\vec f$ une facette de $\vec A$, et $Q$ une famille de parahoriques vérifiant \psph\ et \plienn($\vec m$) pour toute cloison $\vec m$ de $\barre{\vec f^*\cap \vec A}$.
 Alors pour tout $u\in U(\vec f)$, pour tout $a\in A$ tel que $\vec f\subset \vec f_a^*$, ou $u.a=a$, ou bien $u.a\not\in A$.

En terme de groupes, ceci s'écrit $U(\vec f)\cap N.Q(a) \subset Q(a)$.
}
\pv Soit $a\in A$, supposons $u.a\in A$. Soit $\vec m$ une cloison de $\barre{\vec f ^*}\cap \vec A\cap \vec f_a^*$, alors $u\in P(\vec m)$, donc $u.pr_{\vec m}(a) = pr_{\vec m}(u.a)\in A_{\vec m}$. Ainsi $u$ envoie le point $pr_{\vec m}(a)$ sur un autre point de l'appartement $A_{\vec m}$. Sachant que $u$ fixe par ailleurs un cône ouvert de $A_{\vec m}$ (c'est-à-dire une demi-droite, puisque $\dim(A_{\vec m}) =1$), et que $\I_{\vec m}$ est un immeuble affine (en fait un arbre), ceci entraine que $u.pr_{\vec m}(a) = pr_{\vec m}(a)$, d'où $pr_{\vec m}(u.a) = pr_{\vec m}(a)$.\\
 Donc $a$ et $u.a$ sont deux points de $A$ dont les projetés sur tous les $A_{\vec m}$, pour $\vec m$ une cloison de $\barre{\vec f^*}\cap \vec A\cap \vec f_a^*$ coïncident. Ceci entraine que $\fl{a \, u(a)}$ est dans l'intersection de ces $\vec m$. Mais celle-ci est triviale dans $A_{\vec f_a}$.\cqfd

Voici la première conséquence de ce lemme:
\fprop{ Toute famille de parahoriques vérifiant \psph, \pinj, \plien($\vec m$)  et \plienn($\vec m$) pour toute cloison $\vec m$ vérifie aussi \pdec.
}
\rema grâce à \ref{prop:plien et plienn}, on peut remplacer l'hypothèse \plien($\vec m$) par \pinj.\\

\demo

Soit $a\in A$ et $\vec C\in\F(\vec f_a^*\cap \vec A)$. Par la décomposition d'Iwasawa, $G= U(\vec C).N.P(a)$. Prenant l'intersection avec $Q(a)$, il vient $Q(a) = ( U(\vec C).N\cap Q(a)).P(a)$. Mais le lemme entraine facilement (avec \pinj) que $U(\vec C).N\cap Q(a) = (U(\vec C)\cap Q(a)).N(a)$. D'où $Q(a)=(U(\vec C)\cap Q(a)).P(a)$.

 Maintenant, $P$ vérifie \pdec\ par le corollaire \ref{cor:P verifie pdec}, donc $P(a) = (U(\vec C)\cap P(a)).(U(-\vec C)\cap P(a)).N(a)\subset (U(\vec C)\cap Q(a)).(U(-\vec C)\cap Q(a)).N(a)$. D'où le résultat: $Q(a) = (U(\vec C)\cap Q(a)).(U(-\vec C)\cap P(a)).N(a)\subset(U(\vec C)\cap Q(a)).(U(-\vec C)\cap Q(a)).N(a) $.\cqfd
 
 \rema Le résultat obtenu est même un peu plus fort: on a vu que $Q(a) = (U(\vec C)\cap Q(a)).(U(-\vec C)\cap P(a)).N(a)$.\\

 Vu le corollaire \ref{cor:pdec implique plienn} et la proposition \ref{prop:plien implique pinj}, dans le cas où la donnée radicielle est fonctorielle, on a donc pour toute famille de parahoriques $Q$ équivalence entre les assertions suivantes:
\liste{
\item $Q$ est fonctoriellement une bonne famille de parahoriques.
\item $Q$ vérifie fonctoriellement \pdec.
\item $Q$ vérifie fonctoriellement \psph, \plien($\vec m$) et \plienn($\vec m$) pour toute cloison $\vec m$.
\item $Q$ vérifie fonctoriellement \psph, \pinj,  et \plienn($\vec m$) pour toute cloison $\vec m$.
\\
}

Pour une famille $Q$ vérifiant \psph, \pinj, et \plienn\ pour les cloisons, le lemme implique aussi l'égalité $N\cap U(\vec C).Q(a)= N(a)$, pour toute chambre $\vec C$ de $\vec A$ et tout point $a\in A$ tels que $\vec f_a\in \barre{\vec C}$. Et ceci permet, comme en \ref{prop:unicite dans iwasawa}, de prouver l'unicité modulo $N(a)$ du facteur dans $N$ de la décomposition d'Iwasawa:
\fprop{\label{prop:unicite dans iwasawa 2} Soit $Q$ une bonne famille de parahoriques. Soient $\vec C$ une chambre de $\vec A$ et $a\in A$.
 Alors $N\cap U(\vec C).Q(a)= N(a)$ si $\vec f_a\subset \barre{\vec C}$.
 
 En conséquence, pour tous $n,n'\in N$ tels que $\vec f_a\subset n\barre{\vec C}$ ou $\vec f_a\subset n'\barre{\vec C}$, on a:
  $$U(\vec C).n.Q(a)\cap U(\vec C).n'.Q(a)\not=\vide \implique n\inv n'\in N(a) \implique U(\vec C).n.Q(a)= U(\vec C).n'.Q(a).$$
}
\demo
Le premier point découle du lemme \ref{lemme:action de u sur A}, la suite est immédiate (voir la partie \ref{soussoussection:unicie dans iwasawa}).\cqfd

\fprop{\label{prop:points fixes de u1u2...uk} On suppose que $Q$ vérifie \plienn($\vec m$) pour toute cloison $\vec m$.
 Soit $\psi=\{\alpha_1,...,\alpha_k \} \subset \phi^u(\vec C)$ un ensemble réduit de racines,
rangées dans un ordre grignotant. (Ceci signifie que pour tout $i$, $\ker(\alpha_i)$ contient une cloison incluse dans l'intérieur de $\bigcap _{j> i} \D(\alpha_j)$, voir par exemple \cite{remy} 9.1.2). 
 
 Soit $u=u_1...u_k$ avec $u_i\in U_{\alpha_i}$. Alors:
\[\fix_A(u) = \bigcap_i \fix(u_i) = \bigcap_i \D(\alpha_i,\phii_{\alpha_i}(u_i)) \]
}

\demo 

La seconde égalité vient de \ref{prop:points fixes de ualpha}. Pour la première, l'inclusion "$\supset$" est claire.\\

Soit $a\in \fix_A(u)$. Par hypothèse, il existe une cloison $\vec m_1$, incluse dans $\ker(\alpha_1)$ et dans l'intérieur de $\bigcap _{j> 1} \D(\alpha_j)$. Ainsi, $u\in P(\vec m_1)$, et pour tout $j> 1$, $u_j$ fixe $A_{\vec m_1}$.  Donc $u$ agit sur $A_{\vec m_1}$ comme $u_1$. Par la proposition \ref{prop:projection bien def}, $u$ et $u_1$ fixent le point $pr_{\vec m_1}(a)$. Ceci entraine $u_1\in U_{\alpha_1}(pr_{\vec m_1}(a)) = U_{\alpha_1}(a)$.

 Maintenant l'élément $u_1\inv u=u_2...u_k$ fixe $a$ et est le produit de $k-1$ éléments de groupes radiciels dans un ordre grignotant, par récurrence on obtient que pour tout $i\in\ent{2}{k}$, $u_i\in U_{\alpha_i}(a)$.
 
  D'où $u\in U_{\alpha_1}(a)...U_{\alpha_k}(a)$ et $a\in \bigcap _{i=1}^k \fix(u_i)$.\cqfd

\subsubsection{Décomposition de Lévi}

 La proposition suivante donne dans certains cas une décomposition de Lévi pour $Q(\Omega)$:
 
\fprop{\label{prop:levi pour Qa} Soit $\vec f$ une facette sphérique de $\vec A$, soit $Q$ une famille de parahoriques vérifiant \plien($\vec f$) et \plienn($\vec m$) pour chaque cloison $\vec m$ de $\vec f^*\cap \vec A$. Soit $\Omega$ une partie de $A$ telle que $\cl(\vec\Omega) = \barre{\vec f}$ %et $\Omega\cap A_{\vec f} = pr_{\vec f}(\Omega)$ ($\Omega$ contient son projeté sur $A_{\vec f}$)
. Alors:
\calc{
 Q(\Omega) &=& \left( U(\vec f)\cap Q(\Omega)\right) \rtimes \left(N(\Omega). G\left( \phi^m(\vec f),\Omega \right)\right)\\
 &=& \left( U(\vec f)\cap Q(\Omega)\right) \rtimes Q\left(\eng{\Omega,\vec f}_A\right) \point
 }

%Soient $Q$ et $\vec f$ comme dans le lemme, avec $\vec f$ sphérique. On suppose de plus que $Q$ satisfait \plien($\vec f$). Alors pour tout $a\in A_{\vec g}$ avec $\vec g\in\F(\barre{\vec f})$,
%\[ Q(a)\cap P(\vec f) = Q(\eng{a,\vec f}_A) \rtimes \left(U(\vec f)\cap Q(a)\right) \point \]
}
\demo

Il est clair que $  \left( U(\vec f)\cap Q(\Omega)\right) \rtimes \left(N(\Omega). G\left( \phi^m(\vec f),\Omega \right)\right)\subset \left( U(\vec f)\cap Q(\Omega)\right) \rtimes Q\left(\eng{\Omega,\vec f}_A\right) \subset Q(\Omega)$.\\

Soit $g\in Q(\Omega)$. Notons $\Omega_{\vec f}=pr_{\vec f}(\Omega)$, d'après \plien($\vec f$), $g\in Q(\Omega_{\vec f})$. D'après \ref{prop:csq de parasph}, 2, on a $g\in U(\vec f)\rtimes\left( N(\Omega_{\vec f}).G(\phi^m(\vec f), \Omega_{\vec f})\right)$. Mais $G(\phi^m(\vec f), \Omega_{\vec f}) = G(\phi^m(\vec f), \Omega)$, on peut donc supposer $g\in (U(\vec f).N(\Omega_{\vec f})) \cap Q(\Omega)$. Maintenant, le lemme \ref{lemme:action de u sur A} entraine $g\in (U(\vec f)\cap Q(\Omega))\point N(\Omega)$.\cqfd

Cette proposition permet, pour prouver que $Q$ vérifie \plienp($\vec f$), de se ramener à l'étude de $U(\vec f)\cap Q(a)$, comme l'indique le corollaire suivant:

\fcor{Soit $Q$ une famille de parahoriques vérifiant les hypothèses de la proposition. Si $a$ est un point d'une façade $A_{\vec g}$ avec $\vec g\subset \barre{\vec f}$, alors:
\calc{
 Q(a) \cap P(\vec f) &=& \left( U(\vec f)\cap Q(a)\right) \rtimes \left(N(a+\vec f). G\left( \phi^m(\vec f),a \right)\right)\\
 &=& \left( U(\vec f)\cap Q(a)\right) \rtimes Q\left(\eng{a,\vec f}_A\right) \point
 }
}
\demo Ayant remarqué que $Q(a)\cap P(\vec f) = Q(\{a,pr_{\vec f}(a)\})$, on applique la proposition à $\Omega=\{a,pr_{\vec f}(a) \}$.\cqfd

Nous avons ainsi obtenu une décomposition de Lévi pour le groupe $Q(\Omega)$ lorsque $\cl(\vec\Omega)$ est l'adhérence d'une facette $\vec f$ sphérique. Cette décomposition repose sur la décomposition de Lévi vectorielle de $P(\vec f)$. Il y a un autre cas où on dispose d'une décomposition de Lévi vectorielle: c'est pour le fixateur d'une partie $\vec\Omega$ équilibrée (c'est-à-dire une union finie de facettes sphériques dont au moins une est positive et une négative). Ceci fournit naturellement une décomposition des groupes $Q(\Omega)$ correspondants. De plus, dans le cas équilibré on dispose d'une bonne description du facteur unipotent.

\fprop{ \label{prop:levi equilibre}Soit $\Omega$ une partie de $A$ contenant au moins un point sphérique positif et un point sphérique négatif. On suppose que $Q$ est une bonne famille de parahoriques. Alors:
\[ Q(\Omega) =  U(\Omega) \rtimes \left( M(\vec\Omega)\cap Q(\Omega)\right)\]
avec:
\calc{
U(\Omega)&:=& G(\phi^u(\vec\Omega),\Omega)=\prod_{\alpha\in\phi^u_{red}(\vec\Omega)} U_\alpha(\Omega)\\
M(\vec \Omega)\cap Q(\Omega) &=& N(\Omega).G(\phi^m(\vec\Omega),\Omega)\point
}
(Le produit dans $U(\Omega)$ peut s'effectuer dans n'importe quel ordre.)}

\demo

Soient $\vec f^+$ et $\vec f^-$ des facettes maximales de $\cl(\vec\Omega)\cap \vec A^+$ et $\cl(\vec\Omega)\cap \vec A^-$. Alors $\vect(\vec f^+)=\vect(\vec f^-) = \vect(\cl(\vec\Omega))$, donc $M(\vec f^+)= M(\vec f^-)=M(\cl(\vec\Omega))=M(\vec\Omega)$. De plus, $\vec f^+\cup \vec f^-$ est une partie équilibrée, donc son fixateur admet la décomposition de Lévi $P(\vec f^+\cup \vec f^-) = U(\vec f^+\cup \vec f^-) \rtimes M(\vec\Omega)$.

Soit $g\in Q(\Omega)$. En particulier, $g\in P(\vec \Omega)=P(\cl(\vec\Omega)) \subset P(\vec f^+\cup \vec f^-)$, donc $g\in Q(\Omega)\cap (U(\vec f^+\cup \vec f^-) \rtimes M(\vec\Omega))$. Soient $u\in U(\vec f^+\cup \vec f^-)$ et $m\in M(\vec\Omega)$ tels que $g=um$. Montrons que $u$ et $m$ fixent $\Omega$.

 Soit $\omega\in\Omega$. Alors $\vec f^+\cup \vec f^-\cup \vec f_\omega$ est une partie équilibrée, donc 
 $$P(\vec f^+\cup \vec f^-\cup \vec f_\omega) = U(\vec f^+\cup \vec f^-\cup \vec f_\omega)\rtimes M(\vec f^+\cup \vec f^-\cup \vec f_\omega) = U(\vec f^+\cup \vec f^-\cup \vec f_\omega)\rtimes M(\vec\Omega)\point$$
  Le groupe $U(\vec f^+\cup \vec f^-\cup \vec f_\omega)$ fixe toute la façade $A_{\vec f_\omega}$, donc 
  $$g\in Q(\omega)\cap P(\vec f^+\cup \vec f^-\cup \vec f_\omega) = U(\vec f^+\cup \vec f^-\cup \vec f_\omega)\rtimes \left( M(\vec\Omega) \cap Q(\omega) \right)\point$$
  
  Mais $U(\vec f^+\cup \vec f^-\cup \vec f_\omega)\subset U(\vec f^+\cup \vec f^-)$, donc $g=um$ est également l'écriture de $g$ dans la décomposition de $Q(\omega)\cap P(\vec f^+\cup \vec f^-\cup \vec f_\omega)$ qu'on vient d'obtenir. En particulier, $u$ et $m$ fixent $\omega$.
  
  On a ainsi prouvé:
  \calc{
   Q(\Omega) &=& \left( \bigcap_{\vec f\in\F(\vec \Omega)} U(\vec f^+\cup \vec f^-\cup \vec f)\right) \rtimes \left( M(\vec\Omega)\cap Q(\Omega) \right)\\
   &=& \left(U(\vec f^+\cup \vec f^-)\cap Q(\Omega) \right) \rtimes \left( M(\vec\Omega)\cap Q(\Omega) \right)  \point
   }

 Étudions le premier facteur. Comme $\vec f^+\cup \vec f^-$ sont deux facettes sphériques de signes opposés, on sait par \cite{remy}  6.3.1 que $U(\vec f^+\cup \vec f^-)= G(\phi^u(\vec f^+\cup \vec f^-)) = \prod_{\alpha\in\phi^u_{red}(\vec f^+\cup \vec f^-)} U_\alpha$, pour un ordre qu'on peut choisir grignotant sur $\phi^u_{red}(\vec f^+\cup \vec f^-)$ (car $\phi(T)$ est réduit). Alors la proposition \ref{prop:points fixes de u1u2...uk} permet de montrer que:
 $$Q(\Omega)\cap U(\vec f^+\cup \vec f^-)=\prod_{\alpha\in\phi^u_{red}(\vec f^+\cup \vec f^-)} U_\alpha(\Omega) \point $$
 Mais si $\alpha\not\in \phi^u(\vec \Omega)$, alors $U_\alpha(\Omega)=\{e\}$. D'où finalement:
 $$Q(\Omega)\cap U(\vec f^+\cup \vec f^-)=\prod_{\alpha\in\phi^u_{red}(\vec\Omega)} U_\alpha(\Omega) \point $$

Passons à la description de $M(\vec\Omega)\cap Q(\Omega)$. Soit $m\in M(\vec\Omega)\cap Q(\Omega)$. Soit $\omega\in\Omega$. Soit $\vec f=pr_{\vec f_\omega}(\vec f^\epsilon)$, où $\epsilon$ est un signe de $\vec f_\omega$. Il s'agit d'un facette sphérique. Par \plien($\vec f$) puis \ref{prop:csq de parasph}, $m\in Q(pr_{\vec f}(\omega))\cap M(\vec\Omega) = N(pr_{\vec f}(\omega)). G(\phi^m(\vec\Omega),\omega)$. Comme $N(pr_{\vec f}(\omega)) = N(pr_{\vec f^+}(\omega))$, on obtient $m\in Q(pr_{\vec f^+}(\omega))$. Au final, $m$ fixe la partie $pr_{\vec f^+}(\Omega)$, d'où, encore avec \ref{prop:csq de parasph}:
 $$m\in N(pr_{\vec f^+}(\Omega)) . G(\phi^m(\vec\Omega),\Omega)\cap Q(\Omega) = N(\Omega).G(\phi^m(\vec\Omega),\Omega) \point$$ \cqfd

\fcor{\label{cor:enclos equilibre} Dans les conditions de la proposition, $Q(\Omega)=N(\Omega).Q(\cl(\Omega))$.}

\fex{ Lorsque $\Omega=\D(\alpha,k)$ est un demi-appartement dans un appartement $A$, la proposition donne $Q(\Omega) = U_{\alpha,k}\rtimes H$. Soit $\Omega'=(\Omega\cap\inte A)\cup M(\alpha,k)$, alors $\cl(\Omega')=\Omega$, $N(\Omega')=N(\Omega)$, $\vec\Omega'=\vec M(\alpha)$, et la proposition donne alors $Q(\Omega)=Q(\Omega') = \{e\}\ltimes(H.U_{\alpha,k})$. On peut ainsi obtenir plusieurs décompositions de type Lévi d'un même groupe, selon que les facteurs de la forme $U_\alpha(\Omega)$ tels que $U_{-\alpha}(\Omega)=\{e\}$ sont considérés comme inclus dans le facteur unipotent ou dans le facteur de Lévi.

On préfèrera sûrement les placer dans le facteur unipotent, et pour obtenir ce résultat il faut appliquer la proposition à la partie $\cl(\Omega)$ au lieu de $\Omega$.  }

\subsection{La bonne famille de parahoriques maximale}

\label{soussection:famille maximale}

\fdef{ Soit $Q$ une bonne famille de parahoriques. On définit:
\[\bar Q(a) = \ens{g\in P(\vec f_a)}{\forall \vec f\in (\vec f_a^*)\cap \iv\sph,\: g.pr_{\vec f}(a) = pr_{g.\vec f}(a)}\]
}

Ainsi $\bar Q(a)$ est l'ensemble des éléments de $G$ qui permutent les projetés de $a$ dans les façades sphériques.

\flemme{Soit $Q$ une famille de parahoriques vérifiant \psph\ et \plienn(sph) fonctoriellement. Soit $\Delta$ un sous-groupe de $\R$. Alors l'application
\[\psi: \fonc{\I(Q)} {\I(Q^\Delta)} {[g,a]_Q}{[(g,a)]_{Q^\Delta}}\]
est bien définie, $G$-équivariante, et compatible aux projections (c'est-à-dire $\psi(pr_{\vec f}(a)) = pr_{\vec f}(\psi(a))$ pour toute facette $\vec f$ sphérique et tout $a\in A_{\vec g}$ avec $\vec g\subset \barre{\vec f}$). Elle induit une injection de $\I(Q)\sph$ dans $\I(Q^\Delta)\sph$, et plus généralement, l'image inverse d'un point $[(g,a)]_{Q^\Delta}$, $a\in A$ et $g\in G$ est de cardinal 1 dès que $G\cap N^\Delta.Q^\Delta(a) = N.Q(a)$.\\

%En conséquence, la famille $\bar Q$ vérifie \pfonc, autrement dit pour tout $a\in A$, $\bar Q^\Delta(a)\cap G=\bar Q(a)$.
}

\pv Il est clair que $\psi$ est bien définie et $G$-équivariante. Elle est également compatible aux projections dans l'appartement $A$, puis par $G$-équivariance sur $\I(Q)$.

Détaillons un peu l'injectivité: soient $g,h\in G$, $a,b\in A$ tels que $(g,a)\sim_{Q^\Delta} (h,b)$. Alors $g\inv h\in N^\Delta.Q^\Delta(a)\cap G$, et par hypothèse (ou par \ref{prop:csq de fonc} pour un point $a$ sphérique), $N^\Delta.Q^\Delta(a)\cap G=N.Q(a)$. D'où $(g,a)\sim_Q (h,b)$.\cqfd

\fprop{ Pour toute famille de parahoriques $Q$ fonctoriellement bonne, $\bar Q$ est encore une famille de parahoriques fonctoriellement bonne pour $\D$. Elle vérifie en outre \pfonc, et \plien, et non juste \plien(sph).
}

\demo

Pour tout $a\in A$, $\bar Q(a)$ est un sous-groupe de $P(\vec f_a)$ contenant $Q(a)$ et donc $P(a)$ d'après \ref{prop:projection bien def}. De plus, (para 0.4) est claire, donc $\bar Q$ est une famille de parahoriques.

Si $a\in A\sph$, alors $\vec f_a$ est une facette sphérique de $\vec f_a^*$ donc par définition de $\bar Q(a)$, on a pour tout $g\in\bar Q(a)$, $g.a=g.pr_{\vec f_a}(a)=pr_{g\vec f_a}(a)=a$. Donc $\bar Q(a)=Q(a)=P(a)$, car $Q$ vérifie \psph. Donc $\bar Q$ vérifie  \psph.\\

Montrons \pfonc. Soit $a\in A$ et $g\in \bar Q^\Delta(a)\cap G$. Soit $\vec f\in (\vec f_a^*)\sph$. Soit $\psi:\I(Q)\fleche \I(Q^\Delta)$ comme dans le lemme. On a $\psi(g.pr_{\vec f}(a)) = g.\psi(pr_{\vec f}(a)) = g.pr_{\vec f}(\psi(a)) = pr_{g\vec f}(\psi(a)) = \psi(pr_{g\vec f}(a))$. Mais comme $pr_{\vec f}(a)$ est un point sphérique et $\psi$ est injective sur $\I(Q)\sph$, ceci entraine que $g.pr_{\vec f}(a) = pr_{g\vec f}(a)$.

On prouve ainsi que $g\in \bar Q(a)$.\\

 Étudions \plien. Soit $a\in A$, $\vec f\in\vec f_a^*$ et $q\in \bar Q(a)\cap P(\vec f)$. Pour montrer que $q\in \bar Q(pr_{\vec f}(a))$, il suffit de montrer que pour toute facette $\vec h \in (\vec f^*)\sph$, $q.pr_{\vec h}( pr_{\vec f}(a)) = pr_{q.\vec h}(pr_{\vec f}(a))$. Mais $pr_{\vec h}(pr_{\vec f}(a)) = pr_{\vec h}(a)$, et $pr_{q.\vec h}(pr_{\vec f}(a))= pr_{q\vec h}(a)$, le résultat en découle.\\
 
 La condition \pinj\ est alors conséquence de \ref{prop:plien implique pinj}.\\

Il reste à voir \plienn(sph). Soit donc $a\in A$, $\vec f$ une facette sphérique de $\vec f_a^*$, et $g\in N\bar Q(a)\cap P(\vec f)=N(\vec f_a)\bar Q(a)\cap P(\vec f)$. Soit $\Delta\leq \R$ tel que $a$ soit spécial dans $A^\Delta$, il existe alors $t\in T^\Delta$ tel que $tg\in \bar Q^\Delta(a)\cap P^\Delta(\vec f)$, donc par \plien($\vec f$) pour $\bar Q^\Delta$, $tg\in \bar Q^\Delta(\{a,pr_{\vec f}(a)\})$ puis $g\in T^\Delta \bar Q^\Delta(\{a,pr_{\vec f}(a)\}) \cap G \subset  T^\Delta P^\Delta(pr_{\vec f}(a))\cap G$. Alors, comme $pr_{\vec f}(a)$ est sphérique, $P^\Delta$ vérifie \pdec($pr_{\vec f}(a)$) donc par \ref{cor:unicite dans iwasawa}, \ref{prop:csq de fonc} point 5, puis \ref{prop:csq de parasph}, $g\in N(\vec f).P(pr_{\vec f}(a)) = N(\vec f).U(\vec f).G(\phi^m(\vec f),a)$. Quitte à multiplier $g$ à droite par un élément de $G(\phi^m(\vec f),a)$, et à gauche par un élément de $N(\vec f)$, on peut donc supposer $g\in N(\vec f_a).\bar Q(a)\cap  U(\vec f)$.\\

Dans $G^\Delta$, ceci donne $g\in T^\Delta.\bar Q^\Delta(a)\cap U^\Delta(\vec f)$ puisque $a$ est spécial dans $A^\Delta$. Soit $t\in T^\Delta$ tel que $g\in t.\bar Q^\Delta(a)$, alors $t\inv g\in \bar Q^\Delta(a)\cap P(\vec f^*\cap \vec A)$, donc $t\inv g$ fixe tous les projetés de $a$ sur les façades sphériques de direction dans $\barre{\vec f^*\cap \vec A}\cap \vec f_a^*$. En particulier, soit $\vec m$ une cloison dans $\barre{\vec f^*\cap \vec A}\cap \vec f_a^*$, alors $g.pr_{\vec m}(a) = t.pr_{\vec m}(a) \in A_{\vec m}$. Comme $A_{\vec m}$ est une façade sphérique, le lemme \ref{lemme:action de u sur A} y est toujours vrai, et le fait que $g\in U^\Delta(\vec f)$ entraine alors $g.pr_{\vec m}(a) = pr_{\vec m}(a)$. Donc $t$ fixe la façade $A_{\vec m}$. Au total, $t$ induit sur $A_{\vec f_a}$ une translation de vecteur inclus dans chaque $\vect(\vec m)$ pour $\vec m$ une cloison dans $\barre{\vec f^*\cap \vec A}\cap \vec f_a^*$. L'intersection de ces hyperplans et de $\vec A_{\vec f_a}$ est triviale, donc $t$ fixe $A_{\vec f_a}$, et $g\in \bar Q^\Delta(a)\cap U^\Delta(\vec f)\cap G$.

Appliquant \plien($\vec f$), puis \pfonc, on obtient $g\in \bar Q^\Delta(\{a,pr_{\vec f}(a)\}) \cap G = \bar Q(\{a,pr_{\vec f}(a)\})$.\cqfd

%Recommençant alors le raisonnement précédent, on obtient cette fois $g\in N^\Delta.(\bar Q^\Delta(\{a,pr_{\vec f}(a)\})\cap U^\Delta(\vec f))\cap G$. En fonction du signe de $\vec f$, soit  $R$ celle des deux familles de parahoriques vérifiant \pdec\ définies au lemme \ref{lemme:obtenir pdec}, qui est telle que $\bar Q^\Delta(\{a,pr_{\vec f}(a)\})\cap U^\Delta(\vec f) \subset R^\Delta(\{a,pr_{\vec f}(a)\})$. Alors $g\in N^\Delta.\left(U^\Delta(\vec f)\cap R^\Delta(\{a,pr_{\vec f}(a)\})\right)\cap G$.
% Soit $n^\Delta\in N^\Delta$ et $r^\Delta\in U^\Delta(\vec f)\cap R^\Delta(\{a,pr_{\vec f}(a)\})$ tels que $g=n^\Delta u^\Delta$. Par ailleurs, soit $F$ la facette contenant $\gm_a(a+\vec f)$, soit $\vec C$ une chambre de $\vec f^*$ et soit $g=rnu$ une écriture de $g$ dans la décomposition d'Iwasawa $G= R(F).N.U(\vec C)$. C'est aussi, tout comme $g= n^\Delta u^\Delta$ une écriture de $g$ dans la décomposition $G^\Delta= R^\Delta(F).N^\Delta.U(\vec C)$. D'où par la proposition \ref{prop:unicite dans iwasawa}, $n^\Delta\in n.N^\Delta(F)$. Mais $N^\Delta(F)\subset R^\Delta(\{a,pr_{\vec f}(a)\})$ donc $g\in n.R^\Delta(\{a,pr_{\vec f}(a)\})\cap G \subset N.( \bar Q^\Delta(\{a,pr_{\vec f}(a)\})\cap G )$. Enfin, comme $\bar Q$ vérifie \pfonc, $g\in N.\bar Q(\{a,pr_{\vec f}(a)\})$.\cqfd
 
 \flemme{ Soient $Q$ et $R$ deux familles de parahoriques vérifiant \psph\ et \plienn(sph). On suppose que pour tout $a\in A$, $Q(a)\subset R(a)$. Alors le morphisme surjectif $G$-équivariant d'immeubles $\psi:\fonc{\I(Q)}{\I(R)}{[g,a]_Q}{[g,a]_R}$ induit un isomorphisme entre $\I(Q)\sph$ et $\I(R)\sph$, et pour toute facette sphérique $\vec h$, pour tout point $a$ tel que $\vec f_a\subset \barre{\vec h}$, on a $\psi(pr_{\vec h}(a)) = pr_{\vec h}(\psi(a))$.
}
\pv Il est clair que $\psi$ envoie $\I(Q)\sph$ sur $\I(R)\sph$. Montrons qu'elle est injective sur $\I(Q)\sph$. Soient $g,h\in G$, $a,b\in A\sph$ tels que $(g,a)\sim_R (h,b)$. Alors il existe $n\in N$ tel que $b=na$ et $g\inv hn\in R(a)$. Mais $R(a)=P(a)=Q(a)$ car $R$ et $Q$ vérifient \psph. D'où $(g,a)\sim_Q(h,b)$.

Lorsque $Q$ et $R$ vérifient \plienn($\vec h$) pour une certaine facette $\vec h$, la compatibilité entre $\psi$ et $pr_{\vec h}$ est claire sur la définition.\cqfd

\fprop{Pour toute bonne famille de parahoriques $Q$, sa complétion $\bar Q$ est égale à celle de $P$.
}

\demo
Soit $\psi:\I(P)\fleche \I(Q)$ comme dans le lemme. Soit $a\in A$, $g\in\bar P(a)$ et $\vec f\in(\vec f_a^*)\sph$. Alors $g.pr_{\vec f}(a)=pr_{g\vec f}(a)$, d'où, puisque $\psi$ est $G$-équivariante et compatible aux projections, $g.pr_{\vec f}(a)=g.pr_{\vec f}(\psi(a)) = pr_{g\vec f}(\psi(a))=pr_{g\vec f}(a)$. On obtient ainsi $g\in \bar Q(a)$, puis $\bar P\subset \bar Q$.

Pour l'autre inclusion, soit $g\in\bar Q(a)$, $\vec f\in (\vec f_a^*)\sph$. Alors $\psi(g.pr_{\vec f}(a)) = g.pr_{\vec f}(\psi(a)) = pr_{g\vec f}(\psi(a)) = \psi(pr_{g\vec f}(a))$. Alors par injectivité de $\psi$ sur $\I(Q)\sph$, on obtient dans $\I(Q)$, $g.pr_{\vec f}(a) = pr_{g\vec f}(a)$. D'où $\bar Q\subset \bar P$.

Maintenant, pour toute bonne famille $Q$ de parahoriques, on a $Q\subset \bar Q=\bar P$.\cqfd

\fcor{\label{cor:bonnes familles min et max}
Supposons qu'il existe une famille de parahoriques fonctoriellement bonne pour $\D$.

 Alors $\bar P$ est l'unique bonne famille de parahoriques maximale, et $P$ est l'unique bonne famille de parahorique minimale.
 
 Ceci est en particulier le cas dans un groupe de Kac-Moody déployé.
}

\demo
 On a déjà vu que $P$ est l'unique famille de parahoriques minimale, et grâce aux corollaires \ref{cor:P verifie pdec} et \ref{cor:pdec implique plienn} qu'elle est bonne dès qu'il existe une famille de parahoriques fonctoriellement bonne pour $\D$. Concernant $\bar P$, il s'agit des deux propositions précédentes.\cqfd

\subsection{Très bonnes familles de parahoriques}
\label{sousoussection:csq de plienpn}

On donne maintenant quelques propriétés de $\I(Q)$ lorsque $Q$ est une bonne famille de parahoriques vérifiant certaines hypothèses supplémentaires. Ces hypothèses supplémentaires sont vérifiées par exemple par la famille de parahoriques construite par Guy Rousseau pour un groupe de Kac-Moody déployé (\cite{masures2}), elles le seront encore par la famille de parahoriques que nous obtiendrons pour un groupe presque déployé dans la partie \ref{section:KM}.

\subsubsection{Action de $\iv$ sur $\I$}

\fprop{\label{prop:csq de plienp} Soit $\vec f$ une facette de $\vec A$. Si $Q$ vérifie \plienpn$(\vec f)$, alors pour tout $\vec v\in\vec f$ et $a\in A$ tel que $\vec f_a\subset \barre{\vec f}$, pour tout appartement $B$ contenant $a$ et $\vec f$, on a $a+_A\vec v = a+_B\vec v$. Autrement dit, l'opération $a+\vec v$ est bien définie, et indépendante de l'appartement contenant $a$ et $\vec f$ considéré.\\

De plus, tout appartement contenant $a$ et $\vec f$ contient alors $\barre{a+\vec f}^A$, et cet ensemble est indépendant de $A$. On pourra donc le noter juste $\barre{a+\vec f}$.
}
\demo

Soit $g\in G$ tel que $g.B=A$. Alors $g\in N.Q(a)\cap N.P(\vec f)= N.Q(\barre{a+\vec f}^A)$ (adhérence dans $A$). On peut supposer $g\in Q(\barre{a+\vec f}^A)$, en particulier, $g$ fixe $a+_A\vec v$, donc $g\inv(a+_A\vec v) = a+_A\vec v$. Mais par définition de la structure affine sur les appartements, $g\inv(a+_A\vec v) = g\inv(a)+_{g\inv A} g\inv(\vec v) = a+_B\vec v$. D'où $a+_A\vec v=a+_B\vec v$.\\

Ensuite, le fait que $B=g\inv A$, avec $g\inv\in Q(\barre{a+\vec f}^A)$ entraine que $\barre{a+\vec f}^A\subset A\cap B$. Or de manière générale, $g\inv.\barre{a+\vec f}^A = \barre{g\inv a+g\inv \vec f}^{g\inv A}$, on obtient donc ici que $\barre{a+\vec f}^B=\barre{a+\vec f}^A$.\cqfd

\subsubsection{Existence d'isomorphismes entre appartements}

La condition \plienpn(sph) permet également de prouver \pisom\ pour certaines parties $\Omega$. Cependant, une condition un peu plus faible suffit, il s'agit de:\\

\plienpm($\vec f$): $\forall g\in U(\vec f)$, il existe $a\in \inte A$ tel que $g\in Q\left(\barre{a+(\vec f^*\cap \vec A})\right)$.\\

\remas{
\item Comme $U(\vec f)=\bigcap_{\vec C} U(\vec C)$, où $\vec C$ parcourt l'ensemble des chambres de $\vec f^*\cap \vec A$, un élément $u\in U(\vec f)$ fixe toujours un cône dirigé par chacune de ces chambres. Lorsque $\vec f$ est déjà une chambre, \plienpm($\vec f$) est donc toujours vérifiée. Pour une facette générale, cette condition impose que la réunion de ces cônes fixés soit connexe.

\item Cette condition est clairement vérifiée si $Q$ vérifie \plienp($\vec C$) pour toute chambre $\vec C$ de $\vec f^*\cap \vec A$.\\}

%Dans la suite, on adopte la notation suivante: pour $\vec f\in\F(\vec A)$, et pour $\Omega$ une partie d'une façade $A_{\vec g}$ avec $\vec g\subset \barre{\vec f}$, on note $\Omega+\barre{\vec f} := \bigcap_{\omega\in\Omega} \barre{\omega+\vec f}$.\\

\flemmes[\label{lemme:csq de plienpm}Soit $Q$ une bonne famille de parahoriques vérifiant \plienpm(sph).] {  
\item  Soit $A$ un appartement de $\I(Q)$, soit $\vec f\in\F(\vec A\sph)$,  $\Omega\subset A_{\vec f}$ et $g\in P(\Omega)$. Alors il existe une partie $\Omega_0\subset \inte A$ et $n\in N(\Omega)$ tels que $pr_{\vec f}(\Omega_0)=\Omega$ et $ng\in Q(\barre{\Omega_0+\vec f})$.\\

\item Soit $\vec f$ une facette sphérique. Soit $a\in A$ tel que $\vec f_a\subset \vect_{\vec A}(\vec f)$ et telle que soit $\vec f_a$ est sphérique, soit $\vec f_a$ et $\vec f$ sont de même signe. Alors:
\liste{
\item $\forall g\in Q(a)\cap P(\vec f)$, $\exists b\in a+\vec f$ tel que $g\in Q(\{a\}\cup\barre{b+\vec f}) \subset Q(\{a,pr_{\vec f}(a)\} ) $
\item $\forall g\in NQ(a)\cap P(\vec f)$, $\exists b\in a+\vec f$ tel que $g\in NQ(\{a\}\cup \barre{b+\vec f})\subset NQ(\{a,pr_{\vec f}(a)\} )$.\\
}

\item La projection $pr_{\vec f}(a)$ d'un point $a\in I$ sur une façade $\I_{\vec f}$ avec $\vec f$ sphérique est bien définie lorsqu'il existe un appartement $\vec A$ contenant $\vec f_a\cup \vec f$, tel que $\vec f_a\subset \vect_{\vec A}(\vec f)$ et $\{$ $\vec f_a$ est de même signe que $\vec f$ ou $\vec f_a$ est sphérique$\}$. On a alors pour tout $g\in G$, $g.pr_{\vec f}(a) = pr_{g\vec f}(ga)$.\\
}

Le deuxième point constitue un renforcement de \plien(sph) et de \plienn(sph),  non seulement  car il prouve qu'un cône dirigé par $\vec f$ est fixé, mais aussi car il s'applique lorsque $\vec f_a\subset \vect(\vec f)$ et non seulement lorsque $\vec f_a\subset \barre{\vec f}$.\\
 
\pv
On sait par \ref{prop:csq de parasph} que $P(\Omega)=N(\Omega).U(\vec f).G(\phi^m(\vec f),\Omega)$. Soient $n\in N(\Omega)$, $u\in U(\vec f)$ et $m\in G(\phi^m(\vec f),\Omega)$ tels que $ng=um$. Par \plienpm($\vec f$), il existe $a\in A$ tel que $u$ fixe $\barre{a+\vec f^*\cap \vec A}$, et ceci contient une partie de la forme $\barre{\Omega_0+\vec f}$ avec $\Omega_0\subset\inte A$ et $pr_{\vec f}(\Omega_0)=\Omega$. Le facteur $m$ fixe $\eng{\Omega,\vec f}_A$ qui contient $\barre{\Omega_0+\vec f}$. Ceci prouve le premier point.\\

Soient $\vec f$ et $a$ comme dans l'énoncé du second point. Soit $g\in Q(a)\cap P(\vec f)$. En particulier, $g\in P(\vec f_a \cup\vec f)=P(\cl(\vec f_a\cup \vec f))$ par \ref{props:paraboliques}. Soit $\vec h=pr_{\vec f_a}(\vec f)\subset \cl(\vec f_a\cup \vec f)$, alors $\vec h$ est une facette sphérique et $\vect(\vec h)=\vect(\vec f)$. Par \plien($\vec h$), $g\in Q(pr_{\vec h}(a))$. Par le premier point, il existe $n\in N(pr_{\vec h}(a))$ et $a_0\in \inte A$ tels que $pr_{\vec h}(a_0)=pr_{\vec h}(a)$ et $ng\in Q(\barre{a_0+\vec h})$. Le cône $\barre{a_0+\vec h}$ contient un sous cône de la forme $\barre{b+\vec h}$ avec $b\in a+\vec h$. En particulier, $b\in A_{\vec f_a}$ donc $\barre{b+\vec h}=\barre{b+\vec f}$. Finalement, $g\in Q(a)\cap N(pr_{\vec h}(a)).Q(\barre{b+\vec f})$. Modulo $G(\phi^m(\vec h),a)$, on peut supposer $g\in Q(a)\cap N(pr_{\vec h}(a)).(U(\vec h)\cap Q(\barre{b+\vec f}))$. Alors le lemme \ref{lemme:action de u sur A} entraine que le facteur dans $N(pr_{\vec h}(a))$ fixe $a$, et donc $\barre{a+\vec h}$. Ainsi, $g\in Q(a)\cap Q(\barre{b+\vec f})$.\\

Si maintenant $g\in NQ(a)\cap P(\vec f)$, alors en particulier $g\in N.P(\vec f_a)\cap P(\vec f)$ et d'après \ref{props:paraboliques}, ceci vaut $N.P(\cl(\vec f_a\cup \vec f))$. On peut donc supposer $g\in P(\cl(\vec f_a\cup \vec f))$, en particulier, $g\in P(\vec h)$ en notant encore $\vec h=pr_{\vec f_a}(\vec f)$. Alors par \plienn($\vec h$), $g\in N(\vec h).Q(\{a,pr_{\vec h}(a)\})$. Mais nous venons de voir que tout élément de $Q(\{a,pr_{\vec h}(a) \})$ fixe un cône de la forme $\barre{b+\vec f}$ avec $b\in a+\vec f$.\\

Le troisième point se prouve alors tout comme \ref{prop:projection bien def}.\cqfd

\fprop{ \label{prop:isom entre appart} %Soit $Q$ une famille de parahoriques vérifiant \plienpn et \psph.
Soit $Q$ une bonne famille de parahoriques vérifiant \plienpm(sph).
 Soit $\Omega$ une partie de $A$ contenant au moins un point sphérique et incluse dans $A^+$ ou dans $A^-$, ou contenant au moins un point sphérique positif et un point sphérique négatif.

 Alors:
\[\bigcap_{\omega\in\Omega} N.Q(\omega) = N.Q(\Omega) \]
}

\demo

 Soit $g\in \bigcap_{\omega\in\Omega} N.Q(\omega)$. En particulier, $g\in \bigcap_{\vec f\in \vec\Omega}N.P(\vec f) = N.P(\vec \Omega)=N.P(\cl(\vec \Omega))$ par la proposition \ref{props:paraboliques}, on peut donc supposer $g\in P(\cl(\vec \Omega))$, d'où $g\in \bigcap_{\omega\in\Omega} N(\vec f_\omega).Q(\omega)$.\\
 
 On traite le cas où $\vec\Omega$ contient des facettes sphériques positives et négative, le cas où $\vec\Omega\subset \vec A^\pm$ étant plus simple.
 Soient $\vec f^+$ et $\vec f^-$ des facettes maximales de $\cl(\vec\Omega)\cap \vec A^+$ et $\cl(\vec\Omega)\cap \vec A^-$, respectivement. Ces facettes sont sphériques, $pr_{\vec f^+}(\vec f^-)=\vec f^+$, et $pr_{\vec f^-}(\vec f^+)=\vec f^-$, d'où $\vect(\vec f^+)=\vect(\vec f^-)$ et $N(\vec f^+)=N(\vec f^-)=N(\vec\Omega)$. Notons $\Omega^+=\Omega\cap A^+$ et $\Omega^-=\Omega\cap A^-$.
 
 D'après le deuxième point du lemme, $g\in \bigcap_{\omega\in\Omega^+}N(\vec f^+).Q(pr_{\vec f^+}(\omega))$, et ceci vaut $N(\vec f^+).Q(pr_{\vec f^+}(\Omega^+))$, par \ref{prop:csq de parasph}. De même, $g\in N(\vec f^-).Q(pr_{\vec f^-}(\Omega^-))$.
 
 Par le premier point du lemme, il existe $\Theta_0^-\subset \inte A$ telle que $pr_{\vec f^-}(\Theta_0^-) = pr_{\vec f^-}(\Omega^-)$ et $g\in N(\vec f^-).Q(\barre{\Theta_0^- + \vec f^-})$. Appliquant alors \plienn($\vec f^+$) sur la partie $\Theta_0^-$, puis \ref{prop:csq de parasph} 4, on trouve $g\in N(\vec f^+).Q(pr_{\vec f^+}( \Omega^+\cup \Theta_0^-))$. Mais $\vect(\vec f^+)=\vect(\vec f^-)$ donc $pr_{\vec f^+}(\Omega^+\cup \Theta_0^-) = pr_{\vec f^+}(\Omega)$.
 
 Finalement, on peut supposer $g\in Q(pr_{\vec f^+}(\Omega))$, puis encore par le lemme, qu'il existe une partie $\Omega_0^+\subset \inte A$ et $n^+\in N(\vec f^+)$ tels que $pr_{\vec f^+}(\Omega_0^+)= pr_{\vec f^+}(\Omega)$ et  $n^+g$ fixe $\barre{\Omega_0^+ +\vec f^+}$.
 De même, il existe $\Omega_0^- \subset \inte A$ et $n^-\in N(\vec f^-)$ tels que $pr_{\vec f^-}(\Omega_0^-) = pr_{\vec f^-}(\Omega)$ et $n^-.g$ fixe $\barre{\Omega_0^-+\vec f^-}$.\\
 
 On peut supposer $n^+=e$. Soit $\omega\in\Omega_0^-+\vec f^-$. Alors $g\in N.Q(\omega)\cap P(\vec f^+) = N.Q(\{\omega,pr_{\vec f^+}(\omega)\})$, et il existe d'après le lemme un point $\omega'\in\inte A$ tel que $g\in N.Q(\{\omega\}\cup \omega'+\vec f^+)$. Soit $n\in N$ tel que $ng\in Q(\{\omega\}\cup \omega'+\vec f^+)$. L'intersection $(\omega'+\vec f^+)\cap (\Omega_0^+ +\vec f^+)$ contient un scp de $\omega'+\vec f^+$, car $\Omega_0^+$ contient un point dont la projection sur $A_{\vec f^+}$ est la même que celle de $\omega$ et $\omega'$. Comme $g$ fixe $\Omega_0^+ +\vec f^+$, $n$ fixe l'intersection $\omega'+\vec f^+\cap \Omega_0^+ +\vec f^+$, et en particulier ce scp. D'où, puisque $n$ agit par automorphisme affine sur $\inte A$, $n$ fixe $\omega$. Finalement, $g\in Q(\omega)$.
 
 Ainsi, $g$ fixe $(\Omega_0^+ +\vec f^+) \cup (\Omega_0^-+\vec f^-)$.\\
 
 Enfin, soit $\omega\in\Omega$,  soit $\epsilon$ un signe de $\omega$. Soit $\vec h=pr_{\vec f_\omega}(\vec f^\epsilon)$ (bien défini car $\vec f_\omega$ et $\vec f^\epsilon$ sont de même signe), alors $g\in N.Q(\omega)\cap P(\vec h) = N.Q(\{\omega,pr_{\vec h}(\omega)\})$, et il existe $\omega'\in \omega+\vec h$ et $n\in N(\vec f_\omega)$ tels que $ng$ fixe $\barre{\omega'+\vec h}$.
  D'autre part, $g$ fixe $pr_{\vec h}(\Omega_0^\epsilon) = pr_{\vec h}(\Omega)$, par \plien(sph), donc $g\in N(\vec h).U(\vec h).G(\phi^m(\vec\Omega),\Omega)$. Le fait que $g$ fixe $\Omega_0^\epsilon$ entraine avec le lemme \ref{lemme:action de u sur A} que le facteur dans $N(\vec h)$ fixe $\Omega_0^\epsilon$, et donc $\barre{\Omega_0^\epsilon +\vec h}$, donc $g$ fixe un $\barre{\Omega_{\vec h} +\vec h}$ avec $\Omega_{\vec h}\subset \inte A$, $pr_{\vec h}(\Omega_{\vec h}) = pr_{\vec h}(\Omega)$, et ceci contient en particulier un scp de $\omega'+\vec h$. On en déduit alors, comme au paragraphe précédent, $n\in N(\omega)$ d'où $g\in Q(\omega)$.
  
  On a bien prouvé $g\in Q(\Omega)$.\cqfd

\fcor{ \label{cor:isom entre apparts}
%Si $Q$ vérifie \plienpn et \psph,
Si $Q$ est une bonne famille de parahoriques vérifiant \plienpm(sph), et si $A$ et $B$ sont deux appartements de $\I(Q)$ dont l'intersection contient un point sphérique de signe $\epsilon$, alors il existe $g\in G$ tel que $g.A=B$ et $g$ induit de $A$ sur $B$ un isomorphisme fixant $A^\epsilon\cap B^\epsilon$.

Si $A\cap B$ contient un point sphérique de chaque signe, alors il existe $g\in G$ tel que $g.A=B$ et $g$ fixe $A\cap B$.\\

Autrement dit, si $\Omega$ est une partie d'un appartement $A$, contenant un point sphérique et incluse dans $A^\pm$ ou bien contenant un point sphérique de chaque signe, alors $Q(\Omega)$ est transitif sur les appartements contenant $\Omega$.\\
}

\fcor{ \label{cor:enclos bien def} Soit $\Omega$ une partie d'un appartement $A$ contenant au moins un point sphérique positif et un point sphérique négatif.
Alors pour tout appartement $B$ contenant $\Omega$, il existe $g\in G$ qui induit un isomorphisme de $A$ sur $B$ fixant $\cl_A(\Omega)$.

En particulier, la partie $\cl(\Omega)$ est bien définie, indépendamment de l'appartement contenant $\Omega$ considéré.
}

\demo C'est la concaténation de \ref{cor:isom entre apparts} et de \ref{cor:enclos equilibre}.\cqfd

\subsubsection{Intersection d'appartements}
\label{soussection:cloture}

On passe maintenant à l'étude de \pinter.

 On a vu en s'appuyant sur la décomposition de Lévi de $Q(\Omega)$, que cette condition est vérifiée pour des parties $\Omega\subset A$ contenant un point sphérique positif et un point sphérique négatif, dès que $Q$ est une bonne famille de parahoriques (corollaire \ref{cor:enclos equilibre}).
 Nous allons maintenant étudier le cas d'une partie $\Omega$ composée de points sphériques d'un même signe. On ne dispose pas pour le fixateur d'une telle partie de la décomposition de Lévi, mais on prouve tout de même \pinter, sous l'hypothèse \plienpm(sph).

\begin{prop}\label{prop:intersection close, cas non equilibre}
Soit $\Omega\subset A\sph$. Soit $Q$ une bonne famille de parahoriques vérifiant \plienpm. Alors:
 \[Q(\Omega)=N_\Omega.Q(\cl(\Omega)\cap A\sph)\point \]
% De plus, tout $g\in Q(\Omega)$ fixe modulo $N$ une partie de la forme $\bigcup_{\omega\in\Omega_0} \barre{\omega+\cl(\vec\Omega)}$, avec $\Omega_0$ une partie de $\inte A$ telle que pour tout $\vec f\in\F(\vec\Omega)$, $pr_{\vec f}(\Omega_0)= \cl(\Omega)\cap A_{\vec f}$.\\

  %voisinage de $\cl(\Omega)\cap A\sph$ dans $\aff(\cl(\Omega))$.\\
\end{prop}
%\rema C'est en fait, dans le cas d'un groupe de Kac-Moody déployé, la proposition 3.9 de \cite{microaffine}, adaptée au contexte présent.\\

\demo\\
L'inclusion $"\supset"$ est claire. Le cas où $\Omega$ coupe $A\sph^+$ et $A\sph^-$ découle de \ref{prop:levi equilibre}, on suppose donc $\Omega\subset A\sph^+$.\\

Soit $g\in Q(\Omega)$. En particulier, $g\in P(\cl(\vec\Omega))$ par \ref{props:paraboliques}. Par le lemme \ref{lemme:csq de plienpm}, il existe un choix $C\in\maj C$ tel que défini dans le lemme \ref{lemme:trace dun enclos dans une facade spherique}, tel que $g\in\bigcap_{\omega\in\Omega'(\maj C)} N_{\vec f_\omega}.Q(\omega)$, où $\vec f_\omega$ est la direction de la façade contenant $\omega$, et où $\Omega'(\maj C)$ est l'intersection de $\Omega^\infini(\maj C)$ avec l'union des façades sphériques de direction incluse dans $\cl(\vec\Omega)$. Lorsque $\vec f$ est une facette sphérique de $\cl(\vec\Omega)$, on obtient alors $g\in N_{\vec f} . Q(\cl(\Omega^\infini(C)\cap A_{\vec f})\cap A_{\vec f})$. Or la partie $\cl(\Omega^\infini(C)\cap A_{\vec f})\cap A_{\vec f}$ contient $\cl(\Omega)\cap A_{\vec f}$, d'après le lemme \ref{lemme:trace dun enclos dans une facade spherique}, d'où $g\in N_{\vec f}\point Q(\cl(\Omega)\cap A_{\vec f})$.\\

Nous voulons maintenant prouver le résultat similaire lorsque $\vec f$ est une facette sphérique de $\fl{\cl(\Omega)}$. Supposons qu'il existe $\vec f$ une telle facette non incluse dans $\cl(\vec\Omega)$. Quitte à projeter $\vec f$ sur une facette adéquate de $\cl(\vec \Omega)$, on peut supposer qu'une face $\vec m$ de $\vec f$ de codimension 1 est incluse dans $\cl(\vec\Omega)$. La facette $\vec m$ est dans l'intérieur de l'enclos de $\vec\Omega\cup\vec f$ (pour la topologie induite), elle est donc sphérique. Le fait que $\vec f\not\subset\cl(\vec\Omega)$ signifie qu'il existe $\alpha\in \phi(\vec\Omega)$ telle que $\alpha(\vec f)<0$. Tout demi-appartement dirigé par $\alpha$ ne contient aucun point de $A_{\vec f}$, or $\vec f\subset \fl{\cl(\Omega)}$, ceci implique qu'aucun demi-appartement dirigé par $\alpha$ ne contient $\Omega$. Il existe donc $(\omega)\in\Omega^\N$ telle que $\alpha(\omega_n)\fleche_{n\fleche \infini}-\infini$ et $\alpha(\omega_n)\in\R$ pour tout $n$ (autrement dit, $\vec f_{\omega_n}\subset \ker \alpha$ pour tout $n$). La facette $\vec m$ est maximale dans $\cl(\vec\Omega)\cap \ker(\alpha)$, donc tous les $\omega_n$ ont un projeté $\omega'_n$ dans $A_{\vec m}$. Comme $\vec m$ est sphérique, $\barre{A_{\vec m}}$ est compact, et la suite $(\omega')$ a une valeur d'adhérence $\omega_\infini$ dans $\barre{A_{\vec m}}$. Comme $\alpha(\omega_\infini)=-\infini$, $\omega_\infini$ est dans une facette $\vec g\not\subset \cl(\vec\Omega)$, sphérique. Or, par \ref{lemme:csq de plienpm}, $g$ fixe tous les $\omega'_n$, $n\in\N$, puis comme $\vec g$ est sphérique, $g$ fixe $\cl\{\omega'_n|n\in\N\}$. Donc $g$ fixe un cône $c\subset A_{\vec m}$ tel que $\omega_\infini=[c]_{A_{\vec m}}$. Ainsi $g\in Q(\omega_\infini)\subset P(\vec g)$, on peut donc rajouter $\omega_\infini$ à $\Omega$, et $\vec g$ à $\vec\Omega$. On arrive ainsi à prouver que $g\in P(\fl{\cl(\Omega)})$.\\

Alors il existe un choix $C\in\maj C$ et $\Omega^\infini(C)$, comme définis cette fois à la proposition \ref{prop:trace dun enclos dans une facade}, tels que $g\in \bigcap_{\omega\in \Omega^\infini(C)} N_{\vec f_\omega} Q(\omega)$. Donc pour chaque facette sphérique de $\fl{\cl(\Omega)}$, $g\in N_{\vec f}.Q(A_{\vec f}\cap \cl(A_{\vec f}\cap \Omega^\infini(C)))$, or 
$A_{\vec f}\cap \cl(A_{\vec f}\cap \Omega^\infini(C))$ contient $\cl(\Omega)\cap A_{\vec f}$, d'après la proposition \ref{prop:trace dun enclos dans une facade}. On conclut alors grâce à la proposition \ref{prop:isom entre appart} que $g\in N_\Omega.Q(\cl(\Omega)\cap A\sph)$.\cqfd

%Ensuite, par \ref{prop:csq de parasph} et \ref{prop:trace dun enclos dans une facade}, on obtient que $g$ fixe modulo un $n\in N_\Omega$, une partie $V$ contenant $\cl(\Omega)\cap A\sph$ ainsi qu'un voisinage de $\cl(\Omega)\cap A\sph$ dans chaque $\aff_A(\cl(\Omega)\cap A_{\vec f}$ pour $\vec f\subset \fl{\cl(\Omega)}\setminus \vec A\sph$.\cqfd\\

\section{Descente}
\label{section:descente}

Il est facile de généraliser au cas présent le théorème 9.2.10 de \cite{bruhat-tits} qui permet sous certaines hypothèses de descendre la valuation de $\D$ à un sous-groupe. Pour obtenir une masure bordée pour ce sous-groupe, il faudra ensuite descendre la famille de parahoriques.
  Ce sont bien sûr ces résultats qui serviront pour définir une valuation dans un groupe de Kac-Moody presque déployé.

\subsection{Contexte et notations}
\label{soussection:hypotheses de descente}
On fixe un tore maximal $T$ et on note $\phi=\phi(T)$, $\vec V=\vec V(T)$, $\vec A=\vec A(T)$. On fixe également une bonne famille de parahoriques $Q$ et on note $\I=\I(Q)$.\\

On se donne un autre système de racines $\phi\bec\subset (\vec V\bec)^*$ avec $\vec V\bec$ un sous-espace de $\vec V$. On notera $W(\phi\bec)\subset Gl(\vec V\bec)$ le groupe de Weyl et $\vec A\bec\subset \vec V\bec$ l'appartement correspondants. On notera aussi $\vec A\becb= \vec A\cap \vec V\bec$, la trace de l'appartement $\vec A$ sur $\vec V\bec$. On adopte la convention de noter en lettres latines les racines de $\phi\bec$. On suppose:\\
\liste{
\item (DSR): $\phi\bec$ est un système de racines à base libre, et %$\forall a\in\phi\bec$, $\ker a\cap \vec A\sph\not=\vide$.
il existe une chambre $\vec C\bec$ de $\vec A\bec$ telle que toute facette de $\vec C\bec$ rencontre $\vec A$ et toute facette sphérique de $\vec C\bec$ rencontre $\vec A\sph$.
\\}

On notera $\Pi\bec$ la base de $\phi\bec$ correspondant à la chambre $\vec C\bec$ donnée par (DSR).
On ne suppose pas que $\phi\bec$ engendre $(\vec V\bec)^*$, autrement dit $(\phi\bec,\vec V\bec)$ n'est pas forcément essentiel.\\

 On se donne ensuite une donnée radicielle $\D\bec=(G\bec,(U_a\bec)_{a\in\phi\bec})$ de système de racines $\phi\bec$. On notera avec un $\natural$ en exposant tous les objets associés à la donnée radicielle $\D\bec$, par exemple $N\bec$, $T\bec$ ...

 Pour tout $a\in\phi\bec$, $k\in\R$, on note:
 \calc{
 \phi_a &:=& \ens{\alpha\in\phi}{\alpha|_{\vec V\bec}\in\R^{+*}.a}\\
 \phi_0 &:=& \ens{\alpha\in\phi}{\alpha(\vec V\bec) =0}\\
\vec D(a) &:=& \ens{x\in\vec A\becb}{a(x)\geq0}\\
  U_a &:=& U(\vec D(a)) = G(\phi_a) = \engens{U_\alpha}{\alpha\in\phi_a}\\
   U_{a,k} &:=& \engens{ U_{\alpha,rk}}{\alpha\in\phi_a,\, r\in\R^{+*}\et \alpha_{|\vec V\bec} = ra}\\
   Z &:= & P(\vec A\becb) = \eng{T, \ens{U_\alpha}{\alpha\in\phi\et\alpha(\vec V\bec)=0}}\point
   }

  Pour toute racine simple $a\in \Pi\bec$, l'ensemble $\phi_a$ n'est autre que $\phi^u(\vec D(a))$.  
   De plus, $P(\ker(a)\cap \vec A) = M_{\vec A}(\ker(a)\cap \vec A) = \eng{Z,U_a,U_{-a}}$, c'est un groupe muni d'une donnée radicielle de système de racine $\phi_a\cup \phi_{-a}\cup \phi_0$.\\

  On suppose vérifiées les conditions suivantes de compatibilité des données radicielles:\\
  
  \liste{
	\item (DDR 1): $G\bec\subset G$ et  $\forall a\in \phi\bec$, $U_a\bec\subset U_a$.
	\item (DDR 2): $\forall a\in \phi\bec$ tel que $2a\in\phi\bec$, $\card\ens{\alpha_{|\vec V\bec}} {\alpha\in\phi\et \alpha_{|\vec V\bec}\in\R^{+*}.a} \leq 2$.
	\item (DDR 3.1): $T\bec\subset Z$. % et $T\bec$ normalise chaque $U_a$, $a\in\phi\bec$.
	\item (DDR 3.2): Pour tout $a\in\phi\bec$, et $u\in U_a\bec$, $n\bec(u).Z .n\bec(u)\inv = Z$, $\: n\bec(u)U_a n\bec(u)\inv=U_{-a}$ et $\: n\bec(u)U_{-a}n\bec(u)\inv =U_a$.\\
	}
	
	Remarquons que (DDR 1) entraine que pour tout $a$, $\phi_a\not=\vide$. En particulier, il existe $\alpha\in\phi$ tel que $\ker(a)= \ker(\alpha)\cap \vec V\bec$, et même $\ens{x\in\vec V\bec}{a(x)\geq0} = \ens{x\in\vec V\bec}{\alpha(x)\geq 0}$.\\
	
Concernant l'immeuble vectoriel, on suppose:\\

\liste{
\item (DIV): $G\bec.\vec A\becb \cap \vec A = \vec A\becb$.\\}	
	
On se donne encore une partie $\I\die$ de $\I$. On note $A\die=\inte A\cap \I\die$, et on suppose:\\

\liste{
%\item (DIV): $\iv\die$ est stable par $G\bec$ et $\iv\die\cap \vec A=\vec V\bec\cap \vec A$.\\

\item (DM 1): $\I\die$ est stable par $G\bec$ et pour toute facette sphérique $\vec f\subset \iv$, $I\die\cap\I_{\vec f}$ est une partie convexe de $\I_{\vec f}$.
\item (DM 2): $A\die$ est un sous-espace affine de $\inte A(T)$ dirigé par $\vec V\bec$, et $A\cap \I\die$ est l'adhérence dans $A$ de $A\die$.
\item (DM 3): Pour toute facette sphérique $\vec f$ rencontrant $\vec A\bec$, il existe une facette $F$ de $A_{\vec f}$, rencontrant $\I\die$, qui n'est pas contenue dans l'adhérence d'une autre facette $F'$ de $\I_{\vec f}$ rencontrant $\I\die$.
%Il existe une facette $F$ de $\I$ coupant $A\die$ qui n'est contenue dans l'adhérence d'aucune autre facette de $\I$ rencontrant $\I\die$.
\item (DM 4): $\barre{A\die}$ est stable par $N\bec$.\\
	}
	
	La condition (DM 3) est non triviale: par exemple si $A\die$ est un mur, il pourrait exister un autre appartement $Z$ contenant $A\die$ tel que $Z\cap \I\die$ soit une bande dans $Z$, contenant $A\die$ dans son intérieur. Alors si $F$ est une facette de $A$ maximale dans $A\die$ (une cloison), elle est dans l'adhérence d'une chambre de $Z$ coupant $\I\die$.
	
Grâce à (DDR1), on définit pour tout $a\in\phi\bec$, et $u\in U_a\bec$,
\[\phii\bec_a(u):= \sup \ens{k\in\R\cup\{\infini\}}{u\in U_{a,k}}\point \]

On suppose enfin les deux hypothèses suivantes concernant les valuations:\\
\liste{	\item (DV1): "$\phii\in A\die$" autrement dit, le point $o\in A(T)$ tel que $\forall \alpha\in\phi$, $U_{\alpha,0}$ fixe $o$ est dans $A\die$.
	\item (DV 2): $\forall a\in \phi\bec$, $\card(\phii\bec_a(U\bec_a))\geq 3$.\\
	}

Les conditions (DV 1) et (DM 2) imposent donc que $A\cap \I\die = \barre{o+\vec V\bec}$.

Hormis le fait qu'on ne suppose pas $\phi$ fini, ces hypothèses sont plus fortes que celles de \cite{bruhat-tits}: on a rajouté (DSR), (DIV)
 et (DM 4). (Par ailleurs, on a renommé les (DI x) en (DM x).)\\

\subsection{Descente dans l'immeuble vectoriel}

Conformément à nos notations, $\vec A\bec$ est le cône de Tits dans $\vec V\bec$ défini par le système de racines $\phi\bec$. En général, on verra qu'il est plus grand que $\vec A\becb$, ce qui empêche immédiatement de plonger l'immeuble $\iv\bec$ dans $\iv$. Cependant, nous allons voir que $\vec A$ intersecte chaque facette de $\vec A\bec$, ce qui permettra quand même d'identifier à l'intérieur de $\iv$ une réalisation de l'immeuble de $\D\bec$. Cette réalisation sera notée $\vec\I\becb$. %(à distinguer du $\iv\die$ de (DIV)). %pour l'instant on introduit la notation $\vec A\becb:= \vec A\cap \vec A\bec$.

%\flemme{ Le groupe $N\bec$ stabilise $\vec A\cap \iv\die$, et donc $N\bec\subset N.Z$.
%}

\flemme{\label{lemme:descente vect}
Soit $a\in\Pi\bec$ une racine simple. Alors:\liste{
\item Il existe une partie équilibrée $\vec\Omega$ de $\vec A$ telle que $\vec\Omega\subset \vec D(a) \subset \cl_{\vec A}(\vec\Omega)$,
\item $P(\vec D(a)) = Z\ltimes U_a$,
\item Pour tout $u\in U_a\bec\setminus\{e\}$, $n\bec(u)$ permute $\vec D(a)$ et $\vec D(-a)$,
\item $N\bec \subset (N\cap \stab(\vec A\becb)) \point Z$.
}
En particulier, l'action de $N\bec$ sur $\iv$ stabilise $\vec A\becb$.}

\rema Ceci entraine en particulier que $\phi_a$ est fini, et permet de prouver que les groupes $U_{a,\lambda}$ définis plus haut sont les mêmes que ceux définis à partir de $\phii\bec$ comme dans la définition \ref{def:valuation}.\\

\pv 

  Soit $\vec m\bec$ la cloison de $\vec C\bec$ telle que $\ker(a)= \vect_{\vec V\bec}(\vec m\bec)$, alors d'après (DSR), $\vec m\bec$ coupe $\vec A\sph$. Comme $\vec A\sph$ est un cône ouvert dans $\vec V$, il existe un cône $\vec m \subset \vec m\bec$, inclus dans une facette sphérique de $\vec A$ et contenant un ouvert de $\vec m\bec$. Alors l'enveloppe convexe de $\vec m\cup -\vec m$ est $\ker(a)$. Soit encore $\vec x\in\vec C\bec\cap \vec A\sph$, alors $\vec\Omega := \vec m\cup -\vec m\cup \{\vec x\}$ est une partie équilibrée, incluse dans $\vec D(a)$, et $\vec D(a)\subset \cl_{\vec A}(\vec \Omega)$.
  
  Le sous-groupe parabolique $P(\vec D(a))=P(\vec\Omega)$ admet donc une décomposition de Lévi. Son facteur de Lévi est $\fix_G(\vec D(a)\cup - \vec D(a)) = \fix_G(\vec A\becb) = Z$, et son facteur unipotent est $U_a$, par définition. D'où $P(\vec D(a)) = Z\ltimes U_a$.

%La condition (DDR 3.1) entraine directement $T\bec\subset Z\subset (N\cap \stab(A\becb)) .Z$. Il reste à étudier les $n\bec(u)$.\\

Soit $u\in U\bec_a\setminus\{e\}$. Alors (DDR 3.2) entraine que $n\bec(u)$ conjugue $P(\vec D(a))$ en $P(\vec D(-a))$. Mais $P(\vec D(a))$ est un sous-groupe parabolique de $M(\ker(a))$ (correspondant à la facette contenant $\vec D(a)/\ker(a)$). Choisissons un Borel $B^+$ pour $M(\ker(a))$ dans $P(\vec D(a))$, et un autre $B^-$ dans $P(\vec D(-a))$, alors il existe $n\in N\cap M(\ker(a))$ tel que $B^-=nB^+n\inv$. Alors $n.P(\vec D(a))n\inv$ et $P(\vec D(-a)) = n\bec(u).P(\vec D(a)).n\bec(u)\inv$ sont deux paraboliques conjugués dans $M(\ker(a))$ (remarquer que $n\bec(u)\in M(\ker(a))$ ) contenant un même Borel: ils sont égaux. (On a en fait prouvé que les deux facettes contenant $\vec D(a)/\ker(a)$ et $\vec D(-a)/\ker(a)$ sont de même type.)

Dès lors, $nn\bec(u)\inv \in M(\ker(a))$ stabilise les deux paraboliques $P(\vec D(a))$ et $P(\vec D(-a))$. Donc $nn\bec(u)\inv \in P(\vec D(a))\cap P(\vec D(-a)) = Z$ et $n\bec(u)\in N.Z$.

Donc $n\bec(u).\vec A\becb = n.\vec A\becb\subset \vec A$. Avec l'hypothèse (DIV), on obtient en plus $n\bec (u).\vec A\becb \subset G\bec.\vec A\becb \cap \vec A = \vec A\becb$. Donc $n\bec (u)$ stabilise $\vec A\becb$. Comme $Z$ fixe $\vec A\becb$, l'élément $n\in Z.n\bec(u)$ stabilise aussi $\vec A\becb$.

Enfin, le fait que $n\bec(u)$ échange $P(\vec D(a))$ et $P(\vec D(-a))$ entraine que $n\bec (u)$ échange $\cl_{\vec A}(\vec D(a))$ et $\cl_{\vec A}(\vec D(-a))$, mais $\cl_{\vec A}(\vec D(a))\cap \vec A\becb = \vec D(a)$, et similairement pour $-a$, donc $n\bec(u)$ échange $\vec D(a)$ et $\vec D(-a)$.\\

Comme $N\bec =\eng{T\bec,\ens{n\bec(u)}{a\in\Pi\bec,\; u\in U_a\bec}}$, et $T\bec\subset Z$ par (DDR 3.1), l'inclusion $N\bec\subset (N\cap \stab(\vec A\becb)) .Z$ est maintenant claire.
\cqfd

%Le groupe $T\bec$ est inclus par (DDR 3.1) dans $Z$ qui fixe $\vec A\cap \vec V\bec$. Et $\vec A\cap \vec V\bec = \vec A\cap \iv\die$ par (DIV). Il ne reste qu'à étudier les éléments de $N\bec$ de la forme $n\bec(u)$, $u\in U_a\bec$, $a\in\phi\bec$.

%Soient donc $a\in\phi\bec$, $u\in U_a\bec$ et $n=n\bec(u)$.
% Le fixateur dans $G$ du demi-cône $\vec D(a)\subset \vec A\cap \vec V\bec\subset\iv$ est $P(\vec D(a)) = Z.U_a$. Donc l'hypothèse (DDR 3.2), entraine que $n.P(\vec D(a)) .n\inv = P(\vec D(-a))$, et ainsi, dans $\iv$, $n$ échange $\cl(\vec D(a))$ avec $\cl(\vec D(-a))$. Or $\iv\die\cap \vec A \subset \cl(\vec D(a))\cup \cl(\vec D(-a)) \subset \vec A$, donc $n.(\vec A\cap \iv\die) \subset n.\iv\die\cap \vec A=\vec A\cap \iv\die$, par (DIV).\\
 
% Donc $N\bec$ stabilise $\vec A\cap \iv\die$. Mais cette partie de $\iv$ est l'enclos d'une partie équilibrée, car par (DSR) toute chambre de $\vec A\bec$ rencontre $\iv\sph\cap \vec A$. Donc le fixateur $Z$ de $\iv\die\cap \vec A$ est transitif sur les appartements de $\iv$ contenant $\iv\die\cap \vec A$. D'où l'égalité $N\bec\subset N.Z$.
%\cqfd

\fdef{ Pour tout $n\in N\bec$, soient $n'\in N\cap \stab(\vec A\becb)$ et $z\in Z$ tels que $n=n'z$. On pose $\vec\nu\becb(n) = \vec\nu(n')|_{\vec V\bec}$.}

 Ceci est bien défini car $N\cap Z$ fixe $\vec V\bec$. Notons qu'il s'agit juste de la restriction à $\vec A\becb$ de l'action de $n$ sur $\iv$, étendue par linéarité à $\vec V\bec$, ce qui est possible car $\vec A\becb$ engendre $\vec V\bec$, par (DSR).

\fprops{\label{props:action de Nbec}
\item Pour tout $n\in N\bec$, $\vec\nu\becb(n)$ stabilise $\vec V\bec$ et $\vec A\bec$. L'application $\vec\nu\becb$ est une action de groupe, elle stabilise l'ensemble des facettes de $\vec A\bec$, et induit sur cet ensemble la même action que $\vec\nu\bec$.

 \item %$\vec A\cap \vec A\bec= \vec V\bec\cap \vec A$ et 
 Chaque facette $\vec f\bec$ de $\vec A\bec$ rencontre $\vec A$, et $\vec f\bec\cap \vec A$ est l'intersection de $\vec V\bec$ avec une réunion de facettes de $\vec A$. Si $\vec f\bec$ est sphérique, elle rencontre $\vec A\sph$.
}

\rema En conséquence, la condition (DSR) est vérifiée pour n'importe quelle chambre $\vec C\bec$ de $\vec A\bec$, et le lemme précédent est vrai pour n'importe quelle racine $a\in\phi\bec$ (et non seulement  $a\in\Pi\bec$).\\

\demo\\

Soit $t\in T\bec$. Par (DDR 3.1), $t\in Z$ donc $\vec\nu\becb(t)=id= \vec\nu\bec(t)$. Ainsi $\vec\nu\becb$ et $\vec\nu\bec$ coïncident sur $T\bec$.

De plus, la description de $\vec\nu\becb$ comme extension par linéarité de l'action de $N\bec$ sur $\vec A\becb$ montre qu'il s'agit d'une action de groupe.

Ensuite, soit $a\in\Pi\bec$, $u\in U_a\setminus\{e\}$ et $n=n\bec(u)$. Comme $\vec A\becb$ engendre $\vec V\bec$, par la condition (DSR), le lemme montre que $\vec\nu\becb(n)$ stabilise $\vec V\bec$. De plus, il fixe $\ker(a)$ et échange les demi-espaces délimités par cet hyperplan. Enfin, $\vec\nu\becb(n)^2= \vec\nu\becb(n^2)$ mais $n^2\in T\bec$ donc $\vec\nu\becb(n^2)=Id_{\vec V\bec}$, par la première phrase de cette preuve. Ceci prouve que $\vec\nu\becb(n)$ est une réflexion d'hyperplan $\ker(a)$.\\

Montrons maintenant que $\vec\nu\becb(n)$ préserve l'ensemble des murs de $\vec A\bec$. Soit $b\in\Pi\bec$. D'après le lemme, $Z\cap U_b=\{e\}$. Soit $u_b\in U_b\bec\setminus \{e\}$, alors $u_b\not\in Z$. L'élément $u_b$ fixe $\vec D(b)$,
% qui est un demi-cône de $\vec A$ dont le bord rencontre $\vec A\sph$ par (DSR). Comme $\vec A\sph$ est l'intérieur de $\vec A$, le demi-cône $\vec D(-b)$ rencontre aussi $\vec A\sph$, puis par stabilité par $-Id$, $\vec D(b)$ rencontre $\vec A\sph^+$ et $\vec A\sph^-$.
et s'il fixait également une facette $\vec f$ rencontrant $\vec A\becb$ mais pas $\vec D(b)$, alors en appliquant,  grâce au lemme, la proposition \ref{props:paraboliques}, 3, $u_b$ fixerait $\cl(\vec D(b)\cup \vec f)$, ce qui contient $\vec A\becb$. Mais c'est impossible car $u_b\not\in Z$.
% par clôture de $\fix_{\vec A}(u_b)$ (\ref{props:paraboliques}), $u_b$ fixerait en fait $\vec A\becb$, ce qui est impossible car $u_b\not\in Z$. 
Ainsi $\fix_{\vec A\becb}(u_b) = \vec D(b)$.

Alors $\fix_{\vec A\becb}(nu_bn\inv) = \vec\nu\becb(n).\vec D(b)$, et par ailleurs, comme $nu_bn\inv\in U_{r_a.b}$, $nu_bn\inv$ fixe le demi-cône $\vec D(r\bec_a.b)$. Donc $\vec D(r\bec_a.b)\subset \vec\nu\becb(n).\vec D(b)$. De la même manière, $\vec D(r\bec_a(-b))\subset \vec\nu(n).\vec D(-b)$. Et comme $\vec D(r\bec_a(b)) \cup \vec D(r\bec_a(-b))=\vec A\becb = \vec D(b)\cup \vec D(-b)$, on obtient l'égalité $\vec D(r\bec_a.b)= \vec\nu\becb(n).\vec D(b)$. En particulier, le bord de $\vec D(r\bec_a.b)$ rencontre $\vec A\sph$ et engendre un hyperplan de $\vec V\bec$, il s'agit donc de $\ker(r\bec_a.b)$. Ainsi le mur $\ker(b)$ est envoyé par $\vec\nu(n)$ sur $\ker(r\bec_a.b)$.

Nous avons en particulier prouvé que pour toutes racines simples $a,b\in\Pi\bec$, le mur $\ker(r\bec_a.b)$ rencontre $\vec A\sph$. Ceci permet de recommencer le raisonnement du paragraphe précédent, et de prouver que pour toute racine simple $c\in\Pi$, pour tout $u_c\in U_c\bec$, $\vec\nu ( n\bec(u_c) ). \ker(r\bec_a.b) = \ker(r\bec_c.r_a\bec.b)$. Puis finalement, nous obtenons que tout mur de $\vec V\bec$ rencontre $\vec A\sph$, et que $\vec\nu\becb$ permute l'ensemble de ces murs de la même manière que $\vec \nu\bec$.\\

 Soit $a\in\Pi\bec$, $u\in U_a\bec\setminus\{e\}$. Alors l'image $\vec\nu\becb(n).\vec C\bec$ de la chambre $\vec C\bec$ donnée par (DSR) est une autre chambre de $\vec A\bec$, et son adhérence contient $\barre{\vec C\bec}\cap \ker(a)$, il s'agit donc $r\bec_a.\vec C\bec$. Le même résultat est encore vrai pour $-\vec C\bec$, puis par l'argument standard des complexes de chambre minces,  on prouve que $\vec\nu\becb(n)$ agit comme $r\bec_a$ sur l'ensemble des chambres de $\vec A\bec$, et ceci entraine que $\vec\nu\becb(n)$ agit comme $r\bec_a$ sur l'ensemble des facettes de $\vec A\bec$. En particulier, $\vec\nu\becb(n)$ stabilise $\vec A\bec$.
 
 Le premier point est ainsi prouvé.\\
 
 Le second point est maintenant facile. On a $\vec A\bec = \vec\nu\becb(N\bec).\barre{\vec C\bec} \cup \vec\nu\becb(N\bec).\barre{-\vec C\bec}$, donc toute facette de $\vec A\bec$ s'écrit $\vec\nu\becb(n).\vec f\bec$ pour un certain $n\in N\bec $ et $\vec f\bec$ une facette de $\vec C\bec$. Par définition de $\vec\nu\becb$, $\vec\nu\becb(n).\vec f\bec=\vec\nu(n').\vec f\bec$ pour un certain $n'\in N$. Comme $\vec A$ et $\vec A\sph$ sont stables par $\vec\nu$, ceci intersecte $\vec A$, et même $\vec A\sph$ si $\vec f\bec$ est sphérique.
  
De plus, $\vec f\bec\cap \vec A$ est délimitée par des cônes de la forme $\ker(a)\cap \vec A$ qui sont traces de murs de $\vec A$. Ceci entraine l'existence d'une famille de facettes $(\vec f_i)$ telle que $\vec f\bec\cap \vec A = (\bigcup_i\vec f_i)\cap \vec V\bec$.
\cqfd

 On définit les murs de $\vec A\becb$ comme étant les $\ker(a)\cap \vec A\becb$ pour $a\in\phi\bec$. D'après la proposition, $\vec A\becb$ muni de $\vec\nu\becb$ est, tout comme $\vec A\bec$ muni de $\vec\nu\bec$, une réalisation géométrique du complexe de Coxeter associé à $\phi\bec$ muni de son action de $N\bec$.
Cependant, ni l'une ni l'autre ne sont en général l'appartement de référence pour l'immeuble $\iv\bec$ de la donnée radicielle $\D\bec$, au sens de \ref{sousoussection:immeuble vectoriel} car ils ne sont pas forcément essentiels.  Soit $\vec f_0\bec= \bigcap_a \ker(a)$ la plus petite facette de $\vec A\bec$. Il s'agit d'un sous-espace vectoriel de $\vec V\bec$, et l'appartement de référence pour $\iv\bec$ est $\vec A\bec_e:=\vec A\bec/\vec f_0\bec$. Les actions $\vec \nu\bec$ et $\vec\nu\becb$ ne coïncident en général pas sur $\vec V\bec$, puisqu'on peut modifier la famille $(a^\vee)_{a\in\phi\bec}$ par n'importe quelle famille d'éléments dans $\vec f_0$, tout en gardant une famille de coracines pour $\phi\bec$.

 Par contre, sur $\vec V_e\bec:= \vec V\bec/\vec f\bec_0$, on peut montrer que $\vec \nu\bec$ et $\vec\nu\becb$ coïncident. En effet, pour tout $t\in T\bec$, $\vec\nu\bec(t)=Id=\vec\nu\becb(t)$. Ensuite, soit $a\in\phi\bec$ et $n\in n\bec(U_a\setminus\{e\})$. Alors $\vec\nu\bec(n)\vec\nu\becb\inv(n)$ stabilise chaque facette, donc en particulier induit une homothétie de rapport positif sur chaque facette de dimension $1$. En conjuguant par d'autres éléments de $N\bec$, on voit que le rapport d'homothétie est le même dans chaque facette de dimension $1$ d'un même type. Finalement, $\vec\nu\bec(n)\vec\nu\becb(n)\inv$ est une homothétie sur chaque composante irréductible de $\vec V\bec_e$. Alors $\vec\nu\bec(n)\vec\nu\becb(n)\inv$ commute à $\vec\nu\bec(n)$, permettant de voir que $(\vec\nu\bec(n)\vec\nu\becb(n)\inv)^2 = \vec\nu\bec(n^2)\vec\nu\becb(n^2)\inv = Id$ car $n^2\in T\bec$. Donc le rapport d'homothétie est partout $1$.\\

 % Quitte à modifier la famille $(a^\vee)_{a\in\phi\bec}$ par des scalaires, on peut supposer que les actions $\vec\nu\bec$ et $\vec\nu\becb$ coïncident sur $\vec V\bec/\vec f_0\bec$.\\
  
  On pose $\iv\becb=G\bec.\vec A\becb\subset \iv$, on définit ses appartements, ses facettes, ses murs comme étant les images par les éléments de $G\bec$ de l'appartement $\vec A\becb$ et de ses murs et facettes.

  \fprop{\label{prop:descente vect} $\iv\becb$ est une réalisation géométrique de l'immeuble de $\D\bec$, autrement dit le complexe simplicial formé des facettes de $\iv\becb$, avec son action de $G\bec$ et la relation d'ordre "être dans l'adhérence de" est isomorphe à l'immeuble abstrait de $\D\bec$.}

  \demo 
  
  On veut définir une application entre les facettes de $\iv\bec$ et celles de $\iv\becb$. On pose:
  \[ \vec \j:\fonc	{ \F(\iv\bec)} 	{\F(\iv\becb)}
  					{g.\vec f}		{g.((\vec f+\vec f\bec_0)\cap \vec A)} \]
 , pour tout $g\in G\bec$ et $\vec f$ facette de $\vec A\bec$.\\
  
  Vérifions que $\vec\j$ est bien définie. Si $g.\vec f=h.\vec e$ dans $\iv\bec$, avec $g,h\in G\bec$ et $ \vec f,\vec e\in \F(\vec A\bec)$, alors il existe $n\in N\bec$ tel que $\vec e=n.\vec f$ (dans $\vec A\bec_e$) et $g\inv hn\in P\bec(\vec f)$.
  Alors $\vec\nu\bec(n).(\vec e+\vec f_0\bec) =\vec\nu\becb(n).(\vec e+\vec f_0\bec) = \vec f+\vec f\bec_0$.
   Soient $n'\in N$ et $z\in Z$ tels que $n=n'z$, alors par définition de $\vec\nu\becb$, $\vec\nu\becb(n)=\vec\nu(n')|_{\vec V\bec}$.  D'où $\vec e+\vec f\bec_0=\vec\nu(n).(\vec f+\vec f\bec_0)$ dans $\iv$.
  
 Comme $g\inv h n'=g\inv hnz\inv \in P\bec(\vec f).P(\vec f+\vec f\bec_0)$, et comme, d'après (DDR 1) et (DDR 3.1), $P\bec(\vec f)\subset P(\vec f+\vec f\bec_0)$, on prouve que $g.\vec f=h.\vec b$ dans $\iv$.\\

  Ainsi $\vec j$ est bien définie. Elle est clairement $G\bec$-équivariante, son image est $\iv\becb$, et sa restriction à $\F(\vec A\bec_e)$ est un isomorphisme de complexe de Coxeter sur $\F(\vec A\becb)$. Maintenant, le fait que tout couple de facettes de $\iv\bec$ est inclus dans un appartement, qui est image de $\vec A\bec$ par un élément de $G\bec$, prouve l'injectivité.\cqfd
  
  \fdef{Pour toute facette $\vec f\becb$ de $\iv\becb$, on note $P\bec(\vec f\becb) = P\bec(\vec f\becb/\vec f_0\bec)$, $M\bec_{\vec A\becb}(\vec f\becb) = M\bec_{\vec A\bec}(\vec f\becb/\vec f_0\bec)$ et $U\bec(\vec f\becb) = U\bec(\vec f\becb/\vec f_0\bec)$.
  }
  
  Notons que $\vec f\becb/\vec f_0\bec$ est une partie de facette de $\iv\bec$, mais son fixateur et ses facteurs unipotent et de Lévi sont égaux aux fixateur, facteur unipotent et facteur de Lévi de cette facette.
  
 \fprop{ \label{prop:descente vect suite}Soit $\vec f\becb$ une facette de $\iv\becb$, et $\vec f$ une facette de $\iv$ rencontrant $\vec f\becb$. Alors:
 \liste{
 \item $P\bec(\vec f\becb)=P(\vec f\becb)\cap G\bec = P(\vec f)\cap G\bec$,
 \item $U\bec(\vec f\becb) = U(\vec f)\cap G\bec = U(\vec f\becb)\cap G\bec$,
 \item $M\bec_{\vec A\becb}(\vec f\becb)= M_{\vec A}(\vec f\becb)\cap G\bec =M_{\vec A}(\vec f)\cap G\bec$, si $\vec f\becb$ est dans l'appartement $\vec A\becb$,
 \item $T\bec = Z\cap G\bec$,
 \item pour tout $a\in\phi\bec$, $U_a\bec=U_a\cap G\bec$.
 }
  }  
  
  \demo
  
Par la proposition précédente,  $P\bec(\vec f\becb)=\fix_{G\bec}(\vec f\becb/\vec f_0\bec)=P(\vec f\becb\cap \vec A)\cap G\bec \subset P(\vec f)\cap G\bec$. De plus, si $g\in G\bec\cap P(\vec f)$, alors $g$ fixe une partie de la facette $\vec f\becb$, et donc $g$ fixe cette facette. D'où l'inclusion $G\bec\cap P(\vec f)\subset P\bec(\vec f\becb)$.\\

Maintenant, on obtient directement, si $\vec f\becb$ est dans $\vec A\becb$, $M\bec_{\vec A\becb}(\vec f\becb) = P\bec(\vec f\becb)\cap P\bec(\op_{\vec A\becb} (\vec f\becb)) = G\bec\cap P(\vec f\becb)\cap P(\op_{\vec A}(\vec f\becb)) = G\bec\cap M_{\vec A}(\vec f\becb) = G\bec\cap P(\vec f)\cap P(\op_{\vec A}(\vec f)) = G\bec\cap M_{\vec A}(\vec f)$.\\

% Le groupe $M\bec_{\vec A\becb}(\vec f\becb)$ est engendré par $T\bec$ et par les $U\bec_a$, $a\in \phi\bec[m](\vec f\becb)$. Or $T\bec \subset Z = \fix(\vec A\becb)\subset \fix(\vec f\becb\cup -\vec f\becb)=M_{\vec A}(\vec f\becb)$, et pour tout $a\in \phi\bec[m](\vec f\becb)$, $U_a\bec\subset U_a\subset \fix(\vec D(a)) \subset \fix (\vec f\becb\cup -\vec f\becb)$. D'où $M\bec_{\vec A\becb}(\vec f\becb)\subset M_{\vec A}(\vec f\becb)$.

Ensuite, $U\bec(\vec f\becb)$ est le sous-groupe distingué de $P\bec(\vec f\becb)$ engendré par les $U\bec_a$, $a\in\phi\bec[u](\vec f\becb)$. Chacun de ces $U\bec_a$ est bien inclus dans $U(\vec f\becb)$ et dans $U(\vec f)$, et $P\bec(\vec f\becb)\subset P(\vec f\becb)\cap P(\vec f)$, donc $U(\vec f\becb)$ et $U(\vec f)$ sont normalisés par $P\bec(\vec f\becb)$. D'où l'inclusion $U\bec(\vec f\becb)\subset U(\vec f\becb)\cap U(\vec f)\cap G\bec$.  

Une fois choisi un appartement contenant $\vec f\becb$, les décompositions de Lévi $P\bec(\vec f\becb) = M\bec(\vec f\becb)\ltimes U\bec(\vec f\becb)$ et $P(\vec f) = M(\vec f)\ltimes U(\vec f)$ permettent de prouver $U(\vec f)\cap G\bec \subset U\bec(\vec f\becb)$.

 On a alors $U\bec(\vec f\becb) = G\bec\cap U(\vec f) \subset G\bec\cap U(\vec f\becb)$. L'inclusion manquante est évidente car $U(\vec f\becb)\subset U(\vec f)$, car $\vec f\becb$ contient une partie de la facette $\vec f$.\\
 
 Par le premier point, pour toute partie $\vec \Omega$ de $\vec A\becb$, $P\bec(\vec\Omega) = P(\vec\Omega)\cap G\bec$. En particulier, pour tout $a\in\phi\bec$, $P\bec(\vec D(a)) = P(\vec D(a)) \cap G\bec = (Z\ltimes U_a)\cap G\bec$, comme on l'a vu dans le lemme \ref{lemme:descente vect}. Mais $P\bec(\vec D(a))$ admet par ailleurs la décomposition de Lévi $P\bec(\vec D(a)) = T\bec\ltimes U\bec_a$. Les inclusions $T\bec\subset Z$ et $U\bec_a\subset U_a$ permettent alors de prouver que $T\bec= Z\cap G\bec$ et $U_a\bec= U_a\cap G\bec$.\cqfd

%Les égalités $U\bec(\vec f\becb)\rtimes M\bec(\vec f\becb)=P(\vec f\becb)=P(\vec f)\cap G\bec = (U(\vec f)\rtimes M(\vec f))\cap G\bec$ permettent alors, avec les deux inclusions qu'on vient de prouver, d'obtenir les égalités $U(\vec f)\cap G\bec = U\bec(\vec f\becb)$ et $M_{\vec A}(\vec f) \cap G\bec = M\bec_{\vec A\becb}(\vec f\becb)$. Notons enfin que $M_{\vec A}(\vec f) = M_{\vec A}(\vec f\becb)$ car $\vect_{\vec A}(\vec f)= \vect_{\vec A}(\vec f\becb)$ car $\vec f$.\cqfd

  \subsection{Descente de la valuation}

 On suppose désormais que $Q$ vérifie \plienpm(sph).
 
 \flemme{\label{lemme:Nbec et N} $N\bec \subset N.Q(\barre{A\die})$.
 }
\pv Le groupe $N\bec$ stabilise $\barre{A\die}$ par (DM 4). Et le groupe $Q(\barre{A\die})$ est transitif sur les appartements contenant $\barre{A\die}$, par le corollaire \ref{cor:isom entre apparts}, car $\barre{A\die}$ contient des points sphériques positifs et négatifs, et $Q$ vérifie \plienpm(sph).\cqfd

On étudie maintenant l'action de $G\bec$ sur $A$, et on prouve que $\phii\bec$ est une valuation pour $\D\bec$.

Comme $o\in A\die$, on peut identifier chaque $a\in \phi\bec$ à une forme affine sur $A\die$ s'annulant en $o$, et on note pour tout $a\in\phi\bec$ et $\lambda\in\R$, 
$D(a,\lambda)=\ens{x\in A\die}{a(x)+\lambda\geq 0}$ et $M(a,\lambda)=\ens{x\in A\die}{a(x)+\lambda= 0}$.

\fprop{
 Pour tout $a\in\phi\bec$, pour tout $u\in U_a\bec$, $\fix_{A\die}(u) = D(a,\phii\bec_a(u))$, et l'action de $n(u)$ sur $\I$ induit sur $A\die$ une réflexion selon l'hyperplan $M(a,\phii_a\bec(u))$.
 
 La famille $\phii\bec$ forme une valuation pour $\D\bec$.
}

\demo

Soit $a\in\phi\bec$, $k\in\R$. Pour tout $\alpha\in\phi_a$, soit $r\in\R^{+*}$ tel que $\alpha|_{\vec V\bec}=r.a$, alors $U_{\alpha, rk}=U_\alpha(D(a,k))$. 
Or $\prod_{\alpha\in\phi_{a,red}} U_\alpha(D(a,k))$ est égal au groupe $Q(D(a,k))\cap U_a$, grâce à \ref{prop:levi equilibre} car $\phi_{a,red}$ est une partie nilpotente de racines. Le fait qu'il s'agisse d'un groupe prouve que c'est $U_{a,k}$. Donc $U_{a,k} = \prod_{\alpha\in\phi_{a,red}} U_\alpha(D(a,k)) = U_a\cap Q(D(a,k))$.

Soit $u\in U_a\bec$. Sachant que $\fix_A(u)$ est une intersection de demi-appartements dirigés par des $\vec D(\alpha)$, $\alpha\in\phi_a$ (\ref{prop:points fixes de u1u2...uk}), on voit que $\fix_{A\die}(u)$ est un demi-espace dirigé par $\vec D(a)$. La description qu'on vient d'obtenir des $U_{a,k}$ entraine alors que $\fix_{A\die}(u) = D(a,\phii_a\bec(u))$.\\

 Soit $\vec m\becb$ une cloison de $\vec A\becb$ contenue dans $\ker(a)$. Comme $\vec m\becb$ intersecte $\vec A\sph$, il existe une facette $\vec m$ de $\vec A$, sphérique, contenant un ouvert de $\vec m\becb$. On note $\D_{\vec m} = (M_{\vec A}(\vec m), (U_\alpha,\phii_\alpha)_{\alpha\in\phi^m(\vec m)})$. Il s'agit d'une donnée radicielle valuée de type fini. On note également $\D_{\vec m\becb}\bec = (M\bec_{\vec A\becb}(\vec m\becb) , (U_a\bec)_{a\in \phi\bec[m](\vec m\becb)})$, il s'agit d'une donnée radicielle de type fini. Avec la partie d'immeuble $\I\die[\vec m ]=\I\die\cap \I_{\vec m}$, ces deux données radicielles vérifient les hypothèses du théorème \cite{bruhat-tits} 9.2.10.  En conclusion, la valuation $(\phii_\alpha)_{\alpha\in \phi^m(\vec m)}$ se descend à une valuation de $\D\bec_{\vec m\becb}$, qui n'est autre (si on compare la définition de \cite{bruhat-tits} 9.1.6 à celle du présent texte) que $(\phii_a\bec)_{a\in \phi\bec[m](\vec m\becb)}$. Ainsi, $(\phii_a\bec)_{a\in \phi\bec[m](\vec m\becb)}$ est une valuation. Comme $\phi\bec[m](\vec m\becb)$ contient au moins $a$ et $-a$, on obtient directement (V2.2) et (V5) pour la racine $a$.
 
 Soient $u',u''\in U_{-a}$ tels que $n(u) = u'uu''$. Alors par (V5), $\phii_{-a}(u')=\phii_{-a}(u'') = -\phii_a(u)$. Donc $u'$ et $u''$ fixent $D(-a,-\phii_a(u))$, et $n(u)$ fixe $M(u,\phii_a\bec(u))$.
 
 Montrons que $n(u)$ agit sur $A\die$ comme une réflexion. Soient $n'\in N$ et $q\in Q(\barre{A\die})\subset Z$ tels que $n(u)=n'.q$. Alors $\vec\nu\becb(n(u)) = \vec\nu(n')|_{\vec V\bec}$ et c'est une réflexion dans $\vec V\bec$ (\ref{props:action de Nbec}, 1). Donc $n'$ induit une réflexion sur $A\die$. Mais $n(u)$ agit sur $A\die$ comme $n'$ car $n(u)\inv n'\in Q(A\die)$. Donc $n(u)$ induit bien sur $A\die$ une réflexion selon l'hyperplan $M(a,\phii_a\bec(u))$.\\

 Il est maintenant facile de vérifier que $\phii\bec$ est une valuation. On a déjà supposé (V 0): il s'agit de (DV2). Le (V 1) est clair, (V 2) découle du lemme \ref{lemme:pour V2}, valide car $\vec\nu\bec$ et $\vec\nu\becb$ coïncident. La condition (V 3) est facilement vérifiée grâce à (DR 2) ($\D\bec$ est une donnée radicielle) et grâce à la caractérisation de $\phii_a\bec(u)$ par les points fixes de $u$. Enfin (V 4) est évident sur la définition de $\phii\bec$.\cqfd
 
 \subsection{Descente de la famille de parahoriques}
 
 La dernière étape avant d'obtenir une masure bordée pour le groupe $G\bec$ est de définir une famille de parahoriques pour $(\D\bec,\phii\bec)$. Pour commencer, on décrit une réalisation $A\becb$ de l'appartement affine $A\bec$ en utilisant la masure bordée $\I$.
 
 Comme $\phi\bec$ est à base libre dans $(\vec V\bec)^*$, on peut, quitte à remplacer $\phii$ par une valuation équipollente, supposer que $\phii\bec$ est spéciale, comme dans la partie \ref{section:construction generale}.
 
 On construit alors l'appartement $A\becb$ comme dans \ref{soussection:appartement general}: il s'agit de l'ensemble $A\becb$ des cônes dans $A\die$ dirigés par une facette de $\vec A\becb$, quotienté par la relation $f\sim g\ssi f\cap g$ contient un scp de $f$ et de $g$. Pour une facette $\vec f\becb$ de $\vec A\becb$, l'ensemble des classes de cônes dirigés par $\vec f\becb$ est appelé la façade de direction $\vec f\becb$ et noté $A\becb[\vec f\becb]$. Remarquons que pour toute facette $\vec f$ de $\vec A$ contenant un ouvert de $\vec f\becb$, la façade $A\becb[\vec f\becb]$ est isomorphe à $pr_{\vec f}(A\die)$ ou encore à $\barre{A\die}\cap A_{\vec f}$, on pourra noter $pr_{\vec f}$ cet isomorphisme, et $pr_{\vec f\becb}$ sa réciproque.
  La façade principale est $A\becb[\vec f_0\bec]$, où $\vec f_0\bec$ est la plus petite facette de $\vec A\becb$.
  
  Les murs, demi-appartements, facettes sont définis à partir de $\phi$ et de $\phii$ comme dans \ref{soussection:appartement general}.
  L'action de $N\bec$ stabilise $A\die$ et l'ensemble des facettes de $\vec A\becb$, elle induit donc une action sur $A\becb$ par automorphismes.
  \\
  
  Soit $\vec f\becb$ une facette de $\vec A\becb$, soient $(\vec f_i)_{i\in I}$ les facettes de $\vec A$ qui contiennent un ouvert de $\vec f\becb$, autrement dit les facettes maximales parmi celles qui recouvrent $\vec f\becb$. On pose alors, pour tout $a\in A_{\vec f\becb}$:
  \[ Q\bec(a)=\bigcap_{i\in I} Q(pr_{\vec f_i}(a))\cap G\bec \point \]
  
\rema Si $\vec f\becb$ est sphérique, alors les $\vec f_i$ le sont aussi et en utilisant \ref{lemme:csq de plienpm} puis \plien($\vec f_i$) pour chaque $i$, on voit que tous les $Q(pr_{\vec f_i}(a))\cap P\bec(\vec f\becb)$ sont égaux, de sorte que $Q\bec(a)$ est égal à n'importe lequel de ces groupes.\\

\flemme{\label{lemme:action de Nbec} Soit $a\in A\becb$, soit $\vec f\becb$ la direction de la façade de $a$, soient $(\vec f_i)_{i\in I}$ les facettes de $\vec A$ contenant un ouvert de $\vec f\becb$, notons pour tout $i\in I$ $a_i=pr_{\vec f_i}(a)$.

Alors %pour tout $i\in I$, $G\bec\cap N(\vec f_i).Q(a_i) \subset N\bec(\vec f\becb).Q(a_i)$, et 

$$G\bec\cap \stab_N(\vec f\becb).Q(\{a_i\}_{i\in I} ) = N\bec(\vec f\becb).Q\bec(a)\point$$}  
\pv

L'inclusion $\supset$ vient de $N\bec\subset N.Q(\barre{A\die})$.

Réciproquement, soit $g=nq\in G\bec\cap stab_N(\vec f\becb).Q(\{a_i\}_{i\in I})$. Soit $g=u\bec n\bec q\bec$ une écriture de $g$ dans la décomposition d'Iwasawa $G\bec=U\bec(\vec C\becb). N\bec . G\bec(a)$, avec $\vec C\becb$ une chambre dont l'adhérence contient $\vec f\becb$.

Soit $i\in I$.
 Soit $\vec C_i$ une chambre de $\vec A$ dont l'adhérence contient un ouvert de $\vec C\bec$ et $n\inv \vec f_i$. Soient $n'\in N$, $z\in Q(\barre{A\die})$ tels que $n\bec=n'z$. Alors pour tout $i\in I$, $g=e.n.q = u\bec.n'.(z q\bec)$ sont deux écritures de $g$ dans la décomposition d'Iwasawa $G=U(\vec C_i).N.Q(a_i)$. Par la proposition \ref{prop:unicite dans iwasawa 2}, applicable car $\vec f_i\subset n\barre{\vec C_i}$, $n\inv n'\in N(a_i)$. Au total, $n\inv n'=n\inv n\bec z\inv\in N(\{a_i\}_{i\in I})$ d'où $n\in n\bec.Q(\{a_i\}_{i\in I})$. Alors $g=nq\in n\bec . Q(\{a_i\})\cap G\bec \subset N\bec.(Q(\{a_i\}_{i\in I})\cap G\bec) = N\bec Q\bec(a)$.\cqfd

  \fprop{
  La famille $Q\bec$ est une bonne famille de parahoriques pour $(\D\bec,\phii\bec)$. %Elle vérifie \plienpm(sph).
   Si $Q$ vérifie \plienp(sph), alors $Q\bec$ aussi.
  }
  
  \rema La condition \plienpm(sph) ne semble pas se descendre à $Q\bec$, en tout cas pas de manière évidente.\\
  
  \demo 
  
  On fixe une facette $\vec f\becb\in\F(\vec A\becb)$ et un point $x\in A\becb[\vec f\becb]$, on note $(\vec f_i)_{i\in I}$ l'ensemble des facettes de $\vec A$ contenant un ouvert de $\vec f\becb$, et pour tout $i\in I$, on pose $x_i=pr_{\vec f_i}(x)$.\\
  
   Pour tout $i\in I$, on a $U\bec(\vec f\becb)\subset U(\vec f\becb)\subset U(\vec f_i)$, donc $U\becb(\vec f\becb)\subset Q\bec(x)$. D'autre part, pour tout $i\in I$, $Q\bec(x)\subset P(\vec f_i)\cap G\bec=P\bec(\vec f_i)=P\bec(\vec f\becb)$. D'où $Q\bec(x)\subset P\bec(\vec f\becb)$, donc $Q\bec$ vérifie (para 0.1).
   
   Soit $n\in N\bec(x)$. Modulo un élément $z\in Q(\barre{A\die})$, qui fixe donc tous les $x_i$, on peut supposer $n\in N(x)$. Alors il existe un représentant $f\becb=y+\vec f\becb$ de $x$ inclus dans $A\die$, tel que $f\becb\cap n.f\becb$ contient un scp de $f\becb$. Autrement dit, $n\in N(\vec f\becb)$ et $\fl{yn(y)}\in\vect(\vec f\becb)$. Ceci implique directement $n\in N(\vec f_i)$ et $\fl{yn(y)}\in\vect(\vec f_i)$ (car $\vect(\vec f\becb)=\vect(\vec f_i)$). D'où $n\in N(pr_{\vec f_i}(y)) = N(x_i)$. Ceci prouve (para 0.2).
   
   Les deux conditions (para 0.3) et (para 0.4) sont claires, donc $Q\bec$ est une famille de parahoriques.\\
   
Soit $\vec g\becb$ une facette sphérique de $\vec f\becb^*\cap \vec A\becb$. Soit $g\in Q\bec(x)\cap P\bec(\vec g\becb)$. Soit $\vec g$ une facette de $\vec A$ contenant un ouvert de $\vec g\becb$, il existe $i\in I$ tel que $\vec f_i\subset \barre{\vec g}$, alors $g\in Q(x_i)\cap P(\vec g) \subset Q( pr_{\vec g}(x_i) ) = Q(pr_{\vec g}(pr_{\vec g\becb}(x)))$. D'où \plien(sph).\\

Montrons \psph: supposons $\vec f\becb$ sphérique et montrons que $Q\bec(x) = P\bec(x) = U\bec(\vec f\bec). N\bec(x). G\bec(\phi\bec[m](\vec f\bec),x)$. Soit $i\in I$, alors $\vec f_i$ est sphérique (proposition \ref{props:action de Nbec}) et:
% ,  
\calc{
Q\bec(x) &= & Q(x_i)\cap P\bec(\vec f\becb)\\
& =& P(x_i)\cap P\bec(\vec f\becb)\\
&=& U(\vec f_i)\rtimes (N(x_i).G(\phi^m(\vec f_i),x))\cap P\bec(\vec f\becb)}

Avec la décomposition de Lévi $P\bec(\vec f\becb) = U\bec(\vec f\becb)\rtimes M\bec(\vec f\becb) = (U(\vec f_i)\cap G\bec)\rtimes (M(\vec f_i)\cap G\bec)$ (par \ref{prop:descente vect suite}), on obtient:
\calc{
Q\bec(x) &=& \left( U(\vec f_i) \cap G\bec \right) \rtimes\left( N(x_i).G(\phi^m(\vec f_i),x) \cap G\bec \right)\\
&= & U\bec(\vec f\becb) \rtimes  \left( N(x_i).G(\phi^m(\vec f_i),x) \cap G\bec \right)\point
}
 C'est la proposition 9.1.17 de  \cite{bruhat-tits}, appliquée aux données radicielles $\D_{\vec f_i}$ et $\D\bec_{\vec f\becb}$, et avec $S\bec=T\bec$, qui indique alors que $N(x_i).G(\phi^m(\vec f_i),x) \cap G\bec =  N\bec(x_i).G\bec(\phi\bec[m](\vec f_i),x)$. Comme $N\bec(x_i)\subset N\bec(x)$ et $\phi\bec[m](\vec f_i)=\phi\bec[m](\vec f\becb)$, on a bien obtenu $Q\bec(x)\subset P\bec(x)$.\\

Montrons \pinj: soit $n\in N\bec\cap Q\bec(a) = N\bec\cap P(\vec f\becb) \cap Q(\{a_i\}_{i\in I})$. Soit $a_0\in A\die$ tel que $a=[a_0+\vec f\becb]$, alors $n.a= [n.a_0+\vec f\becb]$. De plus, pour un $i\in I$ quelconque, le fait que $n.a_i=n.[a_0+\vec f_i] = a_i$ entraine que $\fl{a_0n(a_0)}\in \vect(\vec f_i) = \vect(\vec f\becb)$. Donc $n.[a_0+\vec f\becb] = [a_0+\vec f\becb]$, autrement dit $n\in N\bec(a)$.\\

Il ne manque à $Q\bec$ plus que \plienn(sph) pour être une bonne famille de parahoriques. Soit donc $\vec g\becb$ une facette de $\vec f\becb^*\cap \vec A\sph$ et $g\in N\bec.Q\bec(x)\cap P(\vec g\becb)= N\bec(\vec f\becb).Q\bec(x)\cap P(\vec g\becb)$. Soient $(\vec g_i)_{i\in I}$ des facettes de $\vec A$ contenant un ouvert de $\vec g\becb$, telles que pour tout $i\in I$ $\vec f_i\subset \barre{\vec g_i}$.

Pour tout $i\in I$, on a $g\in N(\vec f\becb)Q(x_i)\cap P(\vec g_i)\subset N(\vec f_i).Q(x_i)\cap P(\vec g_i) = N(\vec f_i).Q(\{x_i,pr_{\vec g_i}(x_i)\})$. 
Donc par \ref{prop:isom entre appart}, $g\in P(\vec g\becb)\cap \bigcap_{i\in I}N(\vec f_i).Q(\{x_i,pr_{\vec g_i}(x_i)\}) \subset N.Q(\{x_i,pr_{\vec g_i}(x_i)\}_{i\in I})\cap P(\vec g\becb) = N(\vec g\becb).Q(\{x_i,pr_{\vec g_i}(x_i)\}_{i\in I})$. Si $\vec g$ est une facette de $\vec A$ contenant un ouvert de $\vec g\becb$ qui ne figure pas parmi les $\vec g_i$, alors $Q(\{x_i,pr_{\vec g_i}(x_i)\}_{i\in I})\subset Q(pr_{\vec g}(x)$ (voir la remarque précédant le lemme). Alors
comme de plus $g\in G\bec$, on obtient par le lemme $g\in N\bec(\vec g\becb) .Q\bec(pr_{\vec g\bec}(x))$.\\

% = N\bec(\vec f\becb) Q\bec(pr_{\vec g\becb}(x))$. Comme $\vec g$ est sphérique et $Q\bec$ vérifie \psph, on a $Q\bec(pr_{\vec g\becb}(x)) = N\bec(pr_{\vec g\becb}(x))\point U\bec(\vec g\becb)\point G\bec(\phi\bec[m](\vec g\becb),x)$. On peut donc supposer:
%\[g\in N\bec(\vec f\becb).Q\bec(x)\cap U(\vec g\becb) \point \]

%En particulier, $g\in N(\vec f\becb).Q(\{x_i\}_{i\in I})\:\cap\: U(\vec g\becb)$, ce qui vaut $Q(\{x_i\}_{i\in I})\cap U(\vec g\becb)$ par \ref{lemme:action de u sur A}. Ainsi, par les \plien($\vec g_i$), on obtient $g\in G\bec\cap Q(\bigcup_i \{x_i,pr_{\vec g_i}(x_i)\}) = Q\bec(\{x,pr_{\vec g\becb}(x)\})$.\cqfd

Supposons que $Q$ satisfasse à \plienp(sph). Soit $\vec g\becb\in \vec f\becb^*\cap \vec A\becb$, soit $g\in Q\bec(a)\cap P\bec(\vec g\becb)$. Alors pour tout $i$, $g\in Q(a_i)\cap P(\vec g_i) = Q(\barre{a_i+\vec g_i})$. Ceci entraine bien que $g\in Q\bec(\barre{a+\vec g\becb})$.\cqfd

%Montrons enfin \plienpm(sph). Soit $\vec g\becb$ une facette de $\vec f\becb^*\cap \vec A\sph$ et $g\in Q\bec(a)\cap U\bec(\vec g\becb)$.

\subsection{Injection des façades}

Soit $\I\bec=\I(Q\bec)$ la masure bordée pour $G\bec$ qu'on vient d'obtenir. Le but de ce paragraphe est d'identifier certaines façades de $\I\bec$ à des parties de $\I$.
 Soit $\vec f\becb$ une facette de $\vec A\becb$, soit $\vec f$ une facette de $\vec A$ contenant un ouvert de $\vec f\becb$. Comme on l'a déjà dit, l'injection $A\becb[\vec f\becb ] \hookrightarrow A_{\vec f}$ est bien définie par $[a\die+\vec f\becb]\mapsto [a\die+\vec f]$, pour tout $a\die\in A\die$. Nous voulons à présent étudier l'injection de $\I\bec_{\vec f\becb}$ dans $\I_{\vec f}$.

 \fprop{ \label{prop:injection des facades} La fonction $j_{\vec f\becb}: \fonc {\I\bec_{\vec f\becb}} {\I_{\vec f}} {[g,[a\die+\vec f\becb]]} {[g,[a\die+\vec f]]}$ est bien définie et $G\bec$-équivariante.
 
 Si $\vec f$ est sphérique, elle est de plus injective.\\
 }

\demo Soient $g,h\in G\bec$ et $a,b\in A\becb[\vec f\becb ]$ tels que $(g,a) \sim_{Q\bec} (h,b)$. Soient $a\die,b\die\in A\die$ tels que $a=[a\die+\vec f\becb]$ et $b=[b\die+\vec f\becb]$, notons encore $a'=[a\die+\vec f]$ et $b'=[b\die+\vec f]$.

Soit $n\in N\bec$ tel que $n.a=b$ et $g\inv hn\in Q\bec(a)$. Donc $n.(a\die+\vec f\becb)\cap (b\die+\vec f\becb)\not=\vide$. Comme $\vec f$ contient un ouvert de $\vec f\becb$, ceci entraine $n.(a\die+\vec f)\cap (b\die+\vec f)\not=\vide$.
Il existe $n'\in N$ et $z\in Q(\barre{A\die})$ tels que $n=n'z$, alors $n'.(a\die+\vec f) = n.(a\die+\vec f)$, et nous venons de voir que ce cône est équivalent à $b\die+\vec f$. Ainsi, $n'.a'=b'$.

De plus $g\inv hn'=g\inv h nz\inv \in Q\bec(a).Q(a')=Q(a')$. Donc $(g,a) \sim_{Q} (h,b)$, et la fonction $j_{\vec f\becb}$ est bien définie. Elle est clairement équivariante.\\

Gardant les notations précédentes, supposons maintenant $\vec f\becb$ sphérique et $(g,a')\sim_{Q} (h,b')$. Donc $g\inv h.b'=a'$ et $g\inv h \in Q(a').N\cap P\bec(\vec f\becb)$. Comme $Q(a')\subset P(\vec f)$, on a en fait $g\inv h\in Q(a').N(\vec f)$. Mais $N(\vec f)=N(\vec f\becb)$ d'où finalement $g\inv h\in \left( (Q(a')\cap P(\vec f\becb)).N(\vec f\becb)\right) \cap P\bec(\vec f\becb)$.
Soient $(\vec f_i)_{i\in I}$ les facettes de $\vec A$ contenant un ouvert de $\vec f\becb$, et $(a_i)_{i\in I}$ les projetés de $a$ dans les façades $A_{\vec f_i}$. Comme les $\vec f_i$ sont sphériques et comme $Q$ vérifie \plienpm(sph), $Q(a')\cap P(\vec f\becb) = Q(a')\cap P(\bigcup_i \vec f_i) = Q(\{a_i\}_{i\in I})$. Alors le lemme \ref{lemme:action de Nbec} prouve que $g\inv h\in   Q\bec(a).N\bec(\vec f\becb)$. Soit $n\in N\bec(\vec f\becb)$ tel que $g\inv h\in Q\bec(a)n\inv$, il reste à prouver que $n.a=b$. Sachant que $Q\bec(a)\subset Q(a')$ et $g\inv h.b'=a'$, on a déjà $n.a'=b'$. Donc $n.(a\die+\vec f)\cap (b\die+\vec f)\not=\vide$. Donc $\fl{n(a\die) b\die}\in \vect( \vec f)\cap \vec V\bec = \vect(\vec f\becb)$. Donc $n.(a\die+\vec f\becb) \cap (b\die +\vec f\becb)\not=\vide$.\cqfd

\section{Le cas Kac-Moody}
\label{section:KM}

On adopte la définition de J. Tits pour les groupes de Kac-Moody, et on se réfère principalement à \cite{remy}.

\subsection{Rappels et notations}
\subsubsection{Groupes de Kac-Moody déployés}

Soit $G$ un groupe de Kac-Moody déployé sur un corps $\K$, il s'agit donc d'un foncteur des $\K$-algèbres vers les groupes. Comme tout groupe de Kac-Moody, $G$ vient avec une algèbre de Kac-Moody $\g$ et une action $Ad:G\fleche Aut(\g)$ appelée l'action adjointe.

Pour chaque tore maximal $T$, on note $T^*$ son groupe de caractères et $T_*$ son groupe de cocaractères. On définit une forme bilinéaire $\langle . , . \rangle : T^*\times T_*\fleche \Z$  par $\langle \chi,h\rangle =n$ si $\chi\circ h (k)=k^n$ pour tout $k\in \K$.\\

 A chaque tore maximal $T$ correspond un système de racines $\phi^c(T)$ tel que l'algèbre de Kac-Moody $\g(\K)$ est graduée par $\phi^c(T)\cup\{0\}$.
On note toujours $\vec V(T)$ l'espace vectoriel réel tel que $\phi^c(T)\subset \vec V(T)^*$, toute base de $\phi^c(T)$ est une base de $\vec V(T)$. On note aussi $Q(T) = \Z.\phi^c(T)$ le réseau des racines, il existe un morphisme de $Q(T)$ dans $T^*$, noté $\alpha\mapsto \bar\alpha$, qui s'étend à une application linéaire de $\vec V(T)^*$ dans $T^*\tens\R$.

   Une racine est soit réelle, soit imaginaire, on note $\phi(T)$ l'ensemble des racines réelles, et $\phi^{\text{im}}(T)$ celui des racines imaginaires. L'ensemble $\phi^c(T)$ est appelé le système complet de racines. Dans la suite, ce sera le plus souvent $\phi(T)$ qui interviendra, ce qui explique qu'on ait choisi la notation la plus courte pour le désigner. \\
 
  Il existe une base de Chevalley $(e_a)_{a\in\phi\sqcup Im}$ de $\g$. Pour chaque $\alpha\in\phi$, $\g_\alpha$ est de dimension 1 et la base contient un élément noté $e_\alpha$ de $\g_\alpha\setminus\{0\}$. Par contre pour $\alpha$ une racine imaginaire, ou $\alpha=0$, $\g_\alpha$ peut être de dimension supérieure, et la base contiendra plusieurs élément de $\g_\alpha$. A chaque racine réelle $\alpha \in \phi(T)$ correspond un sous-groupe $U_\alpha$ de $G$, isomorphe au groupe additif. Il existe un choix des isomorphismes $(u_\alpha)_{\alpha\in\phi}$ entre les $U_\alpha(\K)$ et $(\K,+)$ tels que pour $\alpha\in\phi$, $Ad(u_\alpha(k)) = exp(ad(k e_\alpha)) = \sum_n k\tens (\frac{ad(e_\alpha)^n}{n!})\in\aut(\g)$ (l'algèbre de Kac-Moody $\g(\K)$ vaut $\K\tens \g_\Z$, où $\g_\Z$ est une $\Z$-algèbre de Lie stable par les $\frac{ad(e_\alpha)^n}{n!}$, $n\in\N$). Le tore $T$ agit diagonalement, par $Ad(t).e_\alpha= \bar\alpha(t).e_\alpha$.\\
  
 Le groupe $U_\alpha$ est normalisé par $T$, plus précisément, la formule suivante est vérifiée pour  $k\in\K$ et $t\in T(\K)$:
\[ t u_\alpha(k) t\inv = u_\alpha(\bar\alpha(t).k)  \]  

On note $N(T)$ le normalisateur de $T$, il est engendré par $T$ et les éléments $n_\alpha(k):= u_{-\alpha}(k\inv) u_\alpha(k) u_{-\alpha}(k\inv)$ pour $\alpha\in \phi$ et $k\in\K$. On note également $W(T)=N(T)/ T$ le groupe de Weyl vectoriel relatif à $T$. La paire $(W(T), (n_\alpha(1).T)_{\alpha\in\pi})$ forme un système de Coxeter pour toute base $\Pi$ de $\phi$.\\

Tout ceci entraine que la famille $(G(\K),(U_\alpha(\K))_{\alpha\in \phi(T)} )$ est une donnée radicielle génératrice.\\

Supposons maintenant $\K$ muni d'une valuation non triviale $\omega:\K\fleche \R\cup\{\infini\}$.  Comme on l'a déjà dit (proposition \ref{prop:existence de valuation}), la famille de fonction $\phii=(\phii_\alpha)_{\alpha\in\phi(T)}$ définie par :
\[\phii_\alpha:\fonc{U_\alpha(\K)}{\R\cup\{\infini\}} {u_\alpha(k)} { \omega(k)}\]
est une valuation de cette donnée radicielle.\\

Dans le cas d'un groupe de Kac-Moody, l'action de $T$ sur $A(T)$ peut être décrite un peu plus directement que dans le cas général d'une donnée radicielle (\ref{prop:action de T}).
\fprop{\label{prop:action de T KM} Soit $T$ un tore maximal de $G$. Alors pour tout $t\in T$, le vecteur $\vec v_t$ est l'unique vecteur de $\vec V(T)$ tel que:
\[ \forall \alpha\in\phi,\: \alpha(\vec v_t) = -\omega(\bar\alpha(t)) \point\]
}
%\rema Dans le terme de droite, $\alpha$ est vu comme un caractère de $T$.\\

\demo

Par la proposition \ref{prop:action de T}, le vecteur $\vec v_t$ est caractérisé par $\alpha(\vec v_t)=\phii_\alpha(u) - \phii_\alpha(tut\inv)$, pour tout $\alpha\in\phi$ et tout $u\in U_\alpha(\K)\setminus\{e\}$.

Soit $\alpha\in\phi$, et prenons $u=u_\alpha(1)$ ($u\not=e$ puisque $u_\alpha\inv(e) = 0$). Alors:
\calc{ \alpha(\vec v_t) &=& \phii_\alpha(u_\alpha(1)) - \phii_\alpha(tu_\alpha(1)t\inv)\\
&=& 0- \phii_\alpha(u_\alpha(\bar\alpha(t).1))\\
&=& -\omega(\bar\alpha(t)).
}

 Donc $\vec v_t$ vérifie les égalités annoncées. L'unicité de $\vec v_t$ est claire car $\phi$ engendre $\vec V^*$.\cqfd

Nous avons vu (\ref{cor:bonnes familles min et max})  que les familles minimale et maximale de parahorique $P$ et $\bar P$ sont fonctoriellement de bonnes familles de parahoriques.

D'après \cite{masures2}, il existe une bonne famille de parahoriques $Q$ pour $\D$, qui vérifie en outre \plienp(sph) et \pdec. En particulier, tous les résultats de \ref{section:construction generale} s'appliquent à $Q$.

%Sauf contre-indication, on note dans la suite $\I_\K=\I(Q)$.

\subsubsection{Groupes de Kac-Moody presque déployés}
\label{soussoussection:rappels presque deployes}

Soit maintenant un groupe de Kac-Moody $G$ presque déployé sur un corps $\K$. On suppose $G$ déployé sur la clôture séparable $\K_s$ de $\K$. Il existe alors une extension galoisienne finie $\L$ de $\K$, incluse dans $\K_s$ qui déploie $G$. On fixe un tore $\K$-déployé $T_\K$. Il existe un tore maximal $T$ $\K$-défini contenant $T_\K$, et quitte à remplacer $\L$ par une autre extension galoisienne un peu plus grande, on peut supposer $T$ $\L$-déployé.

 Le groupe $G(\L)$ est tel que décrit au paragraphe précédent, et il admet une masure bordée $\I_\L$. On notera parfois $G$ pour $G(\L)$, $\iv$ pour $\iv(\L)$ et $\I$ pour $\I_\L$\\

On suppose $\K$ muni d'une valuation non triviale $\omega: \K\rightarrow \R\cup{\infini}$ (en particulier, $\K$ est infini). Soit  $\Gamma = \mathpzc{Gal}(\L|\K)$. On suppose $\K$ complet, ce qui entraîne en particulier que $\omega$ se prolonge de manière unique à $\K_s$, et donc à $\L$, la valuation obtenue est alors nécessairement $\Gamma$-stable.\\% On note $\Lambda=\omega(\L)$ et $\Lambda_\K=\omega(\K)$.\\

Bertrand Rémy a décrit dans \cite{remy} chapitres 11,12 un immeuble vectoriel, et sa donnée radicielle pour $G(\K)$. Voici un résumé:\\
 
 Le groupe de Galois $\Gamma$ agit par définition d'un groupe presque déployé sur $G(\L)$ et sur son algèbre de lie $\mathfrak g(\L)$, en vérifiant $\sigma( Ad(g).x)=Ad(\sigma g).\sigma x$ pour $g\in G$ et $x\in \mathfrak{g}$. Ceci définit une action de $\Gamma$ sur $\iv_\L$, par automorphismes d'immeubles, qui préserve les immeubles positif et négatif, mais qui ne préserve pas le type des facettes.On note $\ikv=\iv^\Gamma$, on prouve qu'il s'agit d'un immeuble pour $G(\K)$.\\

Les appartements de  $\vec\I(\K)$, qu'on appellera les $\K$-appartements vectoriels, sont les parties maximales de $\iv(\K)$ de la forme $\vec E\cap \vec A(T)$ avec $\vec E$ un sous-espace vectoriel de $\vec V(T)$ rencontrant $\vec A\sph(T)$. Ils sont en bijection avec l'ensemble des tores $\K$-déployés maximaux de $G(\K)$. Comme $\Gamma$ ne respecte pas les types, un $\K$-appartement n'est généralement pas une réunion de facettes de $\iv(\L)$. Un appartement de $\ikv$ est toujours inclus dans un appartement de $\iv$ qui est $\Gamma$-stable, ceci correspond à l'inclusion d'un tore $\K$-déployé maximal dans un tore maximal défini sur $\K$. Réciproquement, si $T_d\subset T$ est l'inclusion d'un tore déployé maximal dans un tore maximal défini sur $\K$, alors $\vec A(T)$ est un appartement $\Gamma$-stable, dont le lieu des points fixes sous $\Gamma$ est $\vec A_\K(T_d)$. Le fixateur du $\K$-appartement $\vec A_\K(T_d)$ est le centralisateur de $T_d$, noté $Z(T_d)$. Le fixateur de $\vec A_\K(T_d)$ dans $G(\K)$ est donc $ Z(T_d)(\K) =Z(T_d)^\Gamma$. Le stabilisateur de $\vec A_\K(T_d)$ est $N(T_d)$, il agit donc sur $\vec A_\K$ via $W(\vec A_\K) := N(T_d)/Z(T_d)$.

Soit $\vec A_\K$ un $\K$-appartement vectoriel inclus dans un appartement vectoriel $\vec A$. Ses murs sont les $\vec A_\K\cap M$ pour $M$ un mur d'un appartement $\vec A$ contenant $\vec A_\K$ et tel que $\vec A_\K\cap M\cap \vec A\sph\not=\vide$ (cette dernière condition signifie que $\vec A_\K\cap M$ est un mur \textit{réel}). Ces murs font de $\vec A_\K$ un complexe de Coxeter, dont le groupe de Coxeter est $W(\vec A_\K)$. Les facettes de $\vec A_\K$ sont des réunions de parties de la forme $\vec f^\Gamma$ pour $\vec f$ une facette de $\vec A$ $\Gamma$-stable. Une facette de $\vec A_\K$ est sphérique si et seulement si elle coupe $\vec A\sph$ (et donc contient une $\vec f^\Gamma$, avec $\vec f$ une facette sphérique de $\vec A$).
 Les racines pour $\vec A_\K$ sont les $\alpha|_{\vec A_\K}$ pour $\alpha\in\phi(\vec A)$ telle que $\ker(\alpha)\cap \vec A_\K$ est un mur réel, autrement dit rencontre $\vec A\sph$ et n'est pas égal à $\vec A_\K$ (donc $\alpha|_{\vec A_\K}\not=0$). L'ensemble des racines de $\vec A_\K$ noté $\phi(\vec A_\K)$ ou $\phi(T_d)$ est un système de racines, pas forcément réduit contrairement à $\phi(\vec A)$, et son groupe de Weyl est $W(\vec A_\K)$.
  Les $\K$-racines géométriques sont les demi-appartements de $\vec A_\K$.
  
  Pour tout $a\in\phi(T_d)$, on note, conformément à \ref{soussection:hypotheses de descente}, $\phi_a=\ens{\alpha\in\phi(T)}{\alpha|_{\vec A_\K}\in \R^{+*}.a}$ et $U_a = \eng{U_\alpha}{\alpha\in\phi_a}$. Le sous-groupe radiciel associé à $a$ est alors $U_{a}(\K)$. La famille $\D_\K:=(G(\K), (U_{a}(\K))_{a \in\phi(T_d)})$ est une donnée radicielle pour $G(\K)$ (\cite{remy} 12.6.3). Il y a ici un petit conflit de notation puisque le groupe jouant le rôle de $T$ pour la donnée radicielle $\D_\K$, c'est-à-dire $\bigcap_a N_{G(\K)}(U_a(\K))$, est en fait $Z(T_d)(\K)$.
  
  L'immeuble que définit cette donnée radicielle  n'est pas exactement $\ikv=\iv^\Gamma$, cependant ce dernier en est tout de même une bonne réalisation géométrique (voir  \cite{remy} 12.4.4 et 13.4.2). L'immeuble $\ikv$ correspond en fait à l'immeuble $\iv\becb$ de la partie \ref{section:descente}, et l'immeuble de la donnée radicielle $\D_\K$ correspond à $\iv\bec$.\\
  
  Le cas d'un groupe de Kac-Moody presque déployé comporte une simplification notable par rapport à la situation générale étudiée en \ref{section:descente}: pour toute facette sphérique $\vec f_\K$ de $\vec A_\K$, il n'existe qu'une seule facette de $\vec A$ contenant un ouvert de $\vec f_\K$. En effet, dans le cas contraire il existerait un mur de $\vec A$ coupant l'intérieur de $\vec f_\K$. Comme $\vec f_\K$ est sphérique, ce mur coupe $\vec A\sph\cap \vec A_\K$, et donc induit un mur de $\vec A_\K$ coupant $\vec f_\K$, ce qui est impossible.

\subsection{Action du groupe de Galois}
\label{soussection:action de galois}

On définit dans ce paragraphe une action du groupe $\Gamma$ sur l'immeuble $\I_\L$.
On rapelle qu'on a fixé un tore $\K$-déployé maximal $T_\K$ inclus dans un tore maximal $\K$-défini et $\L$-déployé $T$. Le tore $T$ est donc $\Gamma$-stable, ce qui permettra de définir une action de $\Gamma$ sur l'appartement $A(T)$. L'extension de cette action à $\I$ ne posera ensuite aucun problème.

% Comme l'extension $\K\subset \L$ est finie, $\Gamma$ est également fini, et le théorème de point fixe de Bruhat entraine l'existence d'un point $\Gamma$-fixe dans l'enveloppe convexe de n'importe quelle orbite de $\Gamma$ sur un immeuble.

% Quelques mots sur la topologie de $\Gamma$: On a $\Gamma= \limproj \ens{\gal(\mathbb E | \mathbb K)}{ \mathbb E \subset \L \text{ est une extension galoisienne finie de } \K}$, et la topologie de $\Gamma$ est induite par la topologie produit sur $\prod \ens{\gal(\mathbb E | \mathbb K)}{ \mathbb E\subset \L \text{ extension galoisienne finie de } \K}$. Une base d'ouverts est donc l'ensemble des $\ens{\sigma\in\Gamma}{\sigma_{| \mathbb E}=\sigma_0}$ pour $\mathbb E$ une extension galoisienne finie de $\K$ et $\sigma_0\in \gal(\mathbb E|\K)$. Les sous-groupes ouverts de $\Gamma$ sont donc d'indice fini, et lorsque $\Gamma$ agit continument sur un ensemble discret, les orbites sont finies.\\

\subsubsection{Action de $\Gamma$ sur $A$}

Dans cette partie, on note $\vec A=\vec A(T)$, $A=A(T)$, $\phi=\phi(T)$, $\vec V=\vec V(T)$. Le fait que $T$ est défini sur $\K$ implique que $\vec A$ est $\Gamma$-stable, et cette action de $\Gamma$ sur $\vec A$ s'étend par linéarité à $\vec V$.
 De plus $\Gamma$ permute les racines et les sous-groupes radiciels relatifs à $T$, de manière compatible à son action sur l'algèbre de Lie. On a précisément, pour $\sigma\in\Gamma$ et $\alpha\in \phi$:
\[ \sigma u_\alpha(k) = u_{\sigma \alpha}(\sigma(k).k_\alpha^\sigma)\]

où $k_\alpha^\sigma$ est défini par $\sigma.e_\alpha=k_\alpha^\sigma . e_{\sigma \alpha}$.\\

D'où \[ \sigma U(\alpha,\lambda)= U(\sigma(\alpha), \lambda+\omega_\alpha^\sigma)\]
où $\omega_\alpha^\sigma\in\R$ vaut $\omega(k_\alpha^\sigma)$.

Par conséquent, on veut définir une action de $\Gamma$ sur $A$ telle que $\sigma$ envoie $D(\alpha,\lambda)$ sur 
$D(\sigma(\alpha), \lambda+\omega_\alpha^\sigma)$. Cette action doit être compatible avec l'action vectorielle de $\Gamma$ sur $\vec V$, il ne reste donc qu'à déterminer l'image du point $o$. Voici en quelques mots la justification de la définition qui va suivre:\\
Nous voulons $\sigma M(\alpha,\lambda)=M(\sigma \alpha,\lambda+\omega_\alpha^\sigma)$, c'est-à-dire
\[ \ens{\sigma x\in Y_0}{\alpha(x)+\lambda=0}= \ens{x\in Y_0}{\sigma \alpha(x)+\lambda+\omega_\alpha^\sigma =0}\]

Il faut donc que $\alpha(\sigma\inv(x))=(\sigma\alpha)(x)+\omega_\alpha^\sigma$, et ce pour tout $\alpha\in\phi$. Soit $\vec u\in\vec V$ tel que $x=o+\vec u$, alors $(\sigma \alpha)(x)=(\sigma\alpha)(\vec u)=\alpha(\sigma\inv(\vec u))$. D'autre part, si nous définissons une action affine de $\sigma$ sur $A$, dont la partie vectorielle coïncide avec l'action déjà connue de $\sigma$ sur $\vec V$, nous aurons $\sigma\inv(x)=o+\fl{o\sigma\inv(o)} +\sigma\inv(\vec u)$, puis $\alpha(\sigma\inv(x))=\alpha(\fl{o\sigma\inv(o)}) +\alpha(\sigma\inv(\vec u))$. Finalement, il nous faut faire en sorte que $\alpha(\fl{o\sigma\inv(o)})=\omega_\alpha^\sigma$.\\

\begin{lemme}
Soit $\sigma\in\Gamma$, on note $\omega_\alpha=\omega_\alpha^\sigma$ pour $\alpha\in\phi$. Soit $S\subset W(T)$ un système générateur de réflexions, et $\Pi=(\alpha_s)_{s\in S}\subset \phi$ la base de $\phi$ correspondante.
Soient $\alpha,\beta\in\phi$, $s\in S$ tels que $\alpha=s.\beta=\beta - <\alpha_s,\beta>\alpha_s$. Alors $\omega_\alpha=\omega_\beta-<\alpha_s,\beta>\omega_{\alpha_s}$.
\end{lemme}
\rema On rappelle que par définition, $<\alpha,\beta>=\beta(\alpha^\vee)$.\\

\pv On note pour $\gamma\in \phi$, $k_\gamma=k_\gamma^\sigma$, donc $\omega_\gamma=\omega(k_\gamma)$. On note également $k_s=k_{\alpha_s}^\sigma$, et on reprend les notations de \cite{remy} chapitres 7 et 8.

Puisque $\alpha=s.\beta$, on a $e_\alpha=\pm s^*e_\beta=\pm Ad(n_{\alpha_s}(1)).e_\beta$. (La notation $s^*$ désigne un automorphisme de $\g$ qui relève l'élément du groupe de Weyl $s\in W(\phi)$, en caractéristique nulle, $s^*=\exp(\text{ad } f_s)\exp(\text{ad } e_s)\exp(\text{ad } f_s)$). Appliquant $\sigma$ à cette égalité, on trouve:\[
\begin{array}{rl}
k_\alpha e_{\sigma\alpha}& =\pm Ad(\sigma(n_{\alpha_s}(1))).\sigma(e_\beta)\\
&=\pm k_\beta Ad(n_{\sigma \alpha_s}(k_s))e_{\sigma\beta}\\
&= \pm k_\beta Ad(n_{\sigma \alpha_s}(1).k_s^{-\alpha^\vee_{\sigma s}}) e_{\sigma\beta} \\
&= \pm k_\beta k_s^{ -<\alpha_s,\beta>} Ad(n_{\sigma \alpha_s}(1)).e_{\sigma\beta}\\
&=\pm\pm k_\beta k_s^{ -<\alpha_s,\beta>} e_{\sigma(s).\sigma(\beta)} \\
&=\pm\pm k_\beta k_s^{ -<\alpha_s,\beta>} e_{\sigma \alpha}
\end{array}\]

D'où $\omega_\alpha=\omega_\beta-<\alpha_s,\beta>\omega_{\alpha_s}$.\cqfd\\

Comme $\Pi$ est une base de $\vec V^*$, il existe pour tout $\sigma\in\Gamma$ un vecteur $\vec v_\sigma\in \vec V$ tel que pour tout $s\in S$, $\alpha_s(\vec v_\sigma)=\omega_{\alpha_s}^{\sigma\inv}$. Alors d'après le lemme, on a aussi $\alpha(\vec v_\sigma)=\omega_\alpha^{\sigma\inv}$ pour toute  $\alpha\in\phi$.\\
  
%  Si la famille $(\alpha_s)_{s\in S}$ n'est pas génératrice, $\vec v_\sigma$ n'est pas unique. On choisit dans ce cas un supplémentaire $\vec Y_1$ au sous-espace $\bigcap_{s\in S} \ker(\alpha_s)=\bigcap_{\alpha\in\phi}\ker(\alpha)$, et on prend $\vec v_\sigma$ dans $\vec Y_1$ pour tout $\sigma\in\Gamma$. Le choix est alors unique.  Comme $\bigcap_{\alpha\in \phi}\ker(\alpha)$ est $\Gamma$-stable, on peut choisir $\vec Y_1$ $\Gamma$-stable.\\
%  Comme dans la définition de l'appartement $A(T_0)$ on quotiente par $\bigcap_{\alpha\in\phi}\ker(\alpha)$, le choix de $\vec Y_1$ n'importe pas.\\
 
 \begin{defin}
 Pour $\sigma\in\Gamma$ et $\vec u\in \vec V$, on pose
 \[\sigma(o+\vec u)=o+\vec v_\sigma +\sigma(\vec u)\]
 \end{defin}

\begin{prop}\label{prop:action de Gamma sur l appart}
 La formule ci-avant s'étend à  une action de groupe de $\Gamma$ sur $A$, par isomorphismes d'appartements, compatible avec celle de $N(T)$ au sens où $(\sigma n).x=\sigma(n(\sigma\inv.x))$, et compatible avec l'action de $\Gamma$ sur les sous-groupes radiciels au sens où $\sigma.D(\alpha,\lambda)=D(\sigma\alpha,\mu)$ si $\sigma.U_{\alpha,\lambda}=U_{\sigma\alpha,\mu}$. L'ensemble de ses points fixes est l'adhérence d'un sous-espace affine de $\inte A$ dirigé par $\vec V^\Gamma$.
\end{prop}

\demo\\
Pour vérifier qu'il s'agit d'une action de groupe, on vérifie la condition de cocycle attendue sur les $\omega_\alpha^\sigma$. Soient $\alpha\in\phi$, $\sigma,\gamma\in\Gamma$, on calcule dans l'algèbre de Lie $\mathfrak g$:
\[ \gamma\sigma.e_{\alpha}= \gamma(k_\alpha^\sigma e_{\sigma\alpha})=\gamma (k_\alpha^\sigma)k_{\sigma\alpha}^\gamma e_{\gamma\sigma\alpha} \]
et d'autre part:
\[\gamma\sigma.e_\alpha= k_\alpha^{\gamma\sigma}e_{\gamma\sigma\alpha}\]
D'où la relation $\omega_\alpha^{\gamma\sigma}=\omega(\gamma (k_\alpha^\sigma))+ \omega_{\sigma\alpha}^{\gamma}
=\omega_\alpha^\sigma  + \omega_{\sigma\alpha}^{\gamma} $, la deuxième égalité car $\Gamma$ préserve la valuation $\omega$.\\

 Maintenant, si $\gamma,\sigma\in \Gamma$, on a pour $\vec u\in \vec V$, $\gamma(\sigma(o+\vec u))=o+\vec v_\gamma +\gamma(\vec v_\sigma)+\gamma\sigma(\vec u)$. Par conséquent, il faut vérifier que $\vec v_{\gamma\sigma}=\vec v_\gamma +\gamma(\vec v_{\sigma})$.\\
 Soit $\alpha\in\phi$, alors $\alpha(\vec v_\gamma +\gamma(\vec v_{\sigma})) = \omega_\alpha^{\gamma\inv}+(\gamma\inv\alpha)(\vec v_\sigma)=\omega_\alpha^{\gamma\inv}+\omega_{\gamma\inv\alpha}^{\sigma\inv}$. Par le calcul précédent, ceci vaut $\omega_\alpha^{\sigma\inv\gamma\inv}=\omega_\alpha^{(\gamma\sigma)\inv}=\alpha(\vec v_{\gamma\sigma})$. Comme $\vec v_{\gamma\sigma}$ est l'unique vecteur de $\vec V$ vérifiant cette relation pour tout $\alpha\in\phi$, on obtient bien $\vec v_{\gamma\sigma}=\vec v_\gamma +\gamma(\vec v_{\sigma})$, et l'action de $\Gamma$ sur $A$ est bien une action de groupe.\\

Par construction, $\Gamma$ agit par automorphismes affines et préserve l'ensemble des murs, cette action, a priori sur $\inte A$, s'étend donc à une action sur $A$  par automorphismes d'appartements.\\

 Montrons la compatibilité avec l'action de $N$. Pour commencer, un élément de la forme $n_\alpha(l)$, avec $\alpha\in\phi$ et $l\in\L$, induit sur $A$ la réflexion $r$ selon le mur $M(\alpha,\omega(l))$. L'élément $\sigma(n_\alpha(l))=n_{\sigma\alpha}(\sigma(l)k_\alpha^\sigma)$ induit alors la réflexion selon le mur $M(\sigma \alpha, \omega(l)+\omega_\alpha^\sigma)=\sigma.M(\alpha,\omega(l))$. C'est bien la conjugaison par l'action de $\sigma$ de la réflexion $r$. Comme $N$ est engendré par les $n_\alpha(l)$ et $T$, il reste à vérifier que $(\sigma t).x=\sigma(t(\sigma\inv x))$ pour $t\in T$.\\
  Soit donc $t\in T$, d'après \ref{prop:action de T KM}, $t$ agit sur $\Gamma$ par translation selon le vecteur $\vec v_t$ défini par 
  \[ \alpha(\vec v_t)=-\omega(\bar\alpha(t)),\:\: \forall \alpha\in \phi\subset \vec V^* \]
  
 Or pour $\alpha\in \phi$, $\alpha\circ\sigma=\sigma\inv.\alpha$ est encore dans $\phi$ et $\barre{\sigma\inv.\alpha}  = \sigma\inv\circ\bar\alpha\circ\sigma\in T^*$ d'où: 
   \[ \alpha(\sigma\vec v_t)= -\omega( \sigma\inv\circ\bar\alpha\circ\sigma(t))  =-\omega(\bar\alpha(\sigma t)),\:\: \forall \alpha\in T^*\subset \vec Y^*. \]
   (la deuxième égalité car la valuation $\omega$ est $\Gamma$-stable.)
   
   Ceci entraine $\vec v_{\sigma t}=\sigma(\vec v_t)$ d'où la relation de compatibilité entre les actions de $t$ et de $\sigma$.\\

 Enfin, comme $\Gamma$ est fini, son action sur $\inte A$ fixe un point, disons $o'$. Dès lors, $A^\Gamma=\barre{o'+\vec V^\Gamma}$.\cqfd\\

 %\begin{defin}
%On note $\vec Y(T_d)(\K)= \vec Y_0^\Gamma$ et $Y_0(\K)=Y_0^\Gamma$.
%\end{defin} 

%\begin{prop}\label{prop:caracteres du tore deploye} L'espace $\vec Y(T_d)(\K)$ ne dépend effectivement que de $T_d$ et pas de $T_0$. Il s'agit de $(T_d)_*\otimes \R$.\end{prop}

\subsubsection{Action de $\Gamma$ sur $\I$}

 On dira qu'une famille $Q$ de parahoriques sur $A(T)$ est $\Gamma$-stable si pour tout $a\in A(T)$ et tout $\sigma\in\Gamma$, $\sigma.Q(a) = Q(\sigma a)$. Les familles minimale et maximale de parahoriques sont clairement $\Gamma$-stables, il en est de même de la famille construite dans \cite{masures2}.

Soit $Q$ une bonne famille de parahoriques $\Gamma$-stable, par définition $\I(Q)=G\times A(T)/\sim$, où $\sim$ est la relation d'équivalence définie par 
$(g,x)\sim (h,y)\ssi \exists n\in N(T) \tq y=nx$ et $g\inv hn\in Q(x)$.
Comme $N(T)$ est $\K$-défini, l'action de $\Gamma$ sur $G\times A(T)$ par $\sigma.(g,x)=(\sigma g,\sigma x)$ passe au quotient et définit une action sur $\I(Q)$, cette action stabilise $A(T)$ et prolonge l'action définie précédemment.\\

\begin{lemme}
Le groupe $\Gamma$ agit sur $\I$ par automorphismes de masure bordée, c'est-à-dire qu'il préserve l'ensemble des appartements de $\I$ et induit entre deux appartements un isomorphisme d'appartements.

 De plus, pour tous $x\in\I$, $g\in G$, et $\sigma\in\Gamma$, $\sigma(gx)=\sigma(g)\sigma(x)$. Pour tout tore maximal $T'$, $\sigma(A(T'))=A(\sigma(T'))$. Enfin, pour toute facette $\vec f$ de $\iv$, $\sigma.\I_{\vec f} = \I_{\sigma \vec f}$.
\end{lemme}

\pv\\
Soit $Z$ un appartement, soit $g\in G$ tel que $Z=g.A(T)$. Alors $\sigma Z=\sigma(g).A(T)$, c'est donc un appartement, l'appartement vectoriel lui correspondant est $\sigma(g).\vec A(T)=\sigma(g.\vec A(T))=\sigma(\vec Z)$. De plus, en notant $\sigma_Z$ la restriction de l'action de $\sigma$ à $Z$, et $\sigma_{A(T)}$ sa restriction à $A(T)$, on a $\sigma_Z= \sigma(g)\circ \sigma_{A(T)} \circ g\inv$. Comme $g\inv$, $\sigma(g)$ et $\sigma_{A(T)}$ induisent des isomorphismes entre les appartements concernés, $\sigma_Z$ est bien un isomorphisme d'appartements.

La relation $\sigma(gx)=\sigma(g)\sigma(x)$ vient de la définition de l'action de $\Gamma$.
Si $T'$ est un autre tore maximal, soit $g\in G$ tel que $T'=gTt\inv$, alors $A(T')=g.A(T)$, d'où $\sigma.A(T')=\sigma(g).\sigma(A(T)) = \sigma(g).A(T) = A(\sigma(g)T\sigma(g)\inv)$. Mais $\sigma(g)T\sigma(g)\inv = \sigma(gTg\inv)=\sigma(T')$.\\

Enfin, si $\vec f$ est une facette de $\vec A$, alors $\I_{\vec f}= P(\vec f).A_{\vec f}$ donc $\sigma.\I_{\vec f} = \sigma(P(\vec f)).\sigma(A_{\vec f}) = P(\sigma\vec f).A_{\sigma\vec f} = \I_{\sigma\vec f}$. Le cas où $\vec f\not\subset \vec A$ s'obtient par conjugaison par un élément $g\in G$ tel que $g.\vec f\subset \vec A$.
\cqfd\\

En particulier, si dans un appartement $Z$ $y=x+_Z\vec v$, alors dans $\sigma Z$ $\sigma(y)=\sigma(x)+_{\sigma Z}\sigma(\vec v)$.

\subsection{Action du normalisateur du tore déployé}

%\rema Cette partie n'est plus vraiment utile. Peut-être a-t-elle tout de même un intérêt ?

 Nous avons obtenu pour tout tore maximal $\K$-défini et $\L$-déployé $T$ contenant $T_\K$ une partie $A(T)^\Gamma$ stable par $N(T)(\K)$ .
Nous allons voir qu'on peut, au moins sous l'hypothèse que le corps $\L$ est "maximalement complet", trouver un autre espace affine dirigé par $\vec A_\K$ qui soit stable par $N(T_\K)(\K)$, c'est-à-dire qui permette de vérifier (DM 4).

Cette partie n'est pas utilisée dans la suite. \\

Soit $\vec f_\K$ une facette maximale de $\vec A_\K$, $\vec f$ une facette de $\vec A$ contenant un ouvert de $\vec f_\K$.
 Soit $\mathcal J:=\I_{\vec f}$, c'est un immeuble affine car $\vec f$ est sphérique, c'est en fait l'immeuble de Bruhat-Tits du groupe $M_{\vec A}(\vec f) = \fix_G(\vec A_\K) = Z(T_\K)$.

   Les groupes $Z(T_\K)$ et $\Gamma$ agissent sur $\mathcal J$, et $Z(T_\K)(\K)$ préserve l'ensemble $\mathcal J^\Gamma$, qui contient le singleton $pr_{\vec f}(A(T)^\Gamma)$.\\

 \begin{lemme}
L'ensemble  $\mathcal J^\Gamma$ est une partie non vide, bornée et convexe de $\mathcal J$.
 \end{lemme}
\pv L'immeuble $\mathcal J$ est l'immeuble de Bruhat-Tits de $Z(T_\K)$, c'est donc aussi l'immeuble du semi-simplifié de $Z(T_\K)$, c'est-à-dire de $Z(T_\K)/Z(Z(T_\K))$. Comme $T_\K\subset Z(Z(T_\K))$, et $T_\K$ est un tore $\K$-déployé maximal, ce semi-simplifié n'a pas de tore $\K$-déployé, il est anisotrope. Alors la proposition 5.2.1 de \cite{these-guy} entraine que $\mathcal J^\Gamma$ est borné.

C'est une partie convexe car $\Gamma$ agit sur $\mathcal J$ par automorphismes d'immeuble, et non vide car elle contient $pr_{\vec f}(A(T)^\Gamma)$.\cqfd

\rema Sans supposer que la valuation de $\K$ est discrète, il n'est en général pas clair que l'immeuble de $Z(T_\K)$ contienne un point $\Gamma$-fixe (voir \cite{bruhat-tits2} 5.1.6). Ici, c'est le fait d'avoir choisi $\L$ de manière à assurer l'existence d'un tore maximal $\L$-déployé et $\K$ défini contenant $T_\K$ qui a cette conséquence.\\

\begin{prop}\label{prop:existence de AK}
Si l'immeuble $\mathcal J$ est complet, il existe un sous-espace affine $Y(\K)$ d'un appartement $A$ de $\I$, stable par $N(T_\K)(\K)$, fixe par $\Gamma$, dirigé par $\vect_{\vec A}(\vec A_\K)$.
\end{prop}

\rema L'immeuble $\mathcal J$ est complet si et seulement si $\L$ est "maximalement complet", d'après \cite{bruhat-tits} 7.5.4 et 7.5.5, ce qui est le cas par exemple dès que la valuation de $\L$ est discrète.\\

\demo\\
Le produit $Z(T_\K)(\K)\times \Gamma$ agit sur $\mathcal J$ en stabilisant $\mathcal J^\Gamma$. D'après le lemme et le fait que $\mathcal J$ est complet, le théorème de point fixe de Bruhat-Tits prouve l'existence d'un point $p\in \mathcal J^\Gamma$ fixe par $Z(T_\K)(\K)$. Soit $T'$ un tore maximal de $Z(T_\K)$ tel que $p\in A(T')_{\vec f}$. Soit $Y = pr_{\vec f}\inv(p)\cap A(T')$, il s'agit de l'adhérence d'un espace affine de $A(T')$ stable par $N(T_\K)(\K)$ et par $\Gamma$, et dirigé par $\vect_{\vec A(T')}(\vec f)$. 
 Alors $\Gamma$ agit sur $Y$ par automorphismes affines, avec des orbites finies, donc $Y$ admet lui même un point $\Gamma$-fixe $p$. Finalement, $Y(\K):= Y^\Gamma = \barre{ p+ \vect_{\vec A(T')}(\vec f^\Gamma)} = \barre{ p+ \vect_{\vec A}(\vec A_\K)}$ convient.\cqfd

On obtient de la sorte un espace $Y(\K)\subset \I^\Gamma$ sur lequel agit $N(T_\K)(\K)$. Cependant, on ne sait pas si il existe un appartement $A$ le contenant tel que $A\cap \I^\Gamma = Y(\K)$. Autrement dit, en modifiant $Y(\K)$ pour satisfaire à (DM 4), on a perdu la condition (DM2).\\

Pour résoudre cette difficulté, nous aurons besoin d'hypothèses sur le groupe $G$ ou le corps $\K$, et nous devrons choisir un $\I\die$ plus petit que $\I^\Gamma$.

\subsection{Descente}

\subsubsection{Vérification des premières conditions de descente}

\fprop{\label{prop:verifications des premieres conditions} Soit $G$ un groupe de Kac-Moody presque déployé sur un corps $\K$, déployé sur une extension galoisienne $\L$. On note $\D_\K$ et $\D_\L$ les données radicielles pour $G(\L)$ et $G(\K)$, correspondant à un tore maximal $\L$-déployé $T_\L$ et un tore $\K$-déployé maximal $T_\K$ inclus dans $T_\L$. 
Alors le couple $(\D_\L,\D_\K)$ vérifie les conditions (DSR), (DDR) et (DIV). 

Si $\K$ est muni d'une valuation non triviale $\omega$, alors toute valuation $\phii_\L$ de $\D_\L$ associée comme en \ref{prop:existence de valuation} vérifie (DV 2).\\

 Supposons de plus  $T_\L$ $\K$-défini,  soit $Q_\L$ une bonne famille de parahoriques vérifiant \plienpm(sph) pour $(\D_\L,\phii_\L)$, soit $\I_\L=\I(Q_\L)$, et  $\Gamma=\gal(\L|\M)$, qui agit donc sur $\I_\L$. 
Alors la donnée $(\D_\L,\phii_\L,\D_\K,\I_\L^\Gamma)$ satisfait aussi aux conditions (DM 1), (DM 2).
} 
 
Les conditions manquantes sont donc (DM 3) et (DM 4), ainsi que (DV 1), mais on peut toujours remplacer $\phii_\L$ par une valuation équipollente pour satisfaire à cette dernière.\\

 On rappelle que, contrairement à ce que les notations pourraient faire croire, le groupe jouant le rôle de $T$ dans la donnée radicielle $\D_\K$ est en fait $Z(T_\K)(\K)$, qui contient en général strictement $T_\K$. 
 
 De plus, on a pris $\vec A_\K= \vec A_\L^\Gamma$, de manière à voir $\iv_\K$ comme une partie de $\iv_\L$. C'est donc en général un complexe de Coxeter non essentiel, et ce n'est pas l'appartement obtenu abstraitement à partir de $\phi(T_d)$.
 
 Pour le reste, on notera avec $\K$ en indice tous les objets habituellement obtenus à partir d'une donnée radicielle.
 
Voici le dictionnaire entre les objets $\K$-rationels considérés ici et les objets de la partie \ref{section:descente} que l'on obtient:

\liste{
\item $\D=\D_\L$, $\D\bec=\D_\K$,
\item $G\bec= G(\K)$, $T\bec= Z(T_\K)(\K)$, $N\bec= N(T_\K)(\K)$,
\item $\phi\bec= \phi(T_\K)$, $U_a\bec = U_{a,\K}$ pour tout $a\in\phi(T_\K)$,
\item $\vec V=\vec V(T_\L)$, $\vec V\bec= \vect_{\vec V(T_\L)}(\vec A_\K)$,
\item $\iv\becb = \ikv$, $\vec A\becb=\vec A_\K$,
\item $\phii=\phii_\L$, $\phii\bec = \phii_\K$,
\item pour la deuxième partie de la proposition, on aura $\I\die = \I_\L^\Gamma$, donc $A\die=\inte A(T_\L)^\Gamma$.\\
}

 \demo
 On note $\vec V_\L=\vec V(T_\L)$, $\vec V_\K=\vect_{\vec V_\L}(\vec A_\K)$, $A_\L=A(T_\L)$ et $\phi_\K=\phi(T_\K)$.\\
 
 \liste{
\item (DSR): On a $\phi_\K=\ens{ \alpha|_{\vec A_\K}}{\alpha\in\phi\et \ker(\alpha)\cap \vec A_{\L sph}\not=\vide}$. Par \cite{remy} 12.4.4, toutes les facettes de $\iv_\K$ coupent $\iv$; par 12.6 $\phi_\K$ est à base libre; enfin si $\vec f_\K$ est une facette sphérique de $\vec A_\K$, alors $\phi_{\K}^m(\vec f_\K)$ est fini, mais pour tout $a\in \phi_{\K}^m(\vec f_\K)$, $\phi_{\L}^m(\ker(a))$ est aussi fini (car $a$ est une $\K$-racine réelle), donc au final $\phi_L^m(\vec f_\K)$ est fini, donc $\vec f_\K$ coupe au moins une facette sphérique.\\
 
 \item (DDR1): Pour tout $a\in\phi_\K$, on a $U_{ a,\K} = \engens{U_{\alpha,\L}} {\alpha\in\phi_\L\et \alpha|_{\vec V_\K}\in\{a,2a\}}^\Gamma$. Ceci est clairement inclus dans le groupe $U_a = \engens{U_{\alpha,\L}}{\alpha\in\phi_\L\et \alpha|_{\vec V_\K}\in\R^{+*}.a}$.\\

 \item (DDR2):  Pour tout $a\in\phi_\K$, si $\alpha\in\phi_\L$ est telle que $\alpha|_{\vec V_\K}\in\R^{+*}.a$, alors $\alpha|_{\vec V_\K}\in\phi_\K$. Comme $\phi_\K$ est un système de racines, il ne peut y avoir plus de deux tels $\alpha|_{\vec V_\K}$.\\
 
 \item (DDR3): Le groupe $T\bec$ pour la donnée radicielle $\D_\K$ est en fait $Z(T_\K)(\K)$, avec $Z(T_\K)$ le centralisateur du tore déployé maximal $T_\K$, autrement dit le fixateur de $\vec A_\K$. Or le groupe $Z$ défini dans \ref{soussection:hypotheses de descente} est précisément ce fixateur. D'où l'inclusion $Z(T_\K)(\K)\subset Z(T_\K)=Z$.
 
 Pour tout $a\in\phi_\K$, $u\in U_{\K a}$, $n(u)$ stabilise $\vec A_\K$ et donc normalise $Z$. Il agit sur $\vec A_\K$ comme une réflexion d'hyperplan $\ker a$, d'où $n(u) U_a n(u)\inv = U_{-a}$ et $n(u)U_{-a}n(u)\inv=U_a$.\\
 
\item (DIV): L'immeuble $\ikv = \iv^\Gamma$ est stable par $G(\K)$, et $\ikv\cap \vec A= \vec A_\K$. D'où (DIV).\\

 \item (DV2): Soit $t\in T_\K(\K)$, $t$ induit sur $A$ une translation de vecteur $\vec v_t$ et $t$ stabilise $A_\K$ donc $\vec v_t\in\vec V_\K$. Pour tout  $a\in \phi_\K$ et $k\in\R$, on a $tU_{a,k}t\inv = U_{a,k+a(\vec v_t)}$. Donc l'image de $\phii_{\K a}$ est stable par le groupe $\Z.a(\vec v_t)$, et il nous faut maintenant prouver qu'il existe $t\in T_\K(\K)$ tel que $a(\vec v_t)\not=0$.
 
 Le vecteur $\vec v_t$ est caractérisé par le fait que pour tout $\alpha\in\phi_\L$, $\alpha(\vec v_t)=-\omega(\bar\alpha(t))$ (proposition \ref{prop:action de T KM}). En particulier, si $\alpha\in\phi_\L$ est telle que $\alpha|_{\vec V_\K}=a$, alors $a(\vec v_t) = \alpha(\vec v_t) = -\omega(\bar\alpha(t)) = -\omega(\bar a(t))$, car $\bar a= \bar \alpha|_{T_\K}$.
  Mais $a(T_\K)\supset a(a^\vee(\K^*)) = (\K^*)^2$, car $\langle a,a^\vee\rangle= 2$, et le résultat découle de ce que $\omega$ est une valuation non triviale sur $\K$, donc sur $(\K^*)^2$.\\
}

On suppose maintenant $T_\L$ $\K$-défini, et on prend $\I\die = \I_\L^\Gamma$. On peut alors définir $A_\K=A_\L^\Gamma =\barre{A\die}$.\\
 
 \liste{
 \item (DM1): La partie $\I_\L^\Gamma$ est clairement stable par $G(\K)$. Soit $\vec f$ une facette sphérique. Si $\I_{\vec f}^\Gamma=\vide$, alors il s'agit bien d'une partie convexe et stable par $G(\K)$. Sinon, la facette $\vec f$ est $\Gamma$-stable, et $\Gamma$ agit sur $\I_{\vec f}$ par automorphismes d'immeuble affine, donc $\I_{\vec f}^\Gamma$ est convexe et $G(\K)$-stable.\\
 
 \item (DM2): Conséquence de \ref{prop:action de Gamma sur l appart}.
 
 }
 \cqfd

Dans toute la suite, on fixe une bonne famille de parahoriques $Q$ pour $\D_\L$.

\subsubsection{Corps intermédiaire}

%On suppose maintenant que la valuation de $\K$ est discrète, et que son corps résiduel est parfait.
Soit $T_\K$ un tore $\K$-déployé maximal. On suppose désormais qu'il existe une extension galoisienne $\M$ modérément ramifiée de $\K$, incluse dans $\L$, telle que le groupe réductif $Z(T_\K)$ soit quasi-déployé sur $\M$. Ceci est par exemple le cas dès que la valuation de $\K$ est discrète et que son corps résiduel est parfait.
\\

Il existe alors un sous-groupe de  Borel $B_\M$ du groupe réductif $Z(T_\K)$ défini sur $\M$. Ce Borel contient un tore maximal $T_\L$ défini sur $\M$. Comme $T_\K\subset Z(B_\M)$, $T_\L$ contient $T_\K$. Soit $T_\M$ la partie $\M$-déployée de $T_\L$, c'est-à-dire le groupe engendré par les image des cocaractères de $T_\L$ fixes par $\Gamma_\M:=\gal(\L|\M)$. Cette partie contient encore $T_\K$ puisque tout cocaractère de $T_\K$ est un cocaractère de $T_\L$ fixe par $\Gamma_\M$.

Prouvons que $T_\L= Z_{Z(T_\K)}(T_\M)$. A priori, $Z_{Z(T_\K)}(T_\M)$ est le sous-groupe de Lévi de $Z(T_\K)$ engendré par $T_\L$ et par les groupes radiciels $U_\alpha$ pour $\alpha\in\phi(T_\L)$ tels que $\alpha(T_\M)=\{1\}$ (cette condition entraine automatiquement $\alpha\in\phi^m(T_\K)$, donc $U_\alpha\subset Z(T_\K)$).

 Le système de racines de $Z(T_\K)$ par rapport au tore maximal $T_\L$ est $\phi^m(T_\K)$, où $\phi=\phi(T_\L)$ est le système de racines pour $G$ par rapport  à $T_\L$. Soit $\Pi$ la base de $\phi^m(T_\K)$ correspondant à $B_\M$, et $\Pi^\vee$ sa base duale. Comme $B_\M$ est $\Gamma_\M$-stable, ces bases sont permutées par $\Gamma_\M$. Une base des cocaractères de $T_\M$ est alors l'ensemble des $\sum_{\rho\in \mathcal O} \rho$, pour $\mathcal O$ une $\Gamma_\M$-orbite dans $\Pi^\vee$.
 
 Soit maintenant $\alpha\in\phi^m(T_\K)$ tel que $\alpha(T_\M)=\{1\}$. Alors $\alpha$ s'annule sur tout cocaractère de $T_\M$, donc pour toute orbite $\mathcal O$ dans $\Pi^\vee$, $\alpha(\sum_{\rho\in \mathcal O} \rho) = \sum_{\rho\in \mathcal O} \eng{\alpha,\rho} =0$. Comme $\mathcal O\subset \Pi^\vee$, tous les $\eng{\alpha,\rho}$ sont de même signe, donc finalement, ils sont tous nuls. Au total, $\alpha$ s'annule sur tous les élément de $\Pi^\vee$, ceci est impossible. Il n'existe donc pas de racine $\alpha\in\phi^m(T_\K)$ s'annulant sur $T_\M$. Donc $Z_{Z(T_\K)}(T_\M)=T_\M$.
 
 Vérifions que $T_\M$ est un tore $\M$-déployé maximal. Si $T$ est un tore $\M$-déployé contenant $T_\M$, alors $T\subset Z(T_\M) = T_\L$. Si $\rho$ est un cocaratère de $T$, comme $T$ est $\M$-déployé, $\rho$ est $\Gamma_\M$-fixe. Donc par définition de $T_\M$, son image est dans $T_\M$. On prouve ainsi que $T\subset T_\M$.\\

Maintenant, le fait que $T_\K\subset T_\M$ entraine que $Z_{Z(T\K)}(T_\M) = Z_G(T_\M)$, donc finalement  $Z_G(T_\M)= T_\L$, et le groupe $G$ est lui aussi quasi-déployé sur $\M$.\\

On a finalement trois tores $T_\K\subset T_\M\subset T_\L$. Le premier est $\K$-déployé maximal, le second est $\M$-déployé maximal, le troisième est  maximal (et $\L$-déployé), et $\M$-défini. De plus l'extension de corps $\K\subset \M$ est modérément ramifiée, et $G$ est quasi-déployé sur $\M$.\\

On a déjà introduit $\Gamma_\M=\gal(\L|\M)$, on notera de plus $\Gamma_\K=\gal(\M|\K)$. Le groupe $\Gamma_\M$ est distingué dans $\Gamma$ et $\Gamma_\K\simeq \Gamma/\Gamma_\M$.
Soient $A_\L=A(T_\L)$, $\vec A_\L=\vec A(T_\L)$, $\vec A_\M=\vec A(T_\M)\subset \vec A_\L$ et $\vec A_\K= \vec A(T_\K)\subset \vec A_\M$. Comme $T_\L$ est $\M$-défini, le groupe $\Gamma_\M$ agit sur $A_\L$, et l'ensemble $A_\M= A_\L^{\Gamma_\M}$ est l'adhérence d'un espace affine sous $\vec A_\M$.\\

\rema On ne peut définir un espace $A_\K$ aussi simplement, il faudra attendre \ref{soussection:non ramifie}.

\subsubsection{Descente quasi-déployée}

On commence par appliquer la partie \ref{section:descente} aux données radicielles $\D_\M$ et $\D_\L$. D'après la proposition \ref{prop:verifications des premieres conditions}, toutes les conditions de descente sauf (DM 3) et (DM 4) sont vérifiées, en prenant $\I_\L^{\Gamma_\M}$ pour jouer le rôle de $\I\bec$, et en remplaçant $\phii$ par une valuation équipollente basée en un point de $A_\M$.\\

%Le fait que $G$ est quasi-déployé sur $\mathbb M$ signifie  qu'il vérifie les conditions équivalentes suivantes:
%\liste{
%\item Il existe un sous-groupe de Borel $\M$-défini.

%\item Il existe une chambre de $\iv$ stable sous $\Gamma_\M:=\gal(\L|\M)$.

%\item Les chambres de $\iv(\M)$ rencontrent des chambres de $\iv(\L)$.

%\item Si $T_\M\subset T$ est un tore $\M$-déployé maximal dans un tore maximal, alors $Z(T_\M)=T$.

%\item Si $T_\M$ est un tore $\M$-déployé maximal, alors $T:=Z(T_\M)$ est l'unique tore maximal contenant $T_\M$.
%\\
%}

Nous avons vu que $G$ est quasi-déployé sur $\M$, c'est-à-dire que nous avons trouvé un tore $\M$-déployé maximal $T_\M$ tel que $T_\L=Z(T_\M)$ est un tore maximal de $G$. Le tore $T_\L$ est donc $\M$-défini, et c'est l'unique tore maximal contenant $T_\M$.

On a $N(T_\M) \subset N(T_\L).Z(T_\M) = N(T_\L)$, d'où $N(T_\M)(\M)\subset N(T_\L)(\M)$. Comme la partie $A_\M=A(T_\L)^{\Gamma_\M}$ est en général stable par $N(T_\L)(\M)$, elle l'est ici aussi par $N(T_\M)(\M)$, ainsi la condition (DM 4) est-elle vérifiée.

Concernant (DM 3), soit $\vec C_\M$ une chambre de $\vec A_\M$ et $\vec C$ une chambre de $\vec A_\L$ rencontrant $\vec C_\M$. Soit $x\in A_\M$ et $F$ la facette de $A_\L$ contenant $\gm_x(x+\vec C)$. C'est une chambre de $\I$ qui coupe $A_\M$ et donc $\I^{\Gamma_\M}$.  Pour toute facette sphérique $\vec f\in\F(\vec A_\L)$, la facette contenant $pr_{\vec f}(F)$ coupe $\barre{A\die}$, est une chambre de $\I_{\vec f}$ et donc n'est incluse dans aucune autre facette de $\I_{\vec f}$.\\

On obtient donc par la partie \ref{section:descente} une valuation $\phii_\M$ de la donnée radicielle $\D_\M$, un appartement $A_\M=A_\M(T_\M)$, une bonne famille de parahoriques $Q_\M$ vérifiant \plienp(sph), puis une masure bordée $\I_\M=\I(Q_\M)$ pour $G(\M)$.\\

On notera $\mathcal J_\M = G(\M). A_\M\subset \I^{\Gamma_\M}$, c'est l'analogue des "points invariants ordinaires" de \cite{these-guy} 2.4.13. Par la proposition \ref{prop:injection des facades}, pour toute facette sphérique $\vec f_\M$ de $\iv_\M$, la façade $\I_{\M \vec f_\M}$ de $\I_\M$ s'identifie à $\mathcal J_\M \cap I_{\vec f}$, où $\vec f$ est l'unique facette de $\iv$ contenant un ouvert de $\vec f_\M$ (voir la remarque à la fin de \ref{soussoussection:rappels presque deployes}).\\

Soit $\sigma\in\Gamma$, alors $\sigma.T_\M$ est un autre tore $\M$-déployé maximal de $G$, donc il existe $g\in G(\M)$ tel que $\sigma.T_\M =g.T_\M.g\inv$. Montrons que $\sigma.A_\M = g.A_\M$.

 Pour commencer, $\sigma.A_\M =\sigma. (A_\L^{\Gamma_\M}) = (\sigma A_\L)^{\Gamma_\M}$. En effet, si $x\in A_\L^{\Gamma_\M}$, alors $\sigma.x\in \sigma.A_\L$, et pour tout $\gamma\in\Gamma_\M$, $\gamma.\sigma.x=\sigma \sigma\inv \gamma \sigma x = \sigma x$ car $\sigma\inv \gamma \sigma \in \Gamma_\M$. Donc $\sigma.(A_\L^{\Gamma_\M})\subset (\sigma A_\L)^{\Gamma_\M}$, l'autre inclusion est semblable.
 
 Donc $\sigma.A_\M = (\sigma A_\L)^{\Gamma_\M} = A_\L(\sigma T_\L)^{\Gamma_\M}$. Mais $\sigma.T_\L = gT_\L g\inv$ car c'est l'unique tore maximal contenant $\sigma T_\M$. Donc $\sigma .A_\M = (g.A_\L)^{\Gamma_\M} = g.A_\L^{\Gamma_\M }=g.A_\M$ car $g$ est $\Gamma_\M$-fixe.\\
 
 Ceci prouve que $\mathcal J_\M$ est stable par $\Gamma$.

\subsubsection{Descente modérément ramifiée}
\label{soussection:non ramifie}

On étudie maintenant le groupe $G(\K)$. On va utiliser les résultats de la partie \ref{section:descente}, appliqués aux données radicielles $\D_\L$ et $\D_\K$, le travail sur $\D_\M$ de la partie précédente servira à définir une partie $\I\die[\K]$ de $\I_\L$, et à prouver les conditions (DM x).

On a vu que $\Gamma$ agit sur $\mathcal J_\M$, on peut donc poser $\I\die[\K] =  \mathcal J_\M^\Gamma = \mathcal J_\M^{\Gamma_\K}$. La masure bordée $\I_\L$ admet des points $\Gamma$-fixes car $G$ contient des tores maximaux $\L$-déployés et $\K$-définis. Donc $\I_\L^\Gamma\not=\vide$. Cependant il n'est pas clair que $\mathcal J_\M^\Gamma = \I_\L^\Gamma\cap \mathcal J_\M$ soit non vide. Nous serons en fait obligés de le supposer, mais c'est une hypothèse qui, tout comme l'hypothèse sur l'existence de l'extension $\M$, est vérifiée dès que la valuation de $\K$ est discrète. C'est en fait l'analogue de la condition (DE) de \cite{bruhat-tits2} 5.1.5.

\fprop{\label{prop:dernieres conditions de descente} On suppose que la famille $Q$ vérifie \plienpm(sph), et que l'immeuble de Bruhat-Tits $\I_\M(Z(T_\K))$ du groupe réductif $Z(T_\K)$ sur le corps $\M$ admet un point $\Gamma$-fixe.

Alors il existe un tore maximal $T$ tel que les conditions de descente de la partie \ref{section:descente} sont vérifiées pour les données radicielles $\D_\L(T)$ et $\D_\K(T_\K)$, et pour la partie $\I\die[\K]$ de $\I_\L$.
}

Comme pour \ref{prop:existence de AK}, par le théorème de point fixe de Bruhat, $\I_\M(Z(T_\K))^\Gamma$ est non vide dès que l'immeuble $\I_\M(Z(T_\K))$ est complet, %c'est-à-dire dès que le corps $\L$ est maximalement complet (\cite{bruhat-tits2} 7.5.4 et 7.5.5)
 et ceci est vrai dès que la valuation de $\K$ est discrète.\\

\demo

%On rappelle que la condition \plienp(sph) pour tout point $x\in \I$ et tout vecteur $\vec v\in \iv\sph$, l'addition $x+\vec v$ est bien définie (\ref{prop:csq de plienp}, voir aussi \ref{cor:facette et facette vectorielle dans un appart}). En conséquence, pour tout $g\in G$ ou $\sigma\in\Gamma$, on a $g.(x+\vec v) = g(x)+g(\vec v)$ et $\sigma(x+\vec v) = \sigma(x) +\sigma(\vec v)$. Par exemple, si $x$ est un point $\Gamma$-fixe dans la façade principale d'un appartement $A$, alors $A^\Gamma \supset \barre{x+\vect(\vec A^\Gamma)}$.\\

 Pour toute facette sphérique $\vec f$ de $\iv_\L$, $\I\die[\K]\cap \I_{\L,\vec f}$ est vide si $\vec f\cap \iv_\M=\vide$. Sinon, il s'agit de l'ensemble des points $\Gamma$-fixes dans la façade $\I_{\M \vec f_\M}$ où $\vec f_\M$ est une facette de $\iv_\M$ contenant $\vec f^{\Gamma_\M}$. Cette façade est un immeuble, et $\Gamma$ y agit par automorphismes, donc $\I_{\M \vec f_\M}^\Gamma$ est une partie convexe. Ainsi (DM 1) est vérifié, pour la partie d'immeuble $\I\die[\K]$.

La proposition \ref{prop:verifications des premieres conditions} prouve encore les conditions (DSR), (DDR), (DIV), et (DV 2), pour n'importe quel tore maximal $T$ contenant $T_\K$ et pour n'importe quelle valuation de $\D_\L(T)$. Il reste à voir (DM 2, 3 et 4) ainsi que (DV 1).\\

Le fait que l'extension $\K\subset \M$ soit modérément ramifiée entraine que le lieu $\I_\M(Z(T_\K))^{\Gamma_\K}$ des points $\Gamma_\K$-fixes de l'immeuble $\I_\M(Z(T_\K))$ de $Z(T_\K)$ sur le corps $\M$ est de diamètre nul. En effet le groupe $Z(T_\K)$, ou plutôt son semi-simplifié $Z(T_\K)/Z(Z(T_\K))$ est $\K$-anisotrope, donc sachant que le degré de sauvagerie $s(\M|\K)$ est nul, c'est la proposition 5.2.1 de \cite{these-guy}. Comme nous avons supposé qu'il est non vide, il s'agit d'un singleton $\{p\}$, et le point $p$ est donc fixé par $Z(T_\K)(\K)$.\\

 Le groupe $Z(T_\K)$ est le fixateur dans $G$ du $\K$-appartement $\vec A_\K$.
 Pour toute chambre $\vec f_\K$ de $\vec A_\K$, on notera $\vec f$ la facette de $\vec A$ contenant un ouvert de $\vec f_\K$. La facette $\vec f$ est sphérique, et 
  %Soit $\vec f$ une facette de $\vec A$ contenant un ouvert de $\vec A_\K$, alors $Z(T_\K)= M_{\vec A}(\vec f)$. De plus, $\vec f$ est sphérique, et 
  l'immeuble $\I(Z(T_\K))$ est isomorphe à $\I_{\vec f}$, la façade de $\I=\I_\L(G)$ de type $\vec f$.
%Il est également isomorphe à la façade  $\I_{-\vec f}$, où $-\vec f$ est la facette opposée à $\vec f$ dans $\vec A$.
 L'immeuble $\M$-rationnel $\I_\M(Z(T_\K))$ est alors isomorphe à $\I_{\vec f}\cap \mathcal J_\M$, par la proposition \ref{prop:injection des facades}. Cet immeuble est inclus dans $\I^{\Gamma_\M}$, de sorte que l'action de $\Gamma$ y coïncide avec celle de $\Gamma_\K$.

On note $p_{\vec f}$ l'unique point $\Gamma$-fixe de $\I_{\vec f}\cap \mathcal J_\M$. Soit $\mathcal E$ l'ensemble des points de $\mathcal J_\M$ ainsi obtenus pour toutes les chambres de $\vec A_\K$. L'ensemble $\mathcal E$ est donc fixé par $Z(T_\K)(\K)$, et stabilisé par $N(T_\K)(\K)$.\\

Il existe $z\in Z(T_\K)(\M)$ tel que le point $\Gamma$-fixe de $\I_\M(Z(T_\K))$ est dans l'appartement correspondant au tore $T:= zT_\L z\inv$. Soit $Z=A(T) = z.A_\L$. Alors $Z^{\Gamma_\M}$ est inclus dans $\mathcal J_\M$ et contient $\mathcal E$. Notons que comme $z\in Z(T_\K)$, $\vec A_\K\subset \vec Z$. Nous allons montrer que $Z\cap \I\die[\K]$ est l'adhérence d'un espace affine dirigé par $\vect(\vec A_\K)$, et qu'il est stable par $N(T_\K)(\K)$.\\

  Comme $\mathcal E$ est constitué de points sphériques des deux signes, $\cl_Z(\mathcal E)$ est en fait indépendant de l'appartement le contenant considéré (corollaire \ref{cor:enclos bien def}). En conséquence, cet enclos est stable par $\Gamma$ et par $N(T_\K)(\K)$. Il s'agit d'une partie convexe de $Z$, donc l'action de $\Gamma$ fixe un point $x\in \inte Z\cap \cl(\mathcal E)$.
  % Grâce a \plienp(sph), on obtient que $\barre{ x+\vect(\vec A_\K) }$ est $\Gamma$-fixe.
 Soit $\sigma\in\Gamma$. Par le corollaire \ref{cor:enclos bien def}, il existe $g\in Q(\cl(\mathcal E))$  tel que $\sigma.Z=g.Z$. Alors $g\inv\sigma$ est un automorphisme de $Z$ qui fixe $\mathcal E$, donc sa partie vectorielle fixe $\vec A_\K$. Il fixe de plus $x$, et ceci entraine qu'il fixe $\barre{x+\vect(\vec A_\K)}$. Comme $g$ fixe $\cl(\mathcal E)$ qui contient $\barre{x+\vect(\vec A_\K)}$, on voit que $\sigma$ fixe $\barre{x+\vect(\vec A_\K)}$. Finalement, $\barre{x+\vect(\vec A_\K)}\subset Z^\Gamma$. Comme $Z^\Gamma\subset Z^{\Gamma_\M} \subset \mathcal J_\M$, on obtient $\barre{x+\vect(\vec A_\K)}\subset Z\cap \I\die[\K]$.\\

  Réciproquement, soit $y\in Z\cap \I\die[\K]$, montrons que $y\in \barre{ x+\vect(\vec A_\K) }$. Supposons dans un premier temps $y\in\inte Z$.
 Soit $\vec f$ une facette de $ \vec Z$ contenant un ouvert d'une chambre de $\vec A_\K$. Le point $pr_{\vec f}(y)$ est $\Gamma$-fixe, et comme $y\in\mathcal J_\M$, c'est un point de $\I_{\vec f}\cap  \mathcal J_\M$, c'est donc $p_{\vec f}$. Ceci et le résultat similaire pour $-\vec f$ prouve déjà que $y\in\cl(\mathcal E)$. Ensuite, comme $x$ et $y$ ont la même projection sur $\I_{\vec f}$, quitte à déplacer $y$ selon $\vect(\vec A_\K)$, on peut supposer $y\in x+\vec f$. Soit $\vec v\in\vec f$ tel que $y=x+_Z\vec v$. Alors pour tout $\sigma\in\Gamma$, $y=\sigma(y) = \sigma(x) +_{\sigma Z} \sigma(\vec v) = x +_{\sigma Z} \sigma(\vec v)$.
Comme au paragraphe précédent, soit $q\in Q(\cl(\mathcal E))$ tel que $\sigma.Z=q.Z$. Alors $q\inv \sigma.y=y$ (car $y\in\cl(\mathcal E)$) d'où, dans $Z$,  $x+_Z\vec v = x+_Zq\inv\sigma(\vec v)$, donc $\vec v = q\inv \sigma (\vec v)$, et comme $q$ fixe $\vec f$, $\vec v= \sigma(\vec v)$.

Ainsi, $\vec v\in\vec f^\Gamma \subset \vec A_\K$, et $y\in x+\vect(\vec A_\K)$. De la même manière, on obtient pour toute façade $Z_{\vec g}$ telle que $\vec g\cap \vec A_\K\not=\vide$, $Z_{\vec g}\cap \I\die[\K] = pr_{\vec g}(x) +\vect(\vec A_\K)$.

Ainsi, $Z\cap \I\die[\K]$ est l'adhérence d'un espace affine sous $\vect(\vec A_\K)$. On le note $Z_\K$, son intérieur jouera le rôle du $A\die$ de la partie \ref{section:descente}.

La condition (DM 2) est donc vérifiée pour l'appartement $Z$ (c'est-à-dire le tore maximal $T$ ou la donnée radicielle $\D_\L(T)$). De plus nous avons vu que $Z_\K\subset \cl(\mathcal E)$ donc $Z_\K=\cl(\mathcal E)\cap \I\die[\K]$, et comme ces deux ensembles sont stabilisés par $N(T_\K)(\K)$, $Z_\K$ aussi, d'où (DM 4).\\

Étudions (DM 3). Soit $\vec g$ une facette sphérique de $\vec Z$ coupant $\vec A_\K$.  Soit $F=\gm_x(x+\vec F)$ une facette de $Z_{\vec g}$ coupant $Z_\K$ de dimension maximale, donc $F$ contient un ouvert de $Z_\K\cap Z_{\vec g}$, et il existe une facette $\vec f$ de $\vec Z\cap \vec g^*$ contenant un ouvert de $\vec A_\K$ telle que $\vec F= \vec f/\vect(\vec g)$ (ou plutôt $\vec F=(\vec f+\vect(\vec g)).\vect(\vec g)$.
 Supposons qu'il existe une autre facette $F'$ de $\I_{\vec g}$ rencontrant $\I\die[\K]$ et telle que $F\subset \bar F'$. 
Il existe un appartement $B_{\vec g}$ de $\I_{\vec g}$ contenant $F'$ et tel que $\vec f\subset \vec B_{\vec g}$.
 Soient $x\in F'\cap \I\die[\K]$ et $y\in F\cap \I\die[\K]$. Alors $pr_{\vec f}(x)$ tout comme $pr_{\vec f}(y)$ sont deux points $\Gamma$-fixes dans $\mathcal J_\M\cap \I_{\vec f}$: ils sont égaux. Donc $x\in y+\vect_{\vec B_{\vec g}}(\vec f) = \aff_{B_{\vec g}}(F)$. Mais $x\in F'$ et $F'\cap \aff(F) = \vide$: on obtient une contradiction, et il n'existe pas de telle facette $F'$.\\

Enfin, soit $\phii$ une valuation de $\D_\L(T)$ basée en un point  $o\in Z_\K$, elle vérifie immédiatement (DV 1).

  \cqfd

\subsection{Conclusion}

Résumons les résultats précédents. Soit $G$ un groupe de Kac-Moody presque déployé sur un corps valué $\K$, déployé sur la clôture séparable de $\K$. Soit $T_\K$ un tore $\K$-déployé maximal, il existe une extension galoisienne $\L$ de $\K$ qui déploie $G$ et telle qu'il existe des tore maximaux $\L$-déployés contenant $T_\K$.

 Pour tout tel tore $T$, la famille de parahoriques $Q$ définie dans \cite{masures2} est une bonne famille de parahoriques pour $\D_\L(T)$, et elle vérifie en outre \plienp(sph). Alors la proposition \ref{prop:verifications des premieres conditions} s'applique, permettant de vérifier les conditions de descente (DSR), (DDR) et (DIV), ainsi que (DV 2) pour toute valuation sur $\D_\L(T)$.\\

Si de plus la valuation de $\K$ est discrète, et le corps résiduel parfait, alors il existe une extension intermédiaire $\M\subset \L$ telle que $G$ est quasi-déployé sur $\M$ et telle que l'extension $\K\subset \M$ est non ramifiée. Ceci permet la définition de la partie $\I\die[\K]=\mathcal J_\M\cap \I^\Gamma$. L'hypothèse de discrétion de la valuation de $\K$ permet également d'appliquer la proposition   \ref{prop:dernieres conditions de descente}, prouvant que la partie $\I\die[\K]$ vérifie les conditions (DM).
 La condition (DV 1) est obtenue dès qu'on choisit une valuation basée en un point de $\I\die[\K]$ , alors toutes les conditions de descente de la partie \ref{section:descente} sont vérifiées.\\

On obtient donc une valuation pour la donnée radicielle $\D_\K$, puis un appartement, une bonne famille de parahoriques vérifiant \plienp(sph), et enfin une masure bordée. On sait en outre que les façades sphériques de cette masure bordée sont incluses dans des façades sphériques de la masure bordée $\I_\L$ pour $\D_\L$.

\begin{prout}\label{th:final}
Soit $G$ un groupe de Kac-Moody presque déployé sur un corps $\K$, déployé sur la clôture séparable de $\K$. On suppose $\K$ muni d'une valuation réelle discrète non triviale, telle que son corps résiduel soit parfait.

 Alors il existe une masure bordée $\I_\K$ pour $G(\K)$, qui provient d'une valuation $\phii_\K$ et d'une bonne famille de parahoriques $Q_\K$ vérifiant \plienp(sph). Pour toute facette sphérique $\vec f_\K$ de $\iv(\K)$, la façade $\I_{\K,\vec f_\K}$ s'injecte dans la façade $\I_{\L,\vec f}$, de la masure bordée $\I_\L$ pour $G(\L)$, où $\vec f$ est la facette de $\iv(\L)$ contenant un ouvert de $\vec f_\K$.\\
\end{prout}

%En fait, il est possible d'alléger un peu les hypothèses sur $\K$

\rema Les hypothèses sur le corps $\K$ (valuation discrète et corps résiduel parfait) interviennent pour résoudre deux difficultés: pour assurer l'existence d'une extension $\M$ non ramifiée de $\K$ qui quasi-déploie le groupe réductif $Z(T_\K)$, puis pour assurer l'existence d'un point $\Gamma$-fixe dans l'immeuble de ce dernier. Ces deux difficultés ne font intervenir qu'un groupe réductif, et sont rencontrées de la même manière dans \cite{bruhat-tits2}. Ainsi, si on veut affiner le résultat précédent en affaiblissant les hypothèses sur le corps $\K$, ceci devrait être possible de la même manière que dans \cite{bruhat-tits2}.\\

\subsection{Questions}

Signalons finalement deux points qui restent non résolus.\\

 En premier lieu, on ne sait pas s'il existe en général, pour toute donnée radicielle valuée, une bonne famille de parahoriques. Nous ne disposons a priori que de la famille minimale de parahoriques; même la définition de la famille maximale n'est possible que si l'on suppose l'existence d'au moins une bonne famille. Ce n'est que dans le cas d'une donnée radicielle valuée venant d'un groupe de Kac-Moody que l'on sait, grâce à \cite{masures2} que la famille minimale est bonne, et qu'il existe en outre une bonne famille vérifiant en plus \plienp.\\
 
 Par ailleurs, pour construire la masure bordée d'un groupe de Kac-Moody presque déployé, on définit un appartement, puis une famille de parahoriques, puis on applique la construction générale. Il n'est alors pas clair que la masure obtenue s'injecte (ou au moins que chacune de ses façades s'injecte) dans la masure du groupe déployé. On a seulement prouvé que ses façades sphériques s'injectent dans des façades sphériques de la masure du groupe déployé, et il est facile d'en déduire l'existence d'un plongement pour les immeubles microaffines de $G(\K)$.

\bibliographystyle{alpha}
\bibliography{biblio}

\vspace{.5cm}

\begin{flushright}
\begin{minipage}{7cm}
Cyril Charignon\\
Institut \'Elie Cartan\\
Unité mixte de recherche 7502\\
Nancy-Université, CNRS, INRIA\\
Boulevard des aiguillettes\\
BP 70239\\
54506 Vandoeuvre lès Nancy cedex (France)\\
cyril.charignon@iecn.u-nancy.fr
\end{minipage}\end{flushright}

\end{document}